\documentclass[a4paper]{amsart}
\usepackage{hyperref}
\usepackage{fullpage}
\usepackage{txfonts, amsmath,amstext,amsthm,amscd,amsopn,verbatim,amssymb, amsfonts}

\usepackage[bbgreekl]{mathbbol}
\usepackage{array}
\usepackage{tikz, tikz-cd}
\usetikzlibrary{matrix}
\usetikzlibrary{shapes}
\usetikzlibrary{arrows}
\usetikzlibrary{calc,3d}
\usetikzlibrary{decorations,decorations.pathmorphing}
\usetikzlibrary{through}
\tikzset{ext/.style={circle, draw,inner sep=1pt},int/.style={circle,draw,fill,inner sep=1pt, minimum size=5pt},nil/.style={inner sep=1pt}}
\tikzset{exte/.style={circle, draw,inner sep=3pt},inte/.style={circle,draw,fill,inner sep=3pt}}
\tikzset{diagram/.style={matrix of math nodes, row sep=3em, column sep=2.5em, text height=1.5ex, text depth=0.25ex}}
\tikzset{diagram2/.style={matrix of math nodes, row sep=0.5em, column sep=0.5em, text height=1.5ex, text depth=0.25ex}}
\theoremstyle{plain}
  \newtheorem{thm}{Theorem}  \newtheorem{defi}[thm]{Definition}
  \newtheorem{prop}[thm]{Proposition}
  \newtheorem{defprop}[thm]{Definition/Proposition}
  \newtheorem{cor}[thm]{Corollary}
  \newtheorem{lemma}[thm]{Lemma}
\theoremstyle{definition}
  \newtheorem{ex}[thm]{Example}
  \newtheorem{rem}[thm]{Remark}

\newcommand{\alg}[1]{\mathfrak{{#1}}}
\newcommand{\co}[2]{\left[{#1},{#2}\right]} 
\newcommand{\eref}[1]{\eqref{#1}} 
\newcommand{\pd}[2]{ { \frac{\partial {#1} }{\partial {#2} } } }

\newcommand{\p}{\partial}
\newcommand{\Hom}{\mathop{Hom}}

\newcommand{\C}{{\mathbb{C}}}
\newcommand{\R}{{\mathbb{R}}}
\newcommand{\Z}{{\mathbb{Z}}}

\newcommand{\FF}{{\mathbb{F}}}

 \newcommand{\Graphs}{{\mathsf{Graphs}}}
\newcommand{\fGraphs}{{\mathsf{fGraphs}}}
\newcommand{\pdu}{{}^*} \newcommand{\fml}{{\mathit{fml}}} 
\newcommand{\SGraphs}{{\mathsf{SGraphs}}}

\newcommand{\Gr}{{\mathsf{Gra}}}
\newcommand{\Gra}{{\mathsf{Gra}}}
\newcommand{\dGr}{{\mathsf{dGra}}}
\newcommand{\dGra}{{\mathsf{dGra}}}

\newcommand{\SGr}{{\mathsf{SGra}}}
\newcommand{\SGra}{{\mathsf{SGra}}}
\newcommand{\fSGra}{{\mathsf{fSGra}}}
\newcommand{\PT}{ \mathsf{PT} }

\newcommand{\Tw}{\mathit{Tw}}
\newcommand{\Def}{\mathrm{Def}}

\newcommand{\Fund}{\mathrm{Fund}}

\newcommand{\bigGra}{\mathsf{bigGra}}
\newcommand{\op}{\mathcal}
\newcommand{\hKS}{\mathsf{hKS}}
\newcommand{\KS}{\mathsf{KS}}
\newcommand{\homKS}{\mathsf{homKS}}
\newcommand{\Br}{\mathsf{Br}}
\newcommand{\hBr}{\mathsf{hBr}}
\newcommand{\Lie}{\mathsf{Lie}}
\newcommand{\hLie}{\mathsf{hLie}}
\newcommand{\ELie}{\mathsf{ELie}}
\newcommand{\Ger}{\mathsf{Ger}}
\newcommand{\hooLie}{\mathsf{hoLie}}
\newcommand{\hoELie}{\mathsf{hoELie}}
\newcommand{\fSGraphs}{\mathsf{fSGraphs}}
\newcommand{\SC}{\mathsf{SC}}
\newcommand{\SG}{\mathsf{SG}}

\newcommand{\ESC}{\mathsf{ESC}}
\newcommand{\EESC}{\mathsf{EESC}}

\newcommand{\calc}{\mathsf{calc}}

\newcommand{\bigChains}{\mathsf{bigChains}}
\newcommand{\bigGraphs}{\mathsf{bigGraphs}}
\newcommand{\FM}{\mathsf{FM}}
\newcommand{\EFM}{\mathsf{EFM}}
\newcommand{\cFM}{\mathsf{EFM}}
\newcommand{\hc}{\mathit{hc}}

\newcommand{\vout}{\mathit{out}}
\newcommand{\vin}{\mathit{in}}
\newcommand{\conn}{\mathit{conn}}

\newcommand{\bpm}{\begin{pmatrix}}
\newcommand{\epm}{\end{pmatrix}}
\newcommand{\Tpoly}{T_{\rm poly}}
\newcommand{\Dpoly}{D_{\rm poly}}

\newcommand{\GC}{\mathrm{GC}}

\newcommand{\bigV}{\mathbf{V}}
\newcommand{\mU}{\mathcal{U}}

\DeclareMathOperator{\dv}{div}
\DeclareMathOperator{\End}{End}
\DeclareMathOperator{\sgn}{sgn}

\newcommand{\SymGrp}{S}

\tikzset{col1/.style={dashed}}
\tikzset{col2/.style={very thick, dashed}}
\tikzset{col3/.style={}}
\tikzset{col4/.style={very thick}}

\newcommand{\collegend}
{
\begin{tabular}{|ll|}
 \hline
\begin{tikzpicture}[baseline=-.65ex]
\draw[col1] (0,0) --(1,0);
\end{tikzpicture}
&: Color 1 \\
\begin{tikzpicture}[baseline=-.65ex]
\draw[col2] (0,0) --(1,0);
\end{tikzpicture}
&: Color 2 \\
\begin{tikzpicture}[baseline=-.65ex]
\draw[col3] (0,0) --(1,0);
\end{tikzpicture}
&: Color 3 \\
\begin{tikzpicture}[baseline=-.65ex]
\draw[col4] (0,0) --(1,0);
\end{tikzpicture}
&: Color 4 \\
 \hline
\end{tabular}
}

\begin{document}
\title{The Homotopy Braces Formality Morphism}
\author{Thomas Willwacher}
\address{Institute of Mathematics\\ University of Zurich\\ Winterthurerstrasse 190 \\ 8057 Zurich, Switzerland}
\email{thomas.willwacher@math.uzh.ch}


\begin{abstract}
We extend M. Kontsevich's formality morphism to a homotopy braces morphism and to a homotopy Gerstenhaber morphism. We show that this morphism is homotopic to D. Tamarkin's formality morphism, obtained using formality of the little disks operad, if in the latter construction one uses the Alekseev-Torossian associator. 
Similar statements can also be shown in the ``chains'' case, i.~e., on Hochschild homology instead of cohomology.  
This settles two well known and long standing problems in deformation quantization and unifies the several known graphical constructions of formality morphisms and homotopies by Kontsevich, Shoikhet, Calaque, Rossi, Alm, Cattaneo, Felder and the author.
\end{abstract}
\maketitle
\setcounter{tocdepth}{1}
\tableofcontents

\makeatletter{}\section{Introduction}
The formalism of deformation quantization introduced by Bayen et al. \cite{BFFLS, BFFLS2} is an attempt to understand the transition between classical and quantum mechanics in physics within the framework of algebraic deformation theory in mathematics. In mathematical terms, the central problem is to classify all (formal) deformations of the associative algebra of smooth functions on a manifold $M$. Such deformations are governed in a precise sense by the Lie algebra of multidifferential operators on $M$, which can be considered as a smooth version of the Hochschild complex of $C^\infty(M)$.
M. Kontsevich solved the central problem of deformation quantization by showing that the multidifferential operators are formal as a Lie algebra, i.e., quasi-isomorphic to their cohomology \cite{K1}. This reduced the classification of formal deformations of the algebra structure on $C^\infty(M)$ to the study of formal Poisson structures on $M$.

Let us describe Kontsevich's result in more detail, restricting to the case $M=\R^n$ for now.
Let $\Tpoly$ be the space of multivector fields on $\R^n$, i.e., the space of smooth sections of exterior powers of the tangent bundle $\Lambda T\R^n$. It is equipped with a natural graded commutative (wedge) product, and a compatible Lie bracket (the Schouten-Nijenhuis bracket), endowing $\Tpoly$ with the structure of a Gerstenhaber algebra.
Let $\Dpoly$ the space of normalized\footnote{This means that the multidifferential operators vanish if one inserts a constant function in any of its slots.} multidifferential operators on $\R^n$. The Hochschild differential and the Gerstenhaber bracket endow $\Dpoly[1]$ with a natural differential graded Lie algebra structure.
Furthermore, there is a (non-commutative) product operation on $\Dpoly$ through the cup product. The cup product and the Gerstenhaber bracket are compatible in the sense that these algebraic structures may be extended to an algebra structure over the braces operad, which encodes higher compatibility relations, see \cite{GJ}.

M. Kontsevich's famous Formality Theorem \cite{K1} now asserts the existence of a $\Lie_\infty$ (homotopy Lie) quasi-isomorphism
\[
\mathcal{U} \colon \Tpoly[1] \to \Dpoly[1].
\]
The components of this $\Lie_\infty$-morphism are furthermore given by explicit formulas of the form
\[
 \mathcal{U}_n(\gamma_1,\dots, \gamma_n) = \sum_{\Gamma} \left(\int_{C_\Gamma} \omega_\Gamma\right) D_\Gamma(\gamma_1,\dots,\gamma_n)
\]
for $\gamma_1,\dots,\gamma_n\in \Tpoly$. Here the sum is over a set of graphs (\emph{Kontsevich graphs}), the integral ranges over a compact manifold with corners $C_\Gamma$, the differential form $\omega_\Gamma$ is associated to the graph $\Gamma$ using Feynman rules, and $D_\Gamma(\cdots)$ is a certain multidifferential operator. Kontsevich's construction was very influential and his methods have been applied successfully to many different problems in the field, and to some beyond.

However, Kontsevich's result is unsatisfying in the sense that the algebraic structure on both $\Tpoly$ and $\Dpoly$ is richer than just that of a Lie algebra, namely, there is a Gerstenhaber structure on $\Tpoly$ and a braces structure on $\Dpoly$. A natural question is whether one can find a quasi-isomorphism $\Tpoly\to\Dpoly$ respecting those structures, up to homotopy. 
This question has been answered positively by D. Tamarkin \cite{tamanother, hinich, vasilyPaljug}. Roughly his proof goes as follows, where we restrict ourselves to the algebraic setting for simplicity.
\begin{enumerate}
 \item Find a $\Ger_\infty$ (homotopy Gerstenhaber) structure on the algebraic multidifferential operators $\Dpoly^{alg}\subset \Dpoly$, reducing on cohomology to the same Gerstenhaber structure as that coming from the natural braces structure. This subproblem is called the \emph{Deligne conjecture}. There are by now several solutions to this conjecture, the most relevant for us are due to D. Tamarkin and M. Kontsevich and Y. Soibelman. For this step of the proof one needs to choose a Drinfeld associator or equivalent data.
 \item One can show that the $\Ger_\infty$-structure on the algebraic multivector fields $\Tpoly^{alg}\subset \Tpoly$ is \emph{rigid}, i.~e., it can not be deformed in a homotopy non-trivial way (see, e.g., \cite[Appendix B]{vasilyPaljug}). Using this rigidity result one can construct a $\Ger_\infty$-quasi-isomorphism $\Tpoly^{alg}\to\Dpoly^{alg}$ which is essentially unique up to homotopy. By restriction one of course also obtains an $\Lie_\infty$-morphism $\Tpoly^{alg}[1]\to \Dpoly^{alg}[1]$.

\end{enumerate}

Note that D. Tamarkin's proof is very non-constructive. Explicit formulas for the components of the constructed $\Ger_\infty$-morphism are difficult to attain. Also, the relation to M. Kontsevich's morphism, which naturally restricts to the algebraic setting, is not clear. Hence there remain the following questions, addressed in this paper:
\begin{enumerate}
 \item Are D. Tamarkin's and M. Kontsevich's $\Lie_\infty$-morphisms the same up to homotopy?
 \item Can M. Kontsevich's $\Lie_\infty$-morphism be extended to a $\Ger_\infty$-morphism?
 \item Can this be done with explicit formulas for all components in the form
\[
 \sum_{\Gamma} (\text{some number}) D_\Gamma(\gamma_1,\dots,\gamma_n)
\]
where the (some number)s are given using Feynman rules?
\end{enumerate}

For the first question one needs to be a bit more precise, since there are many ways to solve the Deligne conjecture, and furthermore the solution depends on the choice of a Drinfeld associator. However, we can give the following answer.
\begin{thm}
\label{thm:kontstamequiv}
 D. Tamarkin's and M. Kontsevich's $\Lie_\infty$ morphisms are homotopic, if one uses D. Tamarkin's solution of the Deligne conjecture via the formality of the little disks operad \cite{tamarkin}, and therein one uses the Alekseev-Torossian associator \cite{pavol, ATassoc}.
\end{thm}

\begin{rem}
 We do not know the answer to the first question if one takes Tamarkin's original solution to the Deligne conjecture via Etingof-Kazhdan quantization.
\end{rem}

The answer to the second question is yes.
\begin{thm}
\label{thm:kontsextends}
 M. Kontsevich's $\Lie_\infty$-morphism can be extended to a $\Ger_\infty$-morphism $\Tpoly\to \Dpoly$, for some $\Ger_\infty$-structure on $\Dpoly$ constructed below.
\end{thm}

To answer the third question, we need to change the question a bit. Note that so far we searched for morphisms in the $\Ger_\infty$ setting, retaining the natural Gerstenhaber structure on $\Tpoly$ and changing the natural braces structure on $\Dpoly$ to an unnatural $\Ger_\infty$-structure. In our point of view this is not fortunate and destroys the natural structure of formulas. What we should rather do is search for a $\Br_\infty$ (homotopy braces) morphism, retaining the natural braces structure on $\Dpoly$, and putting some $\Br_\infty$-structure on $\Tpoly$.

\begin{thm}
\label{thm:brinfty}
 There is a $\Br_\infty$-structure $\mu$ on $\Tpoly$ whose components are given by explicit formulas of the form
\begin{equation*}
 \mu_o(\gamma_1,\dots, \gamma_n) = \sum_{\Gamma} \left(\int_{c_o} \omega_\Gamma \right) V_\Gamma(\gamma_1,\dots,\gamma_n)
\end{equation*}
Here the sum is over a certain set of graphs, $c_o$ is a chain in some compact configuration space, depending on the operation $o$, and $V_\Gamma(\gamma_1,\dots,\gamma_n)$ is a multivector field depending combinatorially on $\Gamma$. The induced Gerstenhaber structure on $\Tpoly$ is the usual one.

There is furthermore a $\Br_\infty$-morphism
\[
 \mathcal{U} \colon \Tpoly \to \Dpoly
\]
whose components are given again by explicit formulas of the form
\begin{equation}
\label{equ:Udef}
 \mathcal{U}_o(\gamma_1,\dots, \gamma_n) = \sum_{\Gamma} \left(\int_{\tilde{c}_o\in C_\bullet(C_\Gamma)} \omega_\Gamma \right) D_\Gamma(\gamma_1,\dots,\gamma_n).
\end{equation}
In particular, restricting this morphism to the $\Lie_\infty$ part one recovers M. Kontsevich's formality morphism.
\end{thm}
The formulas above and the notation used will be explained in detail below. The chains $c_o$ and $\tilde{c}_o$ can be specified explicitly, at least up to contractible choices. 
After M. Kontsevich's and D. Tamarkin's seminal papers, there have appeared several extension and variations of their results.
In particular, there is a version for Hochschild homology of the algebra $A=C^\infty(\R^n)$ \cite{shoikhet, dolgushev} and a version for cyclic cohomology \cite{mecyccochains}.
First, let $C_\bullet$ be the (continuous) Hochschild chain complex of $A=C^\infty(\R^n)$ and $\Omega_\bullet$ the differential forms on $\R^n$, with negative grading. $C_\bullet$ is naturally a module over the dg Lie algebra $\Dpoly[1]$, while $\Omega_\bullet$ is naturally a module over $\Tpoly[1]$. By the $\Lie_\infty$-morphism $\Tpoly[1]\to\Dpoly[1]$ one obtains a structure of $\Lie_\infty$-module over $\Tpoly$ on $C_\bullet$.
It has been shown by B. Shoikhet \cite{shoikhet}, and independently by D. Tamarkin and B. Tsygan \cite{tamtsy}, that there is a quasi-isomorphism of $\Lie_\infty$-modules
\[
 C_\bullet \to \Omega_\bullet.
\]
The components of this morphism are constructed by giving explicit formulas resembling those for M. Kontsevich's morphism.
A globalized version has been proven by V. Dolgushev \cite{dolgushev, dolgushev-2005}.
Similarly to the cochains case, there is a much richer structure on $\Omega_\bullet$ than that of a module over $\Tpoly$. Concretely, there is a \emph{calculus} structure on the pair $(\Tpoly,\Omega_\bullet)$. It consists of a Gerstenhaber structure on $\Tpoly$, an operation $d$ of degree -1, which we take as the de Rham differential, and a degree zero operation
\[
 \iota \colon \Tpoly \otimes \Omega_\bullet \to \Omega_\bullet.
\]
In our case the operaton $\iota$ is given by contractions of multivector fields and differential forms.
$d$ and $\iota$ satisfy the following axioms
\begin{itemize}
\item $d^2 = 0$.
\item $\iota$ makes $\Omega_\bullet$ into a module over the graded commutative algebra $\Tpoly$.
\item $\co{\iota_a}{L_b} =\iota_{\co{a}{b}}$ for any $a,b\in \Tpoly$, where $L_b=\co{d}{\iota_b}$.
\end{itemize}

From the first and third axiom it follows that $a\mapsto L_a:=\co{d}{\iota_a}$ defines a Lie algebra action of $\Tpoly[1]$ on $\Omega_\bullet$. From the second axiom it then follows that 
\[
 L_a \iota_b +(-1)^{|a|} \iota_a L_b = L_{a\wedge b}.
\]
We call the 2-colored operad governing calculus structures $\calc$, following \cite{DTT}. 

On the other side, i.~e., on the pair $(\Dpoly, C_\bullet)$, there is similarly a natural algebraic structure extending the braces structure on $\Dpoly$ and the module structure on $C_\bullet$. This structure was first described (in this form) by Kontsevich and Soibelman \cite{KS2}, We call it the structure of a $\KS$-algebra, and the governing 2-colored operad accordingly $\KS$. Concretely, $\KS$ governs pairs $(A,M)$, where 
$A$ is a braces-algebra, that acts on $M$ with various operations extending the operations $\iota$ and $d$ of $\calc$. In this situation we will also call $M$ a \emph{braces-module} over $A$.
For a concrete description of the $\KS$-operad we refer to Section \ref{sec:KS1} or to \cite{KS2}. For now, we just note that the cohomology of $\KS$ is $\calc$.

It is a natural question whether there exists a formality morphism $C_\bullet \to \Omega_\bullet$ that is compatible with these additional structures, up to homotopy. This question has been answered positively by Dolgushev, Tamarkin and Tsygan \cite{DTT}. They construct a morphism of $\calc_\infty$-algebras, however using non-explicit methods. In this paper, we will reformulate the problem slightly and then give more or less explicit formulas for the morphism. Concretely, let $\homKS$ be a 4-colored operad which governs quadruples $(A_1,M_1, A_2,M_2)$ with the following structures:
\begin{enumerate}
 \item A $\KS$-structure on $(A_1,M_1)$.
 \item A $\KS_\infty$-structure on $(A_2,M_2)$.
 \item A $\Br_\infty$-map $A_2\to A_1$.
 \item A map $M_1\to M_2$ of $\Br_\infty$-modules over $A_2$, where the $\Br_\infty$-module structure on $M_1$ is obtained by pullback along the map $A_2\to A_1$.
\end{enumerate}
A more precise description of $\homKS$ will be given below.
We prove the following Theorem 
\begin{thm}
\label{thm:ksinfty}
 There is a representation of $\homKS$ on the quadruple $(\Dpoly, C_\bullet, \Tpoly,\Omega_\bullet)$, extending the usual $\KS$-structure on $(\Dpoly, C_\bullet)$, and the $\Br_\infty$-structure and -map $\Tpoly \to \Dpoly$ from the previous Theorem. The induced calculus structure on $\Tpoly$ and $\Omega_\bullet$ is the standard one. There are explicit integral formulas for the components of the $\homKS$-structure, resembling those of Kontsevich and Shoikhet. Concretely, the $\Lie_\infty$-part of the map $C_\bullet \to \Omega_\bullet$ agrees with Shoikhet's morphism.
\end{thm}

Formulas for some of the homotopies contained in the $\homKS$-structure have been found before: Compatibility of Dolgushev's map with the cyclic differential was shown in \cite{mecycchains}, compatibility with the cap product was shown by D. Calaque and C. Rossi \cite{CR}, and compatibility with another operation occuring in the Gauss-Manin connection by A. Cattaneo, G. Felder and the author in \cite{CFW}.
On the cochains side, some further homotopies have recently been found by Johan Alm \cite{alm}.
The present work gives a unifying framework for those results.

Using standard methods, the morphisms described here can be globalized to smooth manifolds $M$ other then $\R^n$, see section \ref{sec:globalization}, and similarly to complex manifolds and smooth algebraic varieties over $\C$. However, due to a non-explicit inversion in the globalization process the formulas no longer have the simple form \eqref{equ:Udef}.

\begin{rem}[The generalized Delinge conjecture for Hochschild cochains and chains]
The generalized Deligne conjecture for Hochschild cochains and chains states that the action of $\calc$ on the Hochschild cohomology and homology of any algebra may be lifted to a (co)chain level action of a 2-colored operad naturally extending the little disks operad, cf. \cite[section 11.3]{KS2} for more details. A proof of the conjecture was sketched by Kontsevich and Soibelman in loc. cit., but it contained considerable technical gaps. As part of the proof of Theorem \ref{thm:ksinfty} we fill these gaps and hence give the first complete proof of the generalized Deligne conjecture.
\end{rem}

\subsection{Sketch of the proof}
We want to construct a representation of the 4-colored operad $\homKS$ on the (4-colored) vector space 
\[
V=\Dpoly \oplus C_\bullet \oplus \Tpoly \oplus \Omega_\bullet \, .
\]
We do this by constructing colored operad maps
\begin{equation}
\label{equ:proofchain}
\homKS \to \bigChains \to \bigGra \to \End(V)\, .
\end{equation}
Here the operad $\bigChains$ is essentially the operad of chains on a topological (or rather, semi-algebraic) operad, made of configuration spaces of points.
In particular, it contains M. Kontsevich's configuration spaces. The first map is constructed similarly to the map in the Kontsevich-Soibelman paper \cite{KS1}. The operad $\bigGra$ is combinatorial and composed of suitable (``Feynman'') graphs. We will see that there is a natural representation of this operad on $V$. The map $\bigChains \to \bigGra$ is given by integral formulas analogous to the ones occurring in \cite{K1} (``Feynman rules'').
Finally, the proof of Theorem \ref{thm:kontstamequiv} is more or less a standard exercise, using the rigidity of the Gerstenhaber structure on $\Tpoly$.

\subsection{Structure of the paper}
In sections \ref{sec:grops}, \ref{sec:topops} and \ref{sec:KSdef} we define or recall the definitions of the operads $\bigGra$, $\bigChains$ and $\homKS$ occurring in \eqref{equ:proofchain}. In section \ref{sec:ksproof} and \ref{sec:opmaps} the maps \eqref{equ:proofchain} between those operads are constructed and Theorems \ref{thm:brinfty} and \ref{thm:ksinfty} are shown. In particular, in section \ref{sec:ksproof} we recall the relevant constructions of Kontsevich and Soibelman from \cite{KS1,KS2} and fix an oversight in \cite{KS2}. 
Section \ref{sec:twistedops} sketches the globalization of our results from $\R^n$ to an arbitrary smooth manifold.
In section \ref{sec:9} we show how various results in the literature about the homotopy properties of M. Kontsevich's formality morphisms follow immediately from the existence of the $\KS_\infty$ formality morphism constructed in this paper.
Finally in section \ref{sec:10} we show Theorems \ref{thm:kontstamequiv} and \ref{thm:kontsextends}.
The appendices contain several technical constructions that are used in the main text.

\subsection*{Acknowledgements}
I am grateful for helpful discussions with Johan Alm, Alberto Cattaneo, Giovanni Felder, Sergei Merkulov and others. The original motivation for this work came from Sergei Merkulov's article \cite{merkulov_exotic}. Most of this work was written while the author was a Junior Fellow of the Harvard Society of Fellows. I am very grateful for their support.
Furthermore, the author was partially supported by the Swiss National Science Foundation (grants 200020\_105450 and 200021\_150012).

\subsection*{Notation}
We work over the ground field $\R$ unless otherwise stated. If $V$ is a graded or differential graded vector space, then we denote by $V[r]$ its $r$-fold desuspension. The phrase ``differential graded'' will be abbreviated as dg. In general we work in cohomological conventions, i.~e., our differentials have degree $+1$.

We will use the language of operads and colored operads. A good introduction can be found in the textbook \cite{lodayval}, whose conventions we mostly follow. In particular, if $\op P$ is an operad we denote by $\op P\{r\}$ its $r$-fold operadic desuspension. If $\op P$ is a quadratic operad, we denote by $\op P^\vee$ its Koszul dual cooperad, see \cite[7.2]{lodayval}. We will denote by $\Omega(\op C)$ the operadic bar construction of a coaugmented cooperad $\op C$, cf. \cite[7.3.3]{lodayval}. We denote by $\Lie$ the Lie operad and by $\Lie_\infty=\Omega(\Lie^\vee)$ its minimal cofibrant resolution, where $\Lie^\vee$ is the Koszul dual cooperad of $\Lie$. 

A Gerstenhaber algebra is a dg vector space $V$ together with a commutative algebra structure $\wedge$ on $V$ and a Lie algebra structure $[\ ,\ ]$ on $V[1]$, that are compatible in the sense that $[x,y\wedge z]=[x,y]\wedge z\pm y\wedge [x,z]$ for all $x,y,z\in V$.
The operad governing Gerstenhaber algebras will be denoted by $\Ger$ and its minimal resolution by $\Ger_\infty=\Omega(\Ger^\vee)$. We denote by $\Lie^{(k)}_\infty=\Lie_\infty\{k\}$ the minimal resolution of the degree shifted Lie operad. Homotopy morphisms between $\Lie_\infty$ or $\Ger_\infty$ algebras will be called $\Lie_\infty$- and $\Ger_\infty$-morphisms. (In particular, a $\Lie_\infty$-morphism is not required to be a strict morphism of $\Lie_\infty$ algebras.)

There are several versions of the space of multivector fields. The algebraic multivector fields
\[
 \Tpoly^{alg}=\R[x_1,\dots, x_n, \xi_1,\dots,\xi_n]
\]
are polynomials in variables $x_j$ of degree 0 and $\xi_j$ of degree 1, corresponding to $\frac \p {\p x_j}$. The formal multivector fields are jets of multivector fields at the origin
\[
 \Tpoly^{formal}=\R[[x_1,\dots, x_n]][\xi_1, \dots,\xi_n].
\]
Finally the smooth multivector fields are smooth sections
\[
 \Tpoly^{smooth}=C^\infty(\R^n; \wedge T\R^n).
\]
All three versions carry a Gerstenhaber algebra structure, with the obvious product and the Schouten-Nijenhuis Lie bracket, defined such that $[\xi_i,x_j]=\delta_{ij}$.
Similarly, one may define three versions $\Dpoly^{alg}$, $\Dpoly^{formal}$, $\Dpoly^{smooth}$ of the space of multidifferential operators.
They are equipped with the Hochschild differential and a braces algebra structure.
In all three cases there is a Hochschild-Kostant-Rosenberg Theorem stating that the Hochschild-Kostant-Rosenberg map $\Tpoly^{?}\to \Dpoly^{?}$ is a quasi-isomorphism of complexes.

Unless otherwise stated, we will denote by $\Tpoly=\Tpoly^{smooth}$, $\Dpoly=\Dpoly^{smooth}$ the smooth versions for concreteness, however with the implicit understanding that our construction of a $\Br_\infty$ formality morphism and Theorems \ref{thm:kontsextends}, \ref{thm:brinfty} and \ref{thm:ksinfty} go through without any change in all three cases.

\makeatletter{}\section{Our conventions about (colored) operads}
\label{sec:opconventions}
Recall from the introduction the (loose) definition of the 4-colored operad $\homKS$. It governs a $\Br$-algebra, a $\Br_\infty$-algebra, a map between the two, a module (in some sense) for each of the algebras, and a map between those modules. Our eventual goal in this paper is to construct a representation of $\homKS$ on the 4-colored dg vector space $\Dpoly\oplus  C_\bullet(A)\oplus \Tpoly \oplus \Omega_\bullet$.
We will do that by considering the different components (algebra structures, maps, module structures, module maps) separately. Since there are many different color combinations in a 4 colored operad, this may lead to ugly notations. Hence, instead of working with the full colored operad, we will work with smaller sub-structures, which together generate the colored operad, and consider one of those at a time.

Concretely, the operad $\homKS$ is generated by the following parts.
\begin{enumerate}
 \item The braces operad $\Br$, colored with the color corresponding to $\Dpoly$.
 \item Its cofibrant resolution $\Br_\infty$, colored with the color corresponding to $\Tpoly$.
 \item A $\Br$-$\Br_\infty$ operadic bimodule $\hBr_\infty$, governing homotopy maps from a $\Br_\infty$ algebra to a $\Br$ algebra. Elements of $\hBr_\infty$ have all of their inputs colored in the ``$\Tpoly$-color'', and their output in the ``$\Dpoly$-color''.
 \item The component $\KS_1$ governs the module structure of $\Dpoly$ on $C_\bullet(A)$. Its elements have exactly one input and the output colored with the ``$C_{\bullet}(A)$-color'' and zero or more inputs colored with the ``$\Dpoly$-color''.
 \item Similarly, $\KS_{1,\infty}$ governs the (homotopy) action of $T_{poly}$ on $\Omega_\bullet$. Its elements have exactly one input and the output colored with the $\Omega_\bullet$-color, and zero or more inputs colored with the ``$\Tpoly$-color''.
 \item The component $\hKS_{1,\infty}$ governs a homotopy map of modules from $C_\bullet(A)$ to $\Omega_\bullet$. Its elements have exactly one $C_\bullet(A)$-colored input, a $\Omega_\bullet$-colored output, and zero or more $\Tpoly$-colored inputs.
\end{enumerate}
It is easy to check that these 6 components generate $\homKS$.
To our knowledge the algebraic structures on the last three components have no names. Hence we have made up the following.
The pieces of the operad that govern module structures, we call ``moperads''. The piece that governs the map between the modules, we call ``moperadic bimodule''. 
Before introducing the definitions in more detail, let us recall the notion of (colored) operad we use, cf. \cite{leinsterbook}.\footnote{In loc. cit., the name ``symmetric multicategory'' is used instead of ``colored operad''.}

\begin{defi}
Let $C$ be a finite set (of ``colors''). A $C$-colored operad $\op P$ is the following data:
\begin{enumerate}
 \item For each color $c\in C$ and each tuple $(c_1,\dots, c_n)\in C\times\cdots \times C$ a vector space $\op P^c(c_1,\dots, c_n)$.
 \item For each symmetric group element $\sigma \in S_n$ and for all colors $c,c_1,\dots, c_n\in C$ a morphism 
\begin{equation}\label{equ:sigmaacts}
\op P^c(c_1,\dots, c_n) \to \op P^c(c_{\sigma(1)},\dots, c_{\sigma(n)}).
\end{equation}
\item For each color $c\in C$ a unit element $\mathit{id}_c\in \op P^c(c)$.
\item Composition morphisms
\[
\mu_{c_1,\dots, c_n}^c:
 \op P^c(c_1,\dots c_n) \otimes \op P^{c_1}(c_{1,1},\dots,c_{1,m_1}) 
\otimes \cdots \otimes \op P^{c_n}(c_{n,1},\dots,c_{n,m_n})
\to 
\op P^c(c_{1,1},\dots,c_{1,m_1},\dots,c_{n,m_n} ). 
\]
\end{enumerate}
These data have to satisfy the following conditions:
\begin{enumerate}
\item The actions of the permutation group elements $\sigma\in S_n$ assemble to a (right) representation of $S_n$.
\item (Equivariance) The composition morphisms are equivariant under the right $S_n$ action. 
\item (Unit axiom) For all colors $c,c_1,\dots, c_n\in C$ and each $x\in \op P^c(c_1,\dots, c_n)$ we have
 \begin{align*}
  \mu^c_c(\mathit{id}_c, x) = \mu^c_{c_1,\dots, c_n}(x,\mathit{id}_{c_1},\dots,\mathit{id}_{c_n})=x.
 \end{align*}
\item (Associativity) The composition is associative in the natural way.\footnote{We refer to \cite[section 2.1]{leinsterbook} for more details.}
\end{enumerate}
An operad is a $C$-colored operad for a one-element set $C$.
\end{defi}
To simplify the notation, we will often use a finite set of natural numbers $C=\{1,\dots,f\}$ as our set of colors. Then we use the shorthand notation
\[
 \op P^j(n_1,\dots, n_f) := \op P^j(\underbrace{1,\dots, 1}_{n_1\times},\underbrace{2,\dots, 2}_{n_2\times}, \dots, \underbrace{f,\dots, f}_{n_f\times})
\]
for the operations with $n_i$ inputs of color $i$ and the output of color $j$.
If the set of colors $C$ contains only a single color, we will omit the designation of that color and just write $\op P(n)$ for the space of $n$-ary operations, and $\mathit{id}$ for the operadic unit.

For later use, we will also need a slightly weaker modification of the notion of colored operad.
\begin{defi}
 Let $C$ be a finite set and $d\in C$ be a fixed element. A $C$-colored operad \emph{nonsymmetric in color $d$} is the same data as a $C$-colored operad, except that it is only required to provide the action \eqref{equ:sigmaacts} of symmetric group elements $\sigma \in S_n$ on
\[
\op P^c(c_1,\dots, c_n) 
\]
that leave invariant the relative order of the colors $c_j$ which are equal to the fixed color $d$. The axioms to be satisfied by the composition operation are weakened by requiring equivariance only with respect to symmetric group elements for which the action is defined.
\end{defi}

As mentioned above, we now want to give names to commonly encountered pieces of colored operads, so as to simplify the discussion in the following sections.

\begin{defi}
Let $\op P$ be an operad. A $\op P$-\emph{module operad} (or short: $\op P$-\emph{moperad}) $\op P_1$ is the following data:
\begin{enumerate}
 \item A collection of right $S_k$ modules $\op P_1(k)$, $k=0,1,2,\dots$. Here $\op P_1(k)$ will be thought of as a space of operations with $k$ inputs in one color, and one input and the output in another.
 \item A unit element $\mathit{id}_1\in \op P_1(0)$.
 \item Composition morphisms
\begin{equation}\label{equ:moperadiccomp}
 \mu_{1,k} \colon \op P_1(k) \otimes \op P_1(m_0) \otimes \op P(m_1)\otimes \cdots \otimes  \op P(m_k) \to \op P_1(m_0+\cdots +m_k) 
\end{equation}
for each $k=0,1,\dots$ and any $m_1,\dots,m_k\in \Z_{\geq 0}$.
\end{enumerate}
These data are required to satisfy the following axioms
\begin{enumerate}
 \item (Equivariance) The composition is equivariant under the symmetric group action. 
 \item (Unit axiom) 
\[
\mu_{1,0}(\mathit{id}_1, a) = \mu_{1,k}(a, \mathit{id}_1, \mathit{id}, \dots, \mathit{id}) 
=
a 
\]
for each $k=0,1,\dots$ and any $a\in \op P_1(k)$. Here $\mathit{id}\in \op P(1)$ is the unit of $\op P$.
 \item (Associativity) 
\begin{multline*}
 \mu_{1,k}(a, \mu_{1,m}(b, c, x_{0,1},\dots, x_{0,m}), \mu_{m_1}(x_1, x_{1,1}, \dots, x_{1,m_1}), \dots ,  \mu_{m_k}(x_k, x_{k,1}, \dots, x_{k,m_k}))
\\=
\pm 
\mu_{1,m+m_1+\dots+m_k}(\mu_{1,k} (a,b, x_1, \dots, x_k),
 c, x_{0,1},\dots, x_{0,m}, x_{1,1}, \dots, x_{1,m_1},  \\ \dots ,  x_{k,1}, \dots, x_{k,m_k})
\end{multline*}

for all $k=0,1,\dots$, $m,m_1,\dots, m_k\in \Z_{\geq 0}$, $a\in \op P_1(k)$, $b\in \op P_1(m)$, $c\in \op P_1$ and $x_j\in \op P(m_j)$, $x_{i,j}\in \op P$. The sign is that of the permutation bringing the odd elements of $\op P$ occuring on the right hand side into the ordering on the left hand side.
\end{enumerate}

\end{defi}

\begin{rem}
 The equivariance under the symmetric group actions concretely means two statements:
(i) The composition \eqref{equ:moperadiccomp} commutes with the natural right $S_{m_0}\times \cdots \times S_{m_k}$ action on both sides.
(ii) The composition commutes with the $S_k$ action on both sides, where a permutation acts on the left-hand side of \eqref{equ:moperadiccomp} by the given action on $\op P_1(k)$, and by permutation of the other tensor factors, and acts on the right-hand side of \eqref{equ:moperadiccomp} by permuting $k$ ``blocks'' of size $m_1,\dots, m_k$.
\end{rem}

The definition is such that $\op P$ and $\op P_1$ can be assembled into a two-colored operad, with $\op P(k)$ being the space of operations with $k$ inputs and the output of color 1, and $\op P_1(k)$ being the space of operations with the first input and the output of color 2, and the last $k$ inputs of color one. This colored operad we will denote by
\[
 \bpm \op P & \op P_1 \epm.
\]

The notion of operadic bimodule has been introduced in \cite{KapranovManin}. Let us recall it here for completeness.
\begin{defi}
 Let $\op P$, $\op Q$ be two operads. An operadic $\op P$-$\op Q$ bimodule $\op M$ is the following data:
\begin{enumerate}
 \item A collection of right $S_k$ modules $\op M(k)$ for $k=0,1,\dots$. Intuitively, $\op M(k)$ is thought of a space of operations with $k$ inputs of one sort and the output of another.
 \item Left and right composition morphisms
\begin{align*}
 \nu_{k}^{(l)} &\colon \op P(k)\otimes \op M(m_1)\cdots \otimes \op M(m_k) \to \op M(m_1+\dots +m_k) \\
 \nu_{k}^{(r)} &\colon \op M(k)\otimes \op Q(m_1)\cdots \otimes \op Q(m_k) \to \op M(m_1+\dots +m_k).
\end{align*}
\end{enumerate}
These data are required to satisfy the following axioms.
\begin{enumerate}
 \item (Equivariance) The compositions are equivariant under the symmetric group action.
 \item (Unit axiom) 
\[
\nu_{1}^{(l)}(\mathit{id}_{\op P}, a) 
= \nu_{k}^{(r)}(a, \mathit{id}_{\op Q}, \dots, \mathit{id}_{\op Q}) 
=
a 
\]
for each $k=0,1,\dots$ and any $a\in \op M(k)$. Here $\mathit{id}_{\op P}, \mathit{id}_{\op Q}$ are the units in $\op P, \op Q$.
 \item (Associativity)
\begin{gather*}
 \nu_k^{(l)}(x,  \nu_{m_1}^{(l)} (x_1, a_{1,1}, \dots, a_{1,m_1}), \dots , \nu_{m_k}^{(l)} (x_k, a_{k,1}, \dots, a_{k,m_k}))
=
 \pm \nu_{\sum m_j }^{(l)}(\mu_k(x, x_1,\dots, x_k), a_{1,1}, \dots, a_{k,m_k})
\\
 \nu_{\sum m_j}^{(r)}(\nu_{k}^{(r)} (a, y_1, \dots, y_{k}), y_{1,1}, \dots ,y_{k, m_k})
=
 \pm \nu_{k}^{(r)}(a, \mu_{m_1}(y_1, y_{1,1},\dots, y_{1,m_k}), \dots, \mu_{m_k}(y_k, y_{k,1},\dots ,y_{k,m_k}))
 \\
 \nu_{k}^{(l)}(x, \nu_{m_1}^{(r)} (a_1,  y_{1,1},\dots, y_{1,m_1}), \dots, \nu_{m_k}^{(r)} (a_k,  y_{k,1},\dots, y_{1,m_k}))
=
\pm 
\nu_{\sum m_j }^{(r)}(\nu_k^{(l)}(x, a_1,\dots, a_k), y_{1,1}, \dots, y_{k,m_k})
\end{gather*}
for all $k=0,1,\dots$, $m_1,\dots, m_k\in \Z_{\geq 0}$, $a\in \op M(k)$, $a_j\in \op M(m_j)$, $a_{i,j}\in \op M$, $x\in \op P(k)$, $x_j\in \op P(m_j)$, $y_j\in \op Q(m_j)$, $y_{i,j}\in \op Q$.
Here we abuse the symbol $\mu$ to denote both the composition in $\op P$ and that in $\op Q$. The signs are defined similarly to the ones in the last definition.
\end{enumerate}

\end{defi}

For such $\op P$, $\op Q$, $\op M$ there is a two-colored operad, generated by $\op P$, $\op Q$, $\op M$, such that 
the operations of color 1 are given by $\op P$, the operations of color 2 are given by $\op Q$, the operations with $k$ inputs of color 2 and one output of color 1 are given by $\op M(k)$, and the compositions are defined using the compositions of $\op P$, $\op Q$ and the actions $\nu_k^{(l)}$, $\nu_k^{(r)}$ above. This operad we denote by
\[
 \bpm \op P & \op M & \op Q \epm.
\]

Finally let us give the definition of a moperadic bimodule, which governs ``maps between modules''.
\begin{defi}
 Let $\op P$, $\op Q$ be two operads. Let $\op M$ be a $\op P$-$\op Q$ operadic bimodule. Let $\op P_1$ be a $\op P$-moperad and $\op Q_1$ be a $\op Q$-moperad. A $\op P$-$\op P_1$-$\op M$-$\op Q$-$\op Q_1$ moperadic bimodule $\op M_1$ is the following data:
\begin{enumerate}
 \item A collection of $S_k$ modules $\op M_1(k)$ for $k=0,1,\dots$. Intuitively, $\op M_1(k)$ is thought of a space of operations with $k$ inputs of $\op Q$-color, one input of $\op P_1$-color and the output of $\op Q_1$-color.
 \item Left and right composition morphisms
\begin{gather*}
 \nu_{1,k}^{(l)} \colon \op Q_1(k)\otimes \op M_1(m_0)\otimes \op Q(m_1)\otimes \cdots \otimes \op Q(m_k) \to \op M_1(m_0+\dots +m_k) \\
 \nu_{1,k,l}^{(r)} \colon \op M_1(k)\otimes \op P_1(l)\otimes \op M(n_1) \otimes \cdots \otimes \op M(n_l)
\otimes \op Q(m_1)\otimes \cdots \otimes \op Q(m_k)
\\ 
\to \op M_1(n_1+\dots + n_l + m_1 + \dots + m_k).
\end{gather*}
Pictures of these compositions can be found in Figure \ref{fig:moperadicbimodcomp}.
\end{enumerate}
\begin{figure}
 \centering 
\begin{tabular}{cc}
\noindent\parbox[c]{.5\hsize}{
\makeatletter{}\usetikzlibrary{matrix}
\usetikzlibrary{arrows}
\usetikzlibrary{shapes}
\usetikzlibrary{through}
\usetikzlibrary{calc,3d}
\usetikzlibrary{decorations,decorations.pathmorphing}

\begin{tikzpicture}[
int/.style={circle, draw, fill, minimum size=5pt, inner sep=0},
ext/.style={circle, draw, fill=white, minimum size=5pt, inner sep=1pt},
helper/.style={coordinate},point/.style={circle, draw, fill, inner sep =1pt},
de/.style={-triangle 60},
point/.style={circle, draw, fill, minimum size=3pt, inner sep=0pt},
xst/.style={cross out, draw, minimum size=5 },
abox/.style={draw, minimum width=30, minimum height=10}
]

\node[abox] (v1) at (0,0) {};
\draw[col4] (v1.150) -- +(0,.5);
\node[abox, fill=gray] (v2) at (0,-1) {};
\node[ext, fill=gray!50] (e2) at (2.1,-1) {};
\node[ext, fill=gray!50] (e3) at (2.9,-1) {};
\draw (v1.-30) edge[col3] (e2) edge[col3] (e3);
\draw (v1.-150) edge[col2] (v2.150);
\draw[col2] (v2.-150) --+ (0,-1);

\foreach \x in {0,.5}
{
   \draw[col1] (v2.-30) -- +(\x,-.5);
   \draw (v2.-30) +(\x,-.5) node[ext] (v) {};
   \foreach \x in {-.1,0,.1}
     \draw[col3] (v)--+(\x,-.5);
}
\foreach \n in {(e2),(e3)}
  \foreach \x in {-.3,0,.3}
    \draw[col3] \n -- +(\x,-.5);

\node at (4,-.5) {\huge $\to$};
\node[abox] (v3) at (5.5,-.5) {};
\draw[col4] (v3.150) -- +(0,.5);
\draw[col2] (v3.-150) --+ (0,-.5);
\foreach \x in {-.4,-.3,-.2,-.1,0,.1,.2,.3,.4,.5,.6,.7}
   \draw[col3] (v3.-30) -- +(\x,-.5);

\begin{scope}[yshift=-3cm]

\node[abox, fill=gray!50] (v1) at (0,0) {};
\draw[col4] (v1.150) -- +(0,.5);
\node[abox] (v2) at (0,-1) {};
\node[ext, fill=gray!50] (e1) at (1.3,-1) {};
\node[ext, fill=gray!50] (e2) at (2.1,-1) {};
\node[ext, fill=gray!50] (e3) at (2.9,-1) {};
\draw (v1.-30) edge[col3] (e1) edge[col3] (e2) edge[col3] (e3);
\draw (v1.-150) edge[col4] (v2.150);
\draw[col4] (v2.-150) --+ (0,-.5);

\foreach \x in {-.3,0,.3}
   \draw[col3] (v2.-30) -- +(\x,-.5);
\foreach \n in {(e1),(e2),(e3)}
  \foreach \x in {-.3,0,.3}
    \draw[col3] \n -- +(\x,-.5);

\node at (4,-.5) {\huge $\to$};
\node[abox] (v3) at (5.5,-.5) {};
\draw[col4] (v3.150) -- +(0,.5);
\draw[col2] (v3.-150) --+ (0,-.5);
\foreach \x in {-.4,-.3,-.2,-.1,0,.1,.2,.3,.4,.5,.6,.7}
   \draw[col3] (v3.-30) -- +(\x,-.5);

\end{scope}

\end{tikzpicture}
 
}
&
\collegend \\
\end{tabular}

\caption{\label{fig:moperadicbimodcomp} Pictures of the right and left action on a moperadic bimodule $\op M_1$. The white boxes stand for elements of $\op M_1$. The dark gray box on the left stands for an element of $\op P_1$. The white circles stand for elements of $\op M$, and the light gray circles for elements of $\op Q$. Finally the light gray box symbolizes an element of $\op Q_1$. }
\end{figure}
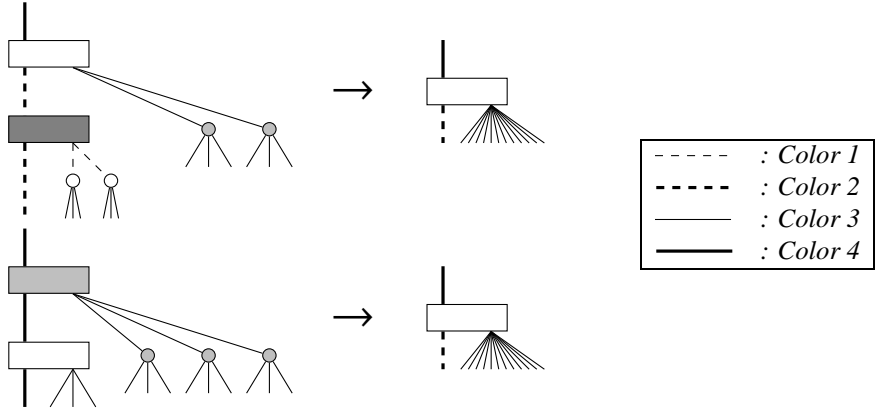

These data are required to satisfy the following axioms.
\begin{enumerate}
 \item (Equivariance) The compositions are equivariant under the symmetric group action.
 \item (Unit axiom) 
\[
\nu_{1,0}^{(l)}(\mathit{id}_{\op Q_1}, a) 
= \nu_{1,k,0}^{(r)}(a, \mathit{id}_{\op P_1}, \mathit{id}_{\op Q}, \dots, \mathit{id}_{\op Q}) 
=
a 
\]
for each $k=0,1,\dots$ and any $a\in \op M_1(k)$. Here $\mathit{id}_{\op P_1}, \mathit{id}_{\op Q}, \mathit{id}_{\op Q_1}$ are the units in $\op P, \op Q, \op Q_1$.
 \item (Associativity) There are three big relations saying that the left action is an action, that the right action is an action and that both actions commute. For brevity, we write them down in a graphical way, see Figure \ref{fig:mopbimodrelations}. 
\end{enumerate}

\end{defi}
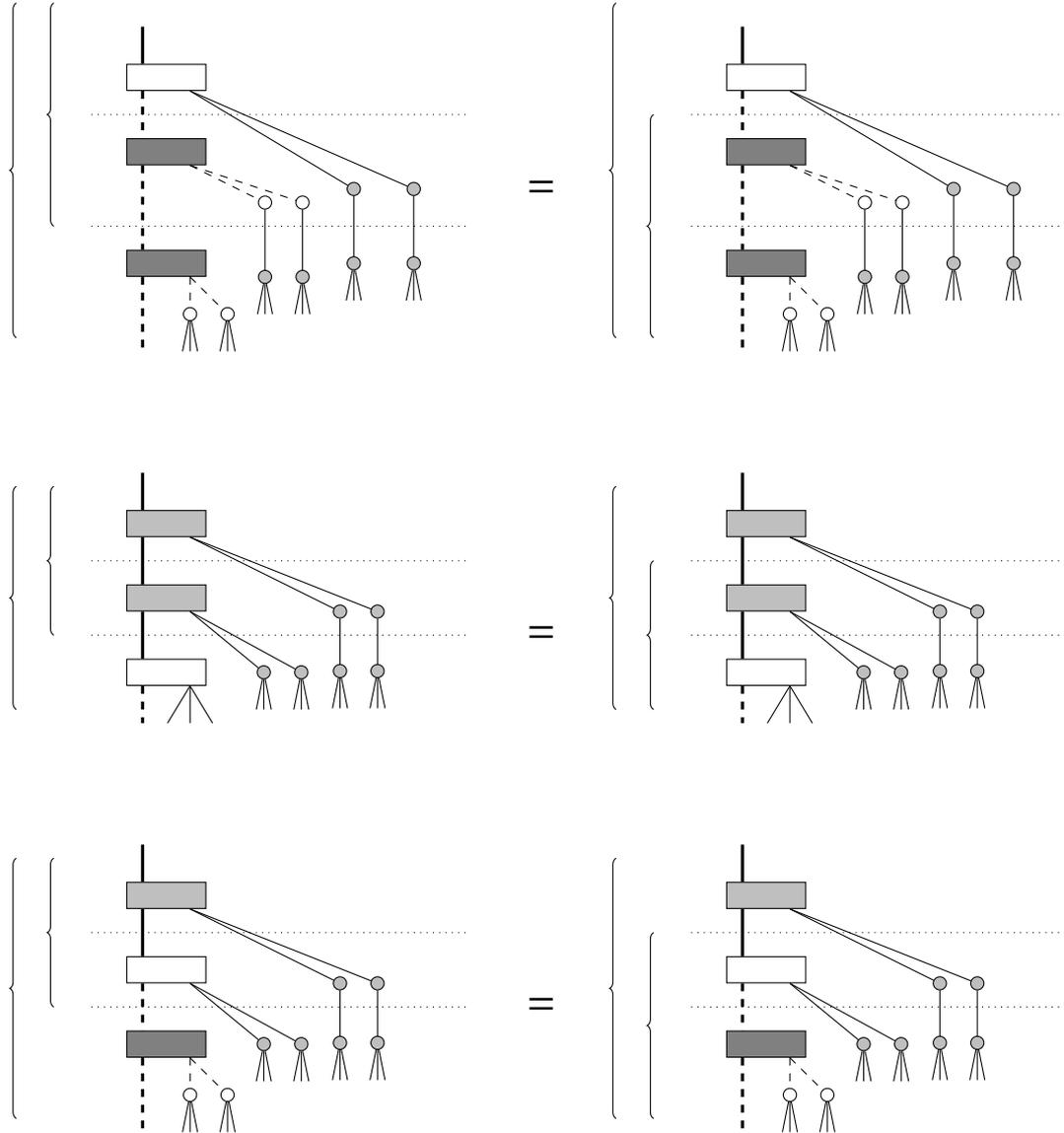
\begin{figure}
 \centering
\makeatletter{}\usetikzlibrary{matrix}
\usetikzlibrary{arrows}
\usetikzlibrary{shapes}
\usetikzlibrary{through}
\usetikzlibrary{calc,3d}
\usetikzlibrary{decorations,decorations.pathmorphing}
\usetikzlibrary{decorations.pathreplacing}

\[
\begin{tikzpicture}[
int/.style={circle, draw, fill, minimum size=5pt, inner sep=0},
ext/.style={circle, draw, fill=white, minimum size=5pt, inner sep=1pt},
helper/.style={coordinate},point/.style={circle, draw, fill, inner sep =1pt},
de/.style={-triangle 60},
point/.style={circle, draw, fill, minimum size=3pt, inner sep=0pt},
xst/.style={cross out, draw, minimum size=5 },
abox/.style={draw, minimum width=30, minimum height=10},
separatorline/.style={dotted}
]

\node[abox] (v1) at (0,0) {};
\draw[col4] (v1.150) -- +(0,.5);
\node[abox, fill=gray] (v2) at (0,-1) {};
\node[abox, fill=gray] (v3) at (0,-2.5) {};

\node[ext, fill=gray!50] (e2) at (2.5,-1.5) {};
\node[ext, fill=gray!50] (e3) at (3.3,-1.5) {};

\draw (v1.-30) edge[col3] (e2) edge[col3] (e3);
\draw (v1.-150) edge[col2] (v2.150);
\draw[col2] (v2.-150) --(v3.150);
\draw[col2] (v3.-150) --+ (0,-1);

\foreach \x in {1,1.5}
{
   \draw[col1] (v2.-30) -- +(\x,-.5);
   \draw (v2.-30) +(\x,-.5) node[ext] (v) {};
   \draw[col3] (v) -- +(0,-1) ;
   \draw (v) +(0,-1) node[ext, fill=gray!50] (vv) {};
   \foreach \x in {-.1,0,.1}
     \draw[col3] (vv)--+(\x,-.5);
}

\foreach \x in {0,.5}
{
   \draw[col1] (v3.-30) -- +(\x,-.5);
   \draw (v3.-30) +(\x,-.5) node[ext] (v) {};
   \foreach \x in {-.1,0,.1}
     \draw[col3] (v)--+(\x,-.5);
}

\foreach \n in {(e2),(e3)}
{
\draw[col3] \n -- +(0,-1);
\draw  \n +(0,-1) node[ext, fill=gray!50] (v) {};
  \foreach \x in {-.1,0,.1}
    \draw[col3] (v) -- +(\x,-.5);
}
\node at (5,-1.5) {\huge $=$};

\draw[separatorline] (-1,-.5)--+(5,0);
\draw[separatorline] (-1,-2)--+(5,0);
\draw[decorate, decoration=brace] (-1.5,-2)--+(0,3);
\draw[decorate, decoration=brace] (-2,-3.5)--+(0,4.5);

\begin{scope}[xshift=8cm]
\node[abox] (v1) at (0,0) {};
\draw[col4] (v1.150) -- +(0,.5);
\node[abox, fill=gray] (v2) at (0,-1) {};
\node[abox, fill=gray] (v3) at (0,-2.5) {};

\node[ext, fill=gray!50] (e2) at (2.5,-1.5) {};
\node[ext, fill=gray!50] (e3) at (3.3,-1.5) {};

\draw (v1.-30) edge[col3] (e2) edge[col3] (e3);
\draw (v1.-150) edge[col2] (v2.150);
\draw[col2] (v2.-150) --(v3.150);
\draw[col2] (v3.-150) --+ (0,-1);

\foreach \x in {1,1.5}
{
   \draw[col1] (v2.-30) -- +(\x,-.5);
   \draw (v2.-30) +(\x,-.5) node[ext] (v) {};
   \draw[col3] (v) -- +(0,-1);
   \draw (v) +(0,-1) node[ext, fill=gray!50] (vv) {};
   \foreach \x in {-.1,0,.1}
     \draw[col3] (vv)--+(\x,-.5);
}

\foreach \x in {0,.5}
{
   \draw[col1] (v3.-30) -- +(\x,-.5);
   \draw (v3.-30) +(\x,-.5) node[ext] (v) {};
   \foreach \x in {-.1,0,.1}
     \draw[col3] (v)--+(\x,-.5);
}

\foreach \n in {(e2),(e3)}
{
\draw[col3] \n -- +(0,-1);
\draw \n +(0,-1) node[ext, fill=gray!50] (v) {};
  \foreach \x in {-.1,0,.1}
    \draw[col3] (v) -- +(\x,-.5);
}

\draw[separatorline] (-1,-.5)--+(5,0);
\draw[separatorline] (-1,-2)--+(5,0);
\draw[decorate, decoration=brace] (-1.5,-3.5)--+(0,3);
\draw[decorate, decoration=brace] (-2,-3.5)--+(0,4.5);

\end{scope}

\begin{scope}[yshift=-7cm]

\node[abox, fill=gray!50] (v1) at (0,0) {};
\node[abox] (v2) at (0,-1) {};
\node[abox, fill=gray!50] (v0) at (0,1) {};
\draw[col4] (v0.150) -- +(0,.5);
\draw[col4] (v0.-150) edge (v1.150);
\node[ext, fill=gray!50] (e1) at (1.3,-1) {};
\node[ext, fill=gray!50] (e2) at (1.8,-1) {};

\draw (v1.-30) edge[col3] (e1) edge[col3] (e2);
\draw (v1.-150) edge[col4] (v2.150);
\draw[col2] (v2.-150) --+ (0,-.5);

\foreach \x in {2,2.5}
{
\draw[col3] (v0.-30)--+(\x, -1);
\draw (v0.-30) +(\x, -1) node[ext, fill=gray!50] (v) {};
\draw[col3] (v) --+(0,-.8);
\draw (v) +(0,-.8) node[ext, fill=gray!50] (vv) {};
   \foreach \x in {-.1,0,.1}
     \draw[col3] (vv)--+(\x,-.5);
}

\foreach \x in {-.3,0,.3}
   \draw[col3] (v2.-30) -- +(\x,-.5);
\foreach \n in {(e1),(e2)}
  \foreach \x in {-.1,0,.1}
    \draw[col3] \n -- +(\x,-.5);

\node at (5,-.5) {\huge $=$};

\draw[separatorline] (-1,.5)--+(5,0);
\draw[separatorline] (-1,-.5)--+(5,0);

\draw[decorate, decoration=brace] (-1.5,-.5)--+(0,2);
\draw[decorate, decoration=brace] (-2,-1.5)--+(0,3);

\end{scope}

\begin{scope}[yshift=-7cm, xshift=8cm]

\node[abox, fill=gray!50] (v1) at (0,0) {};
\node[abox] (v2) at (0,-1) {};
\node[abox, fill=gray!50] (v0) at (0,1) {};
\draw[col4] (v0.150) -- +(0,.5);
\draw[col4] (v0.-150) edge (v1.150);
\node[ext, fill=gray!50] (e1) at (1.3,-1) {};
\node[ext, fill=gray!50] (e2) at (1.8,-1) {};

\draw (v1.-30) edge[col3] (e1) edge[col3] (e2);
\draw (v1.-150) edge[col4] (v2.150);
\draw[col2] (v2.-150) --+ (0,-.5);

\foreach \x in {2,2.5}
{
\draw[col3] (v0.-30)--+(\x, -1);
\draw (v0.-30) +(\x, -1) node[ext, fill=gray!50] (v) {};
\draw[col3] (v) --+(0,-.8);
\draw (v) +(0,-.8) node[ext, fill=gray!50] (vv) {};
   \foreach \x in {-.1,0,.1}
     \draw[col3] (vv)--+(\x,-.5);
}

\foreach \x in {-.3,0,.3}
   \draw[col3] (v2.-30) -- +(\x,-.5);
\foreach \n in {(e1),(e2)}
  \foreach \x in {-.1,0,.1}
    \draw[col3] \n -- +(\x,-.5);

\draw[separatorline] (-1,.5)--+(5,0);
\draw[separatorline] (-1,-.5)--+(5,0);

 \draw[decorate, decoration=brace] (-1.5,-1.5)--+(0,2);
 \draw[decorate, decoration=brace] (-2,-1.5)--+(0,3);

\end{scope}

\begin{scope}[yshift=-12cm]

\node[abox] (v1) at (0,0) {};
\node[abox, fill=gray] (v2) at (0,-1) {};
\node[abox, fill=gray!50] (v0) at (0,1) {};
\draw[col4] (v0.150) -- +(0,.5);
\draw (v0.-150) edge[col4] (v1.150);
\node[ext, fill=gray!50] (e1) at (1.3,-1) {};
\node[ext, fill=gray!50] (e2) at (1.8,-1) {};

\draw (v1.-30) edge[col3] (e1) edge[col3] (e2);
\draw (v1.-150) edge[col2] (v2.150);
\draw[col2] (v2.-150) --+ (0,-1);

\foreach \x in {2,2.5}
{
\draw[col3] (v0.-30)--+(\x, -1);
\draw (v0.-30) +(\x, -1) node[ext, fill=gray!50] (v) {};
\draw[col3] (v) --+(0,-.8);
\draw (v) +(0,-.8) node[ext, fill=gray!50] (vv) {};
   \foreach \x in {-.1,0,.1}
     \draw[col3] (vv)--+(\x,-.5);
}

\foreach \x in {0,.5}
{
   \draw[col1] (v2.-30) -- +(\x,-.5);
   \draw (v2.-30) +(\x,-.5) node[ext] (v) {};
   \foreach \x in {-.1,0,.1}
     \draw[col3] (v)--+(\x,-.5);
}

\foreach \n in {(e1),(e2)}
  \foreach \x in {-.1,0,.1}
    \draw[col3] \n -- +(\x,-.5);

\node at (5,-.5) {\huge $=$};

\draw[separatorline] (-1,.5)--+(5,0);
\draw[separatorline] (-1,-.5)--+(5,0);

\draw[decorate, decoration=brace] (-1.5,-.5)--+(0,2);
\draw[decorate, decoration=brace] (-2,-2)--+(0,3.5);

\end{scope}

\begin{scope}[yshift=-12cm, xshift=8cm]

\node[abox] (v1) at (0,0) {};
\node[abox, fill=gray] (v2) at (0,-1) {};
\node[abox, fill=gray!50] (v0) at (0,1) {};
\draw[col4] (v0.150) -- +(0,.5);
\draw (v0.-150) edge[col4] (v1.150);
\node[ext, fill=gray!50] (e1) at (1.3,-1) {};
\node[ext, fill=gray!50] (e2) at (1.8,-1) {};

\draw (v1.-30) edge[col3] (e1) edge[col3] (e2);
\draw (v1.-150) edge[col2] (v2.150);
\draw[col2] (v2.-150) --+ (0,-1);

\foreach \x in {2,2.5}
{
\draw[col3] (v0.-30)--+(\x, -1);
\draw (v0.-30) +(\x, -1) node[ext, fill=gray!50] (v) {};
\draw[col3] (v) --+(0,-.8);
\draw (v) +(0,-.8) node[ext, fill=gray!50] (vv) {};
   \foreach \x in {-.1,0,.1}
     \draw[col3] (vv)--+(\x,-.5);
}

\foreach \x in {0,.5}
{
   \draw[col1] (v2.-30) -- +(\x,-.5);
   \draw (v2.-30) +(\x,-.5) node[ext] (v) {};
   \foreach \x in {-.1,0,.1}
     \draw[col3] (v)--+(\x,-.5);
}

\foreach \n in {(e1),(e2)}
  \foreach \x in {-.1,0,.1}
    \draw[col3] \n -- +(\x,-.5);

\draw[separatorline] (-1,.5)--+(5,0);
\draw[separatorline] (-1,-.5)--+(5,0);

\draw[decorate, decoration=brace] (-1.5,-2)--+(0,2.5);
\draw[decorate, decoration=brace] (-2,-2)--+(0,3.5);

\end{scope}

\end{tikzpicture}
\] 
\caption{\label{fig:mopbimodrelations} The defining relations for a moperadic bimodule. The color code is as in Figure \ref{fig:moperadicbimodcomp}.
The braces on the left-hand side of drawings shall indicate in which order the corollas are to be composed.
}
\end{figure}

The notion of moperadic bimodule is defined such that $\op P$, $\op Q$, $\op P_1$, $\op Q_1$, $\op M$, $\op M_1$ as above generate a four colored operad, such that 
\begin{enumerate}
 \item $\op P$ lives in color 1, $\op Q$ lives in color 3.
\item $\op P_1(k)$ are the operations with $k$ inputs of color $1$ and one input and the output of color 2, $\op Q_1(k)$ are the operations with $k$ inputs of color $3$ and one input and the output of color 4.
\item $\op M(k)$ are the operations with all $k$ inputs of color $3$ and the output of color $1$.
\item $\op M_1(k)$ are the operations with $k$ inputs of color $3$, one input of color $2$, and the output of color $4$.
\item The compositions between these components agree with the compositions defined above.
\end{enumerate}

We will denote this big colored operad by 
\[
 \bpm \op P & \op M & \op Q \\ \op P_1 & \op M_1 & \op Q_1 \epm.
\]
One can define the notion of \emph{morphism} of moperads, operadic bimodules and moperadic bimodules in a straightforward way. In the following sections we want to construct maps of 4-colored operads of the above form. To give such a map 
\[
 \bpm \op P & \op M & \op Q \\ \op P_1 & \op M_1 & \op Q_1 \epm 
\to
 \bpm \tilde{\op P} & \tilde{\op M} & \tilde{\op Q} \\ \tilde{\op P}_1 & \tilde{\op M}_1 & \tilde{\op Q}_1 \epm
\]
is equivalent to providing the following:
\begin{enumerate}
 \item Operad maps $\op P\to \tilde{\op P}$ and $\op Q\to \tilde{\op Q}$.
 \item Moperad maps $\op P\to \tilde{\op P}_1$ and $\op Q\to \tilde{\op Q}_1$. Here $\tilde{\op P}_1$ is considered a $\op P$-moperad via the map $\op P\to \tilde{\op P}$ and similarly for $\tilde{\op Q}_1$.
 \item A map of operadic bimodules $\op M\to \tilde{\op M}$. Here $\tilde{\op M}$ is considered an $\op P$-$\op Q$ bimodule via the maps $\op P\to \tilde{\op P}$ and $\op Q\to \tilde{\op Q}$.
 \item A map of moperadic bimodules $\op M\to \tilde{\op M}_1$. Here $\tilde{\op M}_1$ is considered a $\op P$-$\op P_1$-$\op M$-$\op Q$-$\op Q_1$ moperadic bimodule via the maps above.
\end{enumerate}

\subsection{Operads of Swiss Cheese type and Extended Swiss Cheese type}

\begin{defi}
\label{def:SCtype}
We say that a two colored operad $\op P$, non-symmetric in color 2, is of \emph{Swiss Cheese type} if all operation with output in color 1 have all its inputs in color 1.
\end{defi}

Denoting the space of operations with output in color $\alpha$ by $\op P^\alpha$ as above, this says that $\op P^1(\cdot, n)=0$ for $n\geq 1$.

\begin{defi}
\label{def:ESCtype}
We say that a three colored operad $\op P$, symmetric in color 2, is of \emph{Extended Swiss Cheese type} if:
\begin{enumerate}
\item All operation with output in color 1 have all its inputs in color 1, i.e., $\op P^1(\cdot, m, n)=0$ if $m+n>0$.
\item All operations with output in color 2 have all its inputs in colors $1$ and $2$, i.e., $\op P^2(\cdot, \cdot, n)=0$ if $n>0$.
\item The operations with output in color 3 have at most one input in color 3. Furthermore those with exactly one input have no inputs of color $2$. In other words $\op P^3(\cdot, m, n)=0$ for $n>1$ and $\op P^3(\cdot, m, 1)=0$ for $m>0$.
\item There are actions of the cyclic groups of order $m$ on the spaces $\op P^3(\cdot, m, 0)$, which are compatible with the operadic compositions.\footnote{This means that the operadic compositions are equivariant in the same sense as the operadic compositions in a (symmetric) operad are $S_m$ equivariant. In other words, one requires only those equivariance relations of a symmetric operad that involve the elements of the cyclic subgroup of $S_m$.}
\end{enumerate}
\end{defi}

Note that in particular $\op P^1$ is an operad and $\op P^3(\cdot,0,1)$ is a $\op P^1$-moperad. Furthermore the operations of output colors 1 and 2 together form an operad of Swiss Cheese type.

Below we will frequently encounter colored operads of Swiss Cheese or Extended Swiss Cheese type. Note also that the definitions are meaningful for operads in any symmetric monoidal category.

\makeatletter{}\section{Several combinatorial operads given by graphs}
\label{sec:grops}
The goal of this section is to define the colored operad $\bigGra$ occurring in the proof of the main theorem, see eqn. \eref{equ:proofchain}. 
Most operads of this section have occurred in the literature in some form or another: The operad $\Gra$ below was introduced in \cite{megrt}. The braces operad goes back to \cite[section 5.2]{GJ}, with the version we are using (sometimes also called the Kontsevich-Soibelman minimal operad) defined in \cite{KS1}. The Kontsevich-Soibelman operad was introduced in \cite{KS2}.
The graphical operad $\bigGra$ combines Feynman diagrammatic objects and algebraic operations thereon that have appeared in some (less general) variant in the proofs of the Kontsevich Formality Theorem \cite{K1} and its extension to Hochschild chains by Shoikhet \cite{shoikhet}.

\subsection{\texorpdfstring{ $\Gr$ and $\dGr$}{Gr and dGr} (``Graphs'' and ``directed Graphs'') operad}
\label{sec:GrdGr}
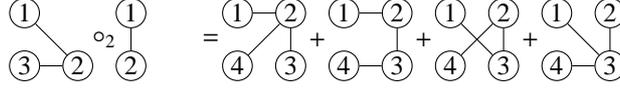
\begin{figure}
\centering
 \begin{tikzpicture}[scale=.7]

\node (e1) at (-1,1) [ext] {$1$};
\node (e2) at (0,0) [ext] {$2$};
\node (e3) at (-1,0) [ext] {$3$};
\draw (e1)--(e2)--(e3);

\node at(.5,.5){$\circ_2$};

\begin{scope}[xshift=1cm]
\node (e1) at (0,1) [ext] {$1$};
\node (e2) at (0,0) [ext] {$2$};
\draw (e1)--(e2);
\end{scope}

\node at(2.5,.5){$=$};

\begin{scope}[xshift=3cm]
\node (e1) at (0,1) [ext] {$1$};
\node (e2) at (1,1) [ext] {$2$};
\node (e3) at (1,0) [ext] {$3$};
\node (e4) at (0,0) [ext] {$4$};
\draw (e2)--(e3) (e1)--(e2)--(e4);
\end{scope}

\node at (4.5,.5) {$+$};

\begin{scope}[xshift=5cm]
\node (e1) at (0,1) [ext] {$1$};
\node (e2) at (1,1) [ext] {$2$};
\node (e3) at (1,0) [ext] {$3$};
\node (e4) at (0,0) [ext] {$4$};
\draw (e2)--(e3) (e1)--(e2) (e3)--(e4);
\end{scope}

\node at (6.5,.5) {$+$};

\begin{scope}[xshift=7cm]
\node (e1) at (0,1) [ext] {$1$};
\node (e2) at (1,1) [ext] {$2$};
\node (e3) at (1,0) [ext] {$3$};
\node (e4) at (0,0) [ext] {$4$};
\draw (e2)--(e3) (e1)--(e3) (e2)--(e4);
\end{scope}

\node at (8.5,.5) {$+$};

\begin{scope}[xshift=9cm]
\node (e1) at (0,1) [ext] {$1$};
\node (e2) at (1,1) [ext] {$2$};
\node (e3) at (1,0) [ext] {$3$};
\node (e4) at (0,0) [ext] {$4$};
\draw (e2)--(e3) (e1)--(e3)--(e4);
\end{scope}
\end{tikzpicture}
\caption{\label{fig:Grcomp}The operadic composition in $\Gr$. Note that the labels on the edges are not shown.}
\end{figure}

Let $\Gr(n)'_k$ be the graded vector space of linear combinations of undirected graphs with vertex set $[n]:=\{1,\dots,n\}$ and edge set $[k]$.\footnote{Concretely, such a graph is an ordered multi-set of $k$ two element subsets of $[n]:=\{1,\dots,n\}$. In particular, note that tadpoles or short cycles, i.e., edges connecting a vertex to itself, are forbidden by this definition.} We consider such a graph as living in degree $-k$, i.e., the edges are considered to be of degree -1. Define 
\[
 \Gr(n) = \prod_{k\geq 0} (\Gr(n)'_k \otimes sgn_k)_{S_k}
\]
where the symmetric group $S_k$ acts on $\Gr(n)'_k$ by permuting the labels on edges and $\sgn_k$ is the sign representation. The spaces $\Gr(n)$ assemble to form an operad $\Gr$. The operadic compositions $\Gamma_1 \circ_j \Gamma_2$ are defined by ``inserting'' graph $\Gamma_2$ into vertex $j$ of graph $\Gamma_1$ and summing over all graphs obtained by reconnecting the edges in $\Gamma_1$ ending at vertex $j$ in all possible ways to vertices of $\Gamma_2$, see Figure \ref{fig:Grcomp}. The labelling on the edges is adjusted such that the edges that came from $\Gamma_2$ have higher labels than those from $\Gamma_1$.

Replacing undirected graphs by directed graphs one can define similarly the operad $\dGr$. There is a map of operads $\Gr\to \dGr$ mapping an undirected graph $\Gamma$ to the sum of all graphs $\Gamma'$ that can be obtained by assigning orientations to the edges.

\begin{rem}
 Note that by the sign convention on edges, some graphs are zero due to odd symmetries. For example, any graph that contains a double edge is zero, because it has an odd symmetry interchanging the two edges forming the double edge.
\end{rem}

\begin{ex}
\label{ex:GraactonTpoly}
 The multivector fields $\Tpoly$ are a $\dGr$- and hence also a $\Gr$-module as follows. Let us use the identification $\Tpoly=C^\infty(T^*[1]\R^d)$. We denote the even coordinates on $T^*[1]\R^d$ by $x_k$ and the odd coordinates by $\xi_k$.
The action of the graph $\Gamma\in \dGr(n)$ on multivector fields $\gamma_1,\dots, \gamma_n$ can then be written as
\[
 \Gamma(\gamma_1,\dots,\gamma_n) = \mu\circ \left( \prod_{(i,j)}\sum_{k=1}^d \pd{}{x_k^{(j)}} \pd{}{\xi_k^{(i)}} \right)\left(\gamma_1\otimes\cdots\otimes\gamma_n \right).
\]
Here $\mu$ is the operation of multiplication of $n$ multivector fields and the product runs over all edges $(i,j)$ in $\Gamma$, in the order given by the numbering of edges. The notation $\pd{}{x_k^{(j)}}$ means that the partial derivative is to be applied to the $j$-th factor of the tensor product, and similarly for $\pd{}{\xi_k^{(i)}}$.
\end{ex}

\begin{rem}
 There is a map of operads $\Ger\to \Gra$, given on generators as follows. The product in $\Ger(2)$ is sent to the graph with two vertices and no edge, and the Lie bracket is sent to the graph with two vertices and one edge between them. In pictures:
\begin{align*}
 \cdot \wedge \cdot &\mapsto 
\raisebox{-.25em}{\tikz{ 
\node[circle,draw, inner sep=1] (v1) at (0,0) {1};
\node[circle,draw, inner sep=1] (v2) at (.8,0) {2};
}}
\\
\co{\cdot}{\cdot} &\mapsto
\raisebox{-.25em}{\tikz{ 
\node[circle,draw, inner sep=1] (v1) at (0,0) {1};
\node[circle,draw, inner sep=1] (v2) at (.8,0) {2};
\draw (v1) edge (v2);
}}
\end{align*}
One checks that the relations in $\Ger$ are respected by that map. 
\end{rem}

\subsection{\texorpdfstring{$\Gra_1$}{Gra1} moperad}
Define $\Gr_1(m)\subset \dGr(m+2)$ as the space of graphs with no incoming edges at the vertex $m+1$ and no outgoing edges at vertex $m+2$. In fact, we will call the vertex $m+1$ the \emph{output vertex}, for short $out$, and vertex $m+2$ the \emph{input vertex}, or $in$. These spaces assemble to form a $\dGr$-moperad $\Gr_1$. 
The operadic right $\dGr$ action is given by insertions at the vertices $1,\dots, m$. It is inherited from $\dGr$, acting on itself from the right. The composition $\Gamma_1\circ \Gamma_2$ of elements $\Gamma_1, \Gamma_2\in \Gr_1$ is given by the following procedure:
\begin{enumerate}
 \item Delete vertex $\vin$ of $\Gamma_1$ and vertex $\vout$ of $\Gamma_2$. This possibly produces several ``dangling edges''.
 \item Reconnect the dangling edges previously attached to $\vin$ of $\Gamma_1$ in an arbitrary manner to vertices of $\Gamma_2$ and reconnect the open edges previously attached to $\vout$ of $\Gamma_2$ in an arbitrary manner to vertices of $\Gamma_1$. (One sums over all graphs thus produced.)
 \item Relabel vertices such that the vertex $\vout$ of $\Gamma_1$ and the vertex $\vin$ of $\Gamma_2$ become the output and input vertices of $\Gamma_1\circ \Gamma_2$.
\end{enumerate}

\begin{figure}
 \centering
 \makeatletter{}\usetikzlibrary{arrows}
\tikzset{ext/.style={circle, draw,inner sep=1pt},int/.style={circle,draw,fill,inner sep=1pt},nil/.style={inner sep=1pt}}
\tikzset{arr/.style={-triangle 60}}
\begin{tikzpicture}
\node [ext, rectangle] (v1) at (-7,2) {$\vout$};
\node [ext] (v4) at (-6,0) {1};

\node [ext, rectangle] (v2) at (-7,-2) {$\vin$};

\draw[arr] (v1) edge (v2);
\draw[arr] (v1) edge (v4);

\begin{scope}[shift={(0,0)}]

\node [ext, rectangle] (v1) at (-4,2) {$\vout$};
\node [ext] (v5) at (-4,0) {1};
\node [ext, rectangle] (v2) at (-4,-2) {$\vin$};
\draw[arr] (v1) edge (v5);
\draw[arr] (v5) edge (v2);
\end{scope}

\node at (-5,0) {$\circ$};
\node at (-2,0) {$=$};
\node at (-7,-3) {$\Gamma_1$};
\node at (-4,-3) {$\Gamma_2$};

\begin{scope}[shift={(7,1)}]
\node [ext, rectangle] (v1) at (-7,2) {$\vout$};
\node [ext] (v4) at (-6,1) {1};
\node [ext, rectangle] (v2) at (-7,-1) {$\vin$};
\node [ext] (v5) at (-6,0) {2};
\draw[arr] (v1) edge (v5);
\draw[arr] (v4) edge (v5);
\draw[arr] (v5) edge (v2);
\draw[arr] (v1) edge (v4);
\end{scope}

\begin{scope}[shift={(9,-3)}]
\node [ext, rectangle] (v1) at (-7,2) {$\vout$};
\node [ext] (v4) at (-6,1) {1};
\node [ext, rectangle] (v2) at (-7,-1) {$\vin$};
\node [ext] (v5) at (-6,0) {2};
\draw[arr] (v1) edge (v5);
\draw[arr] (v1) edge (v2);
\draw[arr] (v1) edge (v4);
\draw[arr] (v5) edge (v2);
\end{scope}

\begin{scope}[shift={(10,1)}]
\node [ext, rectangle] (v1) at (-7,2) {$\vout$};
\node [ext] (v4) at (-6,1) {1};
\node [ext, rectangle] (v2) at (-7,-1) {$\vin$};
\node [ext] (v5) at (-6,0) {2};
\draw[arr] (v1) edge (v2);
\draw[arr] (v1) edge (v4);
\draw[arr] (v5) edge (v2);
\draw[arr] (v4) edge (v5);
\end{scope}

\node at (1,-2.5) {+};
\node at (2,1.5) {+};
\end{tikzpicture}  
 \caption{\label{fig:gra1composition} Illustration of the (m)operadic composition in the moperad $\Gra_1$, of two elements $\Gamma_1,\Gamma_2\in \Gra_1$. One has to delete vertex $\vin$ of $\Gamma_1$ and vertex $\vout$ of $\Gamma_2$ and reconnect the open edges produced in an arbitrary manner. Here we have two open edges. Each can be reconnected in 2 ways. This yields a sum of four graphs, one of which is zero because it contains a double edge. The remaining three graphs are shown.}
\end{figure}
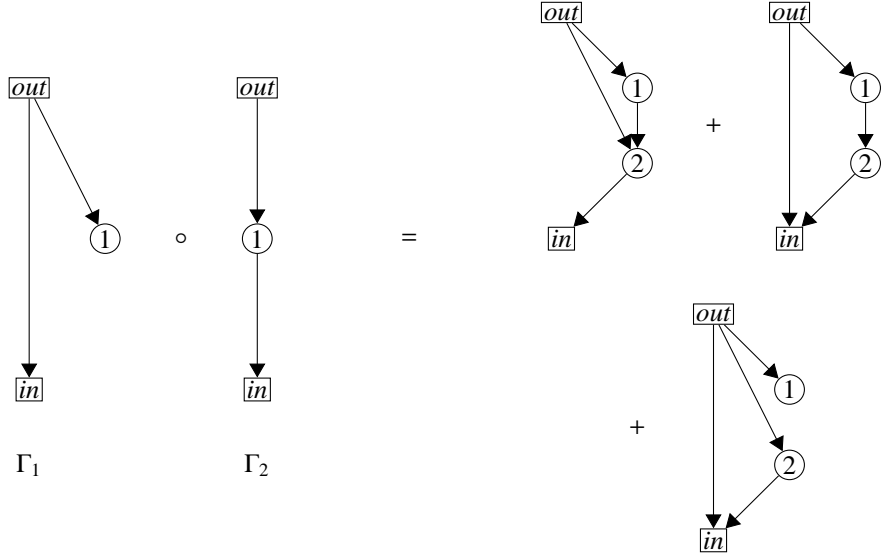

The procedure is depicted in Figure \ref{fig:gra1composition}.

\begin{ex}
\label{ex:Gra1actonOmega}
 The colored operad $\bpm \dGr & \Gr_1\epm$ can be represented on the colored vector space $\Tpoly\oplus \Omega_\bullet$. The action of $\dGr$ on $\Tpoly$ is the one from example \ref{ex:GraactonTpoly}. The action of $\Gr_1$ is defined as follows. Let $\Gamma\in \Gr_1(m)$ and let $\Gamma'$ be the same graph considered as an element of $\dGr(m+2)$. Let $\gamma, \gamma_1,\dots,\gamma_n$ be multivector fields and $\omega$ a differential form. We can assume $\omega = f \omega_0$ with $f$ a function and $\omega_0$ constant. Then we define the action of $\Gamma$ such that 
\[
 \iota_{\gamma}\Gamma(\gamma_1,\dots,\gamma_n; \omega)= (-1)^{|\Gamma||\gamma|} \iota_{\Gamma'(\gamma_1,\dots,\gamma_n, \gamma, f)} \omega_0.
\]
This defines the action uniquely.
\end{ex}

\begin{rem}
 There is map of colored operads $\calc \to \bpm \Gra & \Gra_1 \epm$. The $\Ger$-part was described in the previous subsection. For the $\calc_1$-part one maps the generators $d$ and $\iota$ to the graphs depicted in Figure \ref{fig:dliota}.
\end{rem}
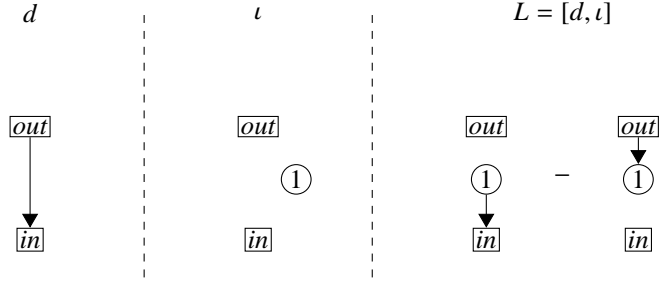
\begin{figure}
 \centering
\makeatletter{}\usetikzlibrary{matrix}
\usetikzlibrary{arrows}
\usetikzlibrary{shapes}
\usetikzlibrary{through}
\usetikzlibrary{calc,3d}
\usetikzlibrary{decorations,decorations.pathmorphing}
\[
\begin{tikzpicture}[
int/.style={circle, draw, fill, minimum size=5pt, inner sep=0},
ext/.style={circle, draw, fill=white, minimum size=5pt, inner sep=1pt},
helper/.style={coordinate},point/.style={circle, draw, fill, inner sep =1pt},
de/.style={-triangle 60},
point/.style={circle, draw, fill, minimum size=3pt, inner sep=0pt},
xst/.style={cross out, draw, minimum size=5 },
scale=.5
]
\begin{scope}[]
\node at (0,5) {$d$};

\node [draw, inner sep=1] (out) at ($(0,0)+(-90:2 and 1)$) {$\vin$};
\node [draw, inner sep=1] (in) at ($(0,3)+(-90:2 and 1)$) {$\vout$};
\draw[-triangle 60] (in)--(out) ;
\end{scope}

\draw[dashed] (3,-2)--(3,5);

\begin{scope}[xshift=6cm]
\node at (0,5) {$\iota$};

\node [draw, inner sep=1] (out) at ($(0,0)+(-90:2 and 1)$) {$\vin$};
\node [draw, inner sep=1] (in) at ($(0,3)+(-90:2 and 1)$) {$\vout$};
\node[ext] at (1,.6) {1};
\end{scope}

\draw[dashed] (9,-2)--(9,5);

\begin{scope}[xshift=12cm]
\node at (2,5) {$L=[d,\iota]$};
\node at (2,.7) {$-$};

\node [draw, inner sep=1] (out) at ($(0,0)+(-90:2 and 1)$) {$\vin$};
\node [draw, inner sep=1] (in) at ($(0,3)+(-90:2 and 1)$) {$\vout$};
\node[ext] (e1) at (0,.6) {1};
\draw[-triangle 60] (e1)--(out) ;
\begin{scope}[xshift=6cm, shift={(-2,0)}]

\node [draw, inner sep=1] (out) at ($(0,0)+(-90:2 and 1)$) {$\vin$};
\node [draw, inner sep=1] (in) at ($(0,3)+(-90:2 and 1)$) {$\vout$};
\node[ext] (e1) at (0,0.6) {1};
\draw[triangle 60-] (e1)--(in) ;
\end{scope}

\end{scope}

\end{tikzpicture}
\] 
\caption{\label{fig:dliota} The images of the $\calc_1$-elements $d$, $\iota$, and $L$ under the map 
$\calc_1 \to \Gra_1$.}
\end{figure}

\subsection{\texorpdfstring{$\PT$}{PT} (``planar trees'') operad}
\label{sec:PT}
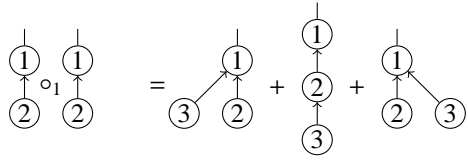
\begin{figure}
\centering
 \begin{tikzpicture}[scale=.7]

\node (e1) at (0,1) [ext] {$1$};
\node (e2) at (0,0) [ext] {$2$};
\draw[->] (e2) edge (e1);
\draw (e1) -- +(0,.6);

\node at(.5,.5){$\circ_1$};

\begin{scope}[xshift=1cm]
\node (e1) at (0,1) [ext] {$1$};
\node (e2) at (0,0) [ext] {$2$};
\draw[->] (e2) edge (e1);
\draw (e1) -- +(0,.6);
\end{scope}

\node at(2.5,.5){$=$};

\begin{scope}[xshift=3cm]
\node (e1) at (1,1) [ext] {$1$};
\node (e2) at (1,0) [ext] {$2$};
\node (e3) at (0,0) [ext] {$3$};
\draw[->] (e2) edge (e1) (e3) edge (e1);
\draw (e1) -- +(0,.6);
\end{scope}

\node at (4.75,.5) {$+$};

\begin{scope}[xshift=5.5cm]
\node (e1) at (0,1.5) [ext] {$1$};
\node (e2) at (0,.5) [ext] {$2$};
\node (e3) at (0,-.5) [ext] {$3$};
\draw[->] (e2) edge (e1) (e3) edge (e2);
\draw (e1) -- +(0,.6);
\end{scope}

\node at (6.25,.5) {$+$};

\begin{scope}[xshift=7cm]
\node (e1) at (0,1) [ext] {$1$};
\node (e2) at (0,0) [ext] {$2$};
\node (e3) at (1,0) [ext] {$3$};
\draw[->] (e2) edge (e1) (e3) edge (e1);
\draw (e1) -- +(0,.6);
\end{scope}
\end{tikzpicture}
\caption{\label{fig:PTcomp}The operadic composition in $\PT$.}
\end{figure}

Let $\PT(n)'$ be the graded vector space of linear combinations of rooted planar trees with vertex set $[n]$ and edge set $[n-1]$.\footnote{In our conventions the trees are not planted.} Such a tree is considered as living in degree $-n+1$, i.e., the edges are considered to have degree $-1$. Define
\[
 \PT(n) = (\PT(n)' \otimes \sgn_{n-1})_{S_{n-1}}
\]
where the symmetric group $S_{n-1}$ acts on $\PT(n)'$ by permuting the labels on edges and $\sgn_{n-1}$ is the sign representation. The spaces  $\PT(n)$ assemble to form an operad $\PT$. The operadic compositions $T_1\circ_j T_2$ are given by inserting the tree $T_2$ into the vertex $j$ of $T_1$ and reconnecting the incoming edges at $j$ in all planar possible ways, cf. Figure \ref{fig:PTcomp} or \cite{MSS}, section I.1.20.

\begin{rem}
 In pictures, we will sometimes draw the edges of the tree as arrows pointing towards the root. Hence the root is the unique vertex without outgoing arrows. Additionally we will draw a small stub at the root, which is however not considered an edge for degree purposes.
\end{rem}

\begin{rem}
 One could define $\PT(n)$ without labelling edges and taking coinvariants, by specifying some ordering on the edges using the planar structure. However, in this manner the signs will be clearer, and also the similarity to the sign rules for $\Gr$.
\end{rem}

\begin{lemma}
\label{lem:PTgenrel}
 The operad $\PT$ has the following presentation in terms of generators and relations: The generators are $\R[\SymGrp_{n+1}] T_n\in \PT(n+1)$ for $n\geq 1$ of degree $-n$, see Figure \ref{fig:KSgen} (left). We denote by $T_0$ the unit of the operad. The relations are given by the ``planar Leibniz rule''
\[
 T_m \circ_1 T_n = \sum_{\substack{ k_1,\dots,k_n \\ J:=\sum_i k_i \leq m}} \sum_{1\leq j_1 <j_2<\cdots <j_n\leq n+m-J}
\pm
T_{n+m-J}\circ_{j_1,\dots, j_n } (T_{k_1},\dots, T_{k_n}) \, .
\]
The symbol $\circ_{j_1,\dots, j_n }$ denotes the operadic insertion into the $j_1$-st, $j_2$-nd, $\dots$, $j_n$-th slots. The signs can be determined by kepping track of the order of edges.
\end{lemma}
\begin{proof}
 Let the operad generated by the above presentation be denoted $\op P$. We want to show $\op P\cong \PT$. There is an obvious surjective map $\op P\to \PT$ since the above relation is satisfied in $\PT$ and $\PT$ is generated by the $T_n$. 
Next let $T'$ be some element of the free operad generated by the $T_n$. We call it ``good'' if it does not contain any insertions into the first (top) slot of any $T_n$. The set of good $T'$ is isomorphic to the set of $\PT$ trees. Hence we are done if we can show that any $T'$ is equivalent in $\op P$ to a good one. But the above relation (from left to right) can be used to eliminate one by one each insertion into a first slot of a generator.
\end{proof}

\begin{ex}
\label{ex:PTmod}
 The total space $\oplus_n \op P(n)[-n]$ of any (non-symmetric) operad $\op P$ over dg vector spaces is a $\PT$-module. In particular $\Dpoly$, considered with zero differential, is a $\PT$-module. Concretely, the action of the generator $T_n$ is given by braces operations
 \[
 T_n(a_0,\dots, a_n) = \pm a_0\{a_1,\dots, a_n\}
=
\pm 
\sum_{ 1\leq j_1<\cdots < j_n \leq |a_0| } 
(-1)^{\sum_i (|a_i|+1)(j_i-1) }
a_0 \circ_{j_1, \dots, j_n}(a_1,\dots, a_n)
 \]
where $a_0,\dots, a_n\in \Dpoly$ and the notation $a_0 \circ_{j_1, \dots, j_n}(a_1,\dots, a_n)$ means that $a_1$ is inserted at the $j_1$-th slot of $a_0$ etc..
The sign is determined similarly to that in the action of $\dGra$. Consider $T_n$ as a non-planar directed graph, say $t_n\in \dGra(n+1)$. As such it acts on $\Tpoly \subset \Dpoly$. We set the sign such that the terms occuring both in the action of $t_n$ and $T_n$ have the same sign.
For example, $T_1$ acts on multidifferential operators $a_0, a_1$ as 
\[
T_1(a_0, a_1) = (-1)^{(|a_0|-1)|a_1| } a_0\{ a_1\}.
\]
Hence the element depicted in figure \ref{fig:PTmc} acts, up to a sign, as the Gerstenhaber bracket. The reader shall not be confused that the latter element of the operad is symmetric, while the Gerstenhaber bracket is antisymmetric. There is a sign hidden inside the Koszul conventions. Concretely, let $\sigma\in S_2$ be the transposition of two elements. Then 
\begin{align*}
(T_1 + T_1\sigma)(a_0, a_1)
&=
T_1(a_0, a_1) + (-1)^{|a_0||a_1|}T_1(a_1, a_0)\\
&=
(-1)^{(|a_0|-1)|a_1|}( a_0\{ a_1\} - (-1)^{(|a_0|-1)(|a_1|-1)} a_1\{ a_0\} ).
\end{align*}
\end{ex}

\begin{ex}
\label{ex:SCtoPTmod}
There is a variant of the previous example we will need below.
Suppose we have a two-colored operad $\op P$ of Swiss Cheese type (see definition \ref{def:SCtype}).
Let $\op P^1(\cdot, 0)$ be the operad of operations with all inputs in color 1 and let $\op P^2$ be the operations with output in color 2. in particular, $\op P^2$ is a non-symmetric operad in the category of right $\op P^1$-modules. By a varying the previous example slightly, the total spaces
\[
\prod_n \op P^2(\cdot, n)[-n]
\]
form an operadic $\PT$-$\op P^1$-module. Furthermore, since we used only natural operations, the map from Swiss Cheese type operads to two colored operads
\[
\op P \mapsto \bpm \PT & \prod_n \op P^2(\cdot, n)[-n] & \op P^1 \epm 
\]
is functorial.
\end{ex}

\subsection{Braces operad -- first definition}
Let us recall the definition of the braces operad, following \cite{GJ}.
Let $V$ be a dg vector space. 
Let $T V[1]$ be the tensor coalgebra on $V[1]$.
A $\op B_\infty$-structure on $V$ is a dg bialgebra structure on $T V[1]$, such that 
\begin{enumerate}
 \item The coproduct is the standard one.
 \item The differential extends the given differential on $V$.
 \item $1\in T V[1]$ is the unit.
\end{enumerate}
Let $D$ be the differential and $m$ be the product.
Since the coalgebra $T V[1]$ is cofree, $D$ and $m$ are uniquely determined by their compositions with the projection onto the generators (i.e., onto $V$). Hence a $\op B_\infty$ structure is given by families of maps
\begin{gather*}
 D_k \colon  V[1]^k \to V[2] \quad\quad\quad \text{for $k\geq 1$}\\
 m_{kl} \colon V[1]^k\otimes V[1]^l \to V[1] \quad\quad\quad \text{for $k,l\geq 0$}
\end{gather*}
satisfying certain compatibility relations. More precisely, $D_1$ is already determined by the second condition, and so are $m_{kl}$ for $k$ or $l$ equal to zero by the third. A \emph{braces algebra} structure on $V$ is a $\op B_\infty$-structure such that $m_{kl}=0$ for $k>1$. The operad $\Br$ is the operad governing braces algebra structures. By the above description, $\Br$ is generated by operations 
\begin{align*}
 D_k \in \Br(k)& \quad\text{of degree $2-k$, for $k\geq 2$} \\
 m_k \in \Br(k+1)& \quad \text{of degree $-k$, for $k\geq 1$}.
\end{align*}

The conditions above lead to following set of relations in $\Br$.
\begin{align*}
dD_n +\sum_{\substack{k,l \\ k+l=n+1}} \sum_{ 1\leq j\leq k } \pm 
D_k \circ_j D_l 
&=0
\\
 d m_n + 
\sum_{\substack{k,l \\ k+l=n+1}} 
\sum_{ 1\leq j\leq k } 
\pm D_k \circ_j m_l 
+ \sum_{\substack{k,l \\ k+l=n+1}} (-1)^{kl} \sum_{1\leq j\leq l } 
\pm m_l \circ_j D_k &= 0
\\
\sum_{\substack{k\geq n', k_1, \dots, k_{n'} \\  k+k_1+\dots+k_{n'}=n+n'}} \sum_{ 1\leq j_1<\cdots < j_{n'} \leq k }
\pm D_k\circ_{j_1, \dots, j_{n'}} (m_{k_1}, \dots,  m_{n'}) 
+ m_n \circ_1 D_{n'} &=0
\\
\sum_{\substack{k\geq n', k_1, \dots, k_{n'} \\ k+k_1+\dots+k_{n'}=n+n'}} \sum_{ 1\leq j_1<\cdots < j_{n'} \leq k }
\pm m_k\circ_{j_1, \dots, j_{n'}} (m_{k_1}, \dots,  m_{k_{n'}}) - m_n \circ_1 m_{n'} &=0
\end{align*}
Here the first equation states that the differential squares to zero. The next two equations come from the compatibility of the differential with the product. The third equation is required to hold for $n'\geq 2$ only.
The fourth equation is the associativity of the product.

\begin{ex}
\label{ex:Bractiononalg}
 The Hochschild complex of any $A_\infty$-algebra is a $\Br$-algebra with the well known formulas
\begin{align*}
 D_k(a_1,\dots, a_k) &=  \mu_k(a_1,\dots, a_k) \\
 m_k(a_0,\dots, a_k) &=  
\sum_{ 1\leq j_1<\cdots < j_k \leq |a_0| } 
(-1)^{\sum_i (|a_i|+1)(j_i-1) }
a_0 \circ_{j_1, \dots, j_k}(a_1,\dots, a_k)
=: a_0\{ a_1,\dots, a_k \},
\end{align*}
where the $\mu_k$ are the $A_\infty$-operations and the notation $a_0 \circ_{j_1, \dots, j_k}(a_1,\dots, a_k)$ means that $a_1$ is inserted at the $j_1$-th slot of $a_0$ etc.
\end{ex}

\begin{rem}
 In the literature the braces operad is often defined as the quotient of our operad $\Br$ by the relation that $m_k=0$ for each $k\geq 3$. This smaller operad then acts on the Hochschild complex of any algebra, but not in general on the Hochschild complexes of $A_\infty$-algebras.
\end{rem}

\subsection{Braces operad -- Kontsevich-Soibelman version}
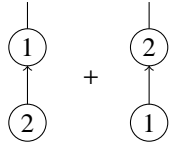
\begin{figure}
\centering
\makeatletter{}\[
 \begin{tikzpicture}[scale=1,
vert/.style={draw,outer sep=0,inner sep=0,minimum size=5,shape=circle,fill},
helper/.style={outer sep=0,inner sep=0,minimum size=5,shape=coordinate},
default_edge/.style={draw,->},
ext/.style={draw,outer sep=0,inner sep=2,minimum size=5,shape=circle},
every loop/.style={}]

\node (v0) at (5,6.5) [ext] {2};
\node (v1) at (5,7.5) [ext] {1};
\node (v3) at (5,8.4) [helper] {};
\node (v10) at (6.6,7.5) [ext] {2};
\node (v9) at (6.6,6.5) [ext] {1};
\node (v12) at (6.6,8.4) [helper] {};
\node (v14) at (5.6,7.1) [helper,label=0:{$+$}] {};

\draw (v1) -- +(0,.6);
\draw (v10) -- +(0,.6);

\draw[default_edge] (v0) to (v1);
\draw[default_edge] (v9) to (v10);
\end{tikzpicture}
\]
 
\caption{\label{fig:PTmc} The Maurer-Cartan element in $\PT$. Note that this element is symmetric. However, the resulting operation on cochains, the Gerstenhaber bracket, is antisymmetric. The additional signs come from the definition of the action.}
\end{figure}

\begin{figure}
\centering
\makeatletter{}\begin{tikzpicture}[scale=1, every edge/.style={draw, <-}]

\node (e1) at (1,1) [ext] {$\phantom{0}$};
\node (e2) at (0,0) [nil] {};
\node (e3) at (.75,0) [nil] {};
\node (ed) at (1.5,0) {$\dots$};
\node (e4) at (2,0) [nil] {};
\node at(-.5,.2){$d$};
\draw (e1) -- +(0,.6);

\draw (e1) edge (e2) edge (e3) edge (e4);

\node at(3,.2) {$=\quad \sum \pm$};

\begin{scope}[xshift=4cm]
\node (etop) at (1,1.5) [ext] {$\phantom{0}$};
\node (emiddle) at (1,.75) [int] {};
\node (e1) at (-0.2,0) [nil] {};
\node at (0.2,0) {$\dots$};
\node (e2) at (.5,0) [nil] {};
\node at (1,0) {$\dots$};
\node (e4) at (1.5,0) [nil] {};
\node at (1.75,0) {$\dots$};
\node (e5) at (2.2,0) [nil] {};
\draw (etop) edge (e1) edge (emiddle) edge (e5);
\draw (emiddle) edge (e2) edge (e4);
\draw (etop) -- +(0,.6);
\end{scope}

\node at(7,.2) {$+\quad \sum \pm$};
\begin{scope}[xshift=8cm]
\node (etop) at (1,1.5) [int] {};
\node (emiddle) at (1,.75) [ext] {$\phantom{0}$};
\node (e1) at (-0.2,0) [nil] {};
\node at (0.2,0) {$\dots$};
\node (e2) at (.5,0) [nil] {};
\node at (1,0) {$\dots$};
\node (e4) at (1.5,0) [nil] {};
\node at (1.75,0) {$\dots$};
\node (e5) at (2.2,0) [nil] {};
\draw (etop) edge (e1) edge (emiddle) edge (e5);
\draw (emiddle) edge (e2) edge (e4);
\draw (etop) -- +(0,.6);
\end{scope}
\end{tikzpicture}
 \begin{tikzpicture}[scale=1, every edge/.style={draw, <-}]

\node (e1) at (1,1) [int] {};
\node (e2) at (0,0) [nil] {};
\node (e3) at (.75,0) [nil] {};
\node (ed) at (1.5,0) {$\dots$};
\node (e4) at (2,0) [nil] {};
\node at(-.5,.2){$d$};
\draw (e1) -- +(0,.6);

\draw (e1) edge (e2) edge (e3) edge (e4);

\node at(3,.2) {$=\quad \sum \pm$};

\begin{scope}[xshift=4cm]
\node (etop) at (1,1.5) [int] {};
\node (emiddle) at (1,.75) [int] {};
\node (e1) at (-0.2,0) [nil] {};
\node at (0.2,0) {$\dots$};
\node (e2) at (.5,0) [nil] {};
\node at (1,0) {$\dots$};
\node (e4) at (1.5,0) [nil] {};
\node at (1.75,0) {$\dots$};
\node (e5) at (2.2,0) [nil] {};
\draw (etop) edge (e1) edge (emiddle) edge (e5);
\draw (emiddle) edge (e2) edge (e4);
\draw (etop) -- +(0,.6);
\end{scope}
\end{tikzpicture}
 
\caption{\label{fig:KSdiff}The differential in the operads $\Tw\PT$ and $\Br'\cong \Br$.}
\end{figure}
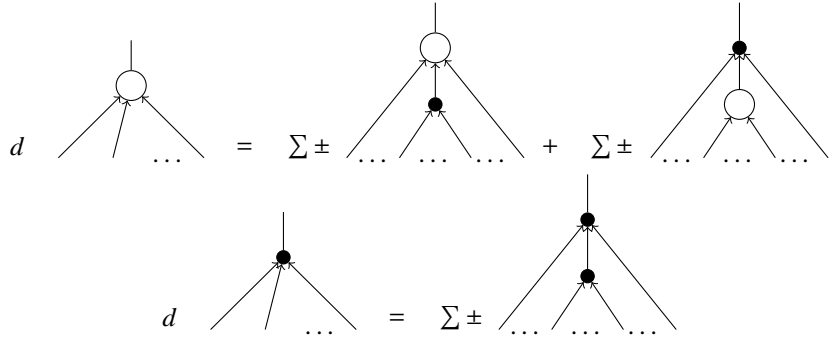

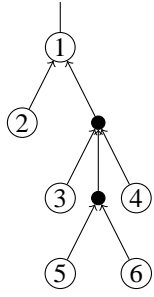
\begin{figure}
\centering
\makeatletter{}\usetikzlibrary{matrix}
\usetikzlibrary{arrows}
\usetikzlibrary{shapes}
\usetikzlibrary{through}
\usetikzlibrary{calc,3d}
\usetikzlibrary{decorations,decorations.pathmorphing}
\[
\begin{tikzpicture}[
yscale=-1,
helper/.style={coordinate},point/.style={circle, draw, fill, inner sep =1pt},
de/.style={-triangle 60},
point/.style={circle, draw, fill, minimum size=3pt, inner sep=0pt},
xst/.style={cross out, draw, minimum size=5 },
every edge/.style={draw, <-}
]

\node[ext] (e1) at (90:1) {$1$};
\draw (e1) +(-.5,1) node[ext] (e2) {$2$} ++(.5,1) node[int] (i1) {} +(-.5,1) node[ext](e3) {$3$}  +(.5,1) node[ext] (e4){$4$} ++(0,1) node[int] (i2) {} +(-.5,1) node[ext](e5) {$5$} +(.5,1) node[ext](e6){$6$} ;
\draw (e1) edge (e2) edge (i1) (i1) edge (e3) edge (i2) edge (e4) (i2) edge (e5) edge (e6);
\draw (e1) -- +(0,-.6);
\end{tikzpicture}
\] 
\caption{\label{fig:Brex} A typical element of $\Br(6)$.}
\end{figure}

Consider again the operad $\PT$. It comes with a natural map $\Lie_1\to \PT$, by sending the generator to the element depicted in Figure \ref{fig:PTmc}. Hence one can apply a general procedure we call \emph{operadic twisting} to $\PT$ and get another operad $\Tw\PT$. For details on this twisting procedure, see Appendix \ref{sec:optwists}, or the more detailed discussion in \cite{vasilydeligne}. Concretely, the operad $\Tw\PT$ is spanned by rooted planar trees with two kinds of vertices, called internal and external. The external vertices are numbered, the internal vertices are not, see Figure \ref{fig:Brex} for an example. There is now a differential, which creates new internal vertices, see Figure \ref{fig:KSdiff}. Let $\Br' \subset \Tw\PT$ be the suboperad spanned by trees all of whose internal vertices have valence at least 2. This operad was introduced by Kontsevich and Soibelman \cite{KS1}.\footnote{To be precise, the way we introduce this operad here by operadic twisting is a bit different from their approach. But the operad is the same.}

\begin{figure}
\centering
 \begin{tikzpicture}[scale=1]

\node (e1) at (1,1) [ext] {$1$};
\node (e2) at (0,0) [ext] {$2$};
\node (e3) at (.75,0) [ext] {$3$};
\node (ed) at (1.5,0) {$\dots$};
\node (e4) at (2,0) [ext] { $\scriptstyle n{+}1$};
\draw (e1) -- +(0,.6);
\draw[->] (e2) edge (e1) (e3) edge (e1) (e4) edge (e1);

\begin{scope}[xshift=3cm]
\node (e1) at (1,1) [int] {};
\node (e2) at (0,0) [ext] {$1$};
\node (e3) at (.75,0) [ext] {$2$};
\node (ed) at (1.5,0) {$\dots$};
\node (e4) at (2,0) [ext] {$n$};
\draw (e1) -- +(0,.6);
\draw[->] (e2) edge (e1) (e3) edge (e1) (e4) edge (e1);
\end{scope}
\end{tikzpicture}
\caption{\label{fig:KSgen}The generators $T_n$ (left) and $T_n'$ (right) of the operad $\Br'\cong\Br$.}
\end{figure}
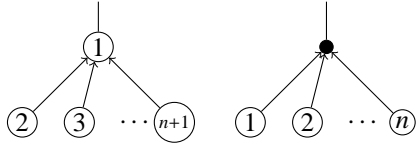
 
Similarly to Lemma \ref{lem:PTgenrel} one proves the following.
\begin{lemma}
 \label{lem:KSgenrel}
 The operad $\Br'$ has the following presentation in terms of generators and relations: The generators are $\R[\Sigma_{n+1}] T_n\in PT(n+1)$ for $n\geq 1$ of degree $-n$, and $\R[\Sigma_n] T_n'\in PT(n)$ for $n\geq 2$ of degree $-n+2$, see Figure \ref{fig:KSgen}.
The differential $d T_n=(\cdots)$, $dT_n'=(\cdots)$ is given pictorially in Figure \ref{fig:KSdiff}.
The relations are given by the ``planar Leibniz rule''
\begin{align*}
 T_m \circ_1 T_n &= 
\sum_{\substack{ j_1,\dots,j_n \\ J:=\sum_k j_k \leq m}} 
\sum_{1\leq i_1 <i_2<\cdots <i_n\leq n+m-J}
\pm
(\cdots(T_{n+m-J}\circ_{i_1}T_{j_1})\cdots) \circ_{i_n}T_{j_n}  \\
T_m \circ_1 T_n' &= \sum_{\substack{ j_1,\dots,j_n \\ J:=\sum_k j_k \leq m}} \sum_{1\leq i_1 <i_2<\cdots <i_n\leq n+m-J}
\pm
(\cdots(T_{n+m-J}'\circ_{i_1}T_{j_1})\cdots) \circ_{i_n}T_{j_n}
\, .
\end{align*}
\end{lemma}

\begin{ex}
$\Br'$ acts on the Hochschild complex of any $A_\infty$-algebra. This follows directly from generalities on operadic twisting (see Appendix \ref{sec:optwists}), from example \ref{ex:PTmod}, and the fact that an $A_\infty$ structure provides a Maurer-Cartan element.
\end{ex}

\begin{cor}
The operads $\Br$ and $\Br'$ are isomorphic.
\end{cor}
In the following we will drop the notation $\Br'$ and call either operad $\Br$.

\begin{proof}
Comparing the generators and relations of $\Br'$ with those of $\Br$ recalled in the previous subsection one sees that $\Br'\cong \Br$, up to signs. To check that the signs can be chosen correctly, it suffices to note that both operads act on the Hochschild complex of an $A_\infty$ algebra by the same universal formulas, up to signs.
\end{proof}

\subsection{\texorpdfstring{$\PT_1$}{PT1} moperad }
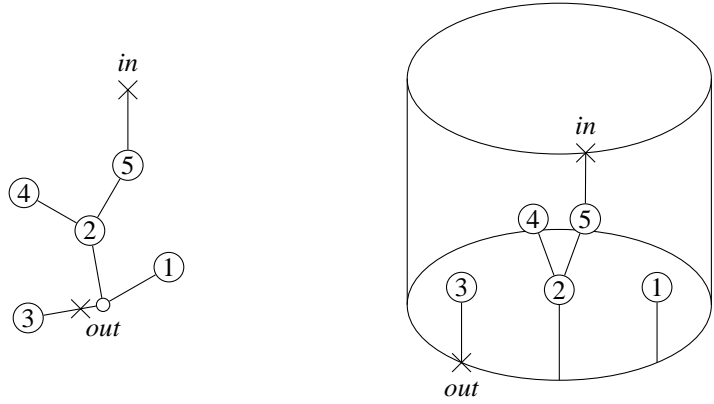
\begin{figure}
 \centering
\makeatletter{}\usetikzlibrary{matrix}
\usetikzlibrary{arrows}
\usetikzlibrary{shapes}
\usetikzlibrary{through}
\usetikzlibrary{calc,3d}
\usetikzlibrary{decorations,decorations.pathmorphing}
\[
\begin{tikzpicture}[
int/.style={circle, draw, fill, minimum size=5pt, inner sep=0},
ext/.style={circle, draw, fill=white, minimum size=5pt, inner sep=1pt},
helper/.style={coordinate},point/.style={circle, draw, fill, inner sep =1pt},
de/.style={-triangle 60},
point/.style={circle, draw, fill, minimum size=3pt, inner sep=0pt},
xst/.style={cross out, draw, minimum size=5pt}
]

\node[ext, label=-90:{$out$}] (out) at (0,0) {};
\node[ext] (e1) at (30:1) {$1$};
\node[ext] (e2) at (100:1) {$2$};
\node[ext] (e3) at (190:1) {$3$};
\draw (e2)+(60:1) node[ext] (e5) {$5$} 
+(150:1) node[ext](e4) {$4$} 
(e5)+(90:1) node[xst, label=90:{$in$}](in) {};

\draw (out) edge (e1) edge (e2) edge (e3) (e2) edge (e5) edge (e4) (e5) edge (in.base);
\node[xst] (mark) at(190:.3){};

\begin{scope}[xshift=6cm]

\draw (0,0) ellipse (2 and 1);
\draw (0,3) ellipse (2 and 1);
\draw (-2,0)--(-2,3) (2,0)--(2,3);

\node [xst, label=-90:{$out$}] (out) at ($(0,0)+(-130:2 and 1)$) {};
\node [xst, label=90:{$in$}] (in) at ($(0,3)+(-80:2 and 1)$) {};

\draw (out)+(0,1) node[ext] (e3) {$3$};
\node[ext] (e2) at ($(0,1.2)+(-90:2 and 1)$) {$2$};
\node[ext] (e1) at ($(0,1)+(-50:2 and 1)$) {$1$};
\draw (e2)+(70:1) node[ext] (e5) {$5$} 
+(110:1) node[ext](e4) {$4$} 
(e5)+(90:1);

\draw (out.base) edge (e3);
\draw (e2) edge (e5) edge (e4) (e5) edge (in.base);
\draw (e2) edge (0,-1) (e1) edge ($(0,0)+(-50:2 and 1)$);

\end{scope}

\end{tikzpicture}
\] 
\caption{\label{fig:PT1ex} An example of a graph in $\PT_1$. Conventionally, one draws it on a cyclinder (right). Here the vertex $\vout$ is extruded and becomes the lower rim of the cylinder, the vertex $\vin$ becomes the upper rim.}
\end{figure}

Let $\PT_1(n)'$ be the graded space spanned by (a priori non-rooted) planar trees $T$ of the following type:
\begin{enumerate}
 \item The vertex set of $T$ consists of two special vertices $\vin$ and $\vout$, and $n$ numbered vertices $1,\dots,n$. We consider $\vout$ as the root vertex of the planar tree. Hence the \emph{children} of some vertex are its neighbors farther away from $\vout$. 
 \item A \emph{framing} is fixed at $\vout$, by which we mean that one edge incident to $\vout$ is marked.
 \item We require that the vertex $\vin$ must not have any children, 
 \item The edges (except the one incident at $\vout$) are labelled by numbers $\{1,\dots,n\}$.
 \item We consider such a graph to live in degree $-n$.
\end{enumerate}
Then we define
\[
 \PT_1(n) = (\PT_1(n)' \otimes \sgn_{n})_{S_{n}}.
\]
where the symmetric group $S_{n}$ acts on $\PT(n)'$ by permuting the labels on edges and $\sgn_{n}$ is the sign representation.
An example of a tree in $\PT_1(n)$ is shown in Figure \ref{fig:PT1ex}. It is conventional to draw such a graph on a cylinder, with the vertex $\vout$ forming the lower rim, and the vertex $\vin$ lying on the upper rim.

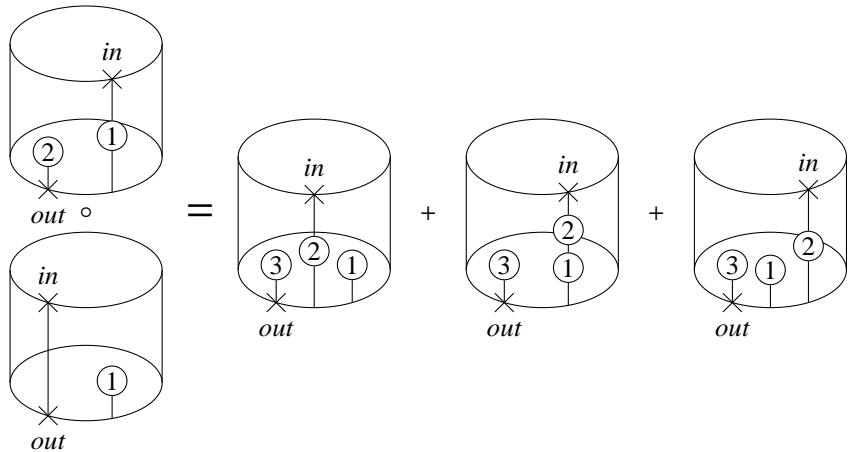
\begin{figure}
 \centering
\makeatletter{}\usetikzlibrary{matrix}
\usetikzlibrary{arrows}
\usetikzlibrary{shapes}
\usetikzlibrary{through}
\usetikzlibrary{calc,3d}
\usetikzlibrary{decorations,decorations.pathmorphing}
\[
\begin{tikzpicture}[
int/.style={circle, draw, fill, minimum size=5pt, inner sep=0},
ext/.style={circle, draw, fill=white, minimum size=5pt, inner sep=1pt},
helper/.style={coordinate},point/.style={circle, draw, fill, inner sep =1pt},
de/.style={-triangle 60},
point/.style={circle, draw, fill, minimum size=3pt, inner sep=0pt},
xst/.style={cross out, draw, minimum size=5 },
scale=.5
]
\begin{scope}[]
\draw (0,0) ellipse (2cm and 1cm);
\draw (0,3) ellipse (2cm and 1cm);
\draw (-2,0)--(-2,3) (2,0)--(2,3);
\node [xst, label=-90:{$out$}] (out) at ($(0,0)+(-120:2 and 1)$) {};
\node [xst, label=90:{$in$}] (in) at ($(0,3)+(-120:2 and 1)$) {};
\draw (in.base)--(out.base) ;
\node[ext] (e1) at ($(0,1)+(-70:2 and 1)$) {1};
\draw (e1) edge  ($(0,0)+(-70:2 and 1)$);
\end{scope}

\begin{scope}[yshift=6cm]
\node at (0,-1.5) {$\circ$};

\draw (0,0) ellipse (2cm and 1cm);
\draw (0,3) ellipse (2cm and 1cm);
\draw (-2,0)--(-2,3) (2,0)--(2,3);
\node [xst, label=-90:{$out$}] (out) at ($(0,0)+(-120:2 and 1)$) {};
\node [xst, label=90:{$in$}] (in) at ($(0,3)+(-70:2 and 1)$) {};
\node[ext] (e1) at ($(0,1.5)+(-70:2 and 1)$) {1};
\node[ext] (one) at ($(0,1)+(-120:2 and 1)$) {$2$};
\draw (in.base) edge (e1) (e1) edge  ($(0,0)+(-70:2 and 1)$) (out.base) edge (one);
\end{scope}

\begin{scope}[xshift=6cm, yshift=3cm]
\node at (-3,1.5) {\huge $=$};

\draw (0,0) ellipse (2cm and 1cm);
\draw (0,3) ellipse (2cm and 1cm);
\draw (-2,0)--(-2,3) (2,0)--(2,3);
\node [xst, label=-90:{$out$}] (out) at ($(0,0)+(-120:2 and 1)$) {};
\node [xst, label=90:{$in$}] (in) at ($(0,3)+(-90:2 and 1)$) {};
\node[ext] (one) at  ($(0,1)+(-120:2 and 1)$)  {$3$};
\draw (one)--(out.base) ;
\node[ext] (e1) at ($(0,1)+(-60:2 and 1)$) {1};
\node[ext](e2) at ($(0,1.5)+(-90:2 and 1)$) {$2$};
\draw (e1) edge  ($(0,0)+(-60:2 and 1)$) (in.base) edge (e2);
\draw (e2) edge  ($(0,0)+(-90:2 and 1)$);
\node at (3,1.5) {$+$};

\begin{scope}[xshift=6cm]

\draw (0,0) ellipse (2cm and 1cm);
\draw (0,3) ellipse (2cm and 1cm);
\draw (-2,0)--(-2,3) (2,0)--(2,3);
\node [xst, label=-90:{$out$}] (out) at ($(0,0)+(-120:2 and 1)$) {};
\node [xst, label=90:{$in$}] (in) at ($(0,3)+(-70:2 and 1)$) {};
\node[ext] (one) at  ($(0,1)+(-120:2 and 1)$)  {$3$};
\draw (one)--(out.base) ;
\node[ext] (e1) at ($(0,1)+(-70:2 and 1)$) {1};
\node[ext](e2) at ($(0,2)+(-70:2 and 1)$) {$2$};
\draw (e1) edge  ($(0,0)+(-70:2 and 1)$) (in.base) edge (e2);
\draw (e2) edge  (e1);
\node at (3,1.5) {$+$};

\end{scope}

\begin{scope}[xshift=12cm]
\draw (0,0) ellipse (2cm and 1cm);
\draw (0,3) ellipse (2cm and 1cm);
\draw (-2,0)--(-2,3) (2,0)--(2,3);
\node [xst, label=-90:{$out$}] (out) at ($(0,0)+(-120:2 and 1)$) {};
\node [xst, label=90:{$in$}] (in) at ($(0,3)+(-60:2 and 1)$) {};
\node[ext] (one) at  ($(0,1)+(-120:2 and 1)$)  {$3$};
\draw (one)--(out.base) ;
\node[ext] (e1) at ($(0,1)+(-90:2 and 1)$) {1};
\node[ext](e2) at ($(0,1.5)+(-60:2 and 1)$) {$2$};
\draw (e1) edge  ($(0,0)+(-90:2 and 1)$) (in.base) edge (e2);
\draw (e2) edge  ($(0,0)+(-60:2 and 1)$);
\end{scope}
\end{scope}

\end{tikzpicture}
\] 
\caption{\label{fig:pt1composition} An example of the (m)operadic composition of two elements in $\PT_1$.}
\end{figure}

The spaces $\PT_1(n)$ assemble to form a $\PT$-moperad. The operadic composition of two graphs $\Gamma_1, \Gamma_2 \in \PT_1$ is computed by the following algorithm:
\begin{enumerate}
 \item Delete the vertex $\vin$ of $\Gamma_1$ and the vertex $\vout$ of $\Gamma_2$. This leaves several open edges. 
 \item Connect the marked edge of $\vout$ of $\Gamma_2$ to the (single) edge of $\vin$ of $\Gamma_1$. 
 \item Reconnect the remaining open edges (previously attached to $\vout$ of $\Gamma_2$) in all planar possible ways to vertices of $\Gamma_1$.
\end{enumerate}
The labels on the edges are adjusted such that the edges stemming from $\Gamma_1$ have lower labels than those from $\Gamma_2$.
An example is shown in Fig. \ref{fig:pt1composition}

\begin{figure}
 \centering
\makeatletter{}\usetikzlibrary{matrix}
\usetikzlibrary{arrows}
\usetikzlibrary{shapes}
\usetikzlibrary{through}
\usetikzlibrary{calc,3d}
\usetikzlibrary{decorations,decorations.pathmorphing}
\begin{tikzpicture}[
int/.style={circle, draw, fill, minimum size=5pt, inner sep=0},
ext/.style={circle, draw, fill=white, minimum size=5pt, inner sep=1pt},
helper/.style={coordinate},point/.style={circle, draw, fill, inner sep =1pt},
de/.style={-triangle 60},
point/.style={circle, draw, fill, minimum size=3pt, inner sep=0pt},
xst/.style={cross out, draw, minimum size=4, inner sep=0 },
scale=.5,
grayfill/.style={fill=gray!30},
]
\begin{scope}[yshift=6cm]
\draw (0,0) ellipse (2cm and 1cm);
\draw (0,3) ellipse (2cm and 1cm);
\draw (-2,0)--(-2,3) (2,0)--(2,3);
\node [xst, label=-90:{$out$}] (out) at ($(0,0)+(-120:2 and 1)$) {};
\node [xst, label=90:{$in$}] (in) at ($(0,3)+(-70:2 and 1)$) {};
\node[ext] (e1) at ($(0,1.5)+(-70:2 and 1)$) {2};
\node[ext] (one) at ($(0,1)+(-120:2 and 1)$) {1};
\draw (in.base) edge (e1) (e1) edge  ($(0,0)+(-70:2 and 1)$) (out.base) edge (one);
\end{scope}

\node at (3.5,7.5) {\huge (};
\node at (10,7.5) {\huge )};
\node [int] at (4.5,7.5) {};
\node [int] at (6,7.5) {};

\node at (7.0235,7.3353) {\large ;};
\node at (5.25,7.25) {\large ,};
\draw(4.5,7.5) -- (4.5,8.5);
\draw(4.5,7.5) -- (4,6.5);
\draw(4.5,7.5) -- (5,6.5);
\draw(6,7.5) -- (6,8.5);
\draw(6,7.5) -- (6,6.5);
\draw(6,7.5) -- (6.5,6.5);
\draw(6,7.5) -- (5.5,6.5);

\draw(8.5,7.5) -- (8.5,8.5);
\draw(8.5,7.5) -- (8,6.5);
\draw(8.5,7.5) -- (9,6.5);
\draw(8.5,7.5) -- (7.5,8);
\draw(8.5,7.5) -- (9.5,8);
\node [int,grayfill] at (8.5,7.5) {};
\node [xst] at (8.5,8) {};

\begin{scope}[shift={(6.5,0)}]
\draw(8.5,7.5) -- (8,6.5);
\draw(8.5,7.5) -- (9,6.5);
\draw(8.5,7.5) -- (9.5,8);
\node [int,grayfill] (v1) at (8.5,7.5) {};
\node [xst] at (8.5,8) {};
\node [int] (v2) at (8.5,8.5) {};
\draw(8,9.5) -- (8.5,8.5) -- (9,9.5);
\node [int] (v3) at (7.5,8) {};
\draw(7,9) -- (7.5,8) -- (6.5,8.5);
\draw(7.5,8) -- (6.5,7.5);
\node [xst] at (7.254,8.5423) {};
\draw (v1) edge (v2);
\draw (v1) edge (v3);
\end{scope}

\begin{scope}[shift={(16.5,0)}]
\draw(8.5,7.5) -- (8,6.5);
\draw(8.5,7.5) -- (9,6.5);
\draw(8.5,7.5) -- (9.5,8);
\node [int,grayfill] (v1) at (8.5,7.5) {};
\node [xst] at (8.5,8) {};
\node [int] (v2) at (8.5,8.5) {};
\draw(8,9.5) -- (8.5,8.5) -- (9,9.5);
\node [int] (v3) at (7.5,8) {};
\draw(7,9) -- (7.5,8) -- (6.5,8.5);
\draw(7.5,8) -- (6.5,7.5);
\node [xst] at (6.9247,7.6956) {};
\draw (v1) edge (v2);
\draw (v1) edge (v3);
\end{scope}

\begin{scope}[shift={(6.5,-4)}]
\draw(8.5,7.5) -- (7.5,8.5);
\draw(8.5,7.5) -- (9,6.5);
\draw(8.5,7.5) -- (9.5,8);
\node [int,grayfill] (v1) at (8.5,7.5) {};
\node [xst] at (8.5,8) {};
\node [int] (v2) at (8.5,8.5) {};
\draw(8,9.5) -- (8.5,8.5) -- (9,9.5);
\node [int] (v3) at (8,7) {};
\draw(7,7.5) -- (8,7) -- (7,6.5);
\draw(8,7) -- (8,6);
\node [xst] at (7.3951,7.301) {};
\draw (v1) edge (v2);
\draw (v1) edge (v3);
\end{scope}

\begin{scope}[shift={(11,-4)}]
\draw(8.5,7.5) -- (7.5,8.5);
\draw(8.5,7.5) -- (9,6.5);
\draw(8.5,7.5) -- (9.5,8);
\node [int,grayfill] (v1) at (8.5,7.5) {};
\node [xst] at (8.5,8) {};
\node [int] (v2) at (8.5,8.5) {};
\draw(8,9.5) -- (8.5,8.5) -- (9,9.5);
\node [int] (v3) at (8,7) {};
\draw(7,7.5) -- (8,7) -- (7,6.5);
\draw(8,7) -- (8,6);
\node [xst] at (7.3951,6.6895) {};
\draw (v1) edge (v2);
\draw (v1) edge (v3);
\end{scope}

\begin{scope}[shift={(11.5,0)}]
\draw(8.5,7.5) -- (8,6.5);
\draw(8.5,7.5) -- (9,6.5);
\draw(8.5,7.5) -- (9.5,8);
\node [int,grayfill] (v1) at (8.5,7.5) {};
\node [xst] at (8.5,8) {};
\node [int] (v2) at (8.5,8.5) {};
\draw(8,9.5) -- (8.5,8.5) -- (9,9.5);
\node [int] (v3) at (7.5,8) {};
\draw(7,9) -- (7.5,8) -- (6.5,8.5);
\draw(7.5,8) -- (6.5,7.5);
\node [xst] at (6.9247,8.2836) {};
\draw (v1) edge (v2);
\draw (v1) edge (v3);
\end{scope}

\node at (11.5,7.5) {=};

\node at (17,8) {$\pm$};
\node at (22,4) {$\pm$};
\node at (17,4) {$\pm$};
\node at (12,4) {$\pm$};
\node at (22,8) {$\pm$};
\node at (23.5,4) {$\cdots$};
\end{tikzpicture}  
\caption{\label{fig:PT1actionexample} Schematic picture of the (right) action of some element of $\PT_1$ on the total space of a cyclic module over an operad. Here the $\PT_1$ element is drawn on the cylinder. It acts on two elements of the operad $\op P$ (notation as in example \ref{ex:PT1mod}) and one element of the module $\op M$. The operad elements are represented by two corollas, with 2 and 3 inputs (black dots). The element of the cyclic module is drawn as a gray corolla, with 5 inputs. The first one is marked by an ``$\times$''. Similarly, on the right hand side the terms occuring in the action are drawn. The operadic composition (or rather module action) is indicated just by connecting the corollas. A cyclic group action is performed so as to make the input indicated by $\times$ the first. }
\end{figure}
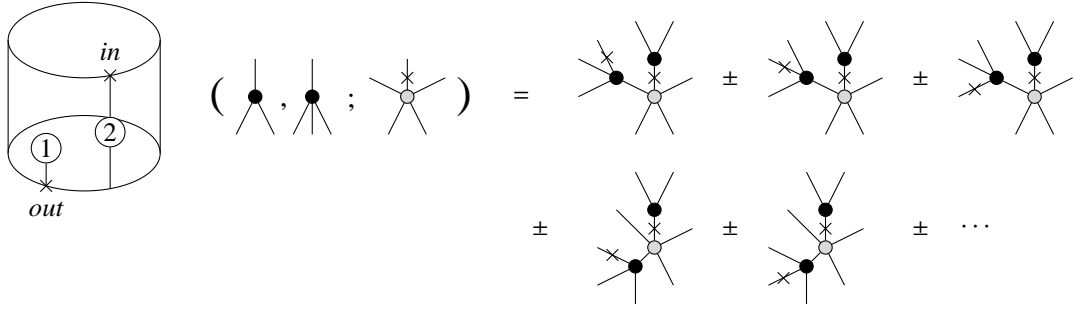

\begin{ex}
\label{ex:PT1mod}
 Recall from example \ref{ex:PTmod} that the total space of any non-symmetric operad $\op P$ is a $\PT$ module. 
Let $\op M$ be an operadic right cyclic module over $\op P$. By this we mean a collection of vector spaces $\op M(n)$, with right actions of the cyclic groups $\Z_{n}$, and with composition morphisms
\[
 \op M(n)\otimes \op P(m_1)\otimes \cdots \otimes \op P(m_n) \to \op M(m_1+\dots + m_n)
\]
satisfying some straightforward axioms.
We claim that in this situation there is a moperadic \emph{right} action of $\PT_1$ on the total space
\[
 M := \prod_n \op M(n).
\]
An example is given in Figure \ref{fig:PT1actionexample}, from which the principle should be clear. Writing down a general formula is left to the reader.
As a special case, consider $\op P$ being the operad whose total space is $\Dpoly$ (with zero differential), and $\op M$ the cyclic module 
\[
 \op M(n) = \Hom(A^{\otimes n+1},\R)
\]
where $A=C^\infty(M)$.
This is the (dual of the) Hochschild chain complex, but considered with zero differential for the moment being. Similarly, one obtains a \emph{left} action of $\PT_1$ on Hochschild chains. We refer to \cite{KS2} for a (slightly) more explicit description. For example, the element
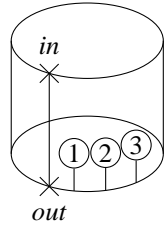
\begin{center}
\begin{tikzpicture}[
ext/.style={circle, draw, fill=white, minimum size=5pt, inner sep=1pt},
xst/.style={cross out, draw, minimum size=5 },
scale=.5
]
\begin{scope}[]
\draw (0,0) ellipse (2cm and 1cm);
\draw (0,3) ellipse (2cm and 1cm);
\draw (-2,0)--(-2,3) (2,0)--(2,3);
\node [xst, label=-90:{$out$}] (out) at ($(0,0)+(-120:2 and 1)$) {};
\node [xst, label=90:{$in$}] (in) at ($(0,3)+(-120:2 and 1)$) {};
\draw (in.base)--(out.base) ;
\node[ext] (e1) at ($(0,1)+(-100:2 and 1)$) {1};
\draw (e1) edge  ($(0,0)+(-100:2 and 1)$);
\node[ext] (e2) at ($(0,1)+(-76:2 and 1)$) {2};
\draw (e2) edge  ($(0,0)+(-76:2 and 1)$);
\node[ext] (e3) at ($(0,1)+(-50:2 and 1)$) {3};
\draw (e3) edge  ($(0,0)+(-50:2 and 1)$);
\end{scope}
\end{tikzpicture}
\end{center}
acts on 3 multidifferential operators $D_1, D_2, D_3\in \Dpoly$ and on a Hochschild chain $a_0\otimes a_1\otimes \cdots \otimes a_n\in C_\bullet$ as
\begin{multline*}
\sum_{\substack{0\leq i_1 , i_2, i_3 \\ 
i_1+i_2+i_3\leq n-k_1-k_2-k_3}}
 \pm  a_0\otimes a_1\otimes \cdots \otimes
D_1(a_{i_1+1},\dots, a_{i_1+k_1})\otimes 
a_{i_1+k_1+1} \otimes \cdots 
\\
\cdots \otimes 
 D_3(a_{i_1+i_2+i_3+k_1+k_2+k_3+1}, \dots ) \otimes  \cdots \otimes a_n.
\end{multline*}
Here $k_j$ is the degree (``number of slots'') of $D_j$, $j=1,2,3$.
\end{ex}

\begin{ex}
\label{ex:ESCtoPT1mod}
Suppose we have a three-colored operad $\op P$ of Extended Swiss Cheese type in differential graded vector spaces, see definition \ref{def:ESCtype}.
Let (as in example \ref{ex:SCtoPTmod}) $\op P^1$ be the operad of operations with output of color 1. Let $\op P^2$ be the space of operations with output in color two. As in example \ref{ex:SCtoPTmod}
we obtain a two colored operad
\[
\bpm 
\PT & 
\prod_n \op P^2(\cdot, n, 0)[-n]
& \op P^1
\epm\, .
\]
By the cyclic action on the color 2 inputs of $\op P^3(\cdot, \cdot,0)$. Then we can make the total space
\[
\prod_n \op P^3(\cdot, n+1,0)[-n]
\] 
 into a moperadic bimodule along the lines of the previous example.
 In other words, we obtain a four colored operad 
 \[
 \bpm
 \PT & \prod_n \op P^2(\cdot, n, 0)[-n]
& \op P^1
\\
 \PT_1 & \prod_n \op P^3(\cdot, n+1, 0)[-n]
& \op P^3(\cdot, 0,1) 
 \epm \, .
 \]
 Since we used only ``intrinsic'' operations to construct the moperadic bimodule structure, the assignment from Extended Swiss Cheese type operads to four-colored operads is functorial.
\end{ex}

\subsubsection{A variation: The moperad \texorpdfstring{$\PT_1^{\mathbb{1}}$}{PT11}}
Consider a $\PT$-moperad $\PT_1^{\mathbb{1}}$ defined in the same manner as $\PT_1$, except that
\begin{enumerate}
 \item The planar trees that generate $\PT_1^{\mathbb{1}}$ may contain another type of vertex, which we call ``unit vertex'' and designate by $\mathbb{1}$ in pictures.
 \item Those unit vertices must not have any children.
 \item The composition law is defined in the same manner as before, except that graphs in which unit vertices acquire children are considered zero.
\end{enumerate}

A typical element of $\PT_1^{\mathbb{1}}$ may look like this:
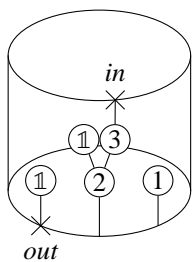
\[
\begin{tikzpicture}[
int/.style={circle, draw, fill, minimum size=5pt, inner sep=0},
ext/.style={circle, draw, fill=white, minimum size=5pt, inner sep=1pt},
helper/.style={coordinate},point/.style={circle, draw, fill, inner sep =1pt},
de/.style={-triangle 60},
point/.style={circle, draw, fill, minimum size=3pt, inner sep=0pt},
xst/.style={cross out, draw, minimum size=5pt},
scale=.6
]
\begin{scope}[xshift=6cm]
\draw (0,0) ellipse (2 and 1);
\draw (0,3) ellipse (2 and 1);
\draw (-2,0)--(-2,3) (2,0)--(2,3);
\node [xst, label=-90:{$out$}] (out) at ($(0,0)+(-130:2 and 1)$) {};
\node [xst, label=90:{$in$}] (in) at ($(0,3)+(-80:2 and 1)$) {};
\draw (out)+(0,1) node[ext] (e3) {$\mathbb{1}$};
\node[ext] (e2) at ($(0,1.2)+(-90:2 and 1)$) {$2$};
\node[ext] (e1) at ($(0,1)+(-50:2 and 1)$) {$1$};
\draw (e2)+(70:1) node[ext] (e5) {$3$} 
+(110:1) node[ext](e4) {$\mathbb{1}$} 
(e5)+(90:1);
\draw (out.base) edge (e3);
\draw (e2) edge (e5) edge (e4) (e5) edge (in.base);
\draw (e2) edge (0,-1) (e1) edge ($(0,0)+(-50:2 and 1)$);
\end{scope}
\end{tikzpicture}
\]

A representation of $\PT_1^{\mathbb{1}}$ is the same as a representation of $\PT_1$, except that there is in addition a singled out ``zero-ary'' element, which we designate $\mathbb{1}$ in formulas. 

\begin{ex}
\label{ex:PT11mod}
Consider the example \ref{ex:PT1mod} of a representation of $\PT_1$. We can make it into a representation of $\PT_1^{\mathbb{1}}$ by specifying some element $\mathbb{1}\in \op P(0)$. It does not have to satisfy any relations for now. The action of any $\PT_1^{\mathbb{1}}$ tree $\Gamma$ is obtained by ``inserting $\mathbb{1}$ for any unit vertex'', or, to be more concrete:
\begin{enumerate}
 \item Consider the $\PT_1$ tree $\Gamma'$ obtained by making the unit vertices into numbered vertices, numbering them in an arbitrary way.
 \item The action of $\Gamma$ is the same as that of $\Gamma'$, with the slots corresponding to unit vertices being filled by copies of $\mathbb{1}\in \op P(0)$.
\end{enumerate}

\end{ex}

\subsection{\texorpdfstring{$\KS_1$}{KS1} moperad and the colored operad \texorpdfstring{$\KS$}{KS}}
\label{sec:KS1}
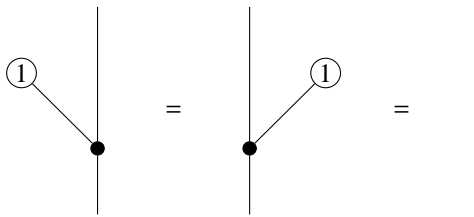
\begin{figure}
 \centering
\makeatletter{}\usetikzlibrary{matrix}
\usetikzlibrary{arrows}
\usetikzlibrary{shapes}
\usetikzlibrary{through}
\usetikzlibrary{calc,3d}
\usetikzlibrary{decorations,decorations.pathmorphing}
\[
\begin{tikzpicture}[
int/.style={circle, draw, fill, minimum size=5pt, inner sep=0},
ext/.style={circle, draw, fill=white, minimum size=5pt, inner sep=1pt},
helper/.style={coordinate},point/.style={circle, draw, fill, inner sep =1pt},
de/.style={-triangle 60},
point/.style={circle, draw, fill, minimum size=3pt, inner sep=0pt},
xst/.style={cross out, draw, minimum size=5 },
]
\begin{scope}[xshift=0cm]
\node (root) at (0,0) {};
\node[int] (i1) at (0,1) {};
\node (nix) at (0,3) {};
\node[ext] (one) at (-1,2) {$\mathrm{1}$};
\draw (root) edge (i1) (i1) edge (one) edge (nix);
\node at (1,1.5) {$=$};
\end{scope}
\begin{scope}[xshift=2cm]
\node (root) at (0,0) {};
\node[int] (i1) at (0,1) {};
\node (nix) at (0,3) {};
\node[ext] (one) at (1,2) {$\mathrm{1}$};
\draw (root) edge (i1) (i1) edge (one) edge (nix);
\node at (2,1.5) {$=$};
\end{scope}
\begin{scope}[xshift=5cm]
\node (root) at (0,0) {};
\node (nix) at (0,3) {};
\draw (root) edge (nix);
\end{scope}
\end{tikzpicture}
\] 
\caption{\label{fig:unitrelation} The unit relation in $\KS_1$.}
\end{figure}

Let us twist the $\PT$-moperad $\PT_1^{\mathbb{1}}$ (see Appendix \ref{sec:optwists} for the definition of moperadic twisting). From the twisting, we get a $\Tw\PT$-moperad $\Tw\PT_1$.
In particular it is a $\Br$-moperad. $\Tw\PT_1$ is spanned by planar trees as in the previous subsection, except that some of the numbered vertices may be replaced by internal vertices. The differential splits vertices, creating a new internal vertex, similar to the operation depicted in \ref{fig:KSdiff}.
It is clear that the subspace $\KS_1'$ of $\Tw\PT_1$ formed by graphs all of whose internal vertices have $\geq 2$ children is a sub-$\Br$-moperad. To define $\KS_1$, we will take its quotient with respect to the following relations:
\begin{enumerate}
 \item Graphs which contain a unit vertex whose parent is an internal vertex with three or more children are set to zero.
 \item If a graph contains an internal vertex with two children, one of which is a unit vertex, this graph is set equal to a graph obtained as follows: Remove the unit and the internal vertex, and connect the two dangling edges remaining, see Figure \ref{fig:unitrelation}.
 \item Graphs which contain a unit vertex whose parent is any numbered vertex are set to zero.
 \item Graphs which contain a unit vertex whose parent is $\vout$ and whose adjacent edge is not marked are considered zero.
\end{enumerate}

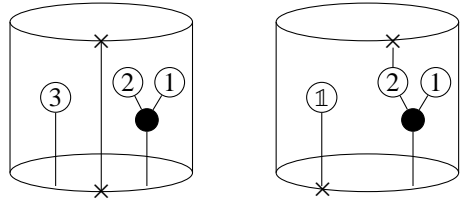
\begin{figure}
\centering
\tikzset{ext/.style={circle, draw,inner sep=1pt},int/.style={circle,draw,fill,inner sep=1pt},nil/.style={inner sep=1pt}}
\tikzset{exte/.style={circle, draw,inner sep=3pt},inte/.style={circle,draw,fill,inner sep=3pt}}
\[
\begin{tikzpicture}[scale=1]
\begin{scope}[xshift=2.7cm]
\draw (0,1) ellipse (1.2cm and .25cm);
\draw (0,-1) ellipse (1.2cm and .25cm);
\draw (-1.2,1)--(-1.2,-1) (1.2,1)--(1.2,-1) (0,.75)--(0,-1.25);
\draw (0,.75) node (in) {$\times$};
\node[outer sep =0, inner sep=0] (ed) at (0,-1.25) {$\times$}; \node (e1) at (.6,-.3) [inte] {};
\node (e1c) at (.6,-1.3) {};
\node (e2) at (.9,.2) [ext] {1};
\node (e3) at (.35,.2) [ext] {2};
\node (e4) at (-.6,0) [ext] {3};
\node (e4c) at (-.6,-1.3) {};
\draw (e1c)--(e1)--(e2) (e1)--(e3) (e4c)--(e4);
\end{scope}
\begin{scope}[xshift=2.7cm, shift={(3.5,0)}]
\draw (0,1) ellipse (1.2cm and .25cm);
\draw (0,-1) ellipse (1.2cm and .25cm);
\draw (-1.2,1)--(-1.2,-1) (1.2,1)--(1.2,-1);
\node[outer sep =0, inner sep=0] (in) at (0.3377,0.7484) {$\times$};
\node[outer sep =0, inner sep=0] (ed) at (-0.5896,-1.2188) {$\times$};
\node [inte] (e1) at (0.6,-0.3) {};
\node (e1c) at (0.6,-1.3) {};
\node [ext] (e2) at (0.9,0.2) {1};
\node [ext] (e3) at (0.3354,0.1999) {2};
\node [ext] (e4) at (-0.6,0) {$\mathbb{1}$};
\node (e4c) at (-0.6,-1.3) {};
\draw (in.base) edge (e3);
\draw (e1c)--(e1)--(e2) (e1)--(e3) (e4c)--(e4);
\end{scope}
\end{tikzpicture}
\]
\caption{\label{fig:KS1ex} Two random elements of $\KS_1$.} 
\end{figure}

One can check that the resulting space, $\KS_1$, is still a $\Br$-moperad. An example of a graph in $\KS_1$ is shown in Figure \ref{fig:KS1ex}
The Kontsevich-Soibelman operad is the colored operad
\[
 \KS = \bpm \Br & \KS_1 \epm.
\]
Here we use the notation from section \ref{sec:opconventions} to denote the colored operad formed by an operad and a moperad.

\begin{ex}
\label{ex:KSactonDpolyCbullet}
Recall from examples \ref{ex:PT1mod}, \ref{ex:PT11mod} the action of $\PT_1^{\mathbb{1}}$ on Hochschild cochains and chains of an $A_\infty$-algebra, both considered with zero differential. Twisting by the Maurer Cartan element provided by the $A_\infty$-structure, we obtain an action of $\Tw\PT_1$ on Hochschild (co-)chains, now considered with the Hochschild differential. The first and second of the above relations now encode that the singled out element $\mathbb{1}$ is a strong unit for the $A_\infty$-product.
The third relation can be satisfied by considering the \emph{normalized} Hochschild cochain complex instead of the ordinary one. The fourth relation can be satisfied by considering the normalized (dual of the) Hochschild chain complex. 
 \end{ex}

\subsection{\texorpdfstring{$\SGra$}{SGra} operadic bimodule}
\label{sec:SGra}
Let the graded space\footnote{The name stands for ``Swiss Cheese Graphs'', in analogy with the ``Swiss Cheese'' operad.} 
$\fSGra(m,n)'_k$ be spanned by directed graphs with vertex set $[m]\sqcup \{\bar 1, \dots, \bar n\}$ and edge set $[k]$, such that none of the vertices in the set $[\bar n]:=\{\bar 1, \dots, \bar n\}$ have any outgoing edges.
Let us define
\[
 \fSGra(m,n) = \prod_{k\geq 0} (\fSGra(m,n)'_k \otimes sgn_k)_{S_k}\, .
\]
We call the vertices in the set $[m]$ ``type I vertices'' and those in the second set $[\bar n]$ ``type II vertices''. The graphs are assigned degree $-k$, so an edge has again degree $-1$. Together with the operad $\dGra$ the spaces $\fSGra(m,n)$ form a two colored operad $\SG$. The operadic composition operations are defined in the same way as those for $\dGra$. We will consider $\SG$ as a partially non-symmetric operad, without an action of the permutation group on the inputs of the second color.\footnote{Although there is an obvious such action.} According to example \ref{ex:SCtoPTmod} the collection of spaces
\[
\fSGra(m) = \prod_n \fSGra(m,n)[-n] 
\]
form a $\PT$-$\dGr$-operadic bimodule. 
Note also that $\fSGra$ is as well a 
$\PT$-$\Gra$-bimodule using the canonical map $\Gra\to \dGra$.

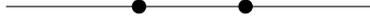
\begin{figure}
 \centering
\makeatletter{}\[
\begin{tikzpicture}[scale=.7,
vert/.style={draw,outer sep=0,inner sep=0,minimum size=5,shape=circle,fill},
helper/.style={outer sep=0,inner sep=0,minimum size=5,shape=coordinate},
default_edge/.style={draw},
ext/.style={draw,outer sep=0,inner sep=0,minimum size=5,shape=circle},
every loop/.style={}]

\node (v0) at (1,6) [helper] {};
\node (v1) at (8,6) [helper] {};
\node (v3) at (3.5,6) [vert] {};
\node (v4) at (5.5,6) [vert] {};

\draw[default_edge] (v0) to (v1);
\end{tikzpicture}
\]
 
\caption{\label{fig:sgramc} The Maurer Cartan element used to twist the left operadic $\PT$-module structure on $\SGra$.}
\end{figure}

We can twist the $\PT$-$\dGr$-operadic bimodule structure to a $\Br$-$\dGr$-operadic bimodule structure on $\fSGra$. According to the general theory of operadic twisting
in Appendix \ref{sec:optwists} this means that one should specify a Maurer-Cartan element (in a sense made precise there) in $\SGr(0)$. 
We take the element given by the graph with two type II vertices, and no edges, depicted in Figure \ref{fig:sgramc}. The differential on the twisted bimodule contains an additional term from the twisting. It is given given by splitting type II vertices. Pictorially it looks like this:
\[
\makeatletter{}\begin{tikzpicture}[
point/.style={circle, draw, fill, minimum size=3pt, inner sep=0pt},
]
\begin{scope}[]
\node[point] (v0) at (-.6, 0) {};
\node[point] (v1) at (-.2, 0) {};
\node[point] (v2) at (.5, 0) {};
\draw (-1,0) -- (1,0);
\foreach \x in {-.3,0,.3}
{
\draw (v1)--+(\x , .7);
\draw (v2)--+(\x , .7);
}
\end{scope}
\begin{scope}[shift={(3,0)}]
\node[draw, dashed, semicircle, minimum size=7] (v1) at (-0.5,0.12) {};
\node[point] (v2) at (0.5,0) {};
\node[point] (v1a) at (-0.6,0) {};
\node[point] (v1b) at (-0.4,0) {};
\draw (-1,0) -- (1,0);
\foreach \x in {-.3,0,.3}
{
\draw (v1)--+(\x , .7);
\draw (v2)--+(\x , .7);
}
\end{scope}
\begin{scope}[shift={(-4,0)}]
\node[point] (v1) at (-0.5,0) {};
\node[point] (v2) at (0.5,0) {};
\draw (-1,0) -- (1,0);
\foreach \x in {-.3,0,.3}
{
\draw (v1)--+(\x , .7);
\draw (v2)--+(\x , .7);
}
\end{scope}
\begin{scope}[shift={(3,-1)}]
\node[point] (v0) at (0.5,0) {};
\node[point] (v1) at (-0.5,0) {};
\node[point] (v2) at (0.2,0) {};
\draw (-1,0) -- (1,0);
\foreach \x in {-.3,0,.3}
{
\draw (v1)--+(\x , .7);
\draw (v2)--+(\x , .7);
}
\end{scope}
\begin{scope}[shift={(0,-1)}]
\node[draw, dashed, semicircle, minimum size=7] (v2) at (0.5,0.12) {};
\node[point] (v1) at (-0.5,0) {};
\node[point] (v2a) at (0.6,0) {};
\node[point] (v2b) at (0.4,0) {};
\draw (-1,0) -- (1,0);
\foreach \x in {-.3,0,.3}
{
\draw (v1)--+(\x , .7);
\draw (v2)--+(\x , .7);
}
\end{scope}
\node at (-5.5,0) {$\delta$};
\node at (-2,0) {$=$};
\node at (-1.5,-1) {$+$};
\node at (1.5,0) {$-$};
\node at (1.5,-1) {$-$};
\end{tikzpicture} \, .
\]
Here the dashed semicircles shall indicate that on sums over all graphs obtained by connecting the incoming edges to vertices inside the semicircle in some way. In our example, each graph with a dashed semicircle hence stands for a sum of $2^3=8$ terms.

\begin{defprop}
The twisted version of the $\Br$-$\dGr$ bimodule $\fSGra$ contains an operadic sub-bimodule $\SGra$ spanned by graphs such that each type II vertex has at least 1 incoming edge.
\end{defprop}
\begin{proof}
One has to show that the $\SGra$ thus defined is closed under (i) the differential (ii) the right $\dGr$ action and (iii) the left $\Br$ action. Since the right $\dGr$ action does not affect the type II vertices, (ii) is clear. Let us consider the actions of generators of $\Br$, see figure \ref{fig:KSgen}. The action of $T_n$ leaves invariant or increases the valence of type II vertices, so it maps $\SGra$ to itself. $T_n'$ acts as zero unless $n=2$. In that case it does not affect the valences of type II vertices.
Hence statement (iii) is shown. 
Note that here it is essential to work with $\Br$ instead of the larger operad $\Tw\PT$. The latter would contain an operation $T_1'$ which would not map $\SGra$ to itself.
Finally consider (i). This statement is similar to the statement that the normalized Hochschild cochain complex of an algebra is a subcomplex of the full (non-normalized) Hochschild complex. More concretely, the differential has the form
\[
\delta \Gamma = T_1(m, \Gamma) \pm T_1(\Gamma, m).
\]
The right hand term splits each type II vertex into two, thus potentially producing two graphs with a valence 0 type II vertex. However, these graphs cancel among neighboring type II vertices. The remaining rightmost and leftmost terms just kill the two terms contributed by $T_1(m, \Gamma)$.
\end{proof}

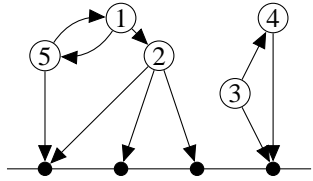
\begin{figure}
 \centering
\makeatletter{}\begin{tikzpicture}[
scale=.5,
int/.style={circle, draw, fill, minimum size=5pt, inner sep=0},
ext/.style={circle, draw, fill=white, minimum size=5pt, inner sep=1pt},
xst/.style={cross out, draw, minimum size=5pt, inner sep=1pt},
point/.style={circle, draw, fill, minimum size=5pt, inner sep=0pt},
arr/.style={-triangle 45},
]
\draw (-5,0)--(3,0) node (v8) {};

\node [ext] (v3) at (-2,4) {1};
\node [ext] (v4) at (-1,3) {2};
\node [ext] (v6) at (1,2) {3};
\node [ext] (v7) at (2,4) {4};
\node [ext] (v1) at (-4,3) {5};
\node [point] (v2) at (-4,0) {};
\node [point] (v5) at (-2,0) {};
\node [point] (v10) at (0,0) {};
\node [point] (v9) at (2,0) {};
\draw[arr] (v1) edge (v2);
\draw[arr] (v1) edge[bend left] (v3);
\draw[arr] (v3) edge[bend left] (v1);
\draw[arr] (v4) edge (v2);
\draw[arr] (v4) edge (v5);
\draw[arr] (v6) edge (v7);
\draw[arr] (v7) edge (v9);
\draw[arr] (v6) edge (v9);
\draw[arr] (v3) edge (v4);
\draw[arr] (v4) edge (v10);
\end{tikzpicture} 
\caption{\label{fig:SGrasample} Some element of the  $\Br-\dGr$-bmodule $\SGra$.}
\end{figure}

An example of an element in $\SGra$ can be found in Figure \ref{fig:SGrasample}.

\begin{rem}
The spaces $\SGr(m)$ are the spaces of Kontsevich graphs. In fact, they have more justly been denoted $\mathsf{KGra}$ in \cite{vasilystable}.
\end{rem}

\begin{ex}
\label{ex:SgraactTpolyDpoly}
 The operadic bimodule $\SGra$ can be represented on the colored vector space $\Tpoly \oplus \Dpoly$. This means that there are maps 
\[
 \SGra(m) \to \Hom(\Tpoly^{\otimes m}, \Dpoly) 
\]
in a way compatible with the left and right actions of $\Br$ and $\Gra$. This map sends a graph $\Gamma\in \SGra(m)$ to the map
\[
(\gamma_1, \dots, \gamma_m)\mapsto D_\Gamma(\gamma_1, \dots, \gamma_n)
\]
with the multidifferential operator $D_\Gamma(\cdots)$ as defined in \cite{K1}. Concretely, we can naturally identify the graph $\Gamma$ (with, say, $n$ type II vertices) with an element of $\dGra(m+n)$. Then the multidifferential operator on the left is defined such that for functions $a_1,\dots, a_n$
\[
D_\Gamma(\gamma_1, \dots, \gamma_n)(a_1,\dots, a_n)
= 
\Gamma(\gamma_1, \dots, \gamma_n, a_1,\dots, a_n).
\]
The action on the right is defined in example \ref{ex:GraactonTpoly} and we consider functions as zero-vector fields.
\end{ex}

\subsection{\texorpdfstring{$\SGra_1$}{SGra1} moperadic bimodule}
\label{sec:sgra1}
Let us define the subspace $\fSGra_1(m,n+1)\subset \fSGra(m+1,n+1)$ spanned by all graphs with no incoming edges to the first type I vertex. We will denote the type $II$ vertices by $\bar 0,\dots, \bar n$.

The spaces $\fSGra_1(m,n+1)$, together with the operad $\dGra$, the moperad $\Gra_1$ and $\fSGra(\cdot,\cdot)$, assemble to form a three colored operad of Extended Swiss Cheese type (see Definition \ref{def:ESCtype}).
Concretely, the action of the cyclic group of order $n+1$ on $\fSGra_1(m,n+1)$ is by cyclically permuting the type II vertices. 
The right $\dGr$-action is by insertion at type I vertices (except the first). The action of $\Gr_1$ is by insertion at the first type I vertex, similar to the composition on $\Gr_1$. The composition of an element in $\fSGra_1(m,n)$ with elements of $\fSGra(\cdot,\cdot)$ is defined by inserting the elements of $\fSGra(\cdot,\cdot)$ into type II vertices and reconnecting the incident edges in all possible ways.

Being part of an Extended Swiss Cheese type operad, we can invoke example \ref{ex:ESCtoPT1mod} to obtain a $\PT$-$\dGr$-$\PT_1$-$\dGr_1$-$\fSGra$ moperadic bimodule structure on 
the total space
\[
 \fSGra_1(m) = \prod_{n\geq 0} \fSGra_1(m,n+1)[-n].
\]
The action of $\PT_1$ can be upgraded to an action of $\PT_1^\mathbb{1}$ along the lines of example \ref{ex:PT11mod}. Here, whenever the element $\mathbb{1}$ is ``inserted'' at a type II vertex $\bar k$ of some graph $\Gamma$, it acts (i) as zero if $\bar k$ has valence $\geq 1$ or (ii) by forgetting the vertex $\bar k$ is it has valence zero, and relabelling the other type II vertices accordingly.

We can twist $\fSGra_1$ together with the operadic bimodule $\fSGra$ and $\bpm \PT & \PT_1^\mathbb{1}\epm$. According to Appendix \ref{sec:optwists} no additional data is needed, on top of the chosen Maurer-Cartan element from figure \ref{fig:sgramc}.
This in particular makes $\fSGra_1$ into a $\Tw\PT$-$\dGra$-$\Tw\PT_1^\mathbb{1}$-(twisted version of) $\fSGra$ moperadic bimodule. 
The twisted differential is the graphical version of the Hochschild differential on the dual space of the Hochschild chain complex. It is given by splitting each type II vertex into two type II vertices, with alternating signs. 

We can restrict the moperadic bimodule structure to $\Br\subset \Tw\PT$ and $\SGra\subset\fSGra$. However, we also want to pass from $\Tw\PT_1^\mathbb{1}$ to its subquotient $\KS_1$. For this we have to check several relations, cf. relations (1)-(4) of section \ref{sec:KS1}. The relations coming from relations (1), (2), (3) are easily checked to hold. However, (4) does not, and requires us to pass to a subspace.

\begin{defprop}
The moperadic bimodule structure on the twisted version of  $\fSGra_1$ descend to a the $\Br$-$\dGr$-$\KS_1$-$\Gra_1$-$\SGra$ moperadic bimodule structure on a subspace $\SGra_1$, which 
is spanned by graphs such that the type II vertices $\bar 1, \dots, \bar n$ have at least 1 incoming edge. (Vertex $\bar 0$ is still allowed to have valence 0.)
\end{defprop}
\begin{proof}
On $\SGra_1$ it is clear that the relation corresponding to relation (4) of section \ref{sec:KS1} holds. However, one still has to show that (i) the subspaces $\SGra_1$ are closed under the differential, (ii) they are closed under the right $\dGr$-action, (iii) they are closed under the left $\Gra_1$ action and (iv) they are closed under the (combined) right $\KS_1$-$\SGra$ action.
Statements (ii) and (iii) are immediate because the corresponding actions cannot decrease the valence of type II vertices. Consider statement (i). The differential splits type II vertices, producing 2 graphs with a valence zero type II vertex other than $\bar 0$. But these graphs cancel in pairs corresponding to neighbouring type II vertices. Hence (i) follows. Finally consider statement (iv). The $\KS_1$ action is built using three ``fundamental'' operations: (a) forgetting vertex $\bar 0$ if it has valence 0, (b) inserting some elements of $\SGra$ into type II vertices and (3) cyclically relabelling type II vertices. Operation (a) clearly does not affect valences of the other 
type II vertices. Since all type II vertices of elements of $\SGra_1$ have valence $\geq 1$ by definition, operation (b) can neither introduce type II vertices of valence 0. Operation (c) could, if the graph it is applied to has a valence $0$ vertex $\bar 0$. However, inspecting the $\KS_1$ operations one sees that whenever the relabelling occurs vertex $\bar 0$ is either forgotten or some $\SGra$ element is inserted. Hence statement (iv) holds as well.
\end{proof}

\begin{figure}
 \centering
\makeatletter{}\tikzset{ext/.style={circle, draw,inner sep=1pt},int/.style={circle,draw,fill,inner sep=1pt},nil/.style={inner sep=1pt}}
\tikzset{exte/.style={circle, draw,inner sep=3pt},inte/.style={circle,draw,fill,inner sep=3pt}, xst/.style={cross out, draw, minimum size=5pt},}
\usetikzlibrary{arrows}
\begin{tikzpicture}[every edge/.style={draw, -triangle 60}]
\draw (0,0) ellipse (2.5 and 2.5);

\node [xst,label=180:{$out$}] (ct) at (0,0) {};
\node [int,label=0:{$in=\bar 0$}] (one) at (0:2.5) {};
\node [int,label=20:{$\bar 1$}] (i1) at (20:2.5) {};
\node [int,label=200:{$\bar 3$}] (i2)at (200:2.5) {};
\node [int,label=-90:{$\bar 4$}] (i3)at (280:2.5) {};
\node [int,label=90:{$\bar 2$}] (i4)at (100:2.5) {};

\node [ext] (e1) at (-20:1.5) {1};
\node [ext] (e2) at (200:1.5) {2};
\node [ext] (e3) at (280:1) {3};
\node [ext] (e4) at (100:1) {4};

\draw (ct.base) edge (e4);
\draw (ct.base) edge (i1);
\draw (ct.base) edge (e1);
\draw (e1) edge (one);
\draw (e1) edge (i1);
\draw (ct) edge (e3);
\draw (e3) edge (i3);
\draw (e3) edge (e2);
\draw (e2) edge (i2);
\draw (e4) edge (i4);
\draw (e4) edge (e2);
\end{tikzpicture} 
\caption{\label{fig:SGra1sample} An example element of $\SGra_1$. Note that by convention we draw the graph in a circle, and label the first type I $\vout$ and the first type II vertex by $\bar 0$ or $\vin$.}
\end{figure}
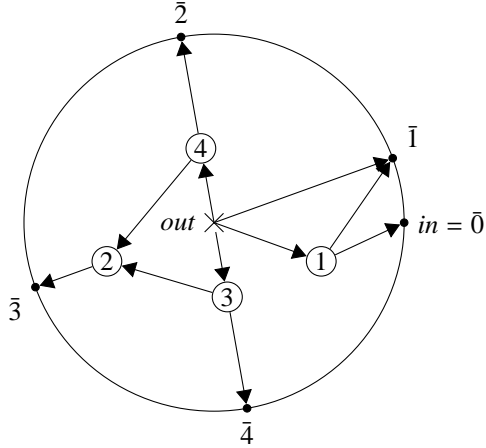

An example of a graph in $\SGra_1$ is shown in Figure \ref{fig:SGra1sample}.

\begin{ex}
\label{ex:Sgra1actCOmega}
 The moperadic bimodule $\SGra_1$ can be represented on the colored vector space $\Dpoly\oplus C_\bullet \oplus \Tpoly\oplus \Omega_\bullet$. Concretely its elements yield operations with 1 input in $C_\bullet$, zero or more inputs in $\Tpoly$ and the output in $\Omega_\bullet$. Concrete formulas for how to associate a graph with such an operation can be found in \cite{shoikhet}.
\end{ex}

\begin{rem}
 The graphs considered above are B. Shoikhet's graphs \cite{shoikhet}.
\end{rem}

\subsection{The four colored operad \texorpdfstring{$\bigGra$}{bigGra}}
\label{sec:bigGra}
Now we can define the four-colored operad $\bigGra$ as follows, using the notation of section \ref{sec:opconventions}.
\[
 \bigGra := \bpm \Br & \SGra & \Gra \\ \KS_1 & \SGra_1 & \Gra_1 \epm  
\]

\begin{ex}
By combining the actions of the various parts of $\bigGra$ from the previous subsections, we obtain an action of $\bigGra$ on the (colored) vector space 
\[
\bigV=\Dpoly\oplus  C_\bullet \oplus \Tpoly \oplus \Omega_\bullet \, .
\]
\end{ex}
 
\makeatletter{}\section{Several topological operads and \texorpdfstring{$\bigChains$}{bigChains}}
\label{sec:topops}
Recall that we want to represent the 4-colored operad $\homKS$ on the 4-colored vector space $\bigV=\Dpoly\oplus  C_\bullet \oplus \Tpoly \oplus \Omega_\bullet$. At the end of the last section we saw that there is a natural action of the 4-colored operad $\bigGra$ on $V$. If we can map $\homKS$ to $\bigGra$ we are hence done. We will do this via an intermediary 4-colored operad $\bigChains$. The goal of the present section is to define this operad.

\begin{rem}
We use the adjective ``topological'' in this section several times, but in fact all our operads will be built of semi-algebraic manifolds \cite{HLTV}, not just topological spaces. This is important for technical reasons since later on we want to integrate differential forms on them.
\end{rem}

\subsection{Topological operad \texorpdfstring{$\FM_2$}{FM2}}
\label{sec:fm2op}
The operad $\FM_2$ introduced by Getzler and Jones \cite{GJ} is a compactification of the space of configurations of points in the plane. For $n\geq2$:
\[
 \FM_2(n) = \left( \{(z_1,\dots z_n)\in \C^n \mid z_i\neq z_j \text{for $i\neq j$} \}/ \R_+ \ltimes \C \right)^-.
\]
Here the final superscript denotes compactification and the group $\R_+ \ltimes \C$ acts by overall scaling and translation. For more details, in particular regarding the compactification, we refer the reader to \cite{GJ,KS2}. For our convenience, we will set $\FM_2(1)=\{pt\}$ to be the space consisting of a single point, so that the operad $\FM_2$ has a unit. 
\begin{ex}
The space $\FM_2(2)\cong S^1$ is a circle. An explicit map $\FM_2(2)\to S^1$ is given by sending a configuration $[(z_1,z_2)]$ to the point
$\frac{z_1-z_2}{|z_1-z_2|}$.
\end{ex}

There will be two series of subspaces of $\FM_2$, which will be important later.
First, for $n\geq 2$, define $U_I^n\subset \FM_2(n)$ as follows.
\[
 U_I^n = \cap_{1\leq i<j\leq n} \pi_{ij}^{-1}( {1} )
\]
Here $\pi_{ij}\colon \FM_2(n) \to \FM_2(2) \cong S^1$ is the forgetful map forgetting all points except the $i$-th and $j$-th. The notation $\pi_{ij}^{-1}( {1} )$ shall mean the preimage of $1\in S^1$.\footnote{We consider $S^1$ as embedded into $\C$, i.~e., as $\{z\in \C \mid |z|=1 \}$. }

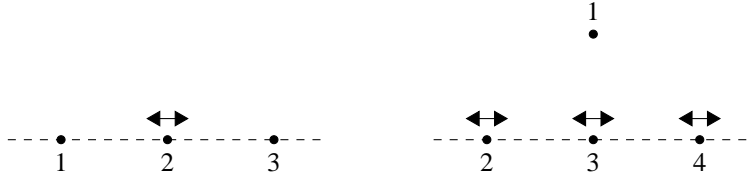
\begin{figure}
 \centering
\makeatletter{}\usetikzlibrary{matrix}
\usetikzlibrary{arrows}
\usetikzlibrary{shapes}
\usetikzlibrary{through}
\usetikzlibrary{calc,3d}
\usetikzlibrary{decorations,decorations.pathmorphing}
\begin{tikzpicture}[
scale=.7,
int/.style={circle, draw, fill, minimum size=5pt, inner sep=0},
ext/.style={circle, draw, fill=white, minimum size=5pt, inner sep=1pt},
helper/.style={coordinate},point/.style={circle, draw, fill, inner sep =1pt},
de/.style={-triangle 60},
point/.style={circle, draw, fill, minimum size=3pt, inner sep=0pt},
]

\draw[dashed] (-3,0)--(3,0);
\node[point, label=-90:1] (i1) at (-2,0) {} node[point, label=-90:2] (i2) at (0,0) {} node[point, label=-90:3] (i3) at (2,0) {}; 
\draw[triangle 60-triangle 60] (i2)+(-.4,.4)-- +(.4,.4);

\begin{scope}[xshift=8cm]
\node[point, label=90:1] (e1) at (0,2){};

\draw[dashed] (-3,0)--(3,0);
\node[point, label=-90:2] (i1) at (-2,0) {} node[point, label=-90:3] (i2) at (0,0) {} node[point, label=-90:4] (i3) at (2,0) {}; 
\draw[triangle 60-triangle 60] (i1)+(-.4,.4)-- +(.4,.4);
\draw[triangle 60-triangle 60] (i2)+(-.4,.4)-- +(.4,.4);
\draw[triangle 60-triangle 60] (i3)+(-.4,.4)-- +(.4,.4);

\end{scope}

\end{tikzpicture} 
\caption{\label{fig:UIUE} The subspaces $U_I^3$ (left) and $U_E^3$ (right).}
\end{figure}

Similarly, for each $n\in \Z_{\geq 0}$ we define a subset $U_E^n\subset \FM_2(n+1)$. For $n=0$, it is given by the single point in $\FM_2(1)$. For $n=2$, it is the (closed) upper semicircle $S^+\subset S^1$. For $n\geq 2$
\[
 U_E^n = \cap_{2\leq k \leq n+1} \pi_{1k}^{-1}(S^+) \cap \cap_{2\leq i<j\leq n+1} \pi_{ij}^{-1}( {1} ).
\]
Pictures of the sets $U_I^n$ and $U_E^n$ are shown in Figure \ref{fig:UIUE}.

\begin{rem}
 The subscripts $I$ and $E$ stand for internal and external. We will see below that these subspaces correspond to internal and external vertices of $\Br$ trees in some sense.
\end{rem}

From the subspaces $U_I^n, U_E^n$ new subspaces can be constructed using the operadic insertion maps. For example, the notation
\[
 U_I^3(U_E^2, U_I^2, U_E^3)
\]
shall denote the image of $U_I^3\times U_E^2 \times  U_I^2 \times U_E^3$ under the operadic composition
\[
\FM_2(3) \times \FM_2(3)\times \FM_2(2)\times \FM_2(4) \to \FM_2(9).
\]

\subsection{Topological moperads \texorpdfstring{$\FM_{2,1}$}{FM2,1}}
\label{sec:FM21}
\begin{figure}
 \centering
\makeatletter{}\usetikzlibrary{matrix}
\usetikzlibrary{arrows}
\usetikzlibrary{shapes}
\usetikzlibrary{through}
\usetikzlibrary{calc,3d}
\usetikzlibrary{decorations,decorations.pathmorphing}
\[
\begin{tikzpicture}[
scale=.7,
int/.style={circle, draw, fill, minimum size=5pt, inner sep=0},
ext/.style={circle, draw, fill=white, minimum size=5pt, inner sep=1pt},
helper/.style={coordinate},point/.style={circle, draw, fill, inner sep =1pt},
de/.style={-triangle 60},
point/.style={circle, draw, fill, minimum size=3pt, inner sep=0pt},
xst/.style={cross out, draw, minimum size=5 },
]

\node[xst, label=-90:{$out$} ] (out) at (0,0) {};
\draw (out.base) --+(130:.8);

\node [point] at (50:2) {};
\node [point] at (130:2.5) {};
\node [point, label=-90:{$j$}] (ej) at (-30:2.5) {};
\node [point] at (177:2) {};
\draw[dashed] (out.base) -- (ej);

\node at (4,0) {\large $\circ_j$};

\begin{scope}[xshift=6cm]
\node[point] at (-.2,-.2) {} node[point] at (-.2,.2) {};
\node[point] at (.3,0) {} node[point] at (1,0) {};
\draw[dashed] (1,0)--(-1,0) (0,1)--(0,-1);
\node  at (1.8,0) {$=$};
\end{scope}

\begin{scope}[xshift=11cm]
\node[xst, label=-90:{$out$} ] (out) at (0,0) {};
\draw (out.base) --+(130:.8);

\node [point] at (50:2) {};
\node [point] at (130:2.5) {};
\node [] (ej) at (-30:2.5) {};
\node [point] at (177:2) {};

\begin{scope}[shift=(ej), scale = 0.5, rotate=60]
\node[point] at (-.2,-.2) {} node[point] at (-.2,.2) {};
\node[point] (c) at (.3,0) {} node[point] at (1,0) {};
\draw[dashed] (c) circle (1);
\end{scope}

\end{scope}

\end{tikzpicture}
\] 
\caption{\label{fig:FM21rightaction} Picture of the right action of $\FM_2$ on $\FM_{2,1}$.}
\end{figure}
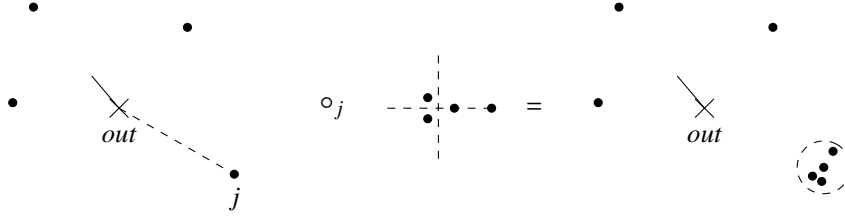

The components of the topological (or rather, semi-algebraic)  $\FM_2$-moperad $\FM_{2,1}$ are defined as 
\[
\FM_{2,1}(n)=\FM_2(n+1) \times S^1\, .
\]
It is composed of configurations of points, with one point distinguished which we call $out$, and some direction at the distinguished point. The (m)operadic compositions are defined as follows. Let $p\in \FM_2(n), q\in \FM_{2,1}(m)$, and assume that we want to determine $q\circ_j p\in \FM_{2,1}(m+n-1)$. First rotate the configuration $p$ such that the positive $y$-axis points into the direction of the ray from the position of $z_j$ in the configuration $q$ to the position of $out$. Second, insert the rotated configuration $p$ at the position of $z_j$. See Figure \ref{fig:FM21rightaction} for an example. Next, let $p\in \FM_{2,1}(n)$, $q\in \FM_{2,1}(m)$. Then $p\circ q$ is obtained by inserting $p$ at the $\vout$-vertex of $q$, after a rotation that aligns the positive $y$-axis with the specified direction at the $\vout$-vertex of $q$.

\begin{rem}
\label{rem:fm21mods1}
Clearly the space $\FM_{2,1}(n)$ can also be identified with 
\[
(\FM_2(n+1) \times S^1 \times S^1)/S^1 
\]
where the action of $S^1$ is by rotations on all three factors.
We think of the additional factor of $S^1$ as ``input direction'', and of the first factor of  $S^1$ as ``output direction''. Typically we will align the input direction with the positive real axis in drawings.

Furthermore, note that $\FM_{2,1}(0)= \FM_2(1) \times S^1= S^1$ is a circle. 
\end{rem}

\begin{ex} The space $\FM_{2,1}(2)$ is depicted in figure \ref{fig:ffm21}. In particular it contains a subspace $H^+\subset \FM_{2,1}(2)$ which is the closure of the set of configurations in which point $1$ is farther away from the origin (i.e., from $out$) than point $2$.
\end{ex}

\begin{figure}
 \centering
\makeatletter{}\usetikzlibrary{arrows}
\begin{tikzpicture}[
point/.style={circle, draw, fill, minimum size=3pt, inner sep=0pt},
xst/.style={cross out, draw, minimum size=5, outer sep =-1, inner sep=0 },
arr/.style={triangle 45-triangle 45},
scale=.7
]
\node[xst, label=-120:{$out$}] at (0,0) {};
\draw[dashed] (0,0) circle (3);
\node[point,label=30:{$k$}] (v0) at (0:3) {};
\node[point,label=140:{$0$}] (v1) at (140:3) {};
\node[point,label=-120:{2}] (v2) at (-100:3) {};
\node[point,label=-180:{1}] (v3) at (-180:3) {};
\node[point,label=80:{$n-1$}] (v4) at (60:3) {};
\node[point,label=100:{$n$}] (v5) at (100:3) {};
\draw[dashed] (0,0) -- (4,0);
\draw[dashed] (0,0) -- (140:3);
\draw[very thick] (0,0) -- (140:.7);
\draw[arr](1.6429,1.9741) -- (0.8618,2.5014);
\draw[arr](-0.056,2.7542) -- (-0.9152,2.5198);
\draw[arr](-2.3016,1.2864) -- (-1.7353,2.0089);
\draw[arr](-2.6726,-0.5644) -- (-2.6726,0.3491);
\draw[arr](-0.9738,-2.5367) -- (0.0026,-2.732);
\node[rotate=45] at (-45:3.3) {$\cdots$};
\node[rotate=-60] at (30:3.3) {$\cdots$};
\end{tikzpicture} 
\caption{\label{fig:VK} A picture of the subspace $V^{n,k}\subset \FM_{2,1}(n+1)$.
All points $0, \dots, n$ are forced on a circle around $\vout$. The radius does not matter since we divide out scalings of $\R^2$ anyway. All points except point $k$ can move. Point $k$ is fixed on the positive real axis, i.e., at the input direction. The output direction (the thick line at $\vout$) is fixed to aim at point $0$. Note that all boundary points are also part of $V^{n,k}$, i.e., some or all points are allowed to come infinitely close together.
}
\end{figure}
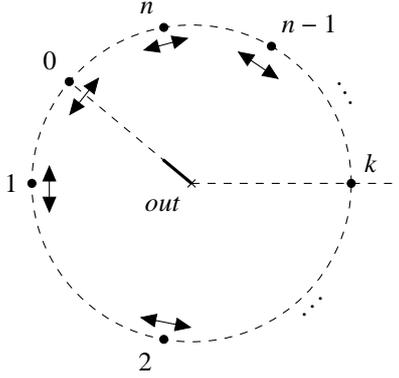

There are important closed subsets $V^{n,k}\subset \FM_{2,1}(n+1)$ that we will be using below, for $n=0,1,2,\dots$ and $k\in \{0,\dots,n\}$. The subset $V^{n,k}$ is (the closure of the) set of configuration in which all $n+1$ points (labelled $0,\dots, n$) are on a circle around $\vout$, in the cyclic order dictated by the labels.
Furthermore the output direction is constrained to point towards point $0$, and the input direction is constrained to point towards point $k$.
Since for us the input direction points along the positive real axis, this means that point $k$ is kept fixed on the positive real axis. A picture of the 
subsets $V^{n,k}$ can be found in figure \ref{fig:VK}. 

The moperad $\FM_{2,1}$, together with the operad $\FM_2$ form a two-colored operad which we call $\EFM_2$ (extended $\FM_2$).

\begin{figure}
 \centering
\makeatletter{}\usetikzlibrary{matrix}
\usetikzlibrary{arrows}
\usetikzlibrary{shapes}
\usetikzlibrary{through}
\usetikzlibrary{calc,3d}
\usetikzlibrary{decorations,decorations.pathmorphing}
\begin{tikzpicture}[scale =2]

\begin{scope}[canvas is xy plane at z=-1]
\draw (-1,-1) rectangle (1,1);
\draw[fill=white, draw=white] (1,0) circle (.2);
\clip (-1,-1) rectangle (1,1);
\draw[fill=white, dashed] (1,0) circle (.2);
\end{scope}
\begin{scope}[canvas is xy plane at z=1]
\draw (-1,-1) rectangle (1,1);
\draw[fill=white, draw=white] (-1,0) circle (.2);
\clip (-1,-1) rectangle (1,1);
\draw[fill=white] (-1,0) circle (.2);

\end{scope}

\begin{scope}[canvas is yz plane at x=-1]
\draw (-1,1) -- (-1,-1) (1,1)--(1,-1);
\clip (-1,-1) rectangle (1,1);
\draw[] (0,1) circle (.2);
\end{scope}

\begin{scope}[canvas is yz plane at x=1]
\draw (-1,1) -- (-1,-1) (1,1)--(1,-1);

\clip (-1,-1) rectangle (1,1);
\draw[] (0,-1) circle (.2);
\end{scope}

\draw[dashed] (-1,0.2,1)--(1,0.2,-1) (-1,-0.2,1)--(1,-0.2,-1);

\begin{scope}[canvas is xz plane at y=0]
\draw[fill=red, opacity=.3] (-1.2,-1.2) rectangle (1.2,1.2);
\end{scope}

\end{tikzpicture}
 
\caption{\label{fig:ffm21} Picture of the space $\FM_{2,1}(2)$, with the factor $S^1$ from the framing omitted. It has the form $S^1\times S^1\times [0,1] - T$ where $T$ is a tube. More concretely, a configuration $(z_1,z_2)$ of two points in $\mathbb{C}\setminus \{0\}$ determines two points in $S^1$, namely $\frac{z_1}{|z_1|}$ and $\frac{z_2}{|z_2|}$, and a relative distance from the origin $r=\frac{|z_1|}{|z_2|}$. The cut-out tube arises because the two points are forbidden to collide, and one performs a real blow-up at this locus. 
In the picture the horizontal directions correspond to the factors of $S^1$, i.e., the cube drawn here should be thought of as periodically repeated in the horizontal plane.
The vertical axis corresponds to $r$.
The space on top of and including the red plane (times the suppressed factor $S^1$ from the framing) is the subspace $H^+$, corresponding to configurations in which point $z_1$ is farther away from the origin than $z_2$. 
}
\end{figure}
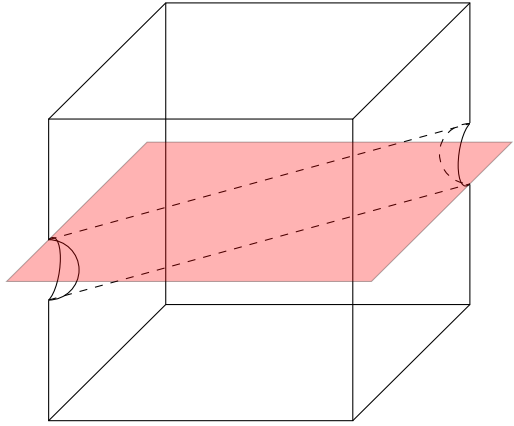

\subsection{Kontsevich's configuration spaces \texorpdfstring{$D_{K}$ and $D_{Ke}$}{DK and DKe}}
\label{sec:DKdef}
The Kontsevich halfspace $D_{Ke}(m,n)$ is the space of configurations of $m$ points in the upper halfplane and $n$ points on $\R$, modulo scaling and translation, suitably compactified.\footnote{See \cite{K1} for more details and the compactification.}
\begin{multline*}
 D_{Ke}(m,n) = \left( \{(z_1,\dots z_m, w_1,\dots,w_n)\in \C^{n+m} \mid 
z_i\neq z_j, w_i\neq w_j \text{for $i\neq j$},
\right. \\ \left.
 \Im(z_j)>0, \Im(w_j)=0\forall j, w_1<w_2<\dots <w_n  \}/ \R_+ \ltimes \R \right)^-
\end{multline*}
For technical convenience we will set $D_{Ke}(0,1)=\{pt\}$.
Together with $\FM_2$ these spaces assemble to form a two-colored operad, which is homotopic to the Swiss Cheese operad, see \cite{Vor} for details. 
We will need this operad later and call it $\SC$. Concretely its component with $m$ inputs of color 1, $n$ inputs of color 2 and the output in color $\alpha\in\{1,2\}$ is\footnote{Here we again use the notation $\op P^\alpha(m_1,m_2,\dots)$ to denote the space of operations with output in color $\alpha$, $m_1$ inputs of color 1, $m_2$ inputs of color 2 etc. of a colored operad $\op P$. }
\[
 \SC^\alpha(m,n) =
\begin{cases}
 \FM_2(m)  &\quad \text{for $n=0,\alpha=1$} \\
 D_{Ke}(m,n)  &\quad \text{for $\alpha=2$}\\
 \emptyset &\quad \text{otherwise} .
\end{cases}
\]
In particular, $D_{Ke}$ is an operadic right $\FM_2$-module. 
The operad $\SC$ is of Swiss Cheese type as defined in definition \ref{def:SCtype}.
The notation $D_K$ is shorthand for $D_K(m) := D_{Ke}(m,0)$. It is an operadic right sub-$\FM_2$-module. We will mainly be interested in the spaces of semi-algebraic chains on these spaces, $C(D_K)$ and $C(D_{Ke})$. The following Proposition is quite important for this paper.

\begin{prop}
\label{prop:KSactonDK}
There is an operadic left $\Br$-action on $C(D_K)$, making $C(D_K)$ an operadic $\Br$-$C(\FM_2)$-bimodule, such that on homology
\[
H\bpm \Br & C(D_K) & C(\FM_2) \epm 
\cong 
\bpm \Ger & \Ger & \Ger \epm \, .
\]
Here the middle $\Ger$ on the right is $\Ger$ considered as operadic $\Ger$-$\Ger$ bimodule.
\end{prop}
\begin{proof}
This action is constructed in Appendix \ref{app:KSactproof}. 
Given this action, consider the induced action on homologies. 
There is a map (of right $C(\FM_2)$-modules)
\begin{gather*}
C(\FM_2) \to C(D_K) \\
c \mapsto c_1 \circ c
\end{gather*}
where $c_1$ is the fundamental chain of $D_K(1)=\{pt\}$. Evidently, this map induces an isomorphism in homology. 
Hence it remains to show that under this map the left action of $\Br$ descends to the usual left action of $H(\Br)\cong \Ger$ on $\Ger$. Clearly it suffices to check this for the generators of $H(\Br)$, namely the ``bracket'' operation (depicted in Figure \ref{fig:PTmc}) and the ``product'' $T_2'$, see figure \ref{fig:KSgen}. 
Given the explicit formulas of the action this is a straightforward verification.
\end{proof}

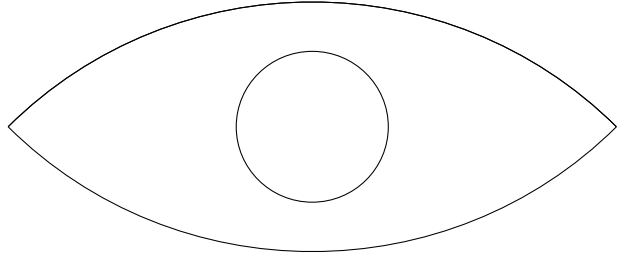
\begin{figure}
 \centering
\makeatletter{}\usetikzlibrary{matrix}
\usetikzlibrary{arrows}
\usetikzlibrary{shapes}
\usetikzlibrary{through}
\usetikzlibrary{calc,3d}
\usetikzlibrary{decorations,decorations.pathmorphing}
\begin{tikzpicture}[scale =2]

\draw[out=-45, in =-135] (-2,0) to (2,0);

\begin{scope}
\draw[clip, out=45, in =135] (-2,0) to (2,0);
\end{scope}
\draw[out=45, in =135] (-2,0) to (2,0);
\draw[fill=white] (0,0) circle (.5);
\end{tikzpicture} 
\caption{\label{fig:kontsevicheye} Picture of the space $D_K(2)$ (``Kontsevich's eye''). }
\end{figure}

A picture of the space $D_K(2)$ is shown in Figure \ref{fig:kontsevicheye}.
The topology of $D_K$ can be understood in terms of that of $\FM_2$:

\begin{lemma}
\label{lem:defretract}
 The spaces $\FM_2(n)$ are strong deformation retracts of $D_K(n)$ for $n=1,2,\dots$. Concretely, the inclusions 
\[
 \iota\colon \FM_2(n) \to D_K(n)
\]
are given by composition with the unique element of $D_K(1)$. The reverse maps are the forgetful maps
\[
 \pi \colon D_K(n) \to\FM_2(n)
\]
forgetting the location of the real line.
\end{lemma}
\begin{proof}
 It is clear that $\iota$ is an inclusion and that $\pi\circ\iota=\mathit{id}$. Hence $\FM_2(n)$ is a retract of $D_K(n)$. To see that it is a deformation retract, one has to specify a homotopy between the identity and $\iota\circ \pi$. It is given by ``moving upwards''. More concretely, if all points $z_1,\dots, z_n$ are finite distance from the real line, this map is $z_j\mapsto z_j+\lambda i$, $\lambda\to\infty$. However, some points may be infinitely close to the real line. General configurations in $D_K(n)$ are given by certain trees decorated with configurations of points at finite distance. The ``moving upwards'' is defined as follows: start at the lowest level of the tree in which points are close to the real axis. Move upwards (as before) all points in the configuration decorating that node. That will produce a tree where points are ``farther away'' from the real axis. By doing this repeatedly one moves away all points from the real axis to $+i\infty$. The projection of the configuration to $\FM_2(n)$ always stays the same.
\end{proof}

\subsection{Shoikhet's disks \texorpdfstring{$D_S$ and $D_{Se}$}{DS and DSe}}
\label{sec:DS}
\begin{figure}
 \centering
\makeatletter{}\usetikzlibrary{matrix}
\usetikzlibrary{arrows}
\usetikzlibrary{shapes}
\usetikzlibrary{through}
\usetikzlibrary{calc,3d}
\usetikzlibrary{decorations,decorations.pathmorphing}
\[
\begin{tikzpicture}[
int/.style={circle, draw, fill, minimum size=5pt, inner sep=0},
ext/.style={circle, draw, fill=white, minimum size=5pt, inner sep=1pt},
helper/.style={coordinate},point/.style={circle, draw, fill, inner sep =1pt},
de/.style={-triangle 60},
point/.style={circle, draw, fill, minimum size=3pt, inner sep=0pt},
xst/.style={cross out, draw, minimum size=5 , outer sep=0},
]
\draw (0,0) circle (3);
\node[xst, label=160:$\vout$] (ctr) at (0,0) {};
\draw (ctr.base) -- +(70:.4);
\node [point, label=0:{$\bar 0$}] at (3,0) {};
\node [point] at (30:3) {};
\node [point] at (130:3) {};
\node (v6) [point] at (230:3) {};
\node  [point] at (234:3) {};
\node  [point] at (226:3) {};
\node [point] at (77:3) {};

\node [point] at (30:1) {};
\node [point] at (130:2.5) {};
\node [point] at (230:.5) {};
\node (v7) [point] at (299:1.7) {};
\node [point] at (7:2) {};
\node [point] at (136:2) {};

\clip (0,0) circle (3);
\draw[dashed] (v6) circle (.5);
\draw (v6) +(30:.2) node[point] {} +(100:.2) node[point] {};

\draw[dashed] (v7) circle (.3);
\draw (v7) +(-30:.2) node[point] {} +(90:.2) node[point] {};

\end{tikzpicture}
\] 
\caption{\label{fig:DSex} Some configuration in $D_{Se}(10,6)$. The dashed circles shall indicate that the points enclosed are infinitesimally close together. The line stub at the center vertex indicates the framing.}
\end{figure}
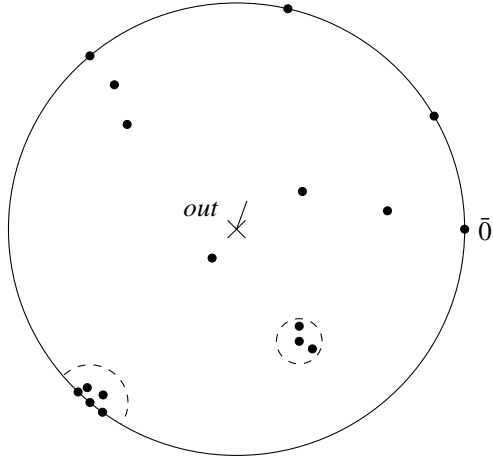

Define 
\[
D_{Se}(m,n+1) := D_{Ke}(m+1,n)\times S^1\, .
\]
It is the same as $D_{Ke}$ except that (i) there is an additional distinguished point, which we call $out$ and (ii) there is a direction specified at $out$. Note that the upper halfplane can be mapped biholomorphically onto a disk, minus one point on the boundary. In fact, we will think of $D_{Se}(m,n+1)$ as the space of configurations of points on a disk, with $m$ points in the interior and $n+1$ points on the boundary. The points on the boundary will be denote by $\bar 0, \dots, \bar n$, in counterclockwise order. The point $\bar 0$ is the ``additional point'', which is the image of $\infty$ under the map from the upper halfplane. 
 By the symmetries we divided out, we can assume that the point $out$ is fixed at the center of the disk, and that the point $\bar 0$ is at $1$. In this way it is clear how to understand the factor $S^1$ as parameterizing directions at $out$. A picture of some configuration in $D_{Se}(m,n)$ is shown in Figure \ref{fig:DSex}.

\begin{rem}
 The letter ``S'' stands for Shoikhet \cite{shoikhet}, who used the spaces $D_{Se}$, without the orientation at $out$, to construct the $\Lie_\infty$ morphism  $C_\bullet \to \Omega_\bullet$ from the introduction.
\end{rem}

The spaces $D_{Se}(m,n)$ are part of a three colored operad $\ESC$ (Extended Swiss Cheese operad) extending the (version of the) Swiss cheese operad $\SC$ from the previous subsection. More concretely, the operad $\ESC$ has colors 1,2,3 and the following color components:
\begin{itemize}
 \item The components with outputs in the first color are the components of $\FM_2$:
\[
 \ESC^1(m,n,r) = 
\begin{cases}
 \FM_2(m) &\quad\text{for $n=r=0$} \\
 \emptyset &\quad\text{otherwise}
\end{cases}
\]
\item The components with output in the second color are the same as in the Swiss cheese operad $\SC$:
\[
 \ESC^2(m,n,r) = 
\begin{cases}
 D_{Ke}(m,n) &\quad\text{for $r=0$} \\
 \emptyset &\quad\text{otherwise}
\end{cases}
\]
\item There are two sorts of components with output in color 3, namely the spaces $\FM_{2,1}$ or the spaces $D_{Se}(m,n)$.
\[
  \ESC^3(m,n,r) = 
\begin{cases}
 \FM_{2,1}(m) &\quad\text{for $n=0$, $r=1$} \\
 D_{Se}(m,n) &\quad\text{for $r=0$} \\
 \emptyset &\quad\text{otherwise}
\end{cases}
\]
\end{itemize}

In addition, the colored operad $\ESC$ is partially cyclic. By this we  mean that on the components $\ESC^3(m,n,0) = D_{Se}(m,n)$ there is an action of the cyclic group of order $n$ by permuting the labels on type II vertices.
The operadic compositions among the various components of $\ESC$ are straightforward ``insertions'', the explicit description in each case is left to the reader. The operad $\ESC$ is of Extended Swiss Cheese type, as defined by definition \ref{def:ESCtype}.

The most important role in this paper will be played by the components $D_{Se}(m,0)$, which we abbreviate by $D_S(m):=D_{Se}(m,0)$.
These spaces are right $\FM_2$-modules. Also, they are moperadic left $\FM_{2,1}$-modules. Similarly to Proposition \ref{prop:KSactonDK} there is the following result.
\begin{prop}
\label{prop:KSactonDS}
There is a moperadic $\Br$-$C(FM_2)$-$\KS_1$-$C(\FM_{2,1})$-$C(D_K)$ bimodule structure on the space of semi-algebraic chains on $D_S$, $C(D_S)$, such that on homology 
\[
H
\bpm
\Br & C(D_K) & C(\FM_2) \\
\KS_1 & C(D_S) & C(\FM_{2,1})
\epm
\cong 
\bpm
\Ger & \Ger & \Ger \\
\calc_1 & \calc_1 & \calc_1
\epm\, .
\]
\end{prop} 
\begin{proof}
The construction of this action can be found in Appendix \ref{app:KSactproof}.
Let us verify the statement about the (co)homology.\footnote{Note that we use cohomological conventions throughout, with the chains $C(\FM_2)$ etc. living in non-positive degrees.} There is an explicit identification of $H(D_S)$ with $\calc_1$ induced by the map
\begin{gather*}
C(\FM_{2,1}) \to  C(D_S) \\
c \mapsto c\circ c_1\, .
\end{gather*}
Here $c_1$ is the chain of a point in $D_S(0)$ and $\circ$ denotes the left action. From this description it is clear that on homology the actions of $H(\FM_2)\cong \Ger$ and $H(\FM_{2,1})\cong \calc_1$ are the standard ones. To check that also the combined action of $H(\KS_1)\cong \calc_1$ and $H(D_K)\cong \Ger$ is the standard one, it suffices to verify the statement on generators. But this is again straightforward given the explicit formulas for the actions.
\end{proof}

\begin{ex}
 Let us consider the simplest cases. $D_S(0)\cong S^1$ is a circle. The space $D_S(1)$ is three-dimensional and of the form $D_S(1)'\times S^1$. The part $D_S(1)'$ is the same as the ``Kontsevich eye'' depicted in Figure \ref{fig:kontsevicheye}. The space $D_S(2)$ is already five-dimensional and hard to depict. \end{ex}

Similar to Lemma \ref{lem:defretract} one proves the following Lemma.
\begin{lemma}
\label{lem:defretract1}
 The spaces $\FM_{2,1}(n)$ are strong deformation retracts of the spaces $D_S(n)$. The embedding 
\[
 \FM_{2,1}(n) \to D_S(n)
\]
is given by the moperadic action on the unit element in $D_S(0)$.\footnote{This element is the point in the configuration space where the framing at $\vout$ aligns with the input framing, i.e., with the positive real axis.} The projection map 
\[
 \pi \colon D_S(n) \to \FM_{2,1}(n)
\]
is the forgetful map sending a configuration of points in the unit disk to its equivalence class under rescalings.
\end{lemma}
Informally speaking, the forgetful map $\pi$ forgets the location of the circle bounding the disk.

\subsubsection{Incorporating the forgetful maps}
\label{sec:DSforgetful}
There is another combinatorial structure on $\ESC$ which we will need later, that is not encoded by the operadic compositions and the additional cyclic structure described in the previous subsection. Namely, there are the forgetful maps
\[
 \pi_{m,n,k} ; \ESC^2(m,n,0) \to \ESC^2(m,n-1,0)
\]
and similarly 
\[
 \pi_{m,n,k}' ; \ESC^3(m,n,0) \to \ESC^3(m,n-1,0)
\]
forgetting the location of the $k$-th type II vertex.
These maps can be conveniently packaged into the colored operad structure by defining $\ESC^2(0,0,0)$ to be a single point, which we call $\mathbb{1}$.\footnote{This notation is not optimal because of the possible confusion with the operadic unit. We will call the latter $\mathit{id}$ if needed. Note also that $\mathit{id}$ is a unary, but $\mathbb{1}$ a zero-ary operation. } 
The new colored operad we call $\EESC$. Its operadic compositions are defined such that an insertion of $\mathbb{1}$ into some vertex simply forgets the location of that vertex. I.e., for $c\in \ESC^2(m,n,0)$, we have
\[
 c \circ_k \mathbb{1} :=  \pi_{m,n,k} (c)
\]
where ``$\circ_k$'' denotes operadic composition at the $k$-th slot of color 2. A similar formula defines the compositions on $\ESC^3(m,n,0)$.

\subsection{The 4 colored operad \texorpdfstring{$\bigChains$}{bigChains}}
We can assemble the pieces introduced above into one big 4 colored operad $\bigChains$. Using the notation of section \ref{sec:opconventions} it is defined as
\[
 \bigChains = \bpm \Br & C(D_K) & C(\FM_2) \\
\KS_1 & C(D_S) & C(\FM_{2,1}) \epm.
\]
Note that we always use semi-algebraic chains as defined in \cite{HLTV}, not the standard singular chains. This is because we later want to integrate certain differential forms over chains, and there is no suitable integration theory for arbitrary continuous chains.

 
\makeatletter{}\section{The \texorpdfstring{$\KS_\infty$}{KS-infinity} operad}
\label{sec:KSdef}
The operad $\Br_\infty$ is by definition the bar-cobar construction of the Braces operad $\Br$.
\[
 \Br_\infty := \Omega(B(\Br)).
\]
For convenience, we give here an explicit combinatorial description, being imprecise with signs however. \begin{defi}
 A $\Br_\infty$-tree $\Gamma$ with $n$ external vertices is a planar rooted tree with four kinds of vertices: internal, external, red and blue, such that:
\begin{enumerate}
\item Every internal vertex has at least two children.
\item Every red and every blue vertex is decorated by a $\Br_\infty$-tree. We will consider the vertices of the decorations also as vertices of $\Gamma$. The reader should think of the decoration as inscribed in this red or blue vertex.
\item There are $n$ external vertices in total (including the decorations), labelled by numbers $\{1,\dots,n\}$.
\item $\Gamma$ is not equal to a tree with one vertex, possibly decorated with some $\Br_\infty$-tree. Nor is any decoration in $\Gamma$ of this form.
\item $\Gamma$ contains only finitely many vertices (including the decorations).
\end{enumerate}
\end{defi}
We leave it to the reader to check that the space of $n$-ary operations $\Br_\infty(n)$ of the operad $\Br_\infty$ may be viewed as the space spanned by all $\Br_\infty$-trees, for $n>1$. For $n=1$, $\Br_\infty(1)$ is by definition a one dimensional space, spanned by the operadic unit.

The operadic composition is as follows. Let $\Gamma_1\in \Br_\infty(n_1), \Gamma_2 \in \Br_\infty(n_2)$ be two trees, for $n_1, n_2\geq 2$. Then $\Gamma_1\circ_j \Gamma_2$ is the tree in $\Br_\infty(n_1+n_2-1)$ obtained by
\begin{itemize}
 \item Making the vertex labelled with $j$ of $\Gamma_1$ into a red vertex with inscribed tree $\Gamma_2$.
 \item Renumbering the external vertices. For example, the label on vertex $j+1$ of $\Gamma_1$ becomes $j+n_2$ etc.
\end{itemize}
It is clear that this operad is free. More concretely, it is the free operad generated by all elements without red vertices.
But the $\Br_\infty$-trees without red vertices may in turn be identified with trees decorated by elements of $\Br$, by ``cutting'' at the blue vertices. 
The space of such decorated trees forms a basis of the operadic bar construction $B(\Br)$. Hence the space of $\Br_\infty$-trees with $n$ external vertices and without red vertices can be identified with the $n$-ary co-operations of the bar construction $B(\Br)$, for $n>1$. In the following we will hence call such a $\Br_\infty$-tree without red vertices a $B(\Br)$-tree.

\begin{figure}
 \centering
\makeatletter{}\usetikzlibrary{through}
\usetikzlibrary{arrows}
\begin{tikzpicture}[
scale=.7,
int/.style={circle, draw, fill, minimum size=5pt, inner sep=0},
ext/.style={circle, draw, fill=white, minimum size=5pt, inner sep=1pt},
helper/.style={coordinate},
]
\node at (-1,1) {$\delta$};
\draw (0,1.5) -- (0,1) node[int, draw=blue, fill=blue] (v1) {} -- +(-.5,-.6);
\draw (v1)-- +(0,-.6)  (v1)--+(.5,-.6);
\draw[dashed] (v1)--+(1,0) +(2,0) node (c) {} circle (1);
\draw (c)+(0,.7)--++(0,.2) node[ext] (a1) {$1$} -- +(0,-.7) node[ext] (a2) {$2$};

\begin{scope}[xshift=5cm]
\node at (-1,1) {$=\, \sum$};
\draw (0,1.5) -- (0,1) node[int, draw=blue, fill=blue] (v1) {} ;
\node[int] (v2) at (0,.5) {};
\draw (v2) -- +(0,-.5) (v2) -- +(-.5,-.5);
\draw (v1) edge (v2)  (v1)--+(.5,-.6);
\draw[dashed] (v1)--+(1,0) +(2,0) node (c) {} circle (1);
\draw (c)+(0,.7)--++(0,.2) node[ext] (a1) {$1$} -- +(0,-.7) node[ext] (a2) {$2$};
\end{scope}

\begin{scope}[xshift=5cm, shift={(5,0)}]
\node at (-1,1) {$+\, \sum$};
\draw (0,1.5) -- (0,1) node [int] (v1) {} ;
\node[int, draw=blue, fill=blue] (v2) at (0,0.5) {};
\draw (v2) -- +(0,-0.5) (v2) -- +(+0.5,-0.5);
\draw (v1) edge (v2)  (v1)--+(-0.5,-0.6);
\draw[dashed] (v2)--+(1,0) +(2,0) node (c) {} circle (1);
\draw (c)+(0,0.7)--++(0,0.2) node [ext] (a1) {$1$} -- +(0,-0.7) node [ext] (a2) {$2$};

\end{scope}

\begin{scope}[xshift=10cm, shift={(-4,-3)}]
\node at (-1,1) {$-$};
\draw (0,1.5) -- (0,1) node [int, draw=blue, fill=blue] (v1) {} -- +(-0.5,-0.6);
\draw (v1)-- +(0,-0.6)  (v1)--+(0.5,-0.6);
\draw[dashed] (v1)--+(1,0) +(2,0) node (c) {} circle (1);
\draw (c)+(0,0.7)--++(0,0.2) node [ext] (a1) {$1$} -- +(0,-0.7) node [ext] (a2) {$2$};
\draw (c)+(-0.5,-0.1) node {$\delta$};
\end{scope}

\begin{scope}[xshift=15cm, shift={(0,0)}]
\node at (-1,1) {$+\, \sum$};
\draw (0,2)--(0,1.5) node [ext] (a1) {$1$} -- +(0,-0.7) node [ext] (a2) {$2$};
\draw (a2)--+(-0.5,-0.5) (a2)--+(0,-0.5) (a1)--+(0.6,-0.5);
\end{scope}

\begin{scope}[xshift=18cm, shift={(-7,-3)}]
\node at (-1,1) {$+$};
\draw (0,1.5) -- (0,1) node [int, draw=red, fill=red] (v1) {} -- +(-0.5,-0.6);
\draw (v1)-- +(0,-0.6)  (v1)--+(0.5,-0.6);
\draw[dashed] (v1)--+(1,0) +(2,0) node (c) {} circle (1);
\draw (c)+(0,0.7)--++(0,0.2) node [ext] (a1) {$1$} -- +(0,-0.7) node [ext] (a2) {$2$};
\node at (3.5,1) {$ $};
\end{scope}
\end{tikzpicture} 
\caption{\label{fig:brinftydiff} An illustration of the differential on $\Br_\infty$ and $B(\Br)$. One blue vertex of the $\Br_\infty$ or $B(\Br)$-tree is shown. One should think of the tree to be continued above and below that vertex. The tree in the dashed circle is the decoration of the blue vertex. The terms in the differential are (from left to right): The splitting of the vertex, the differential applied to the decoration and the insertion at the blue vertex. The last term, coloring the blue vertex red, only occurs for $\Br_\infty$-trees and is absent for $B(\Br)$-trees.}
\end{figure}
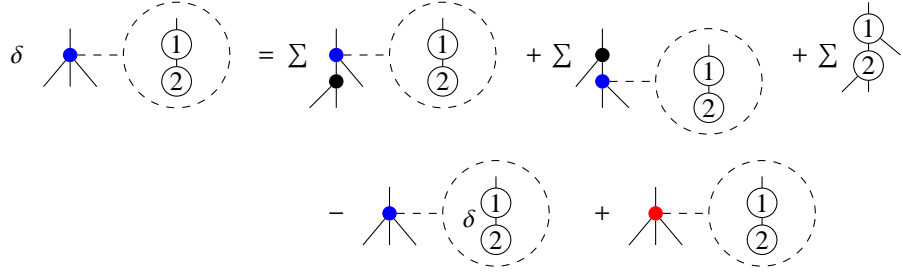

Next define the differential as a sum $d=d_s + d_{br} + d_{bi}$. Here $d_s$ splits off a new internal vertex from any vertex, similarly to the differential in $\Br$. $d_{br}$ has the form
\[
 d_{br}\Gamma = \sum_{v \text{ blue}} \pm \Gamma(v\to \text{red}).
\]
Here $\Gamma(v\to \text{red})$ is the graph obtained by coloring $v$ red. The term $d_{bi}$ acts on each blue vertex as follows: 
\begin{enumerate}
 \item Insert the tree that decorates the blue vertex at that blue vertex.
 \item Reconnect the child edges in all planar possible ways.
\end{enumerate}
A pictorial description of the differential is contained in Figure \ref{fig:brinftydiff}.
The differential on the cooperad $B(\Br)$ (spanned by $\Br_\infty$-trees without red vertices) is given by $d_s+d_{bi}$.
The cooperadic cocompositions are obtained by splitting the graph at some blue vertex, into one graph with the blue vertex made external, and into the decoration.

The degree of a $\Br_\infty$-tree $\Gamma$ can be calculated as follows. 
\begin{enumerate}
 \item Every edge has degree $-1$.
 \item Every red and external vertex has degree 0.
 \item Every blue vertex has degree $-1$.
 \item Every internal vertex has degree $+2$.
\end{enumerate}
The degree of a $B(\Br)$-tree $\Gamma$ (i.e., $\Gamma$ has no red vertices) is its degree as a $\Br_\infty$-tree, minus 1.

\subsection{Notation}
\label{sec:treenotation}
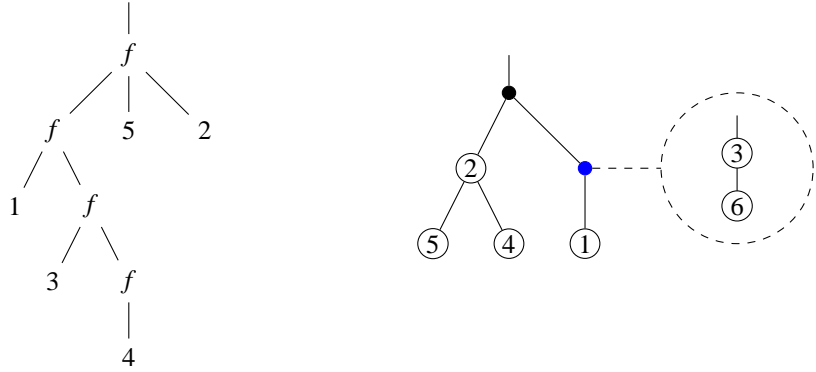
\begin{figure}
 \centering
\makeatletter{}\usetikzlibrary{through}
\usetikzlibrary{arrows}
\begin{tikzpicture}[
int/.style={circle, draw, fill, minimum size=5pt, inner sep=0},
ext/.style={circle, draw, fill=white, minimum size=5pt, inner sep=1pt},
helper/.style={coordinate},
]

\node (v1) at (0,0) {$f$};
\node (v2) at (-1,-1) {$f$};
\node (e5) at (0,-1) {$5$};
\node (e2) at (1,-1) {$2$};
\node (e1) at (-1.5,-2) {$1$};
\node (v3) at (-0.5,-2) {$f$};
\node (v4) at (0,-3) {$f$};
\node (e3) at (-1,-3) {$3$};
\node (e4) at (0,-4) {$4$};

\draw (v1) edge (v2) edge (e5) edge (e2);
\draw (v2) edge (e1) edge (v3) (v3) edge (v4) edge (e3) (v4) edge (e4);
\draw (v1) edge (0,.7);

\begin{scope}[xshift=5cm, yshift=-2cm]
\node[int] (I1) at (0,1.5) {};
\node[ext] (E2) at (-.5,.5) {2};
\node[ext] (E5) at (-1,-.5) {5};
\node[ext] (E4) at (0,-.5) {4};
\node[ext] (E1) at (1,-.5) {1};

\node[int, draw=blue, fill=blue] (v1) at  (1,.5) {};
\draw[dashed] (v1)--+(1,0) +(2,0) node (c) {} circle (1);
\draw (c)+(0,.7)--++(0,.2) node[ext] (a1) {$3$} -- +(0,-.7) node[ext] (a2) {$6$};

\draw (I1) edge (0,2) edge (E2) edge (v1) (E2) edge (E5) edge (E4) (v1) edge (E1); 

\end{scope}
\end{tikzpicture} 
\caption{\label{fig:treeformal} Left: Planar tree corresponding to the formal functional expression $f(f(1,f(3,f(4)),5,2)$. Right: $\Br_\infty$-tree corresponding to the formal functional expression $I(E(2; 5,4 ), B(E(3; 6), 1))$.}
\end{figure}

A planar tree can be defined recursively as an ordered list of planar trees. If the set of terminal (leaf) vertices is $[n]$, then such trees are in 1-1 correspondence to formal expressions that can be built with some formal function $f$, containing each symbol $1,2,\dots,n$ exactly once. For example,
\[
 f(f(1,f(3,f(4)),5,2)
\]
corresponds to the tree shown in Figure \ref{fig:treeformal} (left).
Similarly, a $\Br_\infty$-tree can be written as a formal expression using formal functions $I(\cdot,\cdot,\dots)$, $E(\cdot;\cdot,\dots)$, $R(\cdot;\cdot,\dots)$, $B(\cdot;\cdot,\dots)$ corresponding to internal, external, red and blue vertices.
More precisely, using terms from mathematical logic, one can define a \emph{formal language} whose alphabet consists of 
\begin{itemize}
 \item Functional symbols $I(\cdot,\cdot,\dots)$ (arities $\geq 2$), $E(\cdot;\cdot,\dots)$, $R(\cdot;\cdot,\dots)$, $B(\cdot;\cdot,\dots)$.
 \item Terminal symbols $1,2,\dots$ (corresponding to external vertices)
\end{itemize}
and several formation rules, such that the words of the language are in 1-1 correspondence with the natural basis (given by trees) of $\Br_\infty$. We will however proceed less formally, since the only purpose of this section is to set up some notation that will spare the author from drawing too many pictures. To give some examples, $I(1,3,2)$ denotes an internal vertex with children $1,3,2$ (in this order from left to right). $E(1;2,3)$ denotes an external vertex labelled $1$, with children $2$ and $3$. The first argument of $E(\cdots)$ must be a ``terminal symbol'', i.e., a number $1,2,\dots$.  $B(T;1,2)$ or $R(T;1,2)$ denote a blue or red vertex, decorated with some tree $T$, and having children $1$ and $2$. Here $T$ is just a placeholder, for example we could insert $T=I(1,2,3)$.
To give a more complicated example, the expression
\[
 I(E(2; 5,4 ), B(E(3; 6); 1)) 
\]
corresponds to the $\Br_\infty$-tree shown in Figure \ref{fig:treeformal} (right).
Of course, by restricting to formal functional expressions not containing the ``functions'' $R$ or $B$, we obtain a functional notation to describe $\Br$ elements. 
The tree corresponding to a functional expression can be recovered by replacing each occurrence of $E$ by a generator $T_n$ and each occurence of $I$ by a generator $T_n'$ (see Figure \ref{fig:KSgen}) and interpreting functional composition as operadic composition.
Finally, considering only formal functional expressions not containing the functions $R$ (red vertices) we obtain a functional notation to describe elements of the cooperad $B(\Br)$.
This notation is slightly more economic than drawing pictures, and the author will use it below.

\subsection{The operad \texorpdfstring{$\KS_\infty$}{KS-infinity}}
\label{sec:opKS1}
The operad $\KS_\infty := \Omega(B(\KS))$ is the bar-cobar construction of the colored operad $\KS$, see \cite{KS2}.
It consists of the operad $\Br_\infty$, colored in the first color, and a $\Br_\infty$ moperad $\KS_{1,\infty}$. Elements $\KS_{1,\infty}(n)$ have the output and one input in the second color, and $n$ further inputs in the first color.
Such elements also have a combinatorial description, and one can set up some algebraic notation similar to section \ref{sec:treenotation}. Let us leave the combinatorial description in terms of certain graphs to the reader and jump to the algebraic description.\footnote{The reader should not take the following discussion as a ``definition'' for $\KS_\infty := \Omega(B(\KS))$, just as an eloboration on how its elements look like and set up of some notation. We will be a bit informal language-wise.} First, let us consider the moperad $\KS_1$. Any $\KS_1$-graph can be described by a formal functional expression of the form
\[
 K(T_0 \text{ or $\mathbb{1}$},T_1,\dots, T_n)
\]
where $K$ is some new symbol (formal function), $n=0,1,\dots$ and the $T_j$ are placeholders for some $\Br$ trees, which we think of as functional expressions in formal functions $I$ and $E$ as before, and terminal symbols $1,2,\dots$ and an extra terminal symbol $\vin$. Two legal examples would be 
\[
 K(\mathbb{1}, I(2,I(3,\vin)), 1)
\]
corresponding to the graph 
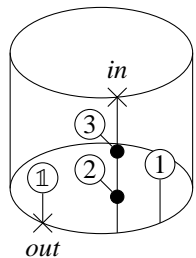
\[
\begin{tikzpicture}[
scale=.6,
int/.style={circle, draw, fill, minimum size=5pt, inner sep=0},
ext/.style={circle, draw, fill=white, minimum size=5pt, inner sep=1pt},
de/.style={-triangle 60},
point/.style={circle, draw, fill, minimum size=3pt, inner sep=0pt},
xst/.style={cross out, draw, minimum size=5pt}
]
\begin{scope}[xshift=6cm]
\draw (0,0) ellipse (2 and 1);
\draw (0,3) ellipse (2 and 1);
\draw (-2,0)--(-2,3) (2,0)--(2,3);
\node [xst, label=-90:{$out$}] (out) at ($(0,0)+(-130:2 and 1)$) {};
\node [xst, label=90:{$in$}] (in) at ($(0,3)+(-80:2 and 1)$) {};
\draw (out)+(0,1) node[ext] (e3) {$\mathbb{1}$};
\draw (out.base) edge (e3);
\draw (in.base) -- ++(0,-1.2) node[int] (i1) {} -- ++(0,-1) node[int] (i2) {} -- +(0,-.8);
\draw (i1)-- +(-.6, .6) node[ext] {3};
\draw (i2)-- +(-.6, .6) node[ext] {2};
\node [ext] (e1) at ($(0,1.3)+(-50:2 and 1)$) {1};
\draw (e1) -- +(0,-1.3);
\end{scope}
\end{tikzpicture}
\]
and
\[
 K(E(1;2,I(4,\vin)),3)
\]
corresponding to the graph 
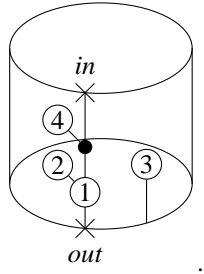
\[
\begin{tikzpicture}[
scale=.6,
int/.style={circle, draw, fill, minimum size=5pt, inner sep=0},
ext/.style={circle, draw, fill=white, minimum size=5pt, inner sep=1pt},
de/.style={-triangle 60},
point/.style={circle, draw, fill, minimum size=3pt, inner sep=0pt},
xst/.style={cross out, draw, minimum size=5pt}
]
\begin{scope}[xshift=6cm]
\draw (0,0) ellipse (2 and 1);
\draw (0,3) ellipse (2 and 1);
\draw (-2,0)--(-2,3) (2,0)--(2,3);
\node [xst, label=-90:{$out$}] (out) at ($(0,0)+(-100:2 and 1)$) {};
\node [xst, label=90:{$in$}] (in) at ($(0,3)+(-100:2 and 1)$) {};
\draw (in.base) -- ++(0,-1.2) node[int] (i1) {} -- ++(0,-1) node[ext] (e1) {1} -- +(0,-.8);
\draw (i1)-- +(-.6, .6) node[ext] {4};
\draw (e1) -- +(-.6, .6) node[ext] {2};
\node [ext] (e1) at ($(0,1.3)+(-60:2 and 1)$) {3};
\draw (e1) -- +(0,-1.3);
\end{scope}
\end{tikzpicture}\, .
\]
This describes a language for specifying $\KS_1$ elements. If we want to describe $B(\KS)$ elements, we have to (i) allow blue vertices in the trees $T_j$ and (ii) introduce another function, say $K_B(\cdot ;\cdots)$. The first slot can be filled with some other $B(\KS)$-tree, the remaining slots can be filled in the same way as the slots of $K(\cdots)$. Some typical element is
\[
 K_B(K(\mathbb{1}, B(E(3;4,5); 6,7), \vin ); 1,2,\vin).
\]
As a graph we could draw this element like this:
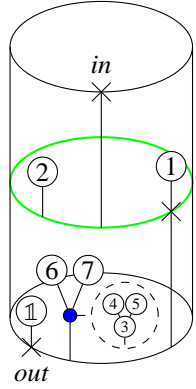
\[
\makeatletter{}\begin{tikzpicture}[
scale=.6,
int/.style={circle, draw, fill, minimum size=5pt, inner sep=0},
ext/.style={circle, draw, fill=white, minimum size=5pt, inner sep=1pt},
de/.style={-triangle 60},
point/.style={circle, draw, fill, minimum size=3pt, inner sep=0pt},
xst/.style={cross out, draw, minimum size=5pt}
]

\begin{scope}[shift={(-1,0)}]

\draw (0,6) ellipse (2cm and 1cm);
\draw (0,0) ellipse (2cm and 1cm);
\draw[draw=green, thick] (0,3) ellipse (2cm and 1cm);
\draw (-2,0)--(-2,6) (2,0)--(2,6);
\node [xst,label=-90:{$out$}] (out) at ($(0,0)+(-140:2 and 1)$) {};
\node [xst] (in) at ($(0,3)+(-40:2 and 1)$) {};
\node [xst,label=90:{$in$}] (inu) at ($(0,6)+(-90:2 and 1)$) {};

\coordinate (I1) at ($(0,3)+(-130:2 and 1)$) {};
\draw (I1) -- +(0,1) node[ext]{2};
\draw (in.base) -- +(0,1) node[ext]{1};

\coordinate (B0) at ($(0,0)+(-110:2 and 1)$) {};
\draw (B0) -- ++(0,1) node[draw, fill=blue, minimum size=5, circle, inner sep=0] (B) {} 
-- +(-.4,1) node[ext] {6}  (B) -- + (.4,1) node[ext] {7};

\draw[dashed] (B) -- +(1.2,0) node[draw, circle, minimum size=25, fill=white] (insert) {};

\begin{scope}[scale=.5]
\draw (insert) +(0,-.5) node[ext] (e3) {\tiny 3};
\draw (e3) -- +(-.5,1) node[ext] {\tiny 4};
\draw (e3) -- +(.5,1) node[ext] {\tiny 5};
\draw (e3) -- +(0,-.8);
\end{scope}

\draw (inu.base)--+(0,-3);
\draw (out.base)--+(0,0.8) node [ext] {$\mathbb{1}$};
\draw (in.base) -- ($(0,0)+(-40:2 and 1)$) ;

\end{scope}

\end{tikzpicture} 
\]

To obtain a functional notation for elements of $\KS_\infty$ one (i) has to allow also red vertices (functional symbols $R(\cdots)$) and (ii) introduce another functional symbol $K_R$.

\makeatletter{}\section{The Kontsevich--Soibelman proof, and an extension}
\label{sec:ksproof}
In this section we will review the construction of Kontsevich and Soibelman \cite{KS1,KS2} of a map of colored operads
\[
 \KS_\infty =\bpm \Br_\infty & \KS_{1,\infty}\epm \to \bpm C(\FM_2) & C(\FM_{2,1}) \epm.
\]
In fact, the author does not understand some part of the Kontsevich-Soibelman construction for the moperad-piece of the map, due to incontractibility of some spaces. Hence we will redo that part with a slightly different argument.  

Finally, in sections \ref{sec:KShomext} and \ref{sec:KShomext2} we extend the arguments to the construction of a map of colored operads
\[
 \homKS = 
\bpm \Br & \hBr_\infty & \Br_\infty \\
\KS_1 & \hKS_{1, \infty} & \KS_{1,\infty}\epm
\to
\bpm \Br & C(D_K) & C(\FM_2) \\
\KS_1 & C(D_S) & C(\FM_{2,1}) \epm = \bigChains.
\]

\subsection{The map \texorpdfstring{$\Br_\infty\to C(\FM_2)$}{Br-infinity to C(FM2)}, following \cite{KS1}}
\label{sec:brtofm2}
In this subsection we review the construction of \cite{KS1}.
The goal is to construct a map of operads
\[
 \Br_\infty=\Omega(B(\Br)) \to C(\FM_2).
\]
Since $\Br_\infty$ is quasi-free, this amounts to constructing a map
\[
 c \colon B(\Br)[-1] \to C(\FM_2)
\]
that is equivariant with respect to the action of the symmetric group and satisfies equations encoding the compatibility with the differential. Concretely, these equations, for $\Gamma$ a tree in $B(\Br)$, have the form\footnote{Here the explicit signs are not important as long as $\p^2=0$.}
\begin{equation}
\label{equ:Brcase}
 \partial c(\Gamma) = c(d_{B(\Br)} \Gamma ) + \sum \pm c(\Gamma') \circ c(\Gamma'').
\end{equation}
Here $d_{B(\Br)}$ is the differential on the bar construction $B(\Br)$. The imprecise notation $\sum \pm c(\Gamma') \circ c(\Gamma'')$ means the following: take the restricted cocomposition in $B(\Br)$ of $\Gamma$, yielding $\Gamma', \Gamma''\in B(\Br)$. Then compute the images of $\Gamma', \Gamma''$ under $c$, and compose again (in $C(\FM_2)$).
Recall from section \ref{sec:KSdef} that $B(\Br)$ is graded in non-positive degrees.\footnote{Note that we also use the non-positive (cohomological) grading on $C(\FM_2)$, so that our differentials always have degree $+1$.} Note that on the right hand side of \eqref{equ:Brcase}, $c$ is applied only to elements of $B(\Br)$ that have strictly larger degree than $\Gamma$. Hence, recursively, one has to solve equations of the form
\[
  \partial c(\Gamma) = (\text{something known}).
\]
Here one can check that $(\text{something known})$ is in fact a cocycle, using \eqref{equ:Brcase} for $\Gamma$'s of higher degrees. The question is hence whether one can make choices such that these cocycles are all exact. 
The ideal, of course, would be to have explicit formulas for all of the $c(\Gamma)$. 
However, an almost as good solution is to define non-empty subsets $U_\Gamma \subset \FM_2(n_\Gamma)$ for each $B(\Br)$-tree $\Gamma$ with $n_\Gamma$ external vertices, satisfying the following requirements.
\begin{itemize}
 \item The $U_\Gamma$ are contractible.
 \item If $\tilde{\Gamma}$ is any tree occuring in the expression $d_{B(\Br)} \Gamma$, then $U_{\tilde{\Gamma}} \subset U_\Gamma$.
 \item For $\Gamma'$, $\Gamma''$ as on the right hand side of \eqref{equ:Brcase}, we have that $U_{\Gamma'}\circ U_{\Gamma''}\subset U_\Gamma$.
 \item The assignment $\Gamma\to U_\Gamma$ is equivariant under the symmetric group action.
\end{itemize}

Given those subsets, one can recursively solve \eqref{equ:Brcase}, in such a way that $c(\Gamma)\in C(U_\Gamma)\subset C(\FM_2(n_\Gamma))$.
Namely, by the second condition we know that the right hand side of \eqref{equ:Brcase} is a cocycle in $C(U_\Gamma)$. Hence, by the first condition the equation can be solved in such a way that $c(\Gamma)\in C(U_\Gamma)$, for all $\Gamma$ of degrees $\leq -2$. For $\Gamma$ of degrees $0$ and $-1$ one must be more careful since there might appear obstructions in $H_0(U_\Gamma)$ on the right hand side of \eqref{equ:Brcase}. Those cases will be treated below.

Kontsevich and Soibelman gave a definition for the spaces $U_\Gamma$. Let us call their spaces $U_\Gamma^{KS}$, because we will define and work with slightly smaller $U_\Gamma$ later. Fix the graph $\Gamma$ with $n=n_\Gamma$ external vertices. On the set of external vertices one can define two half-orderings, the horizontal half-ordering $<_{h}$, and the vertical half-ordering $<_v$, defined as follows:
\begin{itemize}
 \item If vertex $i$ is an ancestor of vertex $j$, then $j <_v i$. Here, if $i$ is contained in a subtree assigned to a blue vertex, it counts as ancestor of all children of the blue vertex.
 \item If $i,j$ are not ancestors of each other, find the ``youngest'' common ancestor $k$ of $i,j$. If the subtree of $i$ stands to the left of the subtree of $j$ in the ordering on $star(k)$, then  $i <_h j$.
\end{itemize}

Consider the configuration space of two points, $\FM_2(2)$. It is the circle. Let $S^+$ be the closed upper semicircle, corresponding to $Im(z_1)\geq Im(z_2)$. Let $\pi_{ij}: \FM_2(n)\to \FM_2(2)$ be the forgetful map, forgetting all but the two points $i,j$.
Then Kontsevich and Soibelman define
\begin{equation}
\label{equ:UKSdef}
 U_\Gamma^{KS} = \cap_{j <_v i} \pi_{ij}^{-1}(S^+) \cap \cap_{i <_h j} \pi_{ij}^{-1}({1}). 
\end{equation}

One can check that the three conditions above are satisfied, and hence one can recursively solve \eqref{equ:Brcase}. 
To actually write down proofs it is helpful to have a recursive definition of the sets $U_\Gamma$. We give such a defintion in the next subsection, and the resulting spaces will satisfy $U_\Gamma\subset U_\Gamma^{KS}$, with equality for ``most'' $\Gamma$. 

\begin{rem}
 We note that operad maps $\Omega(\op C)\to \op P$ from the cobar construction of a coaugmented cooperad $\op C$ into an arbitrary operad $\op P$ correspond to Maurer-Cartan elements in the convolution complex (or rather: convolution dg Lie algebra) of maps of symmetric sequences from $\op C$ to $\op P$. In particular, \eqref{equ:Brcase} may be considered as a Maurer-Cartan equation and our inductive definition of the map $\Br_\infty\to C(\FM_2)$ fits into the general obstruction theoretic framework for defining Maurer-Cartan elements of (filtered) dg Lie algebras. We will, however, not make use of this point below.
\end{rem}

\subsection{Definition of \texorpdfstring{$U_\Gamma$}{UGamma}}
\label{sec:UGammarecursive}
One can use the recursive structure of trees in $B(\Br)$ from Section \ref{sec:KSdef}. Instead of drawing pictures with trees we use the functional notation from Section \ref{sec:treenotation}.  
Moreover, we define the subsets $U_\Gamma$ for \emph{extended} $B(\Br)$-trees to allow for a recursive definition. Here ``extended'' means that we allow the tree or decorations of blue vertices to consist of a single vertex.

Let $\Gamma$ be a $B(\Br)$-tree. Let the associated formal functional expression according to section \ref{sec:treenotation} be $F$, 
\begin{enumerate}
 \item If $F$ is one of the terminal symbols $1,2,\dots$ corresponding to external vertices, we set $U_\Gamma=\{pt\}=: \FM_2(1)$
 \item If $F=I(T_1,\dots,T_n)$, where the functional expressions $T_1,\dots, T_n$ represent trees $\Gamma_1,\dots, \Gamma_n$, we set
\[
U_\Gamma = U_I^n(U_{\Gamma_1}, \dots, U_{\Gamma_n}).
\]
Here we use the subspaces $U_I^n$ from section \ref{sec:fm2op}. Furthermore the notation $U_I^n(U_{\Gamma_1}, \dots, U_{\Gamma_n})$ stands for the image under the operadic composition, inserting configurations in $U_{\Gamma_1}$ into the first slot of configurations in $U_I^n$, configurations in $U_{\Gamma_2}$ into the second etc.
 \item If $F=E(t;T_1,\dots,T_n)$, where the functional expressions $T_1,\dots, T_n$ represent trees $\Gamma_1,\dots, \Gamma_n$, and $t$ is a terminal symbol (i.e., $t=1,2,\cdots$) we set 
\[
U_\Gamma = U_E^n(t; U_{\Gamma_1}, \dots, U_{\Gamma_{n}}).
\]
Here the subspace $U_E^n$ is again defined in section \ref{sec:fm2op} and $U_E^n(t;\cdots)$ means that the first vertex in the configuration in $U_E^n$ gets labelled by $t\in \{1,2,\dots\}$.
 \item Finally we need to consider the case of a blue vertices, i.e., $F=B(T_0; T_1,\dots,T_n)$. Here the functional expressions $T_0, T_1,\dots, T_n$ represent trees $\Gamma_0,\dots, \Gamma_n$. 
Let the number of external vertices in $\Gamma_0$ be $p$. Define an auxiliary subspace $\tilde U\subset \FM(p+n)$ by
\[
 \tilde U = \pi^{-1}(U_{\Gamma_0})
\cap (\pi')^{-1}(U_I^n )
 \cap_{j=1}^n\cap_{\substack{\alpha\in V(\Gamma_0)}} \pi_{\alpha j}^{-1}(S^+).
\]
Here the projections $\pi$, $\pi'$ are the forgetful maps 
\[
 \FM_2(p) \stackrel{\pi}{\leftarrow} \FM_2(p+n) \stackrel{\pi'}{\rightarrow} \FM_2(n).
\]
Similarly $\pi_{\alpha j}:  \FM_2(p+n)\to \FM_2(2)$ forgets all but points $\alpha$ and $j$.
The subspace $U_I^n$ is as before. $V(\Gamma_0)$ is the set of external vertices in $\Gamma_0$. The subspace $S^+\subset \FM_2(2)$ is the upper semicircle, as in section \ref{sec:fm2op}.
In words, $\tilde U$ consists of configurations in $U_{\Gamma_0}$, placed above an additional $n$ points located on a horizontal line.
Then we set 
\[
 U_\Gamma = \tilde U(U_{\Gamma_1}, \dots, U_{\Gamma_n})
\]
where the notation means that we insert configurations from $U_{\Gamma_j}$ at the location of the $j$-th additional point in configurations of $\tilde U$.
 For a picture of the situation look at Figure \ref{fig:ugamma}.
\end{enumerate}
\begin{figure}
 \centering
\makeatletter{}\[
\begin{tikzpicture}[
int/.style={circle, draw, fill, minimum size=5pt, inner sep=0},
ext/.style={circle, draw, fill=white, minimum size=5pt, inner sep=1pt},
helper/.style={coordinate},point/.style={circle, draw, fill, inner sep =1pt}, de/.style={-triangle 60}
]
\draw (0,1.5) -- (0,1) node[int, draw=blue, fill=blue] (v1) {};
\draw (v1) +(-.5,-1) node (g1) {$\Gamma_1$} 
 (v1)+(.5,-1) node (g2){$\Gamma_2$};
\draw[dashed] (v1)--+(1,0) +(1.5,0) node (c) {$\Gamma_0$} circle (.5);
\draw (v1) edge (g1) (v1) edge (g2);

\begin{scope}[xshift=7cm, yshift=-1cm]
\draw[dashed] (-4,0)--(4,0);
\node[point] at (1,1) {};
\node[point] at (1.7,3) {};
\node[point] at (-1,2.3) {};
\node[point] at (-2,3.5) {};
\node[point] at (-2,1.2) {};

\node[semicircle, anchor=south,draw, inner sep=2, dashed] (v1) at (-2,0) {$U_{\Gamma_1}$};
\node[semicircle, anchor=south, draw,inner sep=2, dashed] (v2) at (2,0) {$U_{\Gamma_2}$};
\draw (v1.base)+(1.3,0) node (v3) {};
\draw (v1.base)+(-1.3,0) node (v4) {};
\draw (v2.base)+(-1.3,0) node (v5) {};
\draw (v2.base)+(1.3,0) node (v6) {};
\draw (v1)edge[de] (v3) (v1) edge[de](v4) (v2) edge[de] (v5) (v2) edge[de] (v6);
\node at (0,3) {$U_{\Gamma_0}$};
\end{scope}
\end{tikzpicture}
\] 
\caption{\label{fig:ugamma} Illustration of the recursive definition of $U_\Gamma$, for $\Gamma$ containing a blue vertex. On the left $\Gamma$ is shown. It has a blue vertex, decorated by $\Gamma_0$, with two sub-trees $\Gamma_1$ and $\Gamma_2$. On the right the subspace $U_\Gamma$ is schematically shown. 
}
\end{figure}
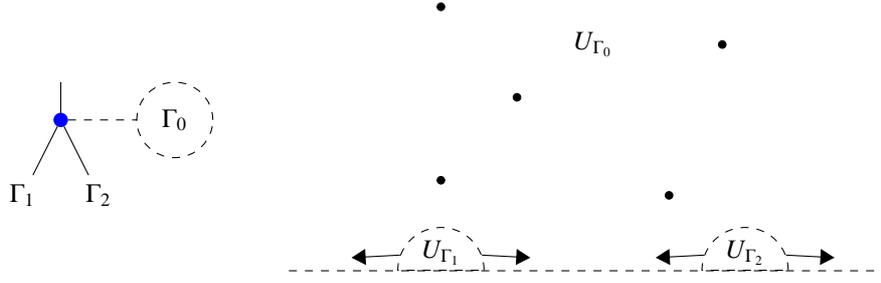
\begin{figure}
 \centering
\makeatletter{}\usetikzlibrary{arrows}
\[
\begin{tikzpicture}[
int/.style={circle, draw, fill, minimum size=5pt, inner sep=0},
ext/.style={circle, draw, fill=white, minimum size=5pt, inner sep=1pt},
helper/.style={coordinate},point/.style={circle, draw, fill, inner sep =1pt},
de/.style={-triangle 60},
point/.style={circle, draw, fill, minimum size=3pt, inner sep=0pt},
]
\begin{scope}[yshift=3cm]
\draw (0,2)--(0,1.5) node[ext] (e1) {1} --++(-1,-1) node[int] (i1){}--+(-.5,-1) node[ext](e2) {2};
\draw (i1)--+(.5,-1) node[ext](e3) {3};
\draw (e1)--++(1,-1) node[int](i2) {} --+(-.5,-1) node[ext](e4){4};
\draw (i2)--+(0,-1) node[ext](e5) {6};
\draw (i2)--+(0.5,-1) node[ext](e6) {5};
\end{scope}

\begin{scope}[xshift=7cm]
\draw[dashed] (-4,2)--(4,2);
\node[point,label=-90:3] at (-1.8,2) {};
\node[point,label=-90:2] at (-2.2,2) {};
\node[point,label=-90:4] at (2,2) {};
\node[point,label=-90:6] (c) at (2.4,2) {};
\node[point,label=-90:5] at (2.8,2) {};
\draw[dashed] (3,2) arc (0:180:.6);
\draw[<->] (c)+(-.2,.2) -- +(.2,.2);
\draw[dashed] (-1.6,2) arc (0:180:.4);
\draw[triangle 60-triangle 60] (-2.4,2.6)--(-1.6,2.6);

\draw[triangle 60-triangle 60] (c)+(-.4,.8)--+(.4,.8);
\node[point, label=0:{$1$}] at (0,5) {};

\end{scope}
\end{tikzpicture}
\] 
\caption{\label{fig:ugammaex} An example of the subspace $U_\Gamma$ corresponding to a $B(\Br)$-tree $\Gamma$, shown on the left. Note that on the right, the point 1 must always be on top of, or touch the dashed line. The point 1 is not allowed to get into the cluster of 2,3 or the cluster of 4,5,6. In other words, the point 2 is always much closer to 3 than to 1. }
\end{figure}
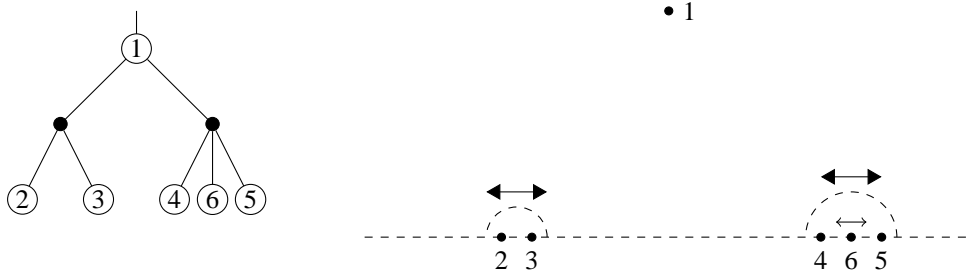

An example of the space $U_\Gamma$ corresponding to some $\Br$-tree $\Gamma$ is shown in Figure \ref{fig:ugammaex}.

\begin{lemma} \label{lem:UGammaprops}
 The spaces $U_\Gamma$ defined above satisfy the following conditions:
\begin{itemize}
 \item The $U_\Gamma$ are contractible.
 \item If $\tilde{\Gamma}$ is any tree occuring in the expression $d_{B(\Br)} \Gamma$, then $U_{\tilde{\Gamma}} \subset U_\Gamma$.
 \item For $\Gamma'$, $\Gamma''$ as on the right hand side of \eqref{equ:Brcase}, we have that $U_{\Gamma'}\circ U_{\Gamma''}\subset U_\Gamma$.
 \item The assignment $\Gamma\to U_\Gamma$ is equivariant with respect to the symmetric group action.
\end{itemize}
\end{lemma}
\begin{proof}[Proof (sketch)]
 To show the first statement, we proceed by induction. We have to consider three cases (i.~e., internal, external and blue vertices). For an internal vertex, the space $U_\Gamma$ is a product of spaces  $U_{\Gamma_1}, \dots, U_{\Gamma_n}$ with the contractible space $U_I^n$ (using the same notation as in the definition of $U_\Gamma$). By the induction hypothesis, $U_{\Gamma_1}, \dots, U_{\Gamma_n}$ are all contractible and hence $U_\Gamma$ is. A similar argument holds for external vertices, by just replacing $U_I^n$ by the space $U_E^n$, which is also contractible.
Next consider the case of a blue vertex. The space $U_\Gamma$ here is again a product
\[
\tilde U \times U_{\Gamma_1}\times \cdots \times U_{\Gamma_n}
\]
using the same notation as in the definition. The spaces $U_{\Gamma_1}, \dots, U_{\Gamma_n}$ are contractible by the induction hypothesis. Hence it suffices to show that $\tilde U$ is contractible as well.
 By ``moving upwards'', similar to the proof of Lemma \ref{lem:defretract}, we can deform $\tilde U$ to a space $\tilde U^\epsilon$, where all points in configurations from $U_{\Gamma_0}$ have finite distance from the $n$ additional points.
This space is in turn homotopic to a product of $U_{\Gamma_0}$ with a contractible space. Then using the induction hypothesis for $U_{\Gamma_0}$ contractibility of $U_\Gamma$ follows.

Next consider the second assertion of the Lemma. Again proceed by induction. Consider first internal vertices, i.~e., trees of the form 
\[
I(T_1, \dots, T_n).
\]
where $T_1, \dots, T_n$ correspond to trees $\Gamma_1,\dots, \Gamma_n$.
The differential has the following form:
\[
d (I(T_1, \dots, T_n)) = 
\sum_j \pm I(T_1,\dots, d T_j, \dots, T_n) 
+
\sum \pm
I(T_1,\dots, I(T_j,\dots, T_{j+k-1}), T_{j+k}, \dots, T_n).
\]
The subspaces corresponding to trees appearing in the first sum have the form 
\[
\bigcup_{\Gamma'} U_I(U_{\Gamma_1},\dots, U_{\Gamma'}, \dots, U_{\Gamma_n})
\]
where the union runs over trees $\Gamma'$ appearing nontrivially in $d \Gamma_j$. But by the induction hypothesis $\bigcup_{\Gamma'}U_{\Gamma'}\subset U_{\Gamma_j}$ and the preceding set is a subset of $U_\Gamma$. The subsets corresponding to terms of the second sum have the form
\[
U_I^{n-k+1}(U_{T_1},\dots, U_I^k(U_{T_j},\dots, U_{T_{j+k-1}}), U_{T_{j+k}}, \dots, U_{T_n}).
\]
It hence suffices to show that 
\[
U_I^{n-k+1}\circ_j U_I^k \subset  U_I^n
\]
but this is immediate from the definition of the $U_I$'s.
For external vertices the proof is similar. Consider next the blue vertices, i.e., a tree of the form
\[
B(T_0; T_1,\dots, T_n). 
\]
Its differential is 
\begin{multline*}
d B(T_0; T_1,\dots, T_n)
=
\sum_{j=0}^n \pm B(T_0; T_1,\dots, d T_j, \dots, T_n) 
\\+
\sum \pm
B(T_0;\dots, I(T_j,\dots, T_{j+k-1}), T_{j+k}, \dots, T_n)
 \\ +
\sum \pm
I(T_1,\dots, B(T_0; T_j,\dots, T_{j+k-1}), T_{j+k}, \dots, T_n)
\\ \pm 
E(0;\dots, I(T_j,\dots, T_{j+k-1}), T_{j+k}, \dots, T_n)
\circ_0 
T_0
.
\end{multline*}
The symbol $\circ_0$ in the last term denotes the operadic insertion of $T_0$ at the first slot (the one labelled by $0$). It is given by a sum of terms according to the operadic composition rules in $\Br$.
We want to show that the subspace associated to any tree occurring in the above expression is contained in $U_\Gamma$.
For the first three sums this is done as in the case of external and internal vertices before.
So let us concentrate on the last term, involving the operadic composition. First note that it is sufficient to consider the case of all subtrees $\Gamma_1, \dots, \Gamma_n$ being single external vertices (using the notation from the definition). Hence $U_\Gamma\cong \tilde U$.
We defined $\tilde U$ as triple intersection of inverse images. So one has to check that the projections of the $U_{\tilde{\Gamma}}$ (for $\tilde{\Gamma}$ some graphs produced by the differential) under each of the forgetful maps is still in the spaces indicated. Start with $U_{\Gamma_0}$. We have to check that $\pi U_{\tilde \Gamma}\subset U_{\Gamma_0}$, where $\pi$ forgets the vertices in $\Gamma_1, \dots, \Gamma_n$ (same notation as in the definition again). This follows from the fact that for each vertex type, forgetting ``downstairs points'' does not take us out of the respective subspace of configurations. For the other two terms of the intersection we leave the proof to the reader. 
Now consider the third assertion. Tracking the definitions it amounts to showing that for a blue vertex, the configurations obtained by downscaling some configuration in $U_{\Gamma_0}$ to a point and inserting is still contained in $U_{\Gamma}$. This is clear.
Finally, the last statement of the Lemma is obvious.
\end{proof}

\begin{figure}
 \centering
\makeatletter{}\usetikzlibrary{arrows}
\[
\begin{tikzpicture}[
scale=.7,
int/.style={circle, draw, fill, minimum size=5pt, inner sep=0},
ext/.style={circle, draw, fill=white, minimum size=5pt, inner sep=1pt},
helper/.style={coordinate},point/.style={circle, draw, fill, inner sep =1pt},
de/.style={-triangle 60},
point/.style={circle, draw, fill, minimum size=3pt, inner sep=0pt},
]
\begin{scope}[yshift=2cm]
\node at (0,4) {$\Gamma$};
\begin{scope}[yshift=.5cm]
\draw (0,2)--(0,1.5) node[ext] (e1) {1} --++(0,-1) node[ext] (e2){2}--+(0,-1) node[ext](e3) {3};
\end{scope}
\end{scope}

\begin{scope}[xshift=8cm, shift={(-3,0)}]
\node at (-3,6) {$U_\Gamma^{\mathit{KS}}$};
\draw[dashed] (-4,2)--(4,2);
\draw[dashed] (-4,3.5)--(4,3.5);
\draw[triangle 60-triangle 60] (-3.8,3)--(-3.8,4);

\node[point,label=-90:2] (e2) at (2,3.5) {};
\node[point,label=-90:3] (e3) at (-2.2,2) {};
\draw[triangle 60-triangle 60] (e3)+(-0.4,0.2)--+(0.4,0.2);
\draw[triangle 60-triangle 60] (e2)+(-0.4,0.2)--+(0.4,0.2);
\node[point, label=0:{$1$}] at (0,5) {};

\end{scope}

\begin{scope}[xshift=20cm, shift={(-6,0)}]
\node at (-3,6) {$U_\Gamma$};
\draw[dashed] (-4,2)--(4,2);

\draw[triangle 60-triangle 60] (-0.7,3)--(-2.1,3);

\node[point,label=-45:2] (e2) at (-1.3,2.5) {};
\node[point,label=-90:3] (e3) at (-1.5,2) {};
\draw[dashed] (-0.6,2) arc (0:180:0.8);
\draw[<->] (e3)+(-0.2,0.2)--+(0.2,0.2);
\node[point, label=0:{$1$}] at (0,5) {};

\end{scope}

\end{tikzpicture}
\] 
\caption{\label{fig:ksbigger} An example of a graph $\Gamma$ for which the Kontsevich-Soibelman subspace $U_\Gamma^{KS}$ is bigger than our $U_\Gamma$.  }
\end{figure}
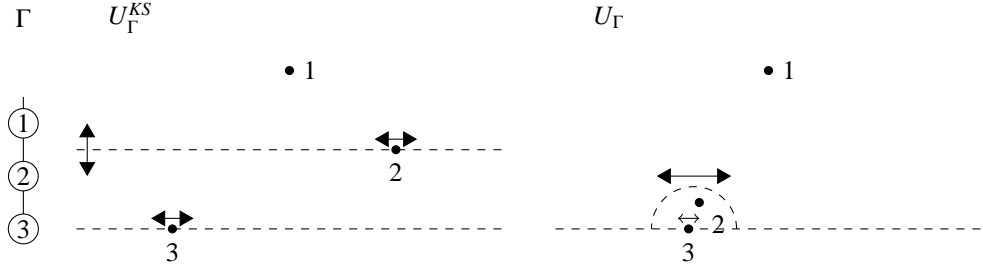

\begin{rem}
 These space are almost the same as the $U_\Gamma$ defined by Kontsevich and Soibelman. However, in some cases they are smaller, see Figure \ref{fig:ksbigger}. 
\end{rem}

\subsection{Starting the recursion and an explicit formula for $\Br$-trees.}
We still have to start the recursion, i.~e., solve equation \eqref{equ:Brcase} for $\Gamma$ of degree $0$ and $-1$. This is not difficult to do explicitly. However, we can show a bit more:
\begin{prop}
\label{prop:explonBr}
 There is a solution of \eqref{equ:Brcase} such that for $B(\Br)$-graphs $\Gamma \in \Br\subset B(\Br)[-1]$ the chain $c(\Gamma)$ is, up to sign, the fundamental chain of $U_\Gamma$.
\end{prop}
\begin{proof}
For graphs $\Gamma \in \Br(n)$ set $c(\Gamma)$ to be the fundamental chain $c(U_\Gamma)$ of $U_\Gamma$. To define the orientation it suffices to define the orientation on the subspaces $U_I^\cdot$ and $U_E^\cdot$ used to build $U_\Gamma$. The spaces $U_I^\cdot$ consist of configurations of points $z_1,z_2,z_3,\dots$ on a line, modulo translation and scaling. If we use the translation and scaling degrees of freedom to fix $z_1=0$, $z_2=1$, then the orientation on $U_I^\cdot$ is such that the form $dz_3\wedge dz_4\wedge\dots$ is positive. Similarly $U_E^\cdot$ consists of configurations of one point $z_0$ above points $z_1,z_2,\dots$ on a line, modulo translation and scaling. The orientation is defined as follows. If we use the translation and scaling to fix $z_0=i$ and $z_1,z_2,\dots$ to be real, 
then $dz_1\wedge dz_2\wedge dz_3\wedge\dots$ shall be positive.

We claim that with this assignment the chain $c(\Gamma)$ solves eqn. \eqref{equ:Brcase}. Since $\Gamma$ does not contain any blue vertices, the second term in that equation vanishes. 
Hence we have to show that
\[
 \p c(\Gamma) = \p c(U_\Gamma) = c(d \Gamma). 
\]
Since $c(\Gamma)$ is defined in a recursive manner, using the operadic insertions, it is enough to check the above statement for generators, i.~e., $\Gamma=T_n$ or $\Gamma=T_n'$, see Figure \ref{fig:KSgen}.\footnote{We want to state clearly that the map $\Br \ni \Gamma\mapsto c(\Gamma) \in C(\FM_2)$ is not a map of operads. However, to define the $c(\Gamma)$ we are using only (some of) the operadic compositions. Hence the compatibility with the differential can be reduced to a verification on generators, by using that the differentials are compatible with the operadic compositions.} In these cases $c(\Gamma)$ is the fundamental chain on some $U_E^\cdot$ or $U_I^\cdot$. Consider the $U_E^\cdot$ case. There are two sorts of boundary strata: either some subset of the points $z_2,z_3,\dots$ comes close together, away from $z_0$, or some subset, together with $z_0$ comes close together. These strata exactly match with the graphs in $d T_n$. In the expression $d T_n$ the two sorts of graphs are those with the internal vertex below or on top the external vertex. It is easy to see that the chains assigned to the two-level trees in $d T_n$ also match with the boundary strata of $U_E$, up to possibly sign. Checking the sign is tedious, and we leave it to the reader.
A similar (slightly simpler) argument goes through for the case of $\Gamma=T_n'$ and $c(\Gamma)$ being the fundamental chain of $U_I^\cdot$. 

So far we have defined $c(\Gamma)$ for trees $\Gamma\in \Br$. 
We want to extend the definition of $c(\cdots)$ by induction on the degree as discussed before. To this end we have to define $c(\Gamma)$ on the degree $0$ elements and check on the degree $-1$ elements that no obstruction appears on the right-hand side of \eqref{equ:Brcase}. Since $\Br$ is concentrated in non-negative degrees all $B(\Br)$-trees of degree 0 are in fact $\Br$-tress. Thus have just defined $c(\Gamma)$ (in particular) for all $\Gamma$ of degree 0.
It remains to check by hand that, for $B(\Br)$-trees $\Gamma$ of degree -1 which are not in $\Br$, i.e., which contain at least one blue vertex, the right-hand side of \eqref{equ:Brcase} vanishes. However, such $\Gamma$ are very simple: they contain exactly one blue vertex, and all external vertices are leaves. There are 2 terms on the right-hand side of \eqref{equ:Brcase}. One of them involves the operadic composition in $\Br$ of two degree zero trees, yielding a single degree zero tree in $\Br$, to which we associate a degree-0-chain, namely a point with coefficient (depending on the orientation) $\pm 1$. The other term similarly involves the operadic composition of two 0-chains (given by single points each, with coefficients $\pm 1$) in $C(\FM_2)$. The result is again a point with coefficient $\pm 1$. The two points cancel and hence setting $c(\Gamma)=0$ in these situations solves \eqref{equ:Brcase}, cf. Example \ref{ex:ksmapainfty} below.
\end{proof}

\begin{ex}
 \label{ex:ksmapainfty}
The braces operad $\Br$ contains the $A_\infty$ (homotopy associative) operad as the suboperad spanned by trees all of whose external vertices are leaves. Hence there is an inclusion
\[
 \Omega(B(A_\infty)) \to \Br_\infty.
\]
One can check that the resulting composition $\Omega(B(A_\infty)) \to \Br_\infty\to C(\FM_2)$ is the same as the natural map
\[
 \Omega(B(A_\infty)) \to A_\infty \to C(\FM_1) \to C(\FM_2).
\]
Here $A_\infty$ is naturally identified with the operad of strata of $C(\FM_1)$. 
\end{ex}

\subsection{The map \texorpdfstring{$\KS_\infty\to C(\cFM_2)$}{KS-infinity to C(EFM2)}}
\label{sec:Ks1map}
Define the topological (or rather, semialgebraic) colored operad
\[
\cFM_2  = \bpm \FM_2 & \FM_{2,1} \epm.
\]
In \cite{KS2} Kontsevich and Soibelman sketch the construction of a map 
\[
 \KS_\infty \to C(\cFM_2) 
\]
extending the map $\Br_\infty\to C(\FM_2)$ defined above. This construction is very similar to the one described in section \ref{sec:brtofm2}.
However, we think there is an oversight in the Kontsevich-Soibelman treatment, which we will fix below.
We want to construct, for each tree $\Gamma_1\in B(\KS_1)(n)$ a chain $c_1(\Gamma_1)\in C(\FM_{2,1})$. It has to satisfy an equation of the form 
\begin{equation}
\label{equ:KScase}
 \p c_1(\Gamma_1) = c_1(d_{B(\KS_1)} \Gamma_1) +\sum \pm c_1(\Gamma_1') \circ c_1(\Gamma_1'') + \sum\pm c_1(\Gamma_1') \circ c(\Gamma'').
\end{equation}
Here the notation on the right is as follows. The differential on $B(\KS_1)$ is denoted by $d_{B(\KS_1)}$. Take the (restricted) cocompositions of $\Gamma_1$. Because the cooperad $B(\KS)$ is colored, there will be two kinds of terms. One kind contains a cocomposition into two graphs $\Gamma_1', \Gamma_1''$ in $B(\KS_1)$. The other kind contains cocomposition into a graph $\Gamma''$ in $B(\Br)$, and a graph $\Gamma_1'$ in $B(\KS_1)$. That is how the elements $\Gamma'$, $\Gamma_1'$ and $\Gamma_1''$ on the right should be understood. The ``$\circ$'' on the right shall denote compositions in $C(\cFM_2)$.
Again, we want to solve \eqref{equ:KScase} by a recursion on the degree of $\Gamma_1$. Note that all arguments to $c_1$ occuring on the right hand side have degree strictly larger than $\Gamma_1$, and are hence known at this stage of the recursion. 
Furthermore the right hand side is closed by the induction assumption. The question is whether the right hand side is exact.
At this point one wants to copy the Kontsevich-Soibelman trick from the previous section, and define certain subsets $V_{\Gamma_1}\subset \FM_{2,1}$, and require that $c_{\Gamma_1}\in C(V_{\Gamma_1})$.
Then eqn. \eqref{equ:KScase} could be solved provided that 
\begin{itemize}
 \item For any graph $\tilde \Gamma$ occuring nontrivially in the expression $d_{B(\KS_1)} \Gamma_1$, we have $V_{\tilde \Gamma}\subset V_{\Gamma_1}$.
 \item For graphs $\Gamma_1', \Gamma_1''$ as in the first sum of \eqref{equ:KScase}, we have $V_{\Gamma_1'}\circ V_{\Gamma_1''} \subset V_{\Gamma_1}$, where $\circ$ is defined similarly to the $\circ$ in \eqref{equ:KScase}.
 \item For graphs $\Gamma_1', \Gamma''$ as in the second sum of \eqref{equ:KScase}, we have $V_{\Gamma_1'}\circ U_{\Gamma} \subset V_{\Gamma_1}$.
 \item $V_{\Gamma_1}$ is contractible.
 \item The assignment $\Gamma_1\mapsto V_{\Gamma_1}$ is equivariant with respect to the symmetric group actions.
\end{itemize}
If we could define such $V_{\Gamma_1}$ we were done (almost) immediately. We will define $V_{\Gamma_1}$ satisfying the first three items below, but the author does not know how to define the $V_{\Gamma_1}$ such that they are also contractible. This is likely an oversight in the constructions of Kontsevich and Soibelman. Hence, we have to live with non-contractible $V_{\Gamma_1}$, and hence there might potentially be obstructions to the exactness of the right hand side of  \eqref{equ:KScase}, parameterized by $H(V_{\Gamma_1})$. Fortunately one can define $V_{\Gamma_1}$ such that the obstructions are present only in low degrees and can be controlled. So we replace the fourth requirement for $V_{\Gamma_1}$ by the following:
\begin{itemize}
 \item $V_{\Gamma_1}$ may be non-contractible, but the potential obstructions generated can be controlled.
\end{itemize}

\subsection{The definition of the \texorpdfstring{$V_{\Gamma}$}{VGamma}}
Let us use a recursive definition similar to the definition of the $U_\Gamma$ from section \ref{sec:UGammarecursive}.\footnote{Again, there is another possible nonrecursive definition, which is shorter to write down, but yields spaces $V_{\Gamma}$ that are slightly bigger and not as simple to handle in proofs. So we prefer the recursive, slightly lengthier definition.}
Let $\Gamma$ be a $B(\KS_1)$-graph.\footnote{We drop the subscript from $\Gamma_1$ here for ease of notation.} Let us use the ``functional'' notation for elements of $B(\KS)$ from Section \ref{sec:opKS1}. So $\Gamma$ corresponds to some functional expression $F$.
More concretely, $F$ can have the form $K(\cdots)$ or $K_B(\cdots)$.
\begin{enumerate}
 \item Suppose $F=K(T_0,T_1,\dots,T_n)$ where the $T_j$ are functional expressions describing sub-trees $\Gamma_0,\dots, \Gamma_n$. 
Note that one of the $T_j$ will contain the terminal symbol $\vin$. Treat it in the same way as the other terminal symbols, i.e., $1,2,\dots$ for now. Let $k\in\{0,1,\dots,n\}$ be such that $\vin$ is contained in $T_k$.
Then we set 
\[
 V_\Gamma' = V^{n,k}(U_{\Gamma_0}, \dots, U_{\Gamma_n} ).
\]
Here the notation $V^{n,k}(U_{\Gamma_0}, \dots, U_{\Gamma_n} )$ denotes the operadic insertion. The subspaces $V^{n,k}\subset \FM_{2,1}(n+1)$ have been defined in section \ref{sec:FM21}, they are depicted in figure \ref{fig:VK}.
Note that configurations in $V_\Gamma'$ thus produced have one point in excess, namely the one labelled by $\vin$. We hence set
\[
 V_\Gamma = \pi_\vin V_\Gamma'
\]
where $\pi_\vin$ is the forgetful map forgetting the location of the point labelled by $\vin$.

\item Suppose $F=K(\mathbb{1},T_1,\dots,T_n)$ where the $T_j$ are functional expressions describing sub-trees $\Gamma_1,\dots, \Gamma_n$. 
So the situation is the same as before, except that instead of the tree $T_0$ there is the special symbol $\mathbb{1}$ in the first slot. We define $V_\Gamma'$ in exactly the same way we did before, treating the symbol $\mathbb{1}$ as one extra terminal symbol (like $1,2,\dots$). So the resulting configurations will contain two points in excess, namely those labelled by $\mathbb{1}$ and by $\vin$. We then define
\[
 V_\Gamma = \pi_{\mathbb{1}} \pi_\vin V_\Gamma'
\]
where $\pi_{\mathbb{1}}$ is the forgetful map forgetting the location of the point labelled by ${\mathbb{1}}$.

\item Suppose $F=K_B(S; T_0,T_1,\dots,T_n)$, where the $T_j$ are functional expressions describing sub-trees $\Gamma_0,\dots, \Gamma_n$ and $S$ is a functional expression describing some $B(\KS_1)$ graph $\tilde \Gamma$. Let $k\in\{0,1,\dots,n\}$ again be such that $\vin$ is contained in $T_k$.
Consider the $\FM_{2,1}(p+n+1)\times S^1$. Think of the $S^1$ term as one extra direction. There are two forgetful maps
\[
  \FM_{2,1}(n+1) \stackrel{\pi'}{\leftarrow}  \FM_{2,1}(p+n+1)\times S^1 \stackrel{\tilde \pi}{\rightarrow} \FM_{2,1}(p).
\]
The map $\pi'$ forgets the position of the first $p$ points and the output direction and takes the direction from the $S^1$ factor as new output direction.
The map $\tilde\pi$ forgets the position of the last $n+1$ points and the output direction and takes the direction from the $S^1$ factor as new output direction.
Here we use the description of $\FM_{2,1}$ from Remark \ref{rem:fm21mods1}.
We define an auxiliary space 
\[
 \tilde V' = (\pi')^{-1}(V^{n,k})
\cap
\tilde\pi^{-1}(V_{\tilde\Gamma})
 \cap_{j=0}^n\cap_{\alpha\in \mathit{Vert}(\tilde \Gamma)} \pi_{j\alpha }^{-1}(H^+).
\]
Here $H^+\subset \FM_{2,1}$ is composed of configuration $[(z_1,z_2,\zeta)]$ with $|z_1|\geq |z_2|$, see section \ref{sec:FM21}. $\mathit{Vert}(\tilde \Gamma)$ is the set of (labels of) external vertices in $\tilde \Gamma$.
The last of the three factors in the intersection defining $\tilde V'$ forces the additional $n+1$ points to be put ``outside of'' configurations in $V_{\tilde\Gamma}$.
At this stage, configurations in $\tilde V'$ include the auxiliary orientation (due to the extra $S^1$ above). We simply set
\[
 \tilde V := \pi_{S^1} \tilde V'
\]
where $\pi_{S^1}: \FM_{2,1}(p+n+1)\times S^1 \to \FM_{2,1}(p+n+1)$ forgets the $S^1$ factor.
Finally we set 
\[
 V_\Gamma := \pi_{\vin} (V_{\Gamma}(U_{\Gamma_0}, \dots, U_{\Gamma_n} ))
\]
where the notation $V_{\Gamma}(\cdots)$ denotes the right $\FM_2$-action. The forgetful map $\pi_\vin$ again forgets the position of the point corresponding to the symbol $\vin$.

\item The remaining case $F=K_B(\mathbb{1},T_1,\dots,T_n; S)$ is handled analogously, by first inserting and then forgetting an auxiliary point labelled $\mathbb{1}$ as before.
\end{enumerate}

\begin{lemma} \label{lem:VGammaprops}
The spaces $V_{\Gamma_1}$ defined above satisfy the following properties.
\begin{enumerate}
 \item For any graph $\tilde \Gamma$ occuring nontrivially in the expression $d_{B(\KS_1)} \Gamma_1$ (see eqn. \eqref{equ:KScase}), we have $V_{\tilde \Gamma}\subset V_{\Gamma_1}$.
 \item For graphs $\Gamma_1', \Gamma_1''$ as in the first sum of \eqref{equ:KScase}, we have $V_{\Gamma_1'}\circ V_{\Gamma_1''} \subset V_{\Gamma_1}$, where $\circ$ is defined similarly to the $\circ$ in \eqref{equ:KScase}.
 \item For graphs $\Gamma_1', \Gamma''$ as in the second sum of \eqref{equ:KScase}, we have $V_{\Gamma_1'}\circ U_{\Gamma''} \subset V_{\Gamma_1}$.
 \item The assignment $\Gamma_1\to V_{\Gamma_1}$ is equivariant with respect to the symmetric group actions.
\end{enumerate}
\end{lemma}

The proof is done by an induction similar to Lemma \ref{lem:UGammaprops}.

However, the spaces $V_\Gamma$ are \emph{not} contractible.\footnote{In fact, they better should not be, since we want to map the graph corresponding to the Connes-Rinehart differential $B$, depicted in Figure \ref{fig:dliota} (left), to the circle.}
This means that solving equation \eqref{equ:KScase} is not as simple as solving equation \eqref{equ:Brcase}. Concretely, there might be obstructions, indexed by homology classes of $V_\Gamma$.
Fortunately, one can get $H(V_\Gamma)$ and the obstructions under control.

\subsection{Studying potential obstructions \texorpdfstring{in \eqref{equ:KScase}}{} }
\begin{prop}
\label{prop:VGammahom}
 Let $\Gamma$ be a $B(\KS_1)$-graph with functional expression $F$. Let $N$ be the number of occurrences of the functional symbols $K(\cdots)$ or $K_B(\cdots)$ in $F$. Then $V_\Gamma$ is homotopic to a product $(S^1)^p\times(S^1\vee S^1)^q$ for some $p.q$ such that $p+q\leq N$, or to a point $\{pt\}$.
More precisely, $p$ and $q$ can be computed as follows
\begin{enumerate}
 \item $p$ is the number of functional symbols $K$ occuring in $F$ in the form $K(T_0,\cdots, T_n)$ or $K_B$ in the form $K_B(S; T_0, \cdots, T_n)$  where either exactly one of the $T_j$ contains no external vertices (i.e., terminal symbols $1,2\dots$), or $n=1$ and neither $T_0$ nor $T_1$ contain external vertices. \item $q$ is the number of functional symbols $K$, $K_B$ as above, for $n\geq 2$, for which two $T_j$'s contain no external vertices.
\end{enumerate}
\end{prop}
Note that any $T_j$ has to contain an external vertex unless $T_j=\vin$ or $j=0$ and $T_j=\mathbb{1}$.
In the proof let us call one occurrence of $K$ or $K_B$ as above a ``layer''. This notation should be natural since in the way we constructed $V_\Gamma$, each such symbol contributed one layer of points to a configuration.
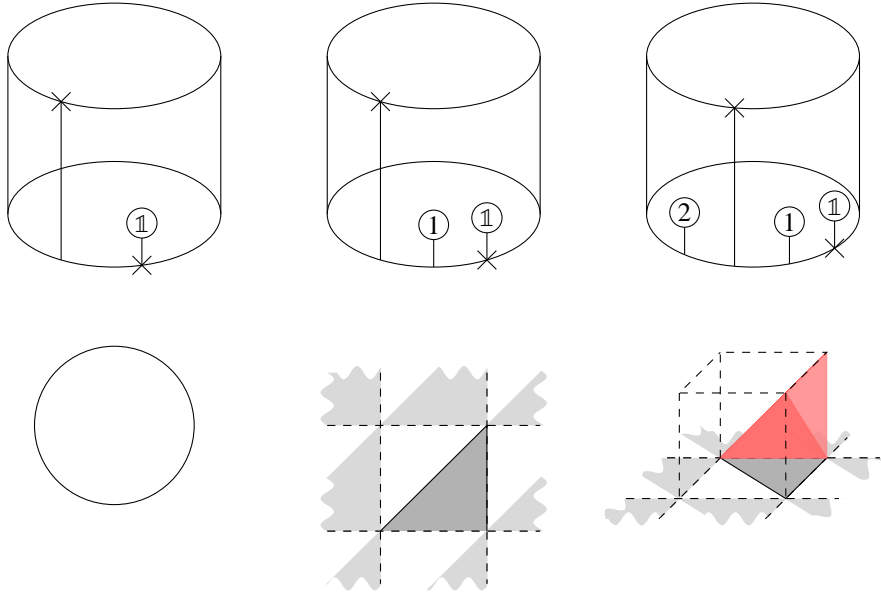
\begin{figure}
 \centering
\makeatletter{}\[
\begin{tikzpicture}[
scale=.7,
int/.style={circle, draw, fill, minimum size=5pt, inner sep=0},
ext/.style={circle, draw, fill=white, minimum size=5pt, inner sep=1pt},
helper/.style={coordinate},point/.style={circle, draw, fill, inner sep =1pt},
de/.style={-triangle 60},
point/.style={circle, draw, fill, minimum size=3pt, inner sep=0pt},
xst/.style={cross out, draw, minimum size=5 },
]
\begin{scope}[yshift=2cm]
\draw (0,0) ellipse (2cm and 1cm);
\draw (0,3) ellipse (2cm and 1cm);
\draw (-2,0)--(-2,3) (2,0)--(2,3);
\node [xst] (out) at ($(0,0)+(-75:2 and 1)$) {};
\node [xst] (in) at ($(0,3)+(-120:2 and 1)$) {};
\draw (out.base)--+(0,.8) node[ext] {$\mathbb{1}$};
\draw (in.base)--($(0,0)+(-120:2 and 1)$) ;
\draw (0,-4) circle (1.5);
\end{scope}
\begin{scope}[xshift =6cm, yshift=2cm]
\draw (0,0) ellipse (2cm and 1cm);
\draw (0,3) ellipse (2cm and 1cm);
\draw (-2,0)--(-2,3) (2,0)--(2,3);
\node [xst] (out) at ($(0,0)+(-60:2 and 1)$) {};
\node [xst] (in) at ($(0,3)+(-120:2 and 1)$) {};
\draw (out.base)--+(0,.8) node[ext] {$\mathbb{1}$};
\draw (in.base)--($(0,0)+(-120:2 and 1)$) ;
\draw ($(0,0)+(-90:2 and 1)$)--+(0,0.8)node[ext] {1};

\end{scope}

\begin{scope}[xshift =12cm, yshift=2cm]
\draw (0,0) ellipse (2cm and 1cm);
\draw (0,3) ellipse (2cm and 1cm);
\draw (-2,0)--(-2,3) (2,0)--(2,3);
\node [xst] (out) at ($(0,0)+(-40:2 and 1)$) {};
\node [xst] (in) at ($(0,3)+(-100:2 and 1)$) {};
\draw (out.base)--+(0,.8) node[ext] {$\mathbb{1}$};
\draw (in.base)--($(0,0)+(-100:2 and 1)$) ;
\draw ($(0,0)+(-70:2 and 1)$)--+(0,0.8)node[ext] {1};
\draw ($(0,0)+(-130:2 and 1)$)--+(0,0.8)node[ext] {2};

\end{scope}

\begin{scope}[xshift=7cm, yshift=-3cm]
\draw[dashed] (-2,2)--(-2,-2) (0,2)--(0,-2) (-3,1)--(1,1) (-3,-1)--(1,-1);
\clip[decorate, decoration={coil, aspect=0}] (-3,2) rectangle (1,-2);
\foreach \x/\y in {-2/-1, -4/-1,0/-1,-2/1,0/1,-2/-3,-4/-3,-4/1}
  \draw[fill=black!15,draw=black!15] (\x,\y)--+(2,2)--+(2,0);
 \draw[fill=black!30] (-2,-1)--+(2,2)--+(2,0);
\draw[dashed] (-2,2)--(-2,-2) (0,2)--(0,-2) (-3,1)--(1,1) (-3,-1)--(1,-1);

\end{scope}

\begin{scope}[xshift=13cm, yshift=-3cm]
 
\begin{scope}[canvas is xz  plane at y=0]
\draw[dashed] (-2,2)--(-2,-2) (0,2)--(0,-2) (-3,1)--(1,1) (-3,-1)--(1,-1);
\clip[decorate, decoration={coil, aspect=0}] (-3,2) rectangle (1,-2);
\foreach \x/\y in {-2/-1, -4/-1,0/-1,-2/1,0/1,-2/-3,-4/-3,-4/1}
  \draw[fill=black!15,draw=black!15] (\x,\y)--+(2,2)--+(2,0);
 \draw[fill=black!30] (-2,-1)--+(2,2)--+(2,0);
\draw[dashed] (-2,2)--(-2,-2) (0,2)--(0,-2) (-3,1)--(1,1) (-3,-1)--(1,-1);
\end{scope}

\draw[fill=red,draw=red,opacity=.8]
 (-2,0,-1) -- (0,2,1) -- (0,2,-1) --cycle
 (-2,-0,-1) -- (0,2,1) -- (0,0,-1) --cycle;
\draw[fill=red!50, draw=red!50, opacity=.8] (-2,0,-1) -- (0,0,-1) -- (0,2,-1) --cycle
;

\draw[dashed]
(-2,0,-1)--(-2,2,-1)--(0,2,-1)--(0,2,1)--(-2,2,1)--(-2,2,-1)
(-2,2,1)--(-2,0,1)  (0,2,1)--(0,0,1);
\end{scope}

\end{tikzpicture}
\]  
\caption{\label{fig:ks1spaces} Several $\KS_1$-graphs $\Gamma$ for which $V_\Gamma$ is non-contractible. For the graph on the left, $V_\Gamma$ is a circle. For the two graphs in the middle, $V_\Gamma$ is a torus segment homotopic to a wedge of two circles. For the graph on the right, $V_\Gamma$ is  three-dimensional subspace of $S^1\times S^1 \times [0,1]$ (the red tetrahedron). Here the drawing should be periodically continued as indicated, However, note that the upper and lower plane, i.e., the boundaries of $[0,1]$, are \emph{not} identified. The space is homotopic to a wedge of two circles. }
\end{figure}

\begin{proof}[Sketch of proof]
 One performs an induction on the number of layers. For one layer, i.e., for $F$ of the form $F=K(T_0, \dots ,T_n)$ one can check that 
$V_\Gamma$ is contractible whenever all $T_0, \dots ,T_n$ have (external) vertices. If exactly one of the $T_j$'s contains no external vertices, then one can retract $V_\Gamma$ to a subspace where all external vertices are infinitesimally close together, in some fixed configuration. This space is then a circle. In case there are two $T_j$'s with no external vertices, they separate the other  $T_j$'s into two sets (possibly empty). One can retract to a subspace where the points belonging to each set are infinitesimally close to each other. There remain three cases: (i) both sets empty (ii) only one set empty (iii) both sets non-empty. Pictures of the resulting spaces are drawn in Figure \ref{fig:ks1spaces}, from which it is clear that in case (i) the space is homotopic to a circle and in cases (ii) and (iii) to a wedge of two circles. Next, assume the Proposition is true for $<N$ layers. Let $F$ have the form 
\[
 F = K_B(S; T_0, \dots, T_n).
\]
Then one can deform the space $V_\Gamma$ to the subspace $V_\Gamma^\epsilon$ consisting of those configurations in which points corresponding to terminal symbols in $S$ have absolute values at least a factor $(1+\epsilon)$ bigger than points from $T_0, \dots, T_n$. This space $V_\Gamma^\epsilon$ has a product structure, $V_\Gamma^\epsilon \cong V_1 \times V_2$, where $V_1$ is homotopic to $V_{\tilde{\Gamma}}$, for $\tilde\Gamma$ the $B(\KS_1)$-graph described by $S$ and $V_2$ is homotopic to $V_{\Gamma'}$ for $\Gamma'$ the $\KS_1$ graph defined by $K(T_0, \dots, T_n)$. For $V_1$ one uses the induction hypothesis to see that it is homotopic to a torus times wedge sums of circles. For $V_2$ one uses the same considerations as in the $N=1$-case.
\end{proof}

Suppose that we want to solve \eqref{equ:KScase} for a given $\Gamma$ of degree $-k$. Then possible obstructions are given by homology classes 
in $H_{-k+1}(V_\Gamma)$. The proposition says that this space is empty, unless $\Gamma$ contains at least $k-1$ layers. In this case $\Gamma$ has degree $-k\leq -2k+3$, since each layer contributes at most $-1$ to the degree. Hence obstructions can only occur if $k\leq 3$. 
These low-degree cases we treat as follows.

\begin{enumerate}
 \item If $\Gamma$ is of degree 0, the space $V_\Gamma$ consists of a single point. We set $c_1(\Gamma)$ to be the chain composed of this point. The sign is determined similarly to the $\Br$-case in the previous subsection.
 \item If $\Gamma$ is of degree -1, it contains either two layers of degree 0, or one of degree -1. In both cases one sees that the obstruction on the right hand side of \eqref{equ:KScase} vanishes. Concretely, the obstruction is a (signed) sum of points, which can be seen to always occur in pairs with different signs. In fact, for the case of one layer, $c_1(\Gamma)$ can be taken to be the fundamental chain of $V_\Gamma$.
 \item If $\Gamma$ is of degree -2, the potential obstructions live in degree -1. Hence it suffices to consider cases in which $V_\Gamma$ has homology in degree -1.\footnote{Note that we use a negative grading for homology to be consistent.} In the case of one layer, we can solve \eqref{equ:KScase} by setting $c_1(\Gamma)$ to be the fundamental chain of $V_\Gamma$. In the case of two layers, one layer must have degree -1 and one degree 0. By the proposition, one can get nontrivial homology in degree -1 only if the degree -1 layer has two subtrees, at least one of which does not contain external vertices. There is a closed 1-form dual to the homology class. Concretely, it is given by $d\phi$, where $\phi$ is the angle ``between the subtrees'', see Figure \ref{fig:homdegm1}.
Integrating this form over the chain on the right hand side of \eqref{equ:KScase}, one sees that one obtains zero. 
\item Similarly one treats the remaining case of $\Gamma$ of degree -3. The potential obstruction lives in degree -2, hence one needs two layers. By the proposition, the layers have to be of degree -1 each. There is again a differential form $d\phi_1\wedge d\phi_2$ that is dual to the top homology class. One can check that the integral over the right hand side of \eqref{equ:KScase} is again zero. Here one needs to know $c_1(\Gamma')$ where $\Gamma'$ is of degree -2 and has one layer. But those $c_1(\Gamma')$ can be taken to be the fundamental class of $V_{\Gamma'}$.
\end{enumerate}
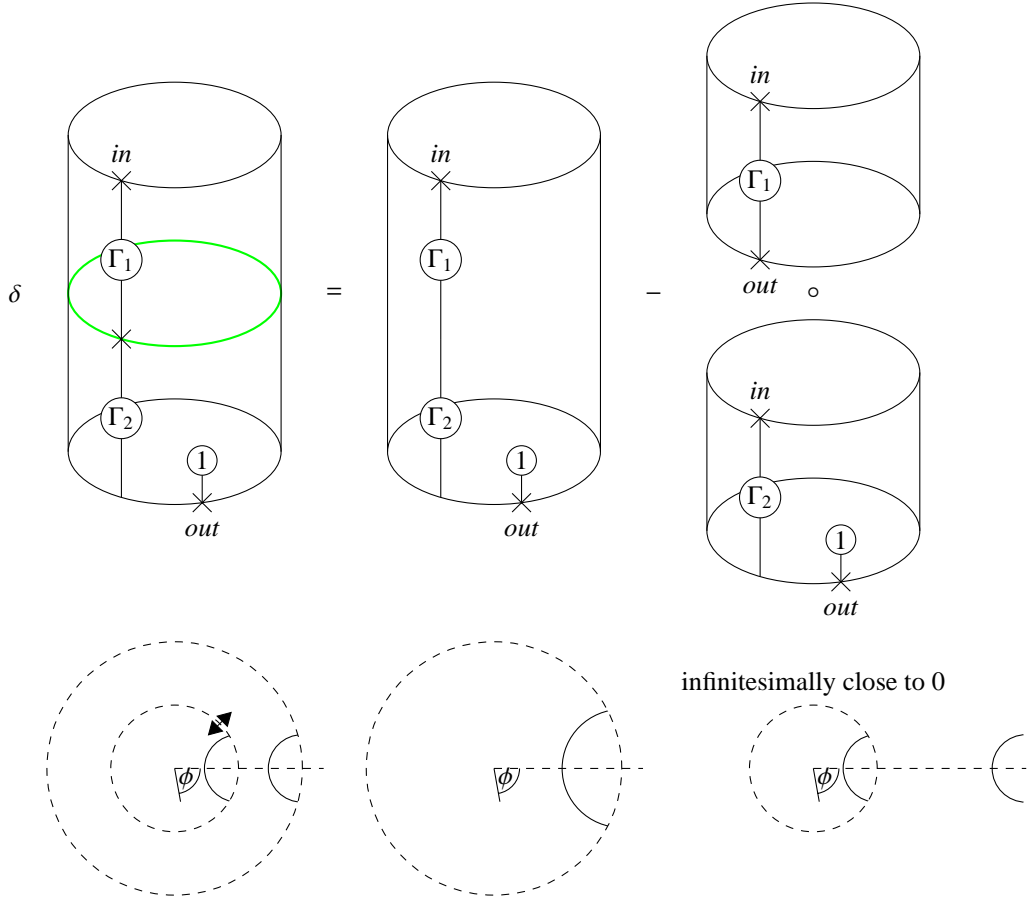
\begin{figure}
 \centering
\makeatletter{}\[
\begin{tikzpicture}[
scale=.7,
int/.style={circle, draw, fill, minimum size=5pt, inner sep=0},
ext/.style={circle, draw, fill=white, minimum size=5pt, inner sep=1pt},
helper/.style={coordinate},point/.style={circle, draw, fill, inner sep =1pt},
de/.style={-triangle 60},
point/.style={circle, draw, fill, minimum size=3pt, inner sep=0pt},
xst/.style={cross out, draw, minimum size=5 },
]
\begin{scope}[]
\node at (-3,3) {$\delta$};
\draw (0,6) ellipse (2cm and 1cm);
\draw (0,0) ellipse (2cm and 1cm);
\draw[draw=green, thick] (0,3) ellipse (2cm and 1cm);
\draw (-2,0)--(-2,6) (2,0)--(2,6);
\node [xst,label=-90:{$out$}] (out) at ($(0,0)+(-75:2 and 1)$) {};
\node [xst] (in) at ($(0,3)+(-120:2 and 1)$) {};
\node [xst,label=90:{$in$}] (inu) at ($(0,6)+(-120:2 and 1)$) {};

\draw (inu.base)--+(0,-1.5) node[ext] {$\Gamma_1$} -- (in.base);
\draw (out.base)--+(0,.8) node[ext] {$\mathrm{1}$};
\draw (in.base)--+(0,-1.5) node[ext] {$\Gamma_2$} -- ($(0,0)+(-120:2 and 1)$) ;
\node at (3,3) {$=$};
\end{scope}

\begin{scope}[xshift=6cm]
\draw (0,6) ellipse (2cm and 1cm);
\draw (0,0) ellipse (2cm and 1cm);
\draw (-2,0)--(-2,6) (2,0)--(2,6);
\node [xst,label=-90:{$out$}] (out) at ($(0,0)+(-75:2 and 1)$) {};
\node [coordinate] (in) at ($(0,3)+(-120:2 and 1)$) {};
\node [xst,label=90:{$in$}] (inu) at ($(0,6)+(-120:2 and 1)$) {};

\draw (inu.base)--+(0,-1.5) node[ext] {$\Gamma_1$} -- (in.base);
\draw (out.base)--+(0,.8) node[ext] {$\mathrm{1}$};
\draw (in.base)--+(0,-1.5) node[ext] {$\Gamma_2$} -- ($(0,0)+(-120:2 and 1)$) ;
\node at (3,3) {$-$};
\end{scope}

\begin{scope}[xshift=12cm, yshift=-1.5cm]
\begin{scope}[yshift=3cm]
\draw (0,6) ellipse (2cm and 1cm);
\draw (0,3) ellipse (2cm and 1cm);
\draw (-2,3)--(-2,6) (2,3)--(2,6);
\node [xst,label=90:{$in$}] (inu) at ($(0,6)+(-120:2 and 1)$) {};
\node [xst,label=-90:{$out$}] (outu) at ($(0,3)+(-120:2 and 1)$) {};
\draw (inu.base)--+(0,-1.5) node[ext] {$\Gamma_1$} -- (outu.base);

\end{scope}

\node at (0,4.5) {$\circ$};

\draw (0,0) ellipse (2cm and 1cm);
\draw (0,3) ellipse (2cm and 1cm);
\draw (-2,0)--(-2,3) (2,0)--(2,3);
\node [xst,label=-90:{$out$}] (out) at ($(0,0)+(-75:2 and 1)$) {};
\node [xst,label=90:{$in$}] (in) at ($(0,3)+(-120:2 and 1)$) {};

\draw (out.base)--+(0,.8) node[ext] {$\mathrm{1}$};
\draw (in.base)--+(0,-1.5) node[ext] {$\Gamma_2$} -- ($(0,0)+(-120:2 and 1)$) ;
\end{scope}

\begin{scope}[yshift=-6cm, scale=.8]
\draw[dashed] (0,0) circle (3) circle(1.5) -- (3.5,0);
\begin{scope}[]
\clip (0,0) circle (3);
\draw (3,0) node[label=180:{$U_{\Gamma_1}$}, outer sep=-8] (c) {} circle (.8);

\clip (0,0) circle (1.5);
\draw (1.5,0) node[label=180:{$U_{\Gamma_2}$}, outer sep=-8] (c) {} circle (.8);
\end{scope}

\draw[triangle 60-triangle 60] (45:1.1)--(45:1.9);

\draw (0,0) node (cc){}--(-80:.8) (-80:.6) arc (-80:0:.6);
\node at (-40:.35) {$\phi$};
\end{scope}

\begin{scope}[xshift=6cm, yshift=-6cm, scale=.8]
\draw[dashed] (0,0) circle (3) -- (3.5,0);
\begin{scope}[]
\clip (0,0) circle (3);
\draw (3,0) node[label=180:{$U_{\Gamma_2\circ \Gamma_1}$}, outer sep=-8] (c) {} circle (1.4);
\node at (-40:.35) {$\phi$};
\end{scope}

\draw (0,0) node (cc){}--(-80:.8) (-80:.6) arc (-80:0:.6);
\end{scope}

\begin{scope}[yshift=-6cm, xshift=12cm, scale=.8]
\draw[dashed] (0,0)  circle(1.5) -- (5,0);
\node at (0,2) {infinitesimally close to 0};
\begin{scope}[]
\clip (0,0) circle (5);
\draw (5,0) node[label=180:{$U_{\Gamma_1}$}, outer sep=-8] (c) {} circle (.8);

\clip (0,0) circle (1.5);
\draw (1.5,0) node[label=180:{$U_{\Gamma_2}$}, outer sep=-8] (c) {} circle (.8);
\end{scope}

\draw (0,0) node(cc){}--(-80:.8) (-80:.6) arc (-80:0:.6);
\node at (-40:.35) {$\phi$};
\end{scope}

\end{tikzpicture}
\]  
\caption{\label{fig:homdegm1} Illustration of the non-occurrence of obstructions on the right hand side of \eqref{equ:KScase} in low orders. The graph on the left is a typical graph of degree $-2$ with two layers. Its differential consists of two terms (right). Their chains are dpeicted schematically below. Note that both terms on the right are assigned nontrivial cycles (namely $S^1$), but they cancel. }
\end{figure}

Similarly to the $\Br$-case one proves the following.
\begin{prop}
\label{prop:KSexpl}
 There is a solution of \eqref{equ:KScase} such that for $B(\KS_1)$-graphs $\Gamma \in \KS_1$ the chain $c_1(\Gamma)$ is the fundamental chain of $V_\Gamma$.
\end{prop}

\subsection{An extension -- \texorpdfstring{$\hBr_\infty$}{hBr-infinity} }
\label{sec:KShomext}
In this section we want to extend the Kontsevich--Soibelman construction to obtain a map of operadic bimodules 
\[
 \hBr_\infty \to C(D_K) \, .
\]
This bimodule is quasi-free, generated by $B(\Br)$. Hence for each graph $\Gamma$ we have to construct a chain $\hc(\Gamma)\in C(D_K)$ satisfying conditions of the following form
\begin{equation}
 \label{equ:hBrcase}
 \p \hc(\Gamma) = \hc(d_{B(\Br)} \Gamma)
+
\sum \pm \hc(\Gamma')\circ c(\Gamma'')
+
\sum \pm \Gamma' \circ (\hc(\Gamma_1''), \hc(\Gamma_2''), \dots).
\end{equation}
Here the notation is as follows. In the first sum one takes the restricted co-compositions of $\Gamma$.
In the second sum one takes the full co-composition, followed by projection of the first factor ($\Gamma'$) onto $\Br$.
In the first sum we allow $\Gamma'$ to be the counit, and in the second we allow $\Gamma_1'',\Gamma_2'', \dots$ to be all counits.
 The $\circ$ in the first sum is the right action of $c(\Gamma'')\in C(\FM_2)$. The $\circ$ in the second sum is the left action
of $\Br$ on $C(D_K)$.
One can apply a variant of the Kontsevich-Soibelman trick again and construct, for each tree $\Gamma\in B(\Br)$ a subspace $W_\Gamma$ such that 
\begin{itemize}
 \item $W_\Gamma$ is contractible.
  \item The spaces $W_{\tilde \Gamma}$ for $\tilde \Gamma$ occuring in $d_{B(\Br)} \Gamma$, are contained in $W_\Gamma$.
 \item The spaces $W_{\Gamma'}\circ U_{\Gamma''}$  are contained in $W_\Gamma$.
 \item The spaces $\Gamma'\circ (W_{\Gamma_1''}, W_{\Gamma_2''},\dots)$ (notation similar to \eqref{equ:hBrcase}) are contained in $W_\Gamma$.\footnote{Here we quietly extended the left action of $\Br$ on chains on $D_K$ to a left ``action'' of $\Br$ trees on subsets of $D_K$. The changes to the definition necessary are marginal. }
 \item The assignment $\Gamma\mapsto W_{\Gamma}$ is equivariant under the symmetric group action.
\end{itemize}

There is a simple definition of the $W_\Gamma$. Recall from Lemma \ref{lem:defretract} the forgetful map $\pi:D_K\to \FM_2$, forgetting the location of the real line. We define
\[
 W_\Gamma = \pi^{-1}U_\Gamma.
\]

\begin{lemma}
 The four assertions above are satisfied.
\end{lemma}
\begin{proof}[Proof sketch]
 The fact that $W_\Gamma$ is contractible follows from (the proof of) Lemma \ref{lem:defretract} and the contractibility of $U_\Gamma$. 
The second assertion follows from the analogous assertion for $U_\Gamma$.
 The third assertion follows since the forgetful map $\pi$ is compatible with the right $\FM_2$ action. 
The fourth assertion is the most difficult. We have to show that 
\[
 \pi (\Gamma'\circ (W_{\Gamma_1''}, W_{\Gamma_2''},\dots)) \subset U_\Gamma.
\]
To see this one can take a small detour and define an ``action'' of $\Br$-trees on subsets of $\FM_2$. It is given by formulas similar to those appearing in the definition of $U_\Gamma$ above. Namely, the generators act as follows: On subsets $U_1,\dots, U_n$ the element $T_n'$ (internal vertex) acts as
\[
 T_n(U_1,\dots, U_n) = U_I(U_1,\dots, U_n).
\]
Here we use the same notation as in the definition of $U_\Gamma$. The other generator, $T_n$ acts in the same way as a blue vertex in the definition of $U_\Gamma$ above, we just replace $U_{\Gamma_j}$ by $U_j$ where $U_0,\dots, U_n$ are again subsets.
Then, more or less by definition 
\[
U_\Gamma = \Gamma'\circ (U_{\Gamma_1''}, U_{\Gamma_2''},\dots). 
\]
However, comparing the $\Br$-tree action on subsets of $D_K$ and on subsets of $\FM_2$, one sees that they are intertwined by $\pi$. Hence the result follows.
\end{proof}

\begin{figure}
 \centering
\makeatletter{}\usetikzlibrary{arrows}
\[
\begin{tikzpicture}[
scale=.6,
int/.style={circle, draw, fill, minimum size=5pt, inner sep=0},
ext/.style={circle, draw, fill=white, minimum size=5pt, inner sep=1pt},
helper/.style={coordinate},point/.style={circle, draw, fill, inner sep =1pt},
de/.style={-triangle 60},
point/.style={circle, draw, fill, minimum size=3pt, inner sep=0pt},
]
\begin{scope}[yshift=3cm, shift={(13,4)}]
\draw (0,2)--(0,1.5) node [ext] (e1) {1} --++(-1,-1) node [int] (i1) {}--+(-0.5,-1) node [ext] (e2) {2};
\draw (i1)--+(0.5,-1) node [ext] (e3) {3};
\draw (e1)--++(1,-1) node [int] (i2) {} --+(-0.5,-1) node [ext] (e4) {4};
\draw (i2)--+(0,-1) node [ext] (e5) {6};
\draw (i2)--+(0.5,-1) node [ext] (e6) {5};
\end{scope}

\begin{scope}[xshift=10cm, shift={(-2,0)}]
\draw (-5,0.5)--(5,0.5);
\draw[dashed] (-4,2)--(4,2);
\node[point,label=-90:3] (e3) at (-1.8,2) {};
\node[point,label=-90:2] (e2) at (-2.2,2) {};
\node[point,label=-90:4] (e4) at (2,2) {};
\node[point,label=-90:6] (c) at (2.4,2) {};
\node[point,label=-90:5] (e5) at (2.8,2) {};
\draw[dashed] (3,2) arc (0:180:0.6);
\draw[<->] (c)+(-0.2,0.2) -- +(0.2,0.2);
\draw[dashed] (-1.6,2) arc (0:180:0.4);

\draw[triangle 60-triangle 60] (-2.4,2.6)--(-1.6,2.6);

\draw[triangle 60-triangle 60] (c)+(-0.4,0.8)--+(0.4,0.8);
\node[point, label=0:{$1$}] at (0,5) {};
\end{scope}

\begin{scope}[xshift=22cm, shift={(-2,0)}]
\node at (-6,3) {$\cup$};
\draw (-5,0.5)--(5,0.5);
\node[point,label=-90:3] (e3) at (-1.8,0.8) {};
\node[point,label=-90:2] (e2) at (-2.2,0.8) {};
\node[point,label=-90:4] (e4) at (2,0.8) {};
\node[point,label=-90:6] (c) at (2.4,0.8) {};
\node[point,label=-90:5] (e5) at (2.8,0.8) {};
\draw[dashed] (3.3,0.5) arc (0:180:0.9);
\draw[<->] (c)+(-0.2,0.2) -- +(0.2,0.2);
\draw[<->] (e5)+(-0.2,0.2) -- +(0.2,0.2);
\draw[<->] (e3)+(-0.2,0.2) -- +(0.2,0.2);
\draw[dashed] (e2)+(-0.3,0) -- +(0.7,0)  (e4)+(-0.3,0) -- +(1.1,0);

\draw[dashed] (-1.3,0.5) arc (0:180:0.7);

\draw[triangle 60-triangle 60] (-2.4,1.4)--+(0.8,0);

\draw[triangle 60-triangle 60] (c)+(-0.4,0.8)--+(0.4,0.8);
\node[point, label=0:{$1$}] at (0,5) {};
\end{scope}
\end{tikzpicture}
\] 
\caption{\label{fig:wgammaex} An example of the space $W_\Gamma$ for a $B(\Br)$-tree $\Gamma$. Note that the dimension of the component on the right-hand side is larger than the dimension of the component on the left-hand side by one.}
\end{figure}
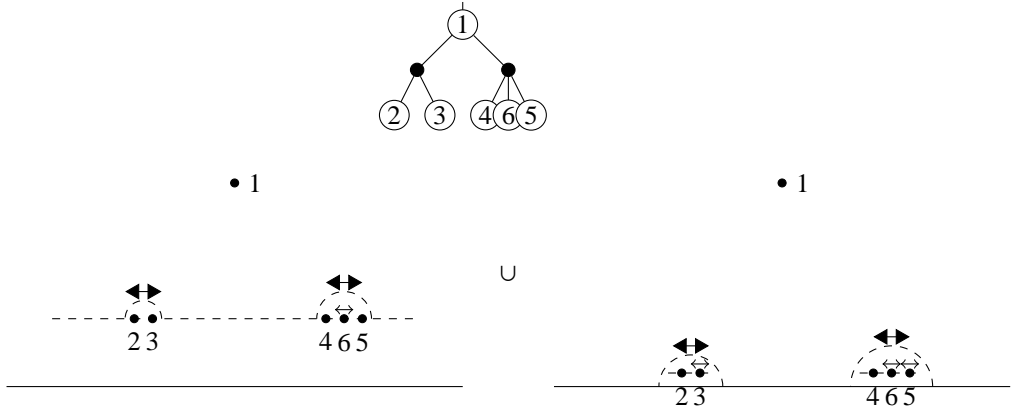

Note that solving \eqref{equ:hBrcase} is ``simpler'' than solving \eqref{equ:Brcase} because there is a smaller problem with obstructions for $\Gamma$ of degree $-1$. In fact, due to the degree shift, the only graph of degree $-1$ are binary trees, with all external vertices being leaves. For those graphs, the right hand side of \eqref{equ:hBrcase} contains exactly two terms (points), which come with opposite signs and hence yield a vanishing degree 0 homology class.

\subsection{Map of moperadic bimodules}
\label{sec:KShomext2}
Finally, we want to construct a map of moperadic bimodules
\[
 \hKS_{1,\infty} \to C(D_{S}) \, .
\]
The moperadic bimodule is (quasi-)freely generated by $B(\KS_1)$. For each graph $\Gamma\in B(\KS_1)$ we want to find a chain $\hc_1(\Gamma)\in C(D_S)$ such that the following equations are satisfied.
\begin{multline}
\label{equ:hKScase}
 \p \hc_1(\Gamma) 
=
 \hc_1(d_{B(\KS_1)} \Gamma)
+\sum \pm  c_1(\Gamma') \circ \hc_1(\Gamma'')
\\  
+\sum \pm \hc_1(\Gamma')\circ \Gamma'' \circ (\hc(\Gamma_1'''),\hc(\Gamma_2'''), \dots )
+ \sum\pm \hc_1(\Gamma') \circ c(\Gamma'').
\end{multline}

Here the notation is similar to the one in equations \eqref{equ:Brcase}, \eqref{equ:KScase} and \eqref{equ:hBrcase}. There will be three different cocompositions, corresponding to the three actions on a moperadic bimodule.
 The ``$\circ$'' in the first sum stands for the left action of $C(\FM_{2,1})$. In the second sum the ``$\circ$''s denotes the moperadic right action of $\KS_1$ on $C(D_K)$. In the third sum on the ``$\circ$'' shall denote the right action of $C(\FM_2)$.
Again, we want to solve \eqref{equ:hKScase} by a recursion on the degree of $\Gamma_1$. Note that all arguments to $c_1$ occuring on the right hand side have degree strictly larger than $\Gamma_1$, and are hence known at this stage of the recursion. 
Furthermore the right hand side is closed by the induction assumption. The question is whether the right hand side is exact.
At this point one wants to copy the Kontsevich-Soibelman trick from the previous sections, define certain subsets $X_{\Gamma}\subset D_S$, and require that $\hc_1(\Gamma)\in C(X_{\Gamma})$. 
The subsets should satisfy
\begin{itemize}
\item For any graph $\tilde \Gamma$ occuring nontrivially in the expression $d_{B(\KS_1)} \Gamma$, we have $X_{\tilde \Gamma}\subset X_{\Gamma}$.
 \item For graphs $\Gamma', \Gamma''$ as in the first sum of \eqref{equ:hKScase}, we have $V_{\Gamma'}\circ X_{\Gamma''} \subset X_{\Gamma}$, where $\circ$ is defined similarly to the $\circ$ in \eqref{equ:hKScase}.
 \item For graphs $\Gamma', \Gamma'', \Gamma_1''',\Gamma_2''',\dots,$ as in the second sum of \eqref{equ:hKScase}, we have $X_{\Gamma'}\circ \Gamma'' \circ (W_{\Gamma_1'''},W_{\Gamma_2'''}, \dots )\subset X_{\Gamma}$.
 \item For graphs $\Gamma', \Gamma''$ as in the third sum of \eqref{equ:hKScase}, we have $X_{\Gamma'}\circ U_{\Gamma''} 
\subset X_{\Gamma}$.
 \item $X_{\Gamma}$ is contractible. (We won't be able to satisfy this requirement.)
  \item The assignment $\Gamma\mapsto X_{\Gamma}$ is equivariant under the symmetric group action.
\end{itemize}

Again, we cannot satisfy the contractibility requirement and replace it with a softer demand.

\begin{itemize}
 \item $X_{\Gamma}$ can be non-contractible, but the homology and possible obstructions to \eqref{equ:hKScase} must be kept under control.
\end{itemize}

We define the spaces $X_\Gamma$ by the formula
\[
 X_\Gamma = \pi^{-1}(V_\Gamma).
\]
Then the above requirements can be verified in a similar manner as in the previous subsections.
Regarding the homology of the $X_\Gamma$ and hence possible obstructions, the following result is immediate
from Lemma \ref{lem:defretract1}:
\begin{lemma}
 $X_\Gamma$ is homotopic to $V_\Gamma$ and hence homotopic to a point or a product of a torus and wedges of circles as in Proposition \ref{prop:VGammahom}.
In particular, the homology groups of $X_\Gamma$ and $V_\Gamma$ are the same.
\end{lemma}

Let us study the possible obstructions that can appear on the right hand side of \eqref{equ:hKScase}, for some fixed $\Gamma$, of degree $-k$ and with $l$ layers.
The obstructions live in $H_{-k+1}(X_\Gamma)$.
By the lemma and proposition \ref{prop:VGammahom} the latter space is zero unless $l'\geq k-1$ of the $l$ layers have one of the forms indicated in the proposition. Each such layer contributes at least $-1$ to the total degree, so
\[
k \geq l + l'.
\]
Inserting, and using that $l\geq l'$, we get 
\[
k\geq 2l'\geq 2 (k-1)
\]
and hence $k\leq 2$ and $l'\leq 1$, $l\leq 2$. For $l=2$ there are no values of $k,l'$ satisfying the above inequalities. 
Let us study the remaining cases:

\begin{itemize}
\item For $k=l=1$, $l'=0$ the possible obstructions live in 
$H_{0}(X_\Gamma)$, i.e., are represented by points. Since $k=1$ the graph $\Gamma$ is closed and contains no blue vertices and hence the right hand side of \eqref{equ:hKScase} simplifies greatly. The first and fourth terms are zero. In the second sum, only $\Gamma'=\Gamma$ and $\Gamma''$ corresponding to the identity survive. This contributes a point.\footnote{In the corresponding configuration of points in the disk, all points are infinitesimally close to the center of the disk, i.e., $\vout$.} 
In the third sum the only remaining term is $\Gamma''=\Gamma$ and $\Gamma'$ and all $\Gamma'''$s the identity elements. Hence the total contribution of the third sum is another point.\footnote{In the corresponding configuration of points in the disk, all points are infinitesimally close to the point $\bar 0$ on the boundary of the disk.} Checking the signs, and since $X_\Gamma$ is connected, one sees that the obstruction vanishes.

\item For $k=2$, $l=l'=1$ possible obstructions live in $H_{1}(X_\Gamma)$.
The possible $\Gamma$'s have either of the forms:
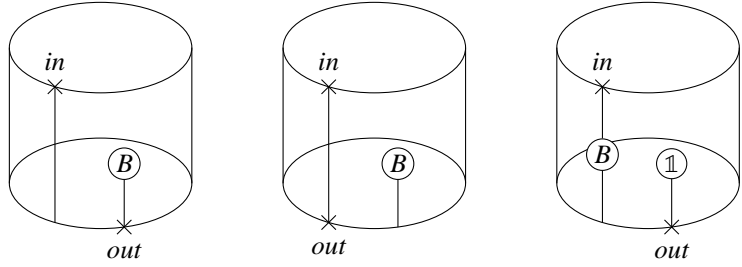
\begin{center}
\makeatletter{}\begin{tikzpicture}[
scale=.6,
int/.style={circle, draw, fill, minimum size=5pt, inner sep=0},
ext/.style={circle, draw, fill=white, minimum size=5pt, inner sep=1pt},
xst/.style={cross out, draw, minimum size=5pt, inner sep=1pt},
helper/.style={coordinate},point/.style={circle, draw, fill, inner sep =1pt},
de/.style={-triangle 60},
point/.style={circle, draw, fill, minimum size=3pt, inner sep=0pt},
]
\begin{scope}[yshift=2cm]
\draw (0,0) ellipse (2cm and 1cm);
\draw (0,3) ellipse (2cm and 1cm);
\draw (-2,0)--(-2,3) (2,0)--(2,3);
\node [xst, label=-90:{$out$}] (out) at ($(0,0)+(-120:2 and 1)$) {};
\node [] (o) at ($(0,0)+(-75:2 and 1)$) {};
\node [xst, label=90:{$in$}] (in) at ($(0,3)+(-120:2 and 1)$) {};
\draw (o.base)--+(0,1.4) node[ext] {$B$};
\draw (in.base)--($(0,0)+(-120:2 and 1)$) ;
\end{scope}
\begin{scope}[yshift=2cm, shift={(-6,0)}]
\draw (0,0) ellipse (2cm and 1cm);
\draw (0,3) ellipse (2cm and 1cm);
\draw (-2,0)--(-2,3) (2,0)--(2,3);
\node [xst, label=-90:{$out$}] (out) at ($(0,0)+(-75:2 and 1)$) {};
\node [xst, label=90:{$in$}] (in) at ($(0,3)+(-120:2 and 1)$) {};
\draw (out.base)--+(0,1.4) node [ext] {$B$};
\draw (in.base)--($(0,0)+(-120:2 and 1)$) ;
\end{scope}
\begin{scope}[yshift=2cm, shift={(6,0)}]
\draw (0,0) ellipse (2cm and 1cm);
\draw (0,3) ellipse (2cm and 1cm);
\draw (-2,0)--(-2,3) (2,0)--(2,3);
\node [xst, label=-90:{$out$}] (out) at ($(0,0)+(-75:2 and 1)$) {};
\node [xst, label=90:{$in$}] (in) at ($(0,3)+(-120:2 and 1)$) {};
\draw (out.base)--+(0,1.4) node [ext] {$\mathbb{1}$};
\draw ($(0,0)+(-120:2 and 1)$) -- +(0,1.5) node[ext] (vB) {$B$};
\draw (in.base) edge (vB) ;
\end{scope}

\end{tikzpicture} 
\end{center}
Here the $B$ symbolizes some (possible empty) binary tree of internal vertices, with external vertices and possibly $\vin$  on the leafs. For these graphs all terms on the right hand side of \eqref{equ:hKScase} can be explicitly computed.
The first term on the right hand side of \eqref{equ:hKScase} (in which the differential of $\Gamma$ appears) is either zero or contributes two line segments sharing one endpoint. The fourth term vanishes in all cases since $\Gamma$ does not contain blue vertices. 
In the second sum only one term survives (since $\Gamma$ has only one layer). This is always a circle. Also in the third sum there is only one term, which contributes either a circle or a line segment that closes the two line segments from the first term to a loop. The loop or circle is homotopic to that coming from the second term, hence the obstruction vanishes. 
\end{itemize}

\subsection{Summary: Map from \texorpdfstring{$\homKS$ to $\bigChains$}{homKS to bigChains}}
Let us summarize the findings of the previous sections.
\begin{thm}\label{thm:homKStobigChains}
The maps constructed above assemble to a map of colored operads
\[
 \homKS \to \bigChains.
\] 
In fact, this map is a quasi-isomorphism.
\end{thm}
\begin{proof}
The only statement not proven in the previous subsections is the fact that the maps are quasi-isomorphisms.
This fact will actually play no role for the present paper, but it is nice to know. For the part $\KS\to C(\cFM_2)$ of the above map this was shown by Kontsevich and Soibelman \cite{KS1,KS2}. Next one checks that the cohomology of $\hBr_\infty$ is the same as that of $\Br$, i.~e., $\Ger$. An explicit isomorphism $H(\Br) \to H(\hBr_\infty)$ is obtained by mapping a cocycle $\Gamma\in \Br$ to the cocycle $\Gamma\circ (f,\dots , f)\in \hBr_\infty$ where $f$ stands for the generator of $\hBr_\infty$ with one input and output. It corresponds to the unit in $B(\Br)$. On the other hand the homology of $D_K$ is the same as that of $\Br$. An explicit isomorphism $H(\Br) \to H(D_K)$ is given by sending the cocycle $\Gamma\in \Br$ to the cocycle $\Gamma\circ (pt,\dots , pt)\in C(D_K)$ where $\mathit{pt}$ stands for the fundamental chain of $D_K(1)$, which is a point. But the map constructed above is compatible with the left $\Br$-module structure and $f$ is mapped to $\mathit{pt}$, hence the induced map of operadic bimodules $H(\hBr_\infty)\to H(D_K)$ is an isomorphism. A similar argument shows that $H(\hKS_{1,\infty})\to H(D_S)$ is an isomorphism.
\end{proof}

\makeatletter{}\section{Maps between the operads, and the proofs of Theorems \ref{thm:brinfty} and \ref{thm:ksinfty}}
\label{sec:opmaps}
Our goal in this section is to prove Theorems \ref{thm:brinfty} and \ref{thm:ksinfty} from the introduction. To achieve this, we need to construct a representation of the big colored operad $\homKS$
on the colored vector space
\[
\bigV :=  \Dpoly \oplus C_\bullet \oplus \Tpoly \oplus \Omega_\bullet 
\]
that reduces on cohomology to the standard representation.
We already saw in Section \ref{sec:grops} that on the colored vector space $\bigV$ there is a natural action of the colored operad 
\[
 \bigGra =  \bpm \Br & \SGra & \Gra \\ \KS_1 & \SGra_1 & \Gra_1 \epm
\]
from section \ref{sec:bigGra}.
Concretely, 
\begin{enumerate}
 \item the action of $\KS$ on $\Dpoly \oplus C_\bullet$ is the standard $\KS$-algebra structure (see examples \ref{ex:Bractiononalg}, \ref{ex:KSactonDpolyCbullet}).
 \item $\bpm \Gra & \Gra_1\epm$ acts on $\Tpoly \oplus \Omega_\bullet$ by examples \ref{ex:GraactonTpoly}, \ref{ex:Gra1actonOmega}.
 \item The $\SGra$-action was described in example \ref{ex:SgraactTpolyDpoly}, the $\SGra_1$-action in example \ref{ex:Sgra1actCOmega}.
\end{enumerate}

Hence it will be sufficient to construct a map of colored operads $\homKS \to \bigGra$, that behaves well on cohomology. 
In Section \ref{sec:ksproof} we constructed a quasi-isomorphism of colored operads
\[
 \homKS \to \bigChains.
\]
Hence it suffices to construct a map 
\[
 \bigChains \to \bigGra.
\]
This section is dedicated to describing that map.
We will split the construction into the different color components of the operads as follows:
\begin{enumerate}
 \item The map of colored operads $\KS\to \KS$ is the identity. Here the first $\KS$ is to be understood as $\KS\subset\bigChains$ and the second $\KS$ as $\KS\subset \bigGra$.
 \item The colored operad map $C(\cFM_2) \to \bpm \Gra & \Gra_1 \epm$ is described in sections \ref{sec:FMtoGra} and \ref{sec:FM1toGra1}.
 \item The operadic bimodule map $C(D_K) \to \SGra$ and the moperadic bimodule map $C(D_S)\to \SGra_1$ are described in sections \ref{sec:DKtoSGra} and \ref{sec:CDStohGra1}. 
\end{enumerate}

\begin{rem}
The map $\bigChains \to \bigGra$ is simple to define using Feynman rules. The verification that it is a map of colored operads can be done combinatorially and by using Stokes' Theorem.
However, checking the signs and prefactors is a very tedious job. For example, both in the works of M. Kontsevich \cite{K1} and of B. Shoikhet \cite{shoikhet} the signs were not displayed explicitly. For Kontsevich's morphism the signs have been verified by hand in \cite{AMM} in a long calculation. 

To circumvent sign calculations, we will proceed as follows:
\begin{itemize}
\item First we define a map of two or three colored operads of Swiss Cheese or Extended Swiss Cheese type. This will involve only very simple sign verifications.
\item Then we extract a map of 2- or 4-colored operads by functoriality of the constructions of examples \ref{ex:SCtoPTmod} and \ref{ex:ESCtoPT1mod}, by Appendix \ref{sec:bimodgeneral} and \ref{sec:bimodgeneralex}, and by operadic twisting. Some things will need to be verified, but the sign calculations are ``hidden''.
\end{itemize} 
\end{rem}

\subsection{\texorpdfstring{ $C(\FM_2) \to \Gra$}{C(FM2) to Gra}}
\label{sec:FMtoGra}
The map $C(\FM_2)\to \Gra$ (or rather, a more complicated one) has been described by M. Kontsevich \cite{K2}. To describe it, it is convenient to introduce the predual\footnote{Here the involved vector spaces are finite dimensional, so the predual is isomorphic to the dual. But later we will encounter similar cases where this is not true.} $\pdu\Gra$ of $\Gra$. It is a cooperad, with the space of $n$-ary cooperations $\pdu\Gra(n)$ having a basis labelled by graphs with $n$ numbered vertices (see the definition of $\Gra$). 
The spaces $\pdu\Gra(n)$ are furthermore free graded commutative algebras and the cooperadic cocompositions respect this structure. The product is given by gluing two graphs together at the external vertices (up to a prefactor). The generators are graphs with a single edge.
Let $\Omega(\FM_2)$ be the collection of PA forms on the operad $\FM_2$ (see \cite{HLTV})
We define a map of collections of dg commutative algebras
\begin{gather*}
\omega: \pdu\Gra \to \Omega(\FM_2) \\
\Gamma \mapsto \omega_\Gamma 
\end{gather*}
which is compatible with the (co)operadic compositions in the sense that the diagrams 
  \begin{equation}\label{equ:opcompat1}
  \begin{tikzcd}
   \Gra(N) \ar{r}{\omega} \ar{dd} & \Omega(\FM_2(N)) \ar{d} \\
   & \Omega(\FM_2(N-k+1)\times \FM_n(k)) \\
   \Gra(N-k+1)\otimes \Gra(k) \ar{r}{\omega\otimes \omega} & \Omega(\FM_2(N-k+1))\otimes \Omega(\FM_2(k)) \ar{u}
  \end{tikzcd}
 \end{equation}
commute, where the left-hand vertical arrow is any cooperadic cocomposition in $\pdu\Gra$ and the right-hand vertical arrows are the pullback along the corresponding operadic composition on $\FM_2$, and the product of PA forms.

To define the map $\omega$, it is sufficient to define its value on commutative algebra generators, i.~e., graphs with one edge. The graph with one edge between vertices $i$ and $j$ is mapped to the PA form $d\arg(z_i-z_j)/2\pi$. Since this form is closed, the map $\omega$ respects the differentials.
It is not hard to see that the map respects the cooperad structure as well (in the sense that the above diagram commutes), since it is sufficient to check the statement on generators.
For graphs $\Gamma$ with more than one edge, $\omega_\Gamma$ is given by a product of 1-forms of the form $d\arg(z_i-z_j)/2\pi$, one for each edge, up to prefactor.
 Now turn to the desired map $C(\FM_2)\to \Gra$. It is given as the following composition, involving the adjoint $\omega^*$ of $\omega$:
\begin{equation}\label{equ:CFM2toGramapdef}
C(\FM_2)\to (\Omega(\FM_2))^* \stackrel{\omega^*}{\to} (\pdu\Gra)^* = \Gra
\end{equation}
Concretely, this map acts on a chain $c\in C(\FM_2(n))$ as follows:
\[
 c \mapsto \sum_\Gamma \Gamma \int_c \omega_\Gamma\, .
\]
Here the sum is over the elements of a basis of $\pdu\Gra(n)$, and we (abusively) use the same symbol $\Gamma$ to denote both the element of the basis of $\pdu\Gra$ (as in $\omega_\Gamma$) and the corresponding element $\Gamma\in \Gra(n)$ of the dual basis. Note that both spaces have a canonical pair of dual basis (up to signs) labelled by (automorphism classes of) graphs.
Note furthermore that commutativity of the diagrams \eqref{equ:opcompat1} implies that the composition \eqref{equ:CFM2toGramapdef} is a map of operads.

\begin{figure}
 \centering
\makeatletter{}\[
\begin{tikzpicture}[scale=1,
vert/.style={draw,outer sep=0,inner sep=0,minimum size=5,shape=circle,fill},
helper/.style={outer sep=0,inner sep=0,minimum size=5,shape=coordinate},
default_edge/.style={draw},
ext/.style={draw,outer sep=0,inner sep=2,minimum size=5,shape=circle},
every loop/.style={}]

\node (v0) at (5,6.5) [ext] {2};
\node (v1) at (5,7.5) [ext] {1};
\node (v3) at (5,8.4) [helper] {};

\draw[default_edge] (v0) to (v1);
\draw[default_edge] (v1) to (v3);
\end{tikzpicture}
\]
 
\caption{\label{fig:gerhomotopy} Gerstenhaber's homotopy. The differential of this $\Br$-tree is the commutator of the cup product. Hence the cup product (on cohomology) is commutative.}
\end{figure}
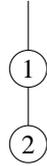

\begin{ex}
\label{ex:BronTpoly}
 The above construction yields an interesting $\Br_\infty$ structure on $\Tpoly$. Let us consider some examples. It can be checked that the tree depicted in Figure \ref{fig:gerhomotopy} acts as $1/2$ the Gerstenhaber bracket. It follows that the Lie bracket (given by the sum of trees in Figure \ref{fig:PTmc}) is the Gerstenhaber bracket. Let us consider the suboperad $A_\infty\subset \Br$ (see also example \ref{ex:ksmapainfty}). The generator $\mu_n$ of $A_\infty$ is mapped to the fundamental chains of $\FM_1(n)\subset \FM_2(n)$. Upon integration, we obtain $0$ for $n\geq 3$ and the usual wedge product for $n=2$. Hence the induced $A_\infty$ structure is the usual one. This statement will however change when we twist the maps by a Poisson structure in section \ref{sec:twistedops} below.
\end{ex}

\subsection{\texorpdfstring{ $C(\FM_{2,1}) \to \Gra_1$}{C(FM2,1) to Gra1}}
\label{sec:FM1toGra1}
\begin{figure}
 \centering
\makeatletter{}\usetikzlibrary{through}
\usetikzlibrary{arrows}
\begin{tikzpicture}[scale=.7,int/.style={circle, draw, fill, minimum size=5pt, inner sep=0, outer sep=0}]
\node[int, label=180:{$out$}] (c) at (0,0) {};

\node[int, label=180:{$\psi$}] (v1) at (-40:2) {};
\node[int] (v2) at (-90:3) {};

\draw[dashed] (c) edge (v1);
\draw (c) to (.5,.5);
\draw (v1) +(130:.7) arc (130:240:.7);
\draw[-triangle 60] (v1)--(v2);

\begin{scope}[xshift=6cm, shift={(-1,0)}]
\node[int, label=90:{$\phi$}, label=-90:{$out$}] (c2) at (0,0) {};

\node[int] (vv1) at (130:2) {};

\draw[-triangle 60] (c2) edge (vv1);
\draw (c2) to (60:1);
\draw (50:0.8) arc (50:140:0.8);

\end{scope}
\draw[dashed] (2.5,3)--+(0,-8);
\draw[dashed] (7,3)--+(0,-8);
\begin{scope}[xshift=11cm, shift={(-2,0)}]
\node[int, label=180:{$out$}] (c) at (0,0) {};

\node[int, outer sep=-2.5pt,label=-130:{$\psi$}] (v1) at (-40:2) {};
\draw (v1)+(5,0) node [int, label=90:{$in$ at $+\infty$}] (v2) {};

\draw[dashed] (c) edge (v1);
\draw (c) to (0.5,0.5);
\draw (v1) +(130:0.8) arc (-230:10:0.8);
\draw[-triangle 60] (v1)--(v2);

\end{scope}
\end{tikzpicture} 
\caption{\label{fig:euclphipsi} The definition of $\psi$ (left) and $\phi$ (middle) for the map $C(\FM_{2,1}) \to \Gra_1$. The input vertex $in$ is thought of as $+\infty$ (right). }
\end{figure}
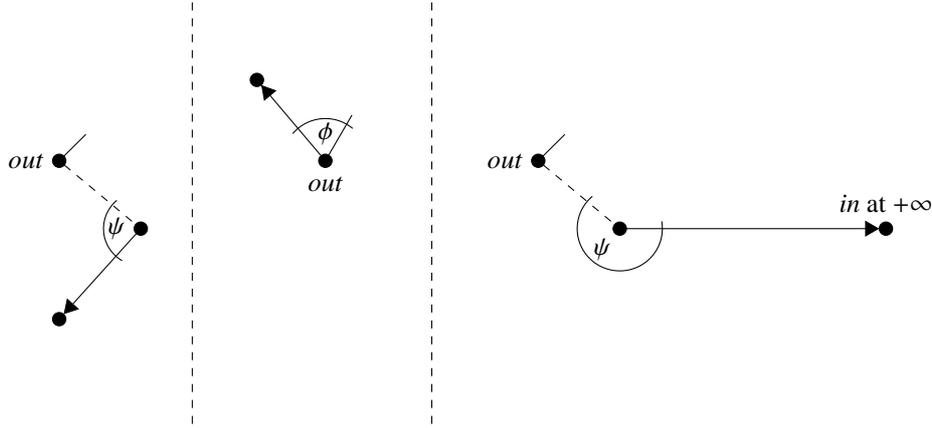

Let us define a map of moperads $C(\FM_{2,1}) \to \Gra_1$.
Again it is easier to first describe the (pre-dual) map of collections of dg commutative algebras
\[
\omega: \pdu\Gra_1 \to \Omega(\FM_{2,1})
\]
where $\pdu\Gra_1$ is the predual of $\Gra_1$ (it is canonically isomorphic to the dual) and $\Omega(\FM_{2,1})$ is the collection of PA forms on $\FM_{2,1}$. Again $\pdu\Gra_1$ is a co-moperad of free graded commutative algebras and the generators are the graphs with exactly one edge. In particular, the algebra structures on the spaces of cooperations are compatible with those on $\pdu\Gra$ and the right $\pdu\Gra$-coaction.
The map $\omega$ will be compatible with the dg commutative algebra structures, so it is sufficient to define it on generators (graphs with a single edge). Here one has to distinguish several kinds of edges.
\begin{enumerate}
 \item For an edge between vertices $i$ and $j$, none of which is the central vertex $out$, one associates the form $\frac{d\psi}{2\pi}$, where $\psi$ is the angle between the lines from $out$ to $z_i$, and from $z_i$ to $z_j$.
\item For an edge between the central vertex $out$ and some other vertex $j\neq in$, one associates the form $\frac{d\phi}{2\pi}$, where $\phi$ is the angle between the framing at $out$ and the line between $out$ and $z_j$.
\item If $j$ above is the vertex $in$, one understands $z_j$ as $+\infty$.
\end{enumerate}
These definitions of $\phi$ and $\psi$ are shown pictorially in Figure \ref{fig:euclphipsi}.
For a graph $\Gamma$ with more than one edge the differential form $\omega_\Gamma$ is a product of 1-forms, one for each edge, up to a (conventional) prefactor. The forms involved are closed, so the map  automatically respects the differentials. We claim that the map is also compatible with the co-moperadic (co)compositions in the sense that the natural diagrams 
\[
 \begin{tikzcd}
  \pdu\Gra_1(N) \ar{r}{\omega}\ar{dd} & \Omega(\FM_{2,1}(N)) \ar{d}\\
& \Omega\left(\FM_{2,1}(N-k+1)\times \FM_2(k)\right) \\
 \pdu\Gra_1(N-k+1) \otimes \pdu\Gra(k) \ar{r}{\omega\otimes \omega} &  \Omega(\FM_{2,1}(N-k+1))\otimes\Omega(\FM_2(k)) \ar{u}
 \end{tikzcd}
\]
\[
 \begin{tikzcd}
  \pdu\Gra_1(N) \ar{r}{\omega}\ar{dd} & \Omega(\FM_{2,1}(N)) \ar{d}\\
& \Omega\left(\FM_{2,1}(N-k)\times \FM_{2,1}(k)\right) \\
 \pdu\Gra_1(N-k) \otimes \pdu\Gra_1(k) \ar{r}{\omega\otimes \omega} &  \Omega(\FM_{2,1}(N-k))\otimes\Omega(\FM_{2,1}(k)) \ar{u}
 \end{tikzcd}
\]
built out of the (co)operadic (co)composition, the $\pdu\Gra$-coaction on $\pdu\Gra_1$ and the $\FM_2$ action on $\FM_{2,1}$ commute. In fact, since it suffices to check the commutativity of those diagrams on algebra generators, our claim is easily verified.

The map $C(\FM_{2,1}) \to \Gra_1$ is then the following composition, involving the adjoint of $\omega$:
\[
C(\FM_{2,1})\to (\Omega(\FM_{2,1}))^* \stackrel{\omega^*}{\to} (\pdu\Gra_1)^* = \Gra_1.
\]
Concretely, for a chain $c\in C(\FM_{2,1}(n))$ we have the ``Feynman rules'' formula 
\[
 c \mapsto \sum_\Gamma \Gamma \int_c \omega_\Gamma\, .
\]
Here the sum is over elements of a basis of $\pdu\Gra_{2,1}(n)$ and we again denote by $\Gamma$ both the basis vector in $\pdu\Gra_{2,1}(n)$ and the dual basis vector in $\Gra_{2,1}(n)$, abusing notation.
It is convenient to take as basis the canonical (up to signs) basis whose elements are labelled by graphs, but of course it does not matter which basis we pick.

\begin{figure}
 \centering
\makeatletter{} \usetikzlibrary{matrix}
\usetikzlibrary{arrows}
\usetikzlibrary{shapes}
\usetikzlibrary{through}
\usetikzlibrary{calc,3d}
\usetikzlibrary{decorations,decorations.pathmorphing}
\[
\begin{tikzpicture}[
scale=.6,
int/.style={circle, draw, fill, minimum size=5pt, inner sep=0},
ext/.style={circle, draw, fill=white, minimum size=5pt, inner sep=1pt},
helper/.style={coordinate},point/.style={circle, draw, fill, inner sep =1pt},
de/.style={-triangle 60},
point/.style={circle, draw, fill, minimum size=3pt, inner sep=0pt},
xst/.style={cross out, draw, minimum size=5pt}
]

\begin{scope}[xshift=6cm]
\draw (0,0) ellipse (2 and 1);
\draw (0,3) ellipse (2 and 1);
\draw (-2,0)--(-2,3) (2,0)--(2,3);
\node [xst, label=-90:{$out$}] (out) at ($(0,0)+(-130:2 and 1)$) {};
\node [xst, label=90:{$in$}] (in) at ($(0,3)+(-80:2 and 1)$) {};
\draw (out)+(0,1) node[ext] (e3) {$\mathbb{1}$};
\draw (out.base) edge (e3);
\draw (in.base) -- +(0,-3);
\end{scope}

\begin{scope}[xshift=6cm, shift={(5,0)}]
\draw (0,0) ellipse (2 and 1);
\draw (0,3) ellipse (2 and 1);
\draw (-2,0)--(-2,3) (2,0)--(2,3);
\node [xst, label=-90:{$out$}] (out) at ($(0,0)+(-90:2 and 1)$) {};
\node [xst, label=90:{$in$}] (in) at ($(0,3)+(-90:2 and 1)$) {};
\draw (out)+(0,1.2) node [int] (e3) {} edge (in.base);
\draw (out.base) edge (e3);
\draw (e3) +(.7,.7) node[ext] (e1) {1};
\draw (e3) edge (e1);
\end{scope}

\begin{scope}[xshift=6cm, shift={(9,0.5)}, scale=0.5]
\draw (0,0) ellipse (2 and 1);
\draw (0,3) ellipse (2 and 1);
\draw (-2,0)--(-2,3) (2,0)--(2,3);
\node [xst, label=-90:{$out$}] (out) at ($(0,0)+(-90:2 and 1)$) {};
\node [xst, label=90:{$in$}] (in) at ($(0,3)+(-90:2 and 1)$) {};
\node [ext] (e1) at ($(0,1.2)+(-60:2 and 1)$) {1};
\draw (out.base) edge (in.base);
\draw (e1)-- +(0,-1.2);
\end{scope}

\begin{scope}[xshift=6cm, shift={(12,0.5)}, scale=0.5]
\draw (0,0) ellipse (2 and 1);
\draw (0,3) ellipse (2 and 1);
\draw (-2,0)--(-2,3) (2,0)--(2,3);
\node [xst, label=-90:{$out$}] (out) at ($(0,0)+(-90:2 and 1)$) {};
\node [xst, label=90:{$in$}] (in) at ($(0,3)+(-90:2 and 1)$) {};
\draw (out)+(0,1.2) node [ext] (e3) {1} edge (in.base) edge (out.base);
\end{scope}

\begin{scope}[xshift=6cm, shift={(16,0)}]
\draw (0,0) ellipse (2 and 1);
\draw (0,3) ellipse (2 and 1);
\draw (-2,0)--(-2,3) (2,0)--(2,3);
\node [xst, label=-90:{$out$}] (out) at ($(0,0)+(-130:2 and 1)$) {};
\node [xst, label=90:{$in$}] (in) at ($(0,3)+(-90:2 and 1)$) {};
\draw (out)+(0,1) node [ext] (e3) {$\mathbb{1}$};
\node [ext] (e1) at ($(0,1.2)+(-60:2 and 1)$) {1};
\draw (out.base) edge (e3);
\draw (in.base) -- +(0,-3);
\draw (e1)-- +(0,-1.2);
\end{scope}

\node at (6,4.5) {$B$};
\node at (22,4.5) {$H$};
\node at (16.5,4.5) {$L$};
\node at (11,4.5) {$I$};
\node at (16.5037,1.2069) {+};
\end{tikzpicture}
\] 
\caption{\label{fig:BLI} Some important operations in the moperad $\KS_1$. On Hochschild chains they act as: (left) the Connes-Rinehart differential $B$, (semi-left) the cap product $I$, (semi-right) the Lie algebra action $L$. The rightmost operation, say $H$, acts as homotopy in the Cartan formula, i.e., $B\circ I + I\circ B -L = \delta H$. }
\end{figure}
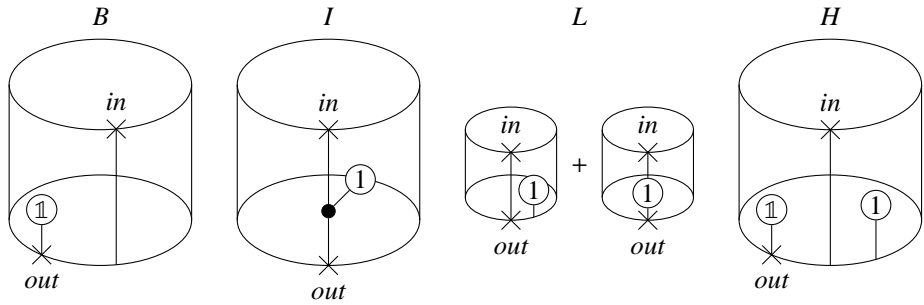

\begin{ex}
 One obtains an interesting ($\KS$-)action of $\Tpoly$ on $\Omega_\bullet$. Let us work it out for the leading order terms, i.~e., for the elements $B,L,I$ depicted in Figure \ref{fig:BLI}.
The element $B$ is mapped to the fundamental chain of $\FM_{2,1}(0)$. There is only one graph giving nonzero value on this chain, namely the one in Figure \ref{fig:dliota} (left). The corresponding operation is the de Rham differential. Next consider $I$. It is mapped to a degree zero chain (point) in $C(\FM_{2,1}(1))$. The only graph that can attain a nonzero value is hence one without any edges, i.~e., the middle graph depicted in Figure \ref{fig:dliota}. This is the standard contraction. The operation $L$ is slightly more complicated. It is mapped to the fundamental chain of $\FM_2(2)$ (i.~e., a circle) inside $\FM_{2,1}$. It is embedded such that the framing is fixed pointing to $+\infty$. There are two graphs that can attain nonzero coefficents, namely the ones depicted in Figure \ref{fig:dliota} (right). The associated operation is the usual Lie derivative. So the lowest degree operations are the standard ones. However, there are other nonzero operations. For example consider the operation $H$ depicted in Figure \ref{fig:BLI} (right). One can check that is is mapped to $\frac{1}{2} d\circ L$. This operation was considered already in \cite{CFW}.
\end{ex}

\subsection{\texorpdfstring{$C(D_{K}) \to \SGra$}{C(DK) to SGra} }
\label{sec:DKtoSGra}
Let us turn to the map of operadic $\Br$-$C(\FM_2)$-bimodules
\[
 \Phi\colon C(D_{K}) \to \SGra\, .
\]

To define it, we will consider first the (version of the) two-colored Swiss cheese operad $\SC$ from section \ref{sec:DKdef}. Its components are the spaces $\FM_2(n)$ in one color and the Kontsevich configuration spaces $D_{Ke}(m,n)$ in mixed colors. 

There is a similar two-colored operad $\SG$, whose components are given by $\Gra(n)$ in one color, and $\fSGra(m,n)$ in mixed colors, as noted in section \ref{sec:SGra}. First we will construct a map of two-colored operads
\[
C(\SC) \to \op SG.
\]
The construction is more or less a copy of the one in the previous two subsections. We construct a map of colored collections of dg commutative algebras
\[
\omega: \pdu \SG \to \Omega(\SC)
\]
where $\pdu \SG$ is the predual and $\Omega(\SC)$ are PA forms.
Furthermore, $\pdu \SG$ is naturally a colored cooperad of free graded commutative algebras, the algebra generators given by graphs with a single edge. The map $\pdu \SG\to \Omega(\SC)$ will respect the dg commutative algebra structures and hence it is sufficient to define it on generators. For the part $\pdu \Gra\subset \pdu \SG$ the map was defined in section \ref{sec:FMtoGra} above. For the part $\pdu \fSGra(m,n) \to \Omega(D_{Ke}(m,n))$, consider a graph with a single edge from vertex $i$ to $j$, where $j$ can be either type I or type II. The 1-form associated to that graph in either case is the differential $d\alpha/2\pi$ of the hyperbolic angle $\alpha$ between the hyperbolic geodesic lines $(z_i,+i\infty)$ and $(z_i,z_j)$, see Figure \ref{fig:kontsangle}. This form was introduced by M. Kontsevich. Again it is easy to check on generators that the map $\omega: \pdu \SG\to \Omega(\SC)$ thus defined is compatible with the colored (co)operadic (co)compositions. The desired map of colored operads $C(\SC) \to \SG$ is then the composition
\[
C(\SC)\to (\Omega(\SC))^* \stackrel{\omega^*}{\to} (\pdu\SG)^* = \SG\, .
\]

Now, given a Swiss Cheese type operad (like $\SG$), one can construct an operadic $\PT$-$\SG^1$-bimodule structure on the total spaces
\[
\prod_n \SG(\cdot,n)[-n]
\]
as in example \ref{ex:SCtoPTmod} in section \ref{sec:PT}. This construction is functorial and from our map $C(\SC)\to \SG$ we hence obtain a map of operadic bimodules
\[
\prod_n C(D_{Ke}(\cdot,n))[-n] \to \prod_n \fSGra(\cdot, n)[-n] = \fSGra(\cdot).
\]
Next we want to twist the right $\PT$ actions to $\Tw\PT$- and hence $\Br$-actions. For this we need to identify a Maurer-Cartan element 
\[
m\in \prod_n C(D_{Ke}(0,n))[-n].
\]
It is given by the sum of the fundamental chains of $D_{Ke}(0,n)$,
\[
m = \sum_{n\geq 2} \mathit{Fund}(D_{Ke}(0,n)).
\]
Note that in our conventions the Maurer-Cartan element has degree 2.
The image of the Maurer-Cartan element in $\fSGra(0)$ is easily checked to be the graph depicted in Figure \ref{fig:sgramc}. I.~e., all fundamental chains are sent to zero, except for $\mathit{Fund}(D_{Ke}(0,2))$, which is a point.
By operadic twisting, we obtain (i) an operadic $\Tw\PT$-$C(\FM_2)$-bimodule structure on $\prod_n C(D_{Ke}(\cdot,n))[-n]$ (with changed differential), (ii) an operadic $\Tw\PT$-$\Gra$-bimodule structure on $\fSGra$ (with changed differential), and (iii) a map between the bimodules. 

Now let us finally construct the map of operadic bimodules $C(D_{K}) \to \SGra$ as planned. First we are (of course) free to restrict the left $\Tw\PT$ actions on the above bimodules to $\Br\subset \Tw\PT$ actions. Secondly, by the very construction of the $\Br$-$C(\FM_2)$ bimodule structure on $C(D_{K})$ we have an embedding of bimodules
\[
\Phi: C(D_{K}) \to \prod_n C(D_{Ke}(0,n))[-n].
\]
By composing this with the map (of bimodules) to $\fSGra$ we obtain the desired map $C(D_{K}) \to \SGra$.

Let us unravel this definition into a concrete formula. One has fibrations $\pi_{m,n}: D_{Ke}(m,n) \to D_{K}(m)=D_{Ke}(m,0)$ by forgetting the positions of the type II vertices. Taking the fibers over chains, one obtains a map $\pi_{m,n}^{-1}:C(D_{K}(m)) \to C(D_{Ke}(m,n))$. For a chain $c\in C(D_{K}(m)$ we then have the ``Feynman rule'' formula:
\[
 c\mapsto \Phi(c) = \sum_\Gamma \Gamma \, \int_{\pi_{m,n}^{-1} (c)} \omega_\Gamma\, .
\]
Here the sum runs over elements of a basis of $\pdu \fSGra(m)$ and $\omega_\Gamma\in \Omega(D_{Ke})$ is the differential form associated to the graph $\Gamma$. (It has first been described in \cite{K1}.)

\begin{figure}
 \centering
\makeatletter{}\usetikzlibrary{through}
\usetikzlibrary{arrows}
\begin{tikzpicture}[int/.style={circle, draw, fill, minimum size=5pt, inner sep=0}]

\draw (-3,0)--(3,0);
\node[int, label=130:{$\alpha$}] (v1) at (1,3) {};
\node[int] (v2) at (-2,1) {};
\draw[dashed] (v1)+(0,-.5) -- +(0,3);
\draw (v1)+(80:.8) arc (80:195:.8);
\clip (-5,-.5) rectangle (2.5,4);
\node [draw, dashed] at (.84,0) [circle through={(v1)}] {};
\draw[-triangle 60] (v1)--(v2);
\end{tikzpicture} 
\caption{\label{fig:kontsangle} The definition of the hyperbolic angle, used by M. Kontsevich.}
\end{figure}
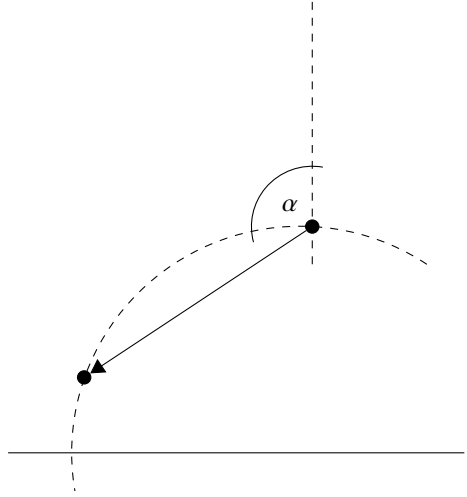

\subsection{\texorpdfstring{$C(D_S) \to \SGra_1$}{C(DS) to SGra1}}
\label{sec:CDStohGra1}
The goal of this subsection is to construct a map of moperadic bimodules $C(D_S) \to \SGra_1$, hence completing the program outlined at the beginning of this section. 
We will arrive at a formula similar to that of B. Shoikhet \cite{shoikhet}.
 The construction of the map will be more a less a copy of the construction of the previous subsection, only the notation is more cumbersome. Let us go through it.

First consider the three-colored version of the Swiss cheese operad $\ESC$ as defined in section \ref{sec:DS}. Its components are $\FM_2$, $\FM_{2,1}$, $D_{Ke}(\cdot, \cdot)$ and the Shoikhet configuration spaces $D_{Se}(\cdot, \cdot)$.
There is a similar Extended Swiss Cheese type operad of graphs, whose components are $\Gra$, $\Gra_1$, $\fSGra(\cdot, \cdot)$ and $\fSGra_1(\cdot, \cdot)$, see section \ref{sec:sgra1}. Let us temporarily call this graphs operad $\op Q$. The first step is to construct a map of three-colored operads
\[
C(\ESC) \to \op Q.
\]
To do this, we construct a map of three colored collections of dg commutative algebras
\begin{gather*}
\omega: \pdu \op Q \to \Omega(\ESC) \\
\Gamma \mapsto \omega_\Gamma
\end{gather*}
compatible with the colored (co)operadic (co)compositions.
Here $\pdu \op Q$ is the predual of $\op Q$ (which is canonically isomorphic to the dual by finite dimensionality) and $\Omega(\cdot)$ again denotes the PA forms. The components of $\pdu \op Q$ are free graded commutative algebras, the generators are graphs with only a single edge. It is hence sufficient to construct the map for such graphs. On the components $\Gra$, $\Gra_1$ and $\fSGra(\cdot, \cdot)$ of $\op Q$ we did that in the previous subsections. Consider the component $\fSGra_1(m, n)$. One has to distinguish two cases, depending on what vertices the edge connects:
\begin{enumerate}
 \item An edge between vertices $i$ and $j$, both not equal to the central vertex $out$, contributes a form $\frac{d\psi}{2\pi}$, where $\psi$ is the angle between the hyperbolic geodesics through $out$ and $z_i$, and through $z_i$ and $z_j$.
 \item An edge from $out$ to a vertex $j$ contributes the 1-form $\frac{d\phi}{2\pi}$, where $\phi$ is the angle between the framing at $out$ and the line from $out$ to $z_i$.
\end{enumerate}
For graphs $\Gamma$ with more than one edge the form $\omega_\Gamma$ is a product of 1-forms, one for each edge. To check that the map $\omega$ respects the moperadic co-bimodule structure it is sufficient to check the statement on generators, which is an easy exercise.
Using the adjoint map $\omega^*$ one constructs a map of 3 colored operads
\[
C(\ESC) \to \op Q
\]
as the composition
\[
C(\ESC)\to (\Omega(\ESC))^* \stackrel{\omega^*}{\to} (\pdu\op Q)^* = \op Q\, .
\]
From a three colored operad like $C(\ESC)$ or $\op Q$ above one can extract by a natural construction (i) an operadic bimodule (with left $\PT$ action) as before and (ii) a moperadic bimodule (with right $\PT_1$ action), see example \ref{ex:ESCtoPT1mod}. Together with the operad pieces one obtains two four colored operads from $C(\ESC)$ and $\op Q$ and by functoriality a map between them:
\[
\bpm
\PT   & \prod_n C(D_{Ke}(\cdot, n))[-n] & C(\FM_2) \\
\PT_1 & \prod_n C(D_{Se}(\cdot, n+1))[-n] & C(\FM_{2,1})
\epm 
\to 
\bpm
\PT   & \fSGra & \Gra \\
\PT_1 & \fSGra_1 & \Gra_1
\epm\, .
\]
Next we we want to extend the actions of $\PT_1$ to actions of $\PT_1^{\mathbb{1}}$ as done in section \ref{sec:sgra1} and Appendix \ref{app:KSactproof}. 
Recall from there that on the left-hand side, the action of the additional element $\mathbb{1}$ is by forgetting the location of points in a configuration in $D_{Se}$. On the right-hand side, the action $\fSGra_1$ is by removing type II vertices from the graph, mapping the graph to zero if the removed vertex has valence $\geq 1$, see section \ref{sec:sgra1}. We have to check by hand that, after extending the operads to include $\PT_1^{\mathbb{1}}$, the above map is still a map of 4-colored operads. Of course, only the part mapping the moperadic bimodules has to be considered. Call this part $F$. Is is sufficient to check that for all chains $c\in C(D_{Se}(m, n))$, and all $m,n$ and $0\leq k \leq m$ we have $\iota_k F(c) = F((\pi_k)_* c)$.
Here $\iota_k$ is the operation of deleting vertex $\bar k$ from the graph, mapping the graph to zero if it does not have valence $0$, and $\pi_k$ is the forgetful map forgetting the location of point $\bar k$ in configurations. Compute:
\[
F((\pi_k)_* c)
=\sum_\Gamma \Gamma \int_{(\pi_k)_*c} \omega_\Gamma
= \sum_\Gamma \Gamma \int_{c} \pi_k^* \omega_\Gamma
= \sum_\Gamma \Gamma \int_{c} \omega_{\Gamma \cup \bar k}
= \sum_\Gamma \iota_k \Gamma \int_{c} \omega_{\Gamma}
\]
Here $\Gamma \cup \bar k$ is the graph obtained from $\Gamma$ by inserting one additional valence zero type II vertex at position $k$.
Hence we have a map of 4 colored operads:
\[
\bpm
\PT   & \prod_n C(D_{Ke}(\cdot, n))[-n] & C(\FM_2) \\
\PT_1^{\mathbb{1}} & \prod_n C(D_{Se}(\cdot, n+1))[-n] & C(\FM_{2,1})
\epm 
\to 
\bpm
\PT   & \fSGra & \Gra \\
\PT_1^{\mathbb{1}} & \fSGra_1 & \Gra_1
\epm\, .
\]
Next we twist the $\PT$ and $\PT_1^{\mathbb{1}}$ actions on the middle (bimodule) part by the Maurer-Cartan element $m$ as in the previous subsection. For this step no additional relations have to be checked, relative to the previous subsection. By naturality of twisting we obtain a map
\[
\bpm
\Tw\PT   & \prod_n C(D_{Ke}(\cdot, n))[-n] & C(\FM_2) \\
\Tw\PT_1^{\mathbb{1}} & \prod_n C(D_{Se}(\cdot, n+1))[-n] & C(\FM_{2,1})
\epm 
\to 
\bpm
\Tw\PT   & \fSGra & \Gra \\
\Tw\PT_1^{\mathbb{1}} & \fSGra_1 & \Gra_1
\epm\, .
\]
Here the (m)operadic bimodules have a modified differential, though we do not reflect this in the notation. We are free to restrict $\Tw\PT$ to the suboperad $\Br$ and $\Tw\PT_1^{\mathbb{1}}$ to the sub-moperad $\KS_1'$ from section \ref{sec:KS1}. We get a map 
\[
\bpm
\Br   & \prod_n C(D_{Ke}(\cdot, n))[-n] & C(\FM_2) \\
\KS_1' & \prod_n C(D_{Se}(\cdot, n+1))[-n] & C(\FM_{2,1})
\epm 
\to 
\bpm
\Br   & \fSGra & \Gra \\
\KS_1' & \fSGra_1 & \Gra_1
\epm\, .
\]
By construction of the actions on the (m)operadic bimodules $C(D_K)$ and $C(D_S)$, we have a map of colored operads
\[
\bpm
\Br   & C(D_K) & C(\FM_2) \\
\KS_1' & C(D_S) & C(\FM_{2,1})
\epm 
\to
\bpm
\Br   & \prod_n C(D_{Ke}(\cdot, n))[-n] & C(\FM_2) \\
\KS_1' & \prod_n C(D_{Se}(\cdot, n+1))[-n] & C(\FM_{2,1})
\epm 
\, .
\]
Composing with the above map, we get a map 
\[
\bpm
\Br   & C(D_K) & C(\FM_2) \\
\KS_1' & C(D_S) & C(\FM_{2,1})
\epm 
\to
\bpm
\Br   & \fSGra & \Gra \\
\KS_1' & \fSGra_1 & \Gra_1
\epm\, .
\]
Next we may pass from $\KS_1'$ to its quotient $\KS_1$, by imposing the relations of section \ref{sec:KS1}. For this step, nothing has to be checked and hence we obtain the desired map of 4 colored operads
\[
\bpm
\Br   & C(D_K) & C(\FM_2) \\
\KS_1 & C(D_S) & C(\FM_{2,1})
\epm 
\to
\bpm
\Br   & \fSGra & \Gra \\
\KS_1 & \fSGra_1 & \Gra_1
\epm\, .
\]

The final step is to check that on the right-hand side we may replace $\fSGra$ and $\fSGra_1$ by their sub-(m)operadic bimodules $\SGra$ and $\SGra_1$. Recall that $\SGra\subset \fSGra$ is spanned by graphs for which all type II vertices have valence $\geq 1$, and similarly $\SGra_1\subset \fSGra_1$ is spanned by graphs for which all type II vertices except possibly the vertex $\bar 0$ have valence $\geq 1$.

To this end, let us first consider the explicit form of the map of operadic bimodules
\[
 C(D_{K}) \to \fSGra \, .
\]
Let $c\in C(D_K(m))$ be a (semi algebraic) chain. We then have
\[
 c \mapsto  
 \sum_{n\geq 0} \sum_\Gamma \Gamma \, \int_{\pi_n^{-1} (c)} \omega_\Gamma.
\]
Here $\pi_n$ is the forgetful map forgetting the locations of all type II vertices.
The sum is over elements of a basis of $\pdu \fSGra(m,n)$. Again we use the same symbol $\Gamma$ for elements of that basis and corresponding elements in the dual basis of $\fSGra(m,n)$.
Note if a graph $\Gamma$ appearing in the sum above contains a type II vertex of valence 0, then the differential form $\omega_\Gamma$ in the integral will have no dependence on the position of the corresponding point in the configuration. Hence the integral is necessarily 0. In other words, all graphs occurring non-trivially in the sum are such that all their type II vertices have valence $\geq 1$.
In yet other words, the map $C(D_{K}) \to \fSGra$ factors through $\SGra$ as we wanted to show.

Similarly, let us study the explicit form of the map of moperadic bimodules
\[
 C(D_{S}) \to \fSGra_1 \, .
\]
Let again $c\in C(D_S(m))$ be a (semi algebraic) chain. We then have
\[
 c\mapsto
 \sum_{n\geq 0} \sum_\Gamma \Gamma \, \int_{\pi_n^{-1} (c)} \omega_\Gamma.
\]
Here $\pi_n$ is the forgetful map forgetting the locations of all but one type II vertices (namely, all but the vertex $\bar 0$). The sum is over elements of a basis of $\pdu \fSGra_1(m,n)$. Again we use the same symbol $\Gamma$ for elements of that basis and corresponding elements in the dual basis of $\fSGra_1(m,n)$.
Note that in case a graph $\Gamma$ is such that a type II vertex $\bar j$ (with $j>0$) has valence 0, then 
the integrand again has no dependence on the position of the corresponding point, and hence the integral vanishes. This then shows that the map $C(D_{S}) \to \fSGra_1$ factors through $\SGra_1$. Hence we finally arrive at the desired map of four colored operads
\begin{equation}\label{equ:bigChainstobigGra}
\bigChains
=
 \bpm
\Br   & C(D_K) & C(\FM_2) \\
\KS_1 & C(D_S) & C(\FM_{2,1})
\epm 
\to
\bpm
\Br   & \SGra & \Gra \\
\KS_1 & \SGra_1 & \Gra_1
\epm
=\bigGra \, .
\end{equation}

\subsection{The proofs of Theorems \ref{thm:brinfty} and \ref{thm:ksinfty}} \label{sec:proofsof4and5}
Composing the map of Theorem \ref{thm:homKStobigChains} with the map \eqref{equ:bigChainstobigGra} just constructed we obtain the map $\homKS\to \bigGra$ whose existence is asserted in Theorem \ref{thm:ksinfty}, and whose part in colors 1 and 3 yields the map whose existence is asserted in Theorem \ref{thm:brinfty}.
There are a few further assertions in these Theorems that we are going to check in this section, thus finishing the proof.

First we need to check that the Gerstenhaber structure on cohomology induced by the $\Br_\infty$ structure on $\Tpoly$ is the standard one.  Of course, it suffices to verify this on cocycles in $\Br_\infty$ generating $\Ger$. Denote the obvious such cocycles by $\co{}{}$ and $\wedge$. They are mapped to chains $c_{\co{}{}}, c_{\wedge}\in C(\FM_2)$. It is easy to see (e.~g. by Proposition \ref{prop:explonBr}) that  $c_{\co{}{}}$ is the fundamental chain while $c_{\wedge}$ is the chain of a point. These chains in turn are mapped to the graph with two vertices and one edge, which acts as the Schouten bracket, and to the graph with two vertices and no edge, which acts as the wedge product.

The final assertion of Theorem 1 is that the $\Br_\infty$ morphism restricts to M. Kontsevich's formality morphism on the $\Lie_\infty$ part. To begin with, note that there is an embedding $\Lie\{1\}\to \Br$. Hence, by functoriality of the bar-cobar construction we get a map of operads $\Omega(B(\Lie\{1\}))\to \Br_\infty$. Composing with the canonical map $\Lie^{(1)}_\infty = \Omega(\Lie\{1\}^{\vee}) \to \Omega(B(\Lie\{1\}))$  we obtain the embedding of $\Lie^{(1)}_\infty\to \Br_\infty$. We claim that under the composition $\Lie^{(1)}_\infty\to \Br_\infty\to C(\FM_2)$ the generator $\mu_n$ is mapped to the fundamental chain $f_n$ of $\FM_2(n)$. This is true for $n=2$. For higher $n$ one shows the statement by induction. Say $\mu_j$ is mapped to some $c_j$. Suppose $c_j=f_j$ for $j<n$. Then we know that $c_n$ and $f_n$ satisfy an equation of the form
\[
\p c_n = (\cdots) = \p f_n
\] 
where $(\cdots)$ are some terms in the $c_j=f_j$ for $j<n$. 
In particular $\p c_n = \p f_n$. Hence $c_n-f_n$ is closed, and thus also exact since $\FM_2$ does not have homology in the relevant degrees. But by dimensional reasons there is no non-zero exact semi algebraic chain in the relevant degree.\footnote{The chain $c_n$ has degree $3-n$, hence the tentative exact element must be the boundary of a chain of degree $2-2n$. But $\dim \FM_2(n) =n-3$ and hence that chain is zero, since there are no semi-algebraic degenerate chains.} Hence $c_n-f_n=0$. By a Kontsevich vanishing Lemma one hence sees that under the composition $\Lie^{(1)}_\infty\to C(\FM_2)\to \Gra$ all generators except $\mu_2$ are sent to zero. It follows that the $\Lie_\infty$ structure induced on $\Tpoly[1]$ is the standard Lie algebra structure.

A similar argument as above shows that the element of $\hBr_\infty$ that correspond to the $j$-th component of the $\Lie^{(1)}_\infty$ map is mapped to the fundamental chain $F_j$ of $D_K(n)$. Concretely, let us assume that that element is mapped to the chain $C_j$. We know that $C_1=F_1$ by construction (note that $D_K(1)=\{pt\}$). Assume inductively that $C_j=F_j$ for $j<n$. Then by definition of the $C_n$ it satisfies an equation
\[
\p C_n = (\cdots) = \p F_n
\]
where $(\cdots)$ are some terms depending on the $c_i=f_i$ and the $C_j=F_j$ for $j<n$.
Hence $C_n -F_n$ is a closed chain of degree $2-2n$. Hence it is exact since $D_K(n)$ does not have homology in this degree for $n>1$. Hence it is zero by dimensionality reasons.
 It follows that restricting the $\Br_\infty$-map $\Tpoly\to \Dpoly$ to its $\Lie^{(1)}_\infty$-part one recovers the Kontsevich formula.\footnote{Strictly speaking M. Kontsevich constructed a $\Lie_\infty$ map $\Tpoly[1]\to \Dpoly[1]$, but this is equivalent to giving a $\Lie^{(1)}_\infty$-map $\Tpoly\to \Dpoly$. }
Hence Theorem \ref{thm:brinfty} is shown.

Finally let us turn to Theorem \ref{thm:ksinfty}. Consider first the induced calculus structure on $\Tpoly$ and $\Omega_\bullet$. We already checked that the induced Gerstenhaber structure on $\Tpoly$ is the standard one. To check the $\calc_1$ structure on $\Omega_\bullet$, it suffices to compute the action of the generators $I$ and $B$. Tracing the construction of the maps (e.~g., using proposition \ref{prop:KSexpl}), we see that $B$ acts as the de Rham differential and $I$ as the contraction operator $\iota$, as was claimed in the Theorem.

To see the final statement of the Theorem (i.~e., the equality of the $\Lie^{(1)}_\infty$ part of the map of modules $C_\bullet\to \Omega_\bullet$ to B. Shoikhet's map), let us make a preliminary remark.
Note that $\KS_1$ contains a $\Br$ sub-moperad $\KS_1'\subset \KS_1$, spanned by graphs such that the subtree of $\vout$ at the marked edge contains $\vin$. For example, the left graph in the following picture is in $\KS_1'$, the right is not.

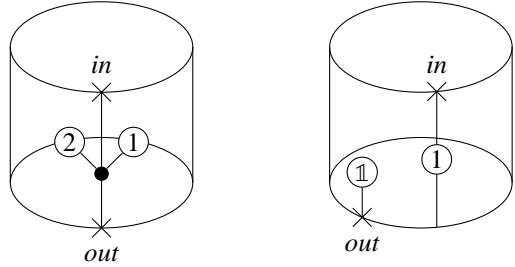
\begin{center}
\begin{tikzpicture}[
scale=.6,
int/.style={circle, draw, fill, minimum size=5pt, inner sep=0},
ext/.style={circle, draw, fill=white, minimum size=5pt, inner sep=1pt},
xst/.style={cross out, draw, minimum size=5pt}
]

\begin{scope}[xshift=6cm, shift={(-2,0)}]
\draw (0,0) ellipse (2 and 1);
\draw (0,3) ellipse (2 and 1);
\draw (-2,0)--(-2,3) (2,0)--(2,3);
\node [xst, label=-90:{$out$}] (out) at ($(0,0)+(-130:2 and 1)$) {};
\node [xst, label=90:{$in$}] (in) at ($(0,3)+(-80:2 and 1)$) {};
\draw (out)+(0,1) node [ext] (e3) {$\mathbb{1}$};
\draw (out.base) edge (e3);
\draw (in.base) -- ++(0,-1.5) node[ext] {1} -- +(0,-1.5);
\end{scope}

\begin{scope}[xshift=6cm, shift={(-9,0)}]
\draw (0,0) ellipse (2 and 1);
\draw (0,3) ellipse (2 and 1);
\draw (-2,0)--(-2,3) (2,0)--(2,3);
\node [xst, label=-90:{$out$}] (out) at ($(0,0)+(-90:2 and 1)$) {};
\node [xst, label=90:{$in$}] (in) at ($(0,3)+(-90:2 and 1)$) {};
\draw (out)+(0,1.2) node [int] (e3) {} edge (in.base);
\draw (out.base) edge (e3);
\draw (e3) +(0.7,0.7) node [ext] (e1) {1};
\draw (e3) +(-0.7,0.7) node [ext] (e2) {2};
\draw (e3) edge (e1) edge (e2);
\end{scope}
\end{tikzpicture}
\end{center}
The cohomology of the colored operad $\KS':=\bpm \Br & \KS_1' \epm$ is the operad governing pre-calculus structures, but that does not play a role here. We can map $\KS'_\infty \to \KS_\infty \to C(\EFM_2)$.
Note that $\FM_{2,1}(j) = \FM_2(j+1)\times S_1$ by definition. We claim that the images of elements of the moperadic part of $\KS'_\infty$ land in the spaces $C(\FM_2(\cdot+1)\times \{1\}\subset \FM_{2,1}(\cdot)$. In other words the input and output directions are always aligned. This is true since the subspaces $V_\Gamma$ for $\Gamma$ graphs in $B(\KS_1')$ are all contained in $\FM_2(\cdot+1)\times \{1\}$. 

Let $\ELie\{1\}$ be the two-colored operad governing a $\Lie\{1\}$-algebra and a module over it. There is a canonical map $\ELie\{1\}\to \KS'\subset \KS$. Hence, we can embed the minimal resolution $\hoELie_1$ into $\KS'_\infty$ by the composition
\[
\hoELie_1 \to \Omega(B(\ELie\{1\})) \to \KS'_\infty.
\]
Let us restrict the map $\KS'_\infty\to \EFM_2$ to $\hoELie_1$. The element in $\hoELie_1$ representing the $j$-th component of the $\Lie^{(1)}_\infty$-module structure is mapped to some chain, say $\tilde c_j\in C(\FM_2(j+1)\times \{1\})$. We claim that $\tilde c_j=\tilde f_j$, where $\tilde f_j$ is the fundamental chain of $\FM_2(j+1)\times \{1\}\subset \FM_{2,1}(j)$. This is true for $j=1$. Assume inductively it is true for $j<n$. Then $\tilde c_n$ satisfies an equation of the form
\[
\p \tilde c_n = (\cdots) = \p \tilde f_n
\]
Hence $\tilde c_n-\tilde f_n$ is closed, hence exact, and hence zero by dimensionality reasons.
Again by a (variant of a) Kontsevich vanishing lemma, it follows that the $\Lie^{(1)}_\infty$-module structure on $\Omega_\bullet$ is in fact a $\Lie\{1\}$-module structure. 

Finally we have to show that the map of $\Lie^{(1)}_\infty$-modules $C_\bullet\to \Omega_\bullet$ agrees with B. Shoikhets map. It suffices to show that the element of $B(\KS_1')$ governing the $j$-th component of that map is sent to the fundamental chain $\tilde F_j$ of
\[
D_K(j+1)\times \{1\}\subset D_S(j)\cong D_K(j+1)\times S^1.
\]
First, by arguments similar to the above, we note that the image of $B(\KS_1')$ is contained in $C(D_K(j+1)\times \{1\})$.
Let the image of the (image of the) $j$-th generator of $\Lie^{(1)}_\infty$ be denoted $\tilde C_j\in C(D_K(j+1)\times \{1\})$.
By construction of the map $\tilde C_0 = \tilde F_0$. (Note that $\tilde F_0$ is the chain of a point.) Next suppose that $\tilde C_j=\tilde F_j$ for $j<n$. Then $C_n$ satisfies an equation of the form 
\[
\p \tilde C_n = (\cdots) = \p \tilde F_n
\]
where $(\cdots)$ is a chain in $C(D_K(j+1)\times \{1\})$ build from the $c_i=f_i$, $C_i=F_i$, $\tilde c_i=\tilde f_i$ and $\tilde C_j=\tilde F_j$ for $j<n$. Hence $\tilde C_n - \tilde F_n$ is closed, hence exact by degree reasons and hence zero by dimensionality reasons. This shows Theorem \ref{thm:ksinfty}.

\makeatletter{}\section{Twisted versions of the operads and operad maps}
\label{sec:twistedops}
In the previous section we constructed maps of colored operads
\[
 \homKS \to \bigChains \to \bigGra \to \End(\bigV)=\End(\Dpoly\oplus C_\bullet\oplus\Tpoly\oplus  \Omega_\bullet ).
\]
This gives us several maps and structures, e.~g., a $\Br_\infty$ map $\Tpoly\to \Dpoly$. In deformation quantization, one usually is given a Poisson structure $\pi$, and then wants to twist the above morphisms and structures. For example, one wants to construct a $\Br_\infty$ map 
\[
 \Tpoly^\pi \to \Dpoly^\pi
\]
where the left hand side are the multivector fields with differential the Schouten bracket with $\pi$ and some $\Br_\infty$-structure to be constructed, and $\Dpoly^\pi$ is the multidifferential Hochschild cochain complex of the quantum (star product) algebra.
In our situation it might not be clear a priori how to twist, in fact we will change (twist) some of the operads to do that. That is the goal of the present section.

\subsection{Short description using the notion of operadic twisting}
Suppose we have a representation of some operad $\op P$ on some vector space $V$ given by a composition
\[
 \op P \to \op Q \to \op R \to \End(V)
\]
where $\op Q, \op R$ are some other auxiliary operads. Suppose further we have some Maurer-Cartan element\footnote{We also tacitly assume that a map $\Lie^{(k)}_\infty\to \op P$ is specified so that we can speak about Maurer-Cartan elements in a $\op P$-algebra.} $\pi\in V$ and want to twist the $\op P$ representation on $V$ to a $\op P$ representation on $V^\pi$. Usually this is not possible, the formalism of operadic twisting merely guarantees a representation of $\Tw \op P$ by:
\[
 \Tw\op P \to \Tw\op Q \to \Tw\op R \to \End(V^\pi).
\]
But now assume further that, say, $\op Q$ (or $\op R$) is natively twistable, i.~e., that there is a map $\op Q\to \Tw\op Q$. Then it is possible to twist the $\op P$ representation via
\[
 \op P \to \op Q \to \Tw\op Q \to \Tw\op R \to \End(V^\pi).
\]
This situation exactly occurs in our case. For us $\op P=\homKS$, $\op Q=\bigChains$ and $\op R=\bigGra$. Their twists are discussed in Appendix \ref{sec:coloredtwists}.
\begin{prop}
 The colored operad $\bigChains$ is natively twistable.
\end{prop}
The proof (sketch) will be given in Appendix \ref{sec:bigchainsnattwist}.
In principle, from this result it is clear how to twist. However, we want to see explicitly how the formulas look like. Furthermore the twisted version $\Tw\bigGra$ contains some very interesting sub-operad $\bigGraphs$.

\subsection{The operad \texorpdfstring{$\Graphs$}{Graphs}}
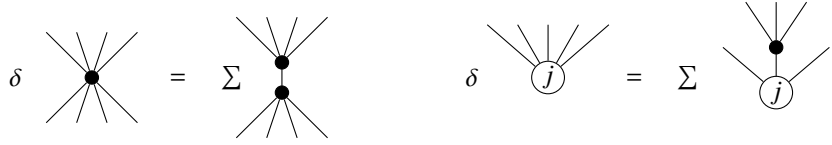
\begin{figure}
 \centering
\makeatletter{}\[
\begin{tikzpicture}[
int/.style={circle, draw, fill, minimum size=5pt, inner sep=0},
ext/.style={circle, draw, fill=white, minimum size=5pt, inner sep=1pt},
helper/.style={coordinate},point/.style={circle, draw, fill, inner sep =1pt},
de/.style={-triangle 60},
point/.style={circle, draw, fill, minimum size=3pt, inner sep=0pt},
xst/.style={cross out, draw, minimum size=5 },
]
\begin{scope}[xshift=0cm]
\node at (0,0) {$\delta$};
\node[int] (i1) at (1,0) {};
\foreach \y in {.6,-.6}
\foreach \x in {-.6,.2,-.2,.6}
\draw (i1)--+(\x,\y) ;

\node at (2.5,0) {$=\quad \sum$};

\draw (3.5,.2) node[int] (i2) {} --+(0,-.4) node[int] (i3) {};

\foreach \x in {-.6,.2,-.2,.6}
{
\draw (i2)--+(\x,.6) ;
\draw (i3)--+(\x,-.6) ;
}
\end{scope}

\begin{scope}[xshift=6cm]
\node at (0,0) {$\delta$};
\node[ext] (e1) at (1,0) {$j$};

\foreach \x in {-.8,.4,0,-.4,.8}
\draw (e1)--+(\x,.7) ;

\node at (2.5,0) {$=\quad \sum$};

\draw (4,.4) node[int] (i1) {} --+(0,-.6) node[ext] (e2) {$j$};

\foreach \x in {.4,0,-.4}
\draw (i1)--+(\x,.6) ;
\foreach \x in {.7,-.7}
\draw (e2)--+(\x,.6) ;

\end{scope}

\end{tikzpicture}
\]  
\caption{\label{fig:graphsdiff} Schematic drawing of the differential on $\Graphs$.}
\end{figure}

Twisting the operad $\Gra$ one obtains an operad $\fGraphs:=\Tw \Gra$.\footnote{cf. the corresponding object in \cite{megrt}.} Elements are (possibly infinite) linear combinations of graphs as in $\Gra$, but with two kinds of vertices: \emph{External} vertices, wich are numbered, and \emph{internal}, ``unidentifiable'' vertices of degree +2. In pictures we will draw external vertices white and internal vertices black.
The differential is given by vertex splitting and depicted schematically in Figure \ref{fig:graphsdiff}.
\begin{defprop}
 The operad $\fGraphs$ contains a suboperad $\Graphs$ spanned by graphs with the following properties:
\begin{enumerate}
 \item All internal vertices are at least trivalent.
 \item There are no connected components containing only internal vertices.
\end{enumerate}
\end{defprop}
\begin{proof}
We have to show that the subspaces spanned by those graphs are closed under the differential and under operadic compositions. The latter statement is easy, since the operadic composition never decrease the valence of vertices, nor does it produce new connected components. The differential also does not produce new connected components, but it might (a priori) create internal vertices of valence 1 or 2. It is shown in \cite{megrt} that in fact it does not. Consider first the graphs with valence 1 vertices potentially occuring in 
\[
\delta \Gamma  
= \pm \Gamma \bullet \multimapdotboth 
\pm \Gamma \circ \multimapdotbothA
\pm \multimapdotbothA \circ \Gamma.
\]
Here the right hand side shall depict the terms occuring in the definition of the twisted differential (see the Appendix in \cite{megrt}).
Graphs with valence 1 internal vertices can be produced by all three terms. However, the contribution of the third term exactly cancels the contributions of the first two.
Valence 2 internal vertices can be produced by both the first and the second term.
However, graphs with such vertices always come in pairs, two for each edge in $\Gamma$. One can check that these two graphs occur with opposite signs and hence cancel.
\end{proof}

The operad $\Graphs$ was introduced by M. Kontsevich \cite{K2}. He also showed the following Theorem.
\begin{thm}[\cite{K2, LV}]
\label{thm:graphscohom}
 $H(\Graphs)=\Ger$.
\end{thm}
There is an explicit quasi-isomorphism $\Ger\to \Graphs$ given by the formulas of the remark in section \ref{sec:GrdGr}.
There is a natural projection $\Graphs\to \Gra$, sending to zero all graphs with internal vertices.
The map $C(\FM_2)\to \Gra$ factors through $\Graphs$. The map $C(\FM_2)\to \Graphs$ has also been constructed by M. Kontsevich \cite{K2}. It is given by the formula:
\[
 \phi(c) = \sum_{n \geq 0} \sum_\Gamma \Gamma \, \int_{\pi_n^{-1}c } \omega_\Gamma.
\]
Here the sum runs over graphs in $\Graphs$ with $n$ internal vertices and $\pi_n : \FM_2(m+n)\to \FM_2(m)$ is the forgetful map. 

\begin{ex}
 In effect, composing the maps $\Br_\infty \to C(\FM_2)\to \Graphs$ one obtains, for example, a $\Br_\infty$-structure on $\Tpoly$ with a specific Poisson structure $\pi$ chosen. 
Let us consider the $A_\infty$-part of that structure. Recall from Example \ref{ex:BronTpoly} that in the untwisted case the $A_\infty$-structure is just the usual commutative algebra structure. In the twisted case this no longer holds, the $A_\infty$-structure is nontrivial. Let me raise the

{\bf Question:} Does this universal $A_\infty$-structure, i.e., the part $A_\infty\to \Graphs$ of a quasi-isomorphism $\Br_\infty\to \Graphs$ already suffice to recover the whole map $\Br_\infty \to \Graphs$ up to homotopy?

This is equivalent to saying that ``all information about the quantization'' is already encoded in this $A_\infty$-structure.
\end{ex}
 
\begin{rem}
One could in fact omit the first condition in the Definition/Proposition, without altering the cohomology of $\Graphs$, as shown in \cite{megrt}.
\end{rem}

The cohomology of the full operad $\Tw\Gra$ is also interesting. It has essentially been computed in \cite{megrt}.

\begin{prop}
\label{prop:twgracohom}
The cohomology of $\fGraphs$ is 
\[
H(\fGraphs(n)) = 
\begin{cases}
\Ger(n) \otimes S( \prod_{k=5,9,\dots} \R[-k]\oplus H(\GC_2)   )
&\quad \text{for $n>0$} \\
S( \prod_{k=5,9,\dots} \R[-k]\oplus H(\GC_2)   )
&\quad \text{for $n=0$} 
\end{cases}
\]
where $H(\GC_2)$ is the cohomology of M. Kontsevich's graph complex and $S(\cdots)$ denotes the completed symmetric product space. The factors $\R[-k]$ correspond to cycles with $k$ vertices and $k$ edges.
\end{prop}

Cohomology classes can be represented by (linear combinations of) graphs that have (i) one or more connected components with external but without internal vertices and (ii) one or more connected components with only internal vertices. Such a class acts on $\Tpoly$ with some chosen Maurer-Cartan element (i.~e., a Poisson structure) $\pi$ in the following manner. First the parts of the graphs with external vertices yield a Gerstenhaber operation. The connected components produce some $\pi$-closed multivector field out of $\pi$. This multivector field gets multiplied to the result of the $\Ger$ operation.

\subsection{\texorpdfstring{$\Graphs_1$}{Graphs1} moperad}
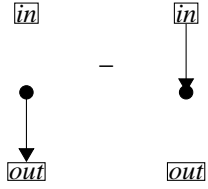
\begin{figure}
 \centering
\makeatletter{}\[
\begin{tikzpicture}[
scale=.7,
int/.style={circle, draw, fill, minimum size=5pt, inner sep=0},
ext/.style={circle, draw, fill=white, minimum size=5pt, inner sep=1pt},
helper/.style={coordinate},point/.style={circle, draw, fill, inner sep =1pt},
de/.style={-triangle 60},
point/.style={circle, draw, fill, minimum size=3pt, inner sep=0pt},
xst/.style={cross out, draw, minimum size=5 },
boxed/.style={ draw, inner sep=0.5 }
]

\begin{scope}[]
\node [boxed] (out) at ($(0,0)+(-90:2 and 1)$) {$\vout$};
\node [boxed] (in) at ($(0,3)+(-90:2 and 1)$) {$\vin$};
\draw[triangle 60-] (out)--+(0,1.5) node[int] {};
\end{scope}
\node at (1.5,1) {$-$};
\begin{scope}[xshift=3cm]
\node [boxed] (out) at ($(0,0)+(-90:2 and 1)$) {$\vout$};
\node [boxed] (in) at ($(0,3)+(-90:2 and 1)$) {$\vin$};
\draw[-triangle 60] (in)--+(0,-1.5) node[int] {};
\end{scope}

\end{tikzpicture}
\]  
\caption{\label{fig:graphs1mc} The Maurer-Cartan element used to twist $\Gra_1$. }
\end{figure}
\begin{figure}
 \centering
\makeatletter{}\[
\begin{tikzpicture}[
scale=.7,
int/.style={circle, draw, fill, minimum size=5pt, inner sep=0},
ext/.style={circle, draw, fill=white, minimum size=5pt, inner sep=1pt},
helper/.style={coordinate},point/.style={circle, draw, fill, inner sep =1pt},
de/.style={-triangle 60},
point/.style={circle, draw, fill, minimum size=3pt, inner sep=0pt},
xst/.style={cross out, draw, minimum size=5 },
boxed/.style={ draw, inner sep=0.5 },
]

\begin{scope}[]
\node [boxed] (out) at ($(0,0)+(-120:2 and 1)$) {$\vout$};
\node [boxed] (in) at ($(0,3)+(-120:2 and 1)$) {$\vin$};
\node[ext]  (e1) at (1,1.5) {$1$};
\draw[triangle 60-] (out)--+(1,1.5) node[int] (i1) {};
\draw[-triangle 60] (in)--(out);
\draw[-triangle 60] (in)--(i1);
\draw[-triangle 60] (i1)--(e1);
\end{scope}
\begin{scope}[xshift=5cm]
\node [boxed] (out) at ($(0,0)+(-120:2 and 1)$) {$\vout$};
\node [boxed] (in) at ($(0,3)+(-120:2 and 1)$) {$\vin$};
\node[ext]  (e1) at (1,1.5) {$1$};
\draw[triangle 60-] (out)--+(1,1.5) node[int] (i1) {};
\draw[-triangle 60] (e1)--(i1);
\end{scope}

\begin{scope}[xshift=10cm]
\node [boxed] (out) at ($(0,0)+(-120:2 and 1)$) {$\vout$};
\node [boxed] (in) at ($(0,3)+(-120:2 and 1)$) {$\vin$};
\node[ext]  (e1) at (1,1.5) {$1$};
\draw[triangle 60-] (out)--+(1,1.5) node[int] (i1) {};
\draw[-triangle 60] (in)--(out);
\draw[-triangle 60] (in)--(i1);
\end{scope}

\end{tikzpicture}
\]  
\caption{\label{fig:graphs1exnonex} The left and middle graphs are in $\Graphs_1$. The right hand graph is not, because after deleting $\vin$ and $\vout$ there is one connected component with only internal vertices. }
\end{figure}
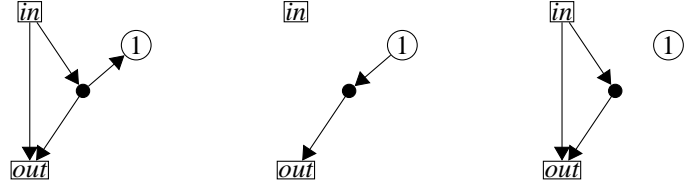
\begin{figure}
 \centering
\makeatletter{}\usetikzlibrary{matrix}
\usetikzlibrary{arrows}
\usetikzlibrary{shapes}
\usetikzlibrary{through}
\usetikzlibrary{calc,3d}
\usetikzlibrary{decorations,decorations.pathmorphing}
\[
\begin{tikzpicture}[
scale =.75,
int/.style={circle, draw, fill, minimum size=5pt, inner sep=0},
ext/.style={circle, draw, fill=white, minimum size=5pt, inner sep=1pt},
helper/.style={coordinate},point/.style={circle, draw, fill, inner sep =1pt},
de/.style={-triangle 60},
point/.style={circle, draw, fill, minimum size=3pt, inner sep=0pt},
xst/.style={cross out, draw, minimum size=5 },
boxed/.style={ draw, inner sep=0.5 },
]
\begin{scope}[xshift=0cm]
\node at (0,0) {$\delta$};
\node[int] (i1) at (1,0) {};
\foreach \y in {.6,-.6}
\foreach \x in {-.6,.2,-.2,.6}
\draw (i1)--+(\x,\y) ;

\node at (2.5,0) {$=\quad \sum$};

\draw (3.5,.3) node[int] (i2) {} --+(0,-.6) node[int] (i3) {};
\draw[-triangle 60] (i2) edge (i3);
\foreach \x in {-.6,.2,-.2,.6}
{
\draw (i2)--+(\x,.6) ;
\draw (i3)--+(\x,-.6) ;
}

\node at (4.3,0) {$+\, \sum$};

\draw (5,.3) node[int] (i2) {} --+(0,-.6) node[int] (i3) {};
\foreach \x in {-.6,.2,-.2,.6}
{
\draw (i2)--+(\x,.6) ;
\draw (i3)--+(\x,-.6) ;
}
\draw[-triangle 60] (i3) edge (i2);
\end{scope}

\begin{scope}[xshift=9cm]
\node at (0,0) {$\delta$};
\node[ext] (e1) at (1,0) {$j$};

\foreach \x in {-.8,.4,0,-.4,.8}
\draw (e1)--+(\x,.7) ;

\node at (2.5,0) {$=\quad \sum$};

\draw (4,.6) node[int] (i1) {} --+(0,-.8) node[ext] (e2) {$j$};
\draw[-triangle 60] (e2) edge (i1);

\foreach \x in {.4,0,-.4}
\draw (i1)--+(\x,.6) ;
\foreach \x in {.7,-.7}
\draw (e2)--+(\x,.6) ;

\node at (5.2,0) {$+ \, \sum$};

\draw (6.4,.6) node[int] (i1) {} --+(0,-.8) node[ext] (e2) {$j$};
\draw[-triangle 60] (i1) edge (e2);
\foreach \x in {.4,0,-.4}
\draw (i1)--+(\x,.6) ;
\foreach \x in {.7,-.7}
\draw (e2)--+(\x,.6) ;

\end{scope}

\begin{scope}[yshift=-4cm]
\node at (0,0) {$\delta$};
\node[boxed] (e1) at (1,0) {$\vin$};

\foreach \x in {-.8,.4,0,-.4,.8}
\draw[triangle 60-] (e1)--+(\x,.7) ;

\node at (2.5,0) {$=\quad \sum$};

\draw (4,.6) node[int] (i1) {} --+(0,-.8) node[boxed] (e2) {$\vin$};
\draw[triangle 60-] (e2) edge (i1);

\foreach \x in {.4,0,-.4}
\draw[triangle 60-] (i1)--+(\x,.7) ;
\foreach \x in {.7,-.7}
\draw[triangle 60-] (e2)--+(\x,.7) ;

\node at (5.2,0) {$+ \, \sum$};

\draw (6.4,.6) node[int] (i1) {} +(0,-.8) node[boxed] (e2) {$\vin$};

\foreach \x in {.4,0,-.4}
\draw[triangle 60-] (i1)--+(\x,.7) ;
\draw[triangle 60-, dotted] (i1)--+(1.4,1.2) ;
\foreach \x in {.7,-.7}
\draw[triangle 60-] (e2)--+(\x,.7) ;
\end{scope}

\begin{scope}[yshift=-4cm, xshift= 9cm]
\node at (0,0) {$\delta$};
\node[boxed] (e1) at (1,0) {$\vout$};

\foreach \x in {-.8,.4,0,-.4,.8}
\draw[-triangle 60] (e1)--+(\x,-.7) ;

\node at (2.5,0) {$=\quad \sum$};

\draw (4,-.6) node[int] (i1) {} --+(0,.8) node[boxed] (e2) {$\vout$};
\draw[-triangle 60] (e2) edge (i1);

\foreach \x in {.4,0,-.4}
\draw[-triangle 60] (i1)--+(\x,-.7) ;
\foreach \x in {.7,-.7}
\draw[-triangle 60] (e2)--+(\x,-.7) ;

\node at (5.2,0) {$+ \, \sum$};

\draw (6.4,-.6) node[int] (i1) {} +(0,.8) node[boxed] (e2) {$\vout$};

\foreach \x in {.4,0,-.4}
\draw[-triangle 60] (i1)--+(\x,-.7) ;
\draw[-triangle 60, dotted] (i1)--+(1.4,-1.2) ;
\foreach \x in {.7,-.7}
\draw[-triangle 60] (e2)--+(\x,-.7) ;
\end{scope}
\end{tikzpicture}
\]  
\caption{\label{fig:graphs1diff} The differential in $\Graphs_1$. The dotted edges are to be reconnected to any other vertex in the graph. }
\end{figure}
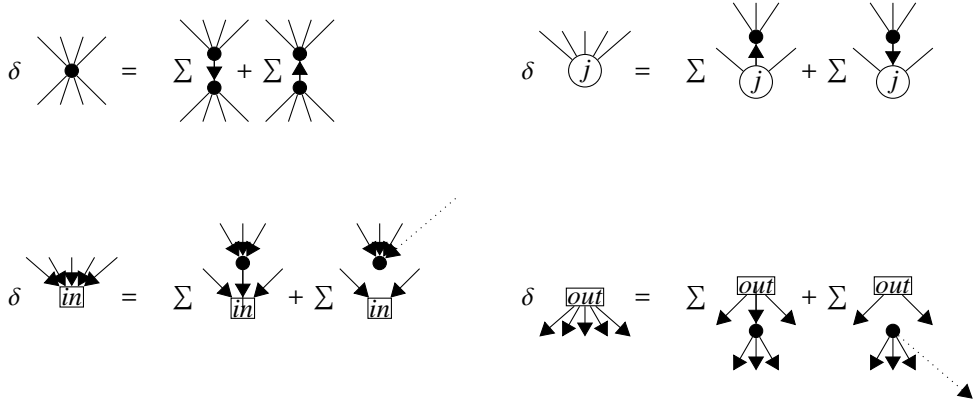
Similarly, the moperad $\Gra_1$ can be twisted to a moperad $\fGraphs_1$. The Maurer-Cartan element necessary for the twist we choose to be the one depicted in Figure \ref{fig:graphs1mc}. On differential forms, it corresponds to taking the Lie derivative by the Poisson structure. 

The differential contains two terms, see Figure \ref{fig:graphs1diff}:
\begin{enumerate}
 \item Splitting of any vertex into that vertex and an internal vertex. The original and newly created vertices are connected by an edge.
 \item Splitting of the vertices $in$ or $out$ into $in$ or $out$ and one internal vertex. The original and newly created vertices are \emph{not} connected by an edge. Instead, the internal vertex is connected to any other vertex in the graph by a new edge.
\end{enumerate}

\begin{defprop}
\label{defprop:Graphs1}
The moperad $\fGraphs_1$ contains a sub-$\Graphs$-moperad $\Graphs_1$ spanned by graphs $\Gamma$ of the following form:
\begin{enumerate}
 \item All internal vertices of $\Gamma$ are at least bivalent.
 \item There are no internal vertices of valence 2 that have 1 incoming and one outgoing edge.
 \item Let $\tilde \Gamma$ be the graph obtained by deleting $\vin$ and $\vout$. Then there are no connected components in $\tilde{\Gamma}$ containing internal, but not external vertices.\footnote{This means, for example, that an edge connecting $\vout$ to $\vin$ is allowed, but not a configuration like $\vout\to\bullet\to\vin$. We use this condition later to compute the cohomology of $\Graphs_1$ to be $\calc_1$.} 
\end{enumerate}
\end{defprop}
\begin{proof}
We have to show that $\Graphs_1$ is closed (i) under the differential (ii) under the operadic right $\Graphs$-action and (iii) under moperadic compositions. Statement (ii) is clear since the $\Graphs$-action can neither produce $<3$-valent internal vertices, nor connected components without external vertices. Similarly, it is easily checked that the moperadic composition of graphs $\Gamma_1, \Gamma_2\in \Graphs_1$ cannot produce new graphs violating either of the three conditions. The most difficult assertion to check is closedness under the differential. Consult Figure \ref{fig:graphs1diff} for a picture of the differential. The operations depicted in the top row can produce internal vertices of valence one, with an incoming or outgoing edge. Those with an incoming edge are killed by the graphs occurring in the second sum of the terms on the bottom left. Those with an outgoing edge are killed by a graph appearing in the second sum on the bottom right. 
Next consider the production of valence two internal vertices, with one incoming and one outgoing edge. Those graphs can potentially be produced by terms in the top row of the figure. However, there are two for every edge:
\begin{center}
\begin{tikzpicture}
[every edge/.style={draw, -triangle 45}, ]
\node[int] (v1) at(0,0) {};
\node[int] (v2) at(.8,0) {};
\node[int] (v3) at(2,0) {};
\draw (v1) edge (v2);
\draw (v2) edge (v3);
\begin{scope}[xshift=4cm]
\node[int] (v1) at(0,0) {};
\node[int] (v2) at(1.2,0) {};
\node[int] (v3) at(2,0) {};
\draw (v1) edge (v2);
\draw (v2) edge (v3);
\end{scope}
\end{tikzpicture}
\end{center}

Here the first term comes from splitting the left vertex and the second from splitting the right. Checking the sign, both terms cancel.
Note that there are no similar terms coming from splittings of other vertices because (and only because) we assume that the graph we started with contains itself no valence 2 internal vertices with one incoming and one outgoing edge, and no valence 1 internal vertices. Finally let us consider condition (3). It is clear that none of the terms of the differential can create connected components with only internal vertices if there were none before. Connected here means connected after deleting $\vin$ and $\vout$. This is important since the terms in the bottom row can increase the number of connected components in the ordinary sense.
\end{proof}

Examples and non-examples of $\Graphs_1$-graphs are shown in Figure \ref{fig:graphs1exnonex}.

\begin{rem}
 Note that by definition $\Graphs_1(0)=\R \oplus \R[1]$. The two operations are the unit and ``$B$'', the graph with one edge from $out$ to $in$. The latter acts as the de Rham differential on differential forms.
\end{rem}

\begin{prop}
\label{prop:graphs1cohom}
 There is a quasi-isomorphism of colored operads 
\[
\calc\to \bpm \Graphs & \Graphs_1\epm. 
\]
The map is given by sending the generators $d$ and $\iota$ to the graphs shown in Figure \ref{fig:dliota}.
\end{prop}
The proof can be found in Appendix \ref{sec:graphs1cohomproof}.
There is a natural projection $\Graphs_1\to \Gra_1$ sending graphs with internal vertices to zero. The map $C(\FM_{2,1})\to \Gra_1$ factors through $\Graphs_1$, with the map $C(\FM_{2,1})\to \Graphs_1$ being described by the formula
\[
  \sum_{n \geq 0} \sum_\Gamma \Gamma \, \int_{\pi_n^{-1}c } \omega_\Gamma.
\]
Here the sum runs over all $\Graphs_1$-graphs and $\pi_n : \FM_{2,1}(m+n)\to \FM_{2,1}(m)$ is the forgetful map. 
The differential form $\omega_\Gamma$ assigned to a graph $\Gamma$ is defined as in section \ref{sec:FM1toGra1}.

\subsection{\texorpdfstring{$\SGraphs$}{SGraphs} operadic bimodule}\label{sec:SGraphs}
We next want to twist the operadic $\PT$-$\Gr$ bimodule $\SGra$ according to the procedure described in appendix \ref{sec:bimodtwist}. To do this, we need a suitable Maurer-Cartan element $m$ in $\Tw\SGra(0)$. Such an element is nothing but (the graph version of) a universal star product. We will take for $m$ the universal star product constructed by M. Kontsevich \cite{K1}.\footnote{In fact, we need to take this particular $m$ because it is the image of the (natural) MC element chosen in $\Tw D_K$ under the map $D_K\to\Tw\SGra$.} Let us call the resulting operadic bimodule $\fSGraphs=\Tw\SGra$.
The differential has three terms:
\begin{itemize}
\item One ``internal vertex splitting'' term which is given by insertions of the graph 
\tikz{ 
\node[int] (v1) at (0,0) {}; 
\node[int] (v2) at (.7,0) {}; 
\draw[-triangle 45] (v1) edge (v2);
}
at all internal vertices.
\item One ``external vertex splitting'' term which is given by insertions of the graph 
$\begin{tikzpicture}
\node[ext] (v1) at (0,0) {}; 
\node[int] (v2) at (0.7,0) {}; 
\draw[-triangle 45] (v1) edge (v2);
\end{tikzpicture}
+
\begin{tikzpicture}
\node[ext] (v1) at (0,0) {}; 
\node[int] (v2) at (0.7,0) {}; 
\draw[-triangle 45] (v2) edge (v1);
\end{tikzpicture}
$
at all external vertices.
\item Terms which together resemble the ``Gerstenhaber bracket with the Hochschild cochain describing the universal star product $m$''.
\end{itemize}

Following the previous section, one would like to identify a certain sub-bimodule of $\fSGraphs$ which is quasi-isomorphic to $C(D_K)$. In particular the cohomology should be isomorphic to $\Ger$. Unfortunately, the author does not know how to define this sub-bimodule. Hence we will have to go with the following definition, which will leave the cohomology ``too big''.

\begin{defprop}
\label{defprop:SGraphs}
The operadic $\Br$-$\Graphs$-sub-bimodule $\SGraphs\subset \fSGraphs$ is spanned by graphs of the following form
\begin{enumerate}
 \item There is at least one external vertex.
 \item There are no internal vertices of valence zero.
 \item There are no internal vertices of valence one, with the incident edge outgoing.
 \item There are no internal vertices of valence two, with one edge incoming and one outgoing.
\end{enumerate}
\end{defprop}

\begin{proof}
We have to check that $\SGraphs$ is closed under (i) the differential, (ii) the right $\Graphs$-action and (iii) the left $\Br$-action. The first two conditions imposed by the definition are immediate since none of the operations can delete external vertices or introduce valence zero vertices. We focus on the other two.

It is evident that the right $\Graphs$-action respects the above conditions, since it does not introduce valence one or two vertices. To check statements (i) and (iii) we need the following properties of M. Kontsevich's universal star product $m$, which we take for granted (see \cite{K1}).
\begin{itemize}
\item Property 1: There is only one graph in $m$ with a valence 1 internal vertex, which is the one with one internal type I vertex, one type II vertex and an edge connecting them. Its prefactor in $m$ is 1.
\item Property 2: No graph in $m$ contains a valence two internal vertex with one incoming and one outgoing edge.
\end{itemize}

Consider the left $\Br$ action. The generator $T_n$ (see figure \ref{fig:KSgen}) cannot produce valence 1 vertices with the edge incoming if there were none before. Similarly, valence two vertices with one edge outgoing and one outgoing can at best be produced from valence one vertices with the edge outgoing, but those we excluded.
Next consider the other generators, $T_n'$. Since $n\geq 2$ the parts of the Maurer-Cartan element $m$ that occur are those with at least two type II vertices. But they do not contain vertices of the excluded types by Property 1 and Property 2.

Finally consider the effect of the differential. Each of the three terms of the differential as recalled above can produce valence 1 and 2 vertices of the forbidden types. One has to check that the contributions cancel. Let us first consider internal valence one vertices with the edge outgoing. They can be produced (a) the first term of the differential as in the list above, this for each internal vertex one graph like this:
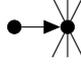
\[
\begin{tikzpicture}
\node[int] (e1) at (0,0) {};
\node[int] (e2) at (0.7,0) {};
\foreach \x in { -.2, 0, .2 }
{
\draw (e2) -- +(\x,.4);
\draw (e2) -- +(\x,-.4);
}
\draw[-triangle 45] (e1) edge (e2);
\end{tikzpicture}
\] 
Here only the relevant internal vertex is shown, not the rest of the graph. Similarly, the second part of the differential produces, for each external vertex, one graph of the form:
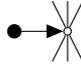
\[
\begin{tikzpicture}
\node[int] (e1) at (0,0) {};
\node[ext] (e2) at (0.7,0) {};
\foreach \x in { -.2, 0, .2 }
{
\draw (e2) -- +(\x,.4);
\draw (e2) -- +(\x,-.4);
}
\draw[-triangle 45] (e1) edge (e2);
\end{tikzpicture}
\] 
Furthermore, there two terms in the third part of the differential which can produce such valence 1 vertices, namely those involving the only graph in $m$ which contains such a vertex, see Property 1 above. Concretely, the first part is:
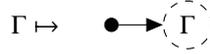
\[
\begin{tikzpicture}
\node at (-1,0) {$\Gamma \mapsto $};
\node[int] (e1) at (0,0) {};
\node[draw, dashed, circle] (e2) at (1,0) {$\Gamma$};
\draw[-triangle 45] (e1) edge (e2);
\end{tikzpicture}
\]
Here the notation shall indicate that one external vertex is added and connected to vertices of $\Gamma$ in all possible ways. In other words, it can be connected to external vertices, internal type I vertices, or type II vertices. The first contributions for type I vertices exactly cancel the corresponding terms we discussed before. (This requires a sign verification.) There remain graphs with forbidden valence one vertices connected to type II vertices as follows:
 
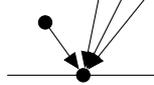
\[
\begin{tikzpicture}
\draw (-1,0)--(1,0);
\node[int] (e2) at (0,0) {};
\node[int] (e1) at (-.5,.7) {};
\draw[-triangle 45] (e1) edge (e2);
\draw[triangle 45-] (e2) -- +(0.2,1);
\draw[triangle 45-] (e2) -- +(.5,1);
\draw[triangle 45-] (e2) -- +(.8,1);
\end{tikzpicture}
\]
These terms cancel with identical ones produced by the following term of the third part of the differential:
 
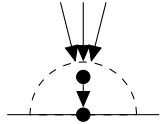
\[
\begin{tikzpicture}
\node[int] (e2) at (0,0) {};
\node[int] (e1) at (0,.5) {};
\node[semicircle,dashed, draw, minimum size=.7cm] (c) at (0,.29) {};

\draw[-triangle 45] (e1) edge (e2);
\draw[triangle 45-] (c) -- +(-0.3,1.2);
\draw[triangle 45-] (c) -- +(0,1.2);
\draw[triangle 45-] (c) -- +(.3,1.2);
\draw (-1,0)--(1,0);
\end{tikzpicture}
\]
Here the edges incoming at the dashed semicircle are to be reconnected in all possible ways. The term relevant to us is when all edges are connected to the type II vertex. This shows that the differential does not produce valence 1 vertices with outgoing edge. Next we need to check that it cannot produce valence 2 vertices with one edge incoming, one outgoing. Graphs with vertices of this form can again be produced by all three parts of the differential. Checking that they cancel is similar to the corresponding step in the proof of Definition/Proposition \ref{defprop:SGraphs} and we omit it here.
\end{proof}

\begin{rem}
Note that it is not possible to forbid valence one internal vertices at all. They may be produced by the differential.
\end{rem}

One has a canonical projection $\SGraphs\to \SGra$ sending all graphs with an internal type I vertex to zero.
The map $C(D_K)\to \SGra$ factors through $\SGraphs$. The map $C(D_K)\to \SGraphs$ takes the form 
\[
 c\mapsto \sum_{n_I,n_{II}\geq 0} \sum_\Gamma \Gamma \, \int_{\pi_{m,n_I, n_{II}}^{-1} c} \omega_\Gamma\, .
\]
Here $c$ is an $m$-chain and the sum runs over all $\SGraphs$-graphs with $n_I$ internal type I and $n_{II}$ internal type II vertices. The map $\pi_{m,n_I, n_{II}}:D_{Ke}(m+n_I, n_{II})$ is the forgetful map. The differential form is defined similarly to section \ref{sec:DKtoSGra}.

\begin{rem}
The cohomology of $\fSGraphs$ and $\SGraphs$ is computed in Appendix \ref{sec:sgraphsproofsuppl}. We will not need the result in this paper.
\end{rem}

\subsection{\texorpdfstring{$\SGraphs_1$}{SGraphs1} moperadic bimodule}
Finally we twist the moperadic bimodule $\SGra_1$ according to Appendix \ref{sec:mopbimodtwist}. No additional choices have to be made, but let us describe the resulting moperadic bimodule $\fSGraphs_1:= \Tw \SGra_1$. It is given by linear combinations of graphs as in $\SGra_1$, but now some (or all) external vertices can be replaced by ``unidentifiable'' \emph{internal} vertices.
For the purpose of globalization, we will identify a certain moperadic sub-bimodule of $\fSGraphs_1$.

\begin{defprop}
We define a $\Br$-$\Graphs$-$\KS_1$-$\Graphs_1$-$\SGraphs$ moperadic sub-bimodule $\SGraphs_1\subset\fSGraphs_1$ by the following constraints on graphs:
\begin{enumerate}
 \item There are no internal type I vertices of valence 0.
 \item There are no internal type I vertices of valence one with the incident edge outgoing.
 \item There are no internal type I vertices of valence 2 with one edge incoming and one outgoing
\end{enumerate}
\end{defprop}
\begin{proof}
We have to show that the indicated subspace is indeed closed under (i) the differential, (ii) the left $\Graphs_1$ action, (iii) the right $\Graphs$-action and (iv) the combined right action of $\KS_1$ and $\SGraphs$. Statements (ii), (iii) and (iv) are immediate since the actions cannot introduce any of the forbidden types of vertices.
The differential can in principal introduce valence one or two vertices of the forbidden type. One has to show that the respective terms cancel. Showing this is done by (almost) a copy of the argument of Definition/Proposition \ref{defprop:SGraphs}.
\end{proof}

There is a canonical projection 
\[
\SGraphs_1 \to \SGra_1
\]
sending to zero all graphs with internal type I vertices.
The map $C(D_S)\to \SGra_1$ constructed in section \ref{sec:opmaps} factors through $\SGraphs_1$.
\[
C(D_S) \to \SGraphs_1 \to \SGra_1.
\]

One can package the various twisted (m)operads and bimodules into a four colored operad:
\[
\bigGraphs := 
\bpm
\Br & \SGraphs & \Graphs \\
\KS_1 & \SGraphs_1 & \Graphs_1
\epm
\]
The map $\bigChains\to \bigGra$ factors through $\bigGraphs$, hence one has maps
\[
\homKS\to \bigChains \to \bigGraphs \to \bigGra.
\]

\subsection{Twisted version of \texorpdfstring{$\KS_\infty$}{KS-infinity} formality}
Suppose now that we are given a Poisson structure $\pi$ in $\Tpoly$. It defines a differential $\epsilon\co{\pi}{\cdot}$ on $\Tpoly[[\epsilon]]$, and a differential $\epsilon L_\pi$ on $\Omega_\bullet[[\epsilon]]$. Let us denote the resulting complexes by $\Tpoly^\pi$ and $\Omega_\bullet^\pi$. By the general theory of twisting, one obtains an action of the colored operad $\bpm \Graphs & \Graphs_1 \epm$ on $\Tpoly^\pi\oplus \Omega_\bullet^\pi$. Concretely, the action is constructed in the same manner as that of $\Gra$ and $\Gra_1$, except that for every internal vertex in graphs in $\Graphs$ or $\Graphs_1$ one inserts one copy of $\epsilon \pi$.

On the other side, one obtains a twisted representation of $\bpm \Br & \KS_1 \epm$ on 
\[
 \Dpoly^\pi \oplus \C_\bullet^\pi =(\Dpoly[[\epsilon]], d_{m_{\epsilon \pi}})\oplus (\C_\bullet[[\epsilon]], L_{m_{\epsilon \pi}}).
\]
Here $m_{\epsilon \pi}$ is Kontsevich's associative product (star product) on $C^\infty(\R^n)$ and $d_{m_{\epsilon \pi}}$ shall denote the Hochschild differential with respect to this star product. This is the same as the Gerstenhaber bracket with  $m_{\epsilon \pi}$.
The $\Br$-action on $\Dpoly^\pi$ is such that each internal vertex with more then two children acts as zero, and each internal vertex with exactly two children is interpreted as a copy of the two-cochain $m_{\epsilon \pi}$.

Twisting the actions of $\SGraphs$ and $\SGraphs_1$ in the same manner, one obtains at the end a representation of the operad $\homKS$ on the twisted (colored) vector space
\[
 \bigV^\pi := \Dpoly^\pi \oplus C_\bullet^\pi  \oplus\Tpoly^\pi \oplus \Omega_\bullet^\pi  .
\]

This amounts to the following data:
\begin{enumerate}
 \item A $\Br_\infty$ structure on $\Tpoly^\pi$.
 \item A $\Br$ structure on $\Dpoly^\pi$. This is just the restriction of the standard $\Br$-structure on the Hochschild complex of the algebra $C^\infty(\R^n)[[\epsilon]]$ with product $m_{\epsilon \pi}$.
 \item A module structure over $\Tpoly^\pi$ on $\Omega_\bullet^\pi$, governed by the moperad $\KS_{1,\infty}$.
 \item A module structure over $\Dpoly^\pi$ on $C_\bullet^\pi$, governed by the moperad $\KS_{1}$. This is again the standard $\KS_1$ structure on the Hochschild chain complex of the algebra $C^\infty(\R^n)[[\epsilon]]$.
 \item A $\Br_\infty$-map $\Tpoly^\pi\to \Dpoly^\pi$.
 \item A compatible map of modules $C_\bullet^\pi \to \Omega_\bullet^\pi$.
\end{enumerate}

\subsection{Globalization (sketch)}
\label{sec:globalization}
The ability to twist makes it possible to globalize the $\Br_\infty$ morphism $\mU$. For globalization, we use the framework of V. Dolgushev \cite{dolgushev}. It will allow us to construct a formality morphism for smooth manifolds and complex manifolds. With some changes, one can also target smooth algebraic varieties over $\C$. We will focus on the smooth case for simplicity.
The machinery of globalization is technically and notationally heavy. We will only sketch the procedure. If you are a non-expert reader, skip this section.

Let $M$ be a smooth manifold. Let $\Tpoly$ be the multivector fields on $M$ and $\Dpoly$ the multidifferential operators. Similarly $\Omega_\bullet$ will denote the differential forms on $M$ and $C_\bullet$ the continuous Hochschild chains on $A=C^\infty(M)$. We will also need Dolgushev-Fedosov resolutions (see \cite{dolgushev}) of these objects, which we denote by $\Dpoly^\fml, C_\bullet^\fml, \Tpoly^\fml, \Omega_\bullet^\fml$. We assume in the following that the reader knows these objects.
There are natural injective quasi-isomorphisms $\Dpoly\to \Dpoly^\fml$, $C_\bullet\to C_\bullet^\fml$, $\Tpoly\to \Tpoly^\fml$, $\Omega_\bullet\to\Omega_\bullet^\fml$ between the considered algebraic objects and their resolutions. We want to construct a representation of $\homKS$ on the 4 colored space
\[
V := \Dpoly\oplus C_\bullet \oplus \Tpoly \oplus \Omega_\bullet.
\]
We will construct it implicitly only, by giving zigzags of quasi-isomorphisms. First we need a $\KS$ structure on $\Tpoly\oplus \Omega_\bullet$. There is a natural $$\bpm \Graphs & \Graphs_1 \epm$$-algebra structure on $\Tpoly^\fml\oplus \Omega_\bullet^\fml$. 
It is defined by acting ``fiberwise'', with internal vertices representing copies of the Maurer-Cartan form (some one-form valued formal fiberwise vector field). To show that this defines an action, one needs to use properties of the $\bpm \Graphs & \Graphs_1 \epm$ action in the local case similar to properties P1), P3), P4), P5) written down by M. Kontsevich (\cite{K1}, section 7). For example, for the action to be well defined (gauge invariant), it is important that there are no internal vertices of valence $2$ with one incoming and one outgoing edge in $\Graphs$ and $\Graphs_1$, corresponding to property P5) of M. Kontsevich.
 By the map $\KS_\infty \to\bpm \Graphs & \Graphs_1 \epm$ just constructed, one hence obtains a $\KS_\infty$ algebra structure on $\Tpoly^\fml\oplus \Omega_\bullet^\fml$. By homotopy transfer along the quasi-isomorphism 
\[
 \Tpoly\oplus \Omega_\bullet\to \Tpoly^\fml\oplus \Omega_\bullet^\fml
\]
one hence obtains a $\KS_\infty$ algebra structure on $\Tpoly\oplus \Omega_\bullet$, together with a $\KS_\infty$ quasi-isomorphism to $\Tpoly^\fml\oplus \Omega_\bullet^\fml$. On the other side, we have natural $\KS$ structures on $\Dpoly\oplus C_\bullet$ and $\Dpoly^\fml\oplus C_\bullet^\fml$ and the quasi-isomorphism 
\[
 \Dpoly\oplus C_\bullet \to \Dpoly^\fml\oplus C_\bullet^\fml
\]
is compatible with the $\KS$ structures.
Hence it suffices to represent $\homKS$ on the 4 colored vector space
\[
\Dpoly^\fml\oplus  C_\bullet^\fml\oplus\Tpoly^\fml\oplus \Omega_\bullet^\fml.
\]
There is a natural action of $\bigGraphs$ on the above 4-colored vector space, by using the formulas from the local case fiberwise.
To check that the action is well defined and is an action, one again needs to check technical properties similar to P1)-P5) of M. Kontsevich \cite{K1}.
Hence one can pull back the $\bigGraphs$-algebra structure to obtain the desired $\homKS$-algebra structure.

\section{Recovery of several results in the literature}
\label{sec:9}
It has been observed before that the $\Lie_\infty$-formality morphisms by M. Kontsevich and B. Shoikhet on homology preserve more algebraic structure than the Lie bracket (respectively, the Lie action), see \cite{K1,mecycchains, CR, CFW}. The present work generalizes and unifies these results. As an illustration, let us show how to recover the statements in those references from Theorems \ref{thm:brinfty} and \ref{thm:ksinfty}.

\subsection{Compatibility with cup products}
M. Kontsevich \cite{K1} observed that his formality morphism respects the cup product on cohomology. He used this statement to re-prove M. Duflo's theorem. Kontsevich also gave an explicit description of the homotopy. Let us see how to recover it. The relevant $\Br$-operation, say $o$, is the product:
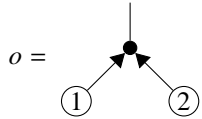
\[
 o = \raisebox{-2em}{
 \tikz{
\usetikzlibrary{arrows}
\tikzset{every edge/.style={draw,-triangle 60}}
\node[int, minimum size=5] (e1) at (0,0) {};
\node[ext] (e2) at (-135:1) {1};
\node[ext] (e3) at (-45:1) {2};
\draw (e2) edge (e1) (e3) edge (e1);
\draw (e1) -- +(0,.6);
} }
\]
We want to find a formula for the homotopy protecting the product on homology. This is given by the component $\mU_o$ of the $\Br_\infty$ map of Theorem \ref{thm:brinfty}. This component is given by eqn. \eqref{equ:Udef}, but we need to know the chain $\tilde{c}_o\in C(D_K(2))$. This chain is defined in eqn. \eqref{equ:hBrcase}, as a solution of the equation
\[
 \p \tilde{c}_o = p_1 - p_2
\]
 where $p_1, p_2$ are chains given by single points, shown in Figure \ref{fig:Kontshomotopy} (left). The obvious solution for $\tilde{c}_o$ is a degree 1 chain connecting these two points, as shown in Figure \ref{fig:Kontshomotopy} (right). This recovers M. Kontsevich's homotopy.
 
\begin{figure}
 \centering
\makeatletter{} \begin{tikzpicture}[
scale=.6,
int/.style={circle, draw, fill, minimum size=5pt, inner sep=0},
ext/.style={circle, draw, fill=white, minimum size=5pt, inner sep=1pt},
helper/.style={coordinate},point/.style={circle, draw, fill, inner sep =1pt},
de/.style={-triangle 60},
point/.style={circle, draw, fill, minimum size=3pt, inner sep=0pt},
xst/.style={cross out, draw, minimum size=5 },
boxed/.style={ draw, inner sep=0.5 },
]

\draw(-5,-0.5) -- (0,-0.5);
\draw[dashed] (-3.5,-0.5) arc (0:180:0.5);
\draw[dashed] (-0.5,-0.5) arc (0:180:0.5);
\node [int] at (-4,-0.25) {};
\node [int] at (-1,-0.25) {};
\draw(1,-0.5) -- (5.5,-0.5);
\node [int] at (3,2) {};
\node [int] at (3.5,2) {};
\draw[dashed] (3.25,2) ellipse (0.5 and 0.5);
\node at (-2.5,3.5) {$p_1$};
\node at (3.25,3.5) {$p_2$};

\begin{scope}[xshift=12cm]

\draw(-5,-0.5) -- (0,-0.5);
\draw[dashed] (-3.5,-0.5) arc (0:180:0.5);
\draw[dashed] (-0.5,-0.5) arc (0:180:0.5);
\node [int] at (-4,-0.25) {};
\node [int] at (-1,-0.25) {};
\begin{scope}[xshift=-5.5cm]
\node [int] at (3,2) {};
\node [int] at (3.5,2) {};
\draw[dashed] (3.25,2) ellipse (0.5 and 0.5);
\end{scope}

\node at (-2.5,3.5) {$\tilde c_o$};
\end{scope}

\draw plot[smooth, tension=.7] coordinates {(8,-0.25) (8.5,0.75) (9.25,1.25) (9.5,2)};
\draw plot[smooth, tension=.7] coordinates {(11,-0.25) (10.75,0.75) (10.25,1.25) (10,2)};
\end{tikzpicture} 
\caption{\label{fig:Kontshomotopy} M. Kontsevich's homotopy protecting the product on cohomology is defined by a chain (right) whose boundary is the difference of the two chains (points) on the left.}
\end{figure}
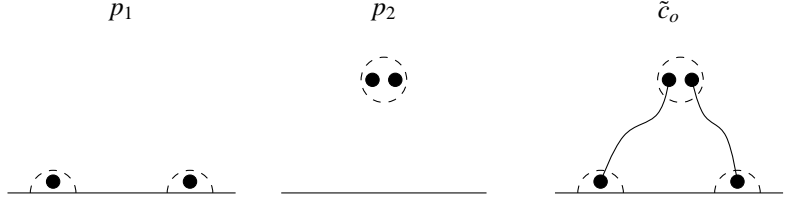

\subsection{Compatibility with the Connes-Rinehart differential}
Let us switch to Hochschild chains. It has been observed by the author in \cite{mecycchains} that B. Shoikhet's $\Lie_\infty$ morphism is compatible with the Connes-Rinehart differential $B$. The $\KS$ operation corresponding to $B$ is depicted in Figure \ref{fig:BLI}. What we can deduce from Theorem \ref{thm:ksinfty} is that this operation is preserved up to homotopy. The homotopy is given by a formula similar to \eqref{equ:Udef}, once we know the chain $c_B\in C(D_S)$ governing the homotopy. Writing down the defining equation (eqn. \eqref{equ:hKScase}), we get 
\[
 \p c_B = 0.
\]
Hence we can choose $c_B= 0$. Hence the operation $B$ is actually preserved on the nose by B. Shoikhet's morphism. This is the main result of \cite{mecycchains}.

\subsection{Compatibility with the cap product}
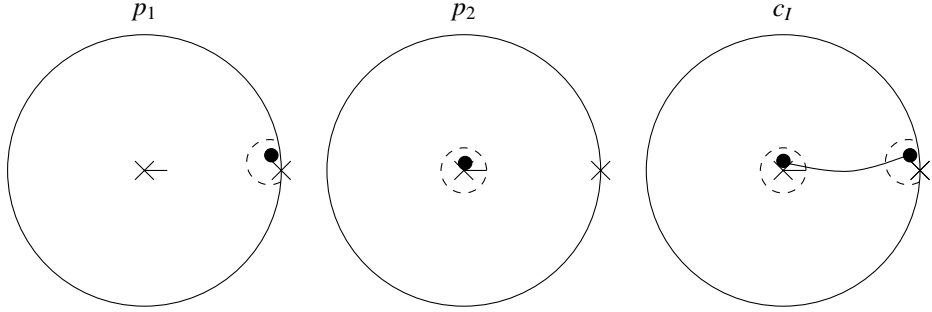
\begin{figure}
 \centering
\makeatletter{}\begin{tikzpicture}[
scale =.6,
int/.style={circle, draw, fill, minimum size=5pt, inner sep=0},
ext/.style={circle, draw, fill=white, minimum size=5pt, inner sep=1pt},
helper/.style={coordinate},point/.style={circle, draw, fill, inner sep =1pt},
de/.style={-triangle 60},
point/.style={circle, draw, fill, minimum size=3pt, inner sep=0pt},
xst/.style={cross out, draw, minimum size=5 },
boxed/.style={ draw, inner sep=0.5 },
]

\draw (0,0) node[xst] {} ellipse (3 and 3) -- (0:.5);
\node [int] at (2.7813,0.3354) {};
\node [xst] at (3,0) {};
\draw[dashed](2.991,-0.2497) arc (-60:-294:.5);

\begin{scope}[xshift=7cm]
\draw (0,0) node[xst] {} ellipse (3 and 3) -- (0:.5);
\node [int] at (0.0296,0.1646) {};
\node [xst] at (3,0) {};
\draw[dashed] (0,0) ellipse (0.5 and 0.5);

\end{scope}

\begin{scope}[xshift=14cm]
\draw (0,0) node[xst] {} ellipse (3 and 3) -- (0:.5);
\node [int] at (0.0061,0.2117) {};
\node [xst] at (3,0) {};
\draw[dashed] (0,0) ellipse (0.5 and 0.5);

\node [int] at (2.7813,0.3354) {};
\node [xst] at (3,0) {};
\draw[dashed](2.991,-0.2497) arc (-60:-294:.5);
\end{scope}

\draw plot[smooth, tension=.7] coordinates {(16.7883,0.3589) (15.453,-0.0174) (14.0531,0.1707)};
\node at (0,3.5) {$p_1$};
\node at (7,3.5) {$p_2$};
\node at (14,3.5) {$c_I$};
\end{tikzpicture} 
\caption{\label{fig:CRhomotopy} The homotopy of D. Calaque and C. Rossi is governed by a chain $c_I\in C(D_S)$ (right) whose boundary is the difference of the two points on the left.}
\end{figure}

It has been shown by D. Calaque and C. Rossi \cite{CR} and used in \cite{CRvdB} that the Shoikhet morphism is compatible with the cap product. The cap product corresponds to the $\KS$ operation $I$ depicted in Figure \ref{fig:BLI}. Again the relevant homotopy will be governed by some chain $c_I\in C(D_S)$. Eqn. \eqref{equ:hKScase} tells us to choose it to connect two points:
\[
 \p c_I = p_1-p_2.
\]
The points and the solution $c_I$ are shown in Figure \ref{fig:CRhomotopy}. This recovers the homotopy of Calaque and Rossi.

\subsection{Compatibility with the Gauss-Manin connection}
A. Cattaneo, G. Felder and the author showed in \cite{CFW} that the Shoikhet morphism also respects the Gauss-Manin connection on homology. For this it is important that the operation $H$ depicted in \ref{fig:BLI} is respected (in the correct sense). Solving equation \eqref{equ:hKScase}, one again find a certain chain $c_H\in C(D_S)$ producing the relevant homotopy. We will not depict the chain $c_H$ here but refer to \cite{CFW} instead, where a picture can be found.

\makeatletter{}\section{\texorpdfstring{$\Ger_\infty$}{G-infinity}-morphism and relation to D. Tamarkin's quantization}
\label{sec:10}
In this section we show how to convert the $\Br_\infty$ formality morphism $\Tpoly\to \Dpoly$ constructed above into a $\Ger_\infty$ morphism and prove Theorems \ref{thm:kontstamequiv} and \ref{thm:kontsextends} in the introduction. The extension to a $\Ger_\infty$ formality morphism, i.~e., the proof of Theorem \ref{thm:kontsextends}, will go through without changes in the algebraic, formal and smooth settings. For the comparison to D. Tamarkin's morphism, i.~e., Theorem \ref{thm:kontstamequiv}, we will restrict to the algebraic case as in Tamarkin's original construction.

\subsection{Review of D. Tamarkin's construction of a Formality morphism}
\label{sec:tamreview}
D. Tamarkin's construction \cite{tamanother, hinich} of formality morphisms proceeds in two steps:
\begin{enumerate}
 \item Produce a solution of Deligne's conjecture, i.~e., construct a $\Ger_\infty$-structure on the multidifferential operators, inducing the standard Gerstenhaber structure on cohomology, and whose $\Lie_\infty$-part coincides with the usual Lie algebra structure.  \item Produce a $\Ger_\infty$-formality morphism $\Tpoly^{alg}\to \Dpoly^{alg}$.
\end{enumerate}

Let us begin with the first step, i.~e.,  the construction of a $\Ger_\infty$ structure on the multidifferential operators. We will actually consider only such structures which factor through the braces operad $\Br$. Note that by cofibrancy of $\Ger_\infty$ it is sufficient to construct a chain of quasi-isomorphisms 
\[
 \Ger_\infty \to \cdots \leftarrow \Br_\infty.
\]
Then by lifting up to homotopy one can construct a quasi-isomorphism 
\[
 \Ger_\infty \to \Br_\infty \to \Br
\]
and hence obtains an action of $\Ger_\infty$ on $\Dpoly$ by pulling back the original action of $\Br$. Note that by the Kontsevich-Soibelman construction (see section \ref{sec:brtofm2} or \cite{KS1}) there is a natural (up to contractible choices) map $\Br_\infty\to C(\FM_2)$. Furthermore there is a natural quasi-isomorphism $\Ger_\infty\to \Ger$. Hence, equivalently, one may construct a zig-zag of quasi-isomorphisms
\[
 C(\FM_2) \to \cdots \leftarrow \Ger.
\]
Such a zig-zag is called a \emph{formality morphisms of the little disks operad}. Up to homotopy, one has a 1:1 correspondence
\[
 (\text{Maps $\Ger_\infty\to \Br$}) \leftrightarrow (\text{formality morphisms of the little disks operad})
\]
where the maps on the left hand side are required to induce the usual isomorphism on cohomology.
In this paper two procedures of constructing a formality of the little disks operad will be important:
\begin{itemize}
 \item The formality morphisms of the little disks operad obtained by Kontsevich \cite{K2}, defined by the zigzag
\[
 \Ger \to \Graphs \leftarrow C(\FM_2).
\]
We will denote the corresponding map from $\Ger_\infty$ to the braces operad via lifting by $K:\Ger_\infty\to \Br$.
 \item The formality morphisms of the little disks operad obtained by D. Tamarkin, using a Drinfeld associator \cite{tamarkin}, see also \cite[section 6.2]{pavol}.
Of particular interest will be the formality morphism thus obtained by using the Alekseev-Torossian Drinfeld associator \cite{ATassoc,pavol}. We will denote the corresponding map from $\Ger_\infty$ to the braces operad via lifting by $T:\Ger_\infty\to \Br$.
\end{itemize}

It was shown in \cite{pavol}, that M. Kontsevich's formality morphism of the little disks operad is homotopic to the formality morphism constructed by D. Tamarkin \cite{tamarkin}, using the Alekseev-Torossian associator. In particular, it means that the maps $K,T:\Ger_\infty\to \Br$ are homotopic.

Let us turn to the second step of D. Tamarkin's construction, assuming that a quasi-isomorphism $\Ger_\infty\to \Br$ is given.
In the original construction of D. Tamarkin one works with the algebraic version of multivector fields and multidifferential operators $\Tpoly^{alg}$ and $\Dpoly^{alg}$. 
We know that $H(\Dpoly^{alg})=\Tpoly^{alg}$, so by general homotopy transfer arguments, there is some $\Ger_\infty$ structure on $\Tpoly^{alg}$, such that there is a $\Ger_\infty$ quasi-isomorphism
\[
 (\Tpoly^{alg})' \to \Dpoly^{alg}.
\]
Here we denote by $(\Tpoly^{alg})'$ the space $\Tpoly^{alg}$ with the non-standard $\Ger_\infty$-structure. Furthermore, the induced Gerstenhaber structure on $\Tpoly^{alg}$ is the standard one. The following result of D. Tamarkin then finishes the construction.

\vspace{.2cm}

{\bf Rigidity of the Gerstenhaber algebra of multivector fields:}
For any $\Ger_\infty$ structure on the algebraic multivector fields $\Tpoly^{alg}$ inducing the standard Gerstenhaber structure, there is an essentially unique $\Ger_\infty$ quasi-isomorphism
\[
 \Tpoly^{alg} \to (\Tpoly^{alg})'.
\]

\vspace{.2cm}

For an accessible review of the proof, see \cite{hinich}.

More recently, a more explicit construction of D. Tamarkin's formality morphism has been found by Dolgushev, Tamarkin and Tsygan \cite{DTT1}, that also allows for globalization. Concretely one can write an explicit zig-zag of $\Ger_\infty$ morphisms of $\Ger_\infty$ algebras, cf. \cite[eqn. 3.8]{DTT}
\begin{equation}\label{equ:DTTdef}
 \Tpoly^{alg} \leftarrow \Omega_\Ger(B_{\Ger^\vee}(\Tpoly^{alg}))
\leftarrow \Omega_\Ger(\Xi) \to \Omega_\Ger(\FF_{\Ger^\vee}(\Dpoly^{alg})) \leftarrow \Dpoly^{alg}.
\end{equation}
Here $\Omega_\Ger(\dots)$ denotes the cobar construction of a $\Ger^\vee$-coalgebra, and $B_{\Ger^\vee}$ denotes the bar construction of a $\Ger$ algebra. $\FF_{\Ger^\vee}(\Dpoly^{alg})$ is the quasi-free $\Ger^\vee$ coalgebra cogenerated by $\Dpoly^{alg}$, with the differential canonically defined by the $\Ger_\infty$ structure on $\Dpoly^{alg}$.
More concretely, elements of $\FF_{\Ger^\vee}(\Dpoly^{alg})$ can be identified as linear combinations of co-Gerstenhaber words
\[
 \underline{A_1^1\cdots A_1^{n_1}}\wedge \cdots \wedge  \underline{A_k^1\cdots A_k^{n_k}}
\]
where the $A_{i}^j\in \Dpoly^{alg}$ and the underline indicates that one considers the sub-words $\underline{A_i^1\cdots A_i^{n_i}}$ modulo shuffle products.
The subspace $\Xi\subset \FF_{\Ger^\vee}(\Dpoly^{alg})$ consists of those words for which all $A_i^j$'s are functions, except possibly for one element of each subword $\underline{A_i^1\cdots A_i^{n_i}}$, which is allowed to be a vector field.
The space $\Xi$ also naturally embeds into $B_{\Ger^\vee}(\Tpoly^{alg})$, hence the second and third map of the zig-zag \eqref{equ:DTTdef} are defined. The first map in \eqref{equ:DTTdef} is the natural projection, while the last map is a canonically defined $\Ger_\infty$-morphism.
For more details we refer the reader to \cite{DTT1}, in particular to eqn 2.11 therein for a more detailed definition of $\Xi$.

We will take the zig-zag \eqref{equ:DTTdef} as the definition of the Tamarkin formality morphism induced by a quasi-isomorphism $\Ger_\infty\to\Br$.

\subsection{The \texorpdfstring{$\Lie_\infty$}{L-infinity}-part of the map \texorpdfstring{$\Br_\infty \to \Graphs$}{Br-infinity to Graphs}}
One has the natural inclusion of operads $\Lie_\infty^{(1)}\to \Ger_\infty$, where $\Lie_\infty^{(1)}=(\Lie\{1\})_\infty$ is the minimal cofibrant resolution of the degree shifted Lie operad.
Also, from the map $\Lie\{1\}\to \Br$ one obtains the embedding $\Omega(B(\Lie\{1\})) \to \Br_\infty$. Precomposing with $\Lie_\infty^{(1)}\to \Omega(B(\Lie\{1\}))$ one obtains a natural map $\Lie_\infty^{(1)}\to \Br_\infty$, and by postcomposing with the Kontsevich-Soibelman map $\Lie_\infty \to C_\bullet(\FM_2)$.

\begin{lemma}
 The maps $\Lie_\infty^{(1)}\to \Br_\infty$ and $\Lie_\infty^{(1)} \to C_\bullet(\FM_2)$ are up to homotopy uniquely determined by the image of the binary generator.
\end{lemma}
\begin{proof}
 It follows from degree reasons; there are no classes of appropriate degrees in $H(\Br_\infty)\cong H(\FM_2)\cong \Ger$.
\end{proof}
The following statement has already been verified in section \ref{sec:proofsof4and5}, but let us re-state it here.
\begin{cor}
The $n$-ary generator is mapped under $\Lie_\infty^{(1)} \to C_\bullet(\FM_2)$ to the fundamental chain of $\FM_2(n)$.
\end{cor}
\begin{proof}
 It is known that the map sending the $n$-ary generator to the fundamental chain of $\FM_2(n)$ is a map of operads. We also know that our map sends the binary generator to the fundamental chain of $\FM_2(2)$. Hence both maps must agree by the lemma, possibly up to homotopy. But there are no degenerate semi-algebraic chains and hence no elements in $C_\bullet(\FM_2)$ of degrees exceeding the dimension of $\FM_2(n)$. Hence the homotopy must necessarily be zero, and hence the generators are indeed mapped to the fundamental chains as claimed.
\end{proof}

A similar argument shows that the $n$-ary part of $\hLie_\infty\subset \hBr_\infty$\footnote{Here we denote the operadic bimodule governing $\Lie_\infty^{(1)}$ morphisms by $\hLie_\infty$.} is mapped to the fundamental chain of $D_K(n)$. From the previous statements it follows that the $\Lie_\infty^{(1)}$ part of the $\Br_\infty$-structure on $\Tpoly$ is the standard one, and that the $\Lie_\infty^{(1)}$ part of the $\Br_\infty$ formality morphism constructed above agrees with Kontsevich's $\Lie_\infty^{(1)}$ morphism $\Tpoly\to \Dpoly$ \cite{K1}.

\subsection{Proof of Theorem \ref{thm:kontsextends}}
\label{sec:extendsproof}
First construct a map $\Ger_\infty \to \Br_\infty$ by lifting up to homotopy:
\[
 \begin{tikzpicture}
  \matrix(m)[diagram]{\Lie_\infty^{(1)} &  & \Br_\infty \\
\Ger_\infty & \Ger & \Graphs \\};
\draw[->]
 (m-1-1) edge (m-1-3)
 (m-1-1) edge (m-2-1)
 (m-1-3) edge (m-2-3)
 (m-2-1) edge (m-2-2) edge[dashed] (m-1-3)
 (m-2-2) edge (m-2-3);
 \end{tikzpicture}
\]

This is possible since $\Br_\infty\to \Graphs$ is a quasi-isomorphism and $\Lie_\infty^{(1)} \to \Ger_\infty$ is a cofibration. 
Note also that the homotopy may be chosen to vanish on $\Lie_\infty^{(1)}\subset \Ger_\infty$.
Using the map $\Ger_\infty \to \Br_\infty$ one can pull back our $\Br_\infty$ morphism to a $\Ger_\infty$ morphism
\[
 \Tpoly' \to \Dpoly.
\]
Here the $\Ger_\infty$ structure on the left is not the standard one. We use the prime to distinguish it from the space $\Tpoly$ with the standard Gerstenhaber structure. However, the $\Lie_\infty^{(1)}$-structure on $\Tpoly'$ is the same as that on $\Tpoly$.
The homotopy in the above diagram, composed with the representation $\Graphs\to \End(\Tpoly)$, gives us (-by integration-) a $\Ger_\infty$ map 
\[
\phi : \Tpoly\to \Tpoly'.
\]
The $\Lie_\infty$-part of this map is the identity. Composing $\phi$ with the above formality morphism, one obtains a $\Ger_\infty$ formality morphism
\[
 \mU : \Tpoly \to \Dpoly.
\]
Here the Gerstenhaber structure on the left is the usual one, and the $\Ger_\infty$ structure on the right (solution to Deligne's conjecture) comes from pulling back the natural $\Br_\infty$ structure via the map $\Ger_\infty \to \Br_\infty$ constructed above. 
Since the $\Lie_\infty$-part of $\phi$ was trivial, the $\Lie_\infty$-part of $\mU$ is exactly M. Kontsevich's $\Lie_\infty$-morphism. This proves Theorem \ref{thm:kontsextends}.
Note also that by constructions all the components of the $\Ger_\infty$ morphism we constructed are expressible through graphical formulas, i.~e., using only operations from the operad $\SGra$.

\subsection{Relation to D. Tamarkin's morphism and proof of Theorem \ref{thm:kontstamequiv} }
Given Theorem \ref{thm:kontsextends}, the proof of Theorem \ref{thm:kontstamequiv} is more or less a standard argument using the rigidity of the Gerstenhaber algebra $\Tpoly^{alg}$. The statements of this section are ``well known'', but is hard to cite any reference.

\begin{lemma}\label{ref:lemhomotopy}
 Fix a $\Ger_\infty$ structure on $\Dpoly^{alg}$ obtained by pull-back along a quasi-isomorphism $\Ger_\infty\to\Br$ as above. Let $\phi$ be any $\Ger_\infty$ quasi-isomorphism $\phi: \Tpoly^{alg}\to \Dpoly^{alg}$ whose components are given by graphical formulas (i.~e., $\SGra$ operations). Then $\phi$ is homotopic to the $\Ger_\infty$ morphism defined by the zig-zag \eqref{equ:DTTdef}, provided both morphisms induce the same map in cohomology.
\end{lemma}
\begin{proof}
We give a short proof which is a simple adaptation of the proof of \cite[Theorem 3]{DTT1}. 
By definition, a $\Ger_\infty$ morphism $\Tpoly^{alg}\to \Dpoly^{alg}$ is the same as a map of dg $\Ger^\vee$ coalgebras
\[
 B_{\Ger^\vee}(\Tpoly^{alg}) \to \FF_{\Ger^\vee}(\Dpoly^{alg})
\]
where the right hand side is the quasi-free $\Ger^\vee$ coalgebra generated by $\Dpoly^{alg}$, with the differential defined by the $\Ger_\infty$ structure on $\Dpoly^{alg}$.
Consider the following diagram of $\Ger^\vee$ coalgebras
\begin{equation}\label{equ:cogercomm}
 \begin{tikzpicture}
  \matrix(m)[diagram]{ B_{\Ger^{\vee}}(\Tpoly^{alg}) & \Xi & \FF_{\Ger^\vee}(\Dpoly^{alg}) \\};
\draw[->]
 (m-1-1) edge[bend left] node[auto]{$\scriptstyle \phi$} (m-1-3)
  (m-1-2) edge node[auto, swap]{$\scriptstyle \iota$} (m-1-1)
 (m-1-2) edge node[auto]{$\scriptstyle \sigma$} (m-1-3);
 \end{tikzpicture}
\end{equation}
where $\Xi$ is as in \eqref{equ:DTTdef}. We claim that this diagram commutes. Indeed, since $\FF_{\Ger^\vee}(\Dpoly^{alg})$ is cofree as a graded $\Ger^\vee$ coalgebra it suffices to check that the two compositions with the projection to cogenerators
\[
\Xi
\mathrel{\mathop{\rightrightarrows}^{\phi\circ \iota}_{\sigma}} 
\FF_{\Ger^\vee}(\Dpoly^{alg}) \to 
\Dpoly
\]
agree. By degree reasons one can see that the only elements of $\Xi$ that can possibly be mapped to nonzero values are:
\begin{enumerate}
 \item Expressions $u$, for $u$ a function or vector field, which are mapped to $u$ since $\phi$ induce the identity map in cohomology by assumption.
 \item Expressions $u\wedge v$ for $u,v$ vector fields. By antisymmetry in $u,v$, and since the operations are defined by graphical formulas, this basis element can only be mapped to $\lambda \dv \co{u}{v}$, for some constant $\lambda$. Here $\dv$ is the divergence operator, which in local coordinates $\{x_i\}$ maps a vector field $w=w^i \frac{\p}{\p x_i}$ to $\frac{\p}{\p x_i}w^i$.
 \item Expressions $\underline{fu}$ for $f$ a function and $u$ a vector field. These basis elements can be mapped to $\mu f \dv u + \nu u\cdot f$.
\end{enumerate}
 If we denote by $d$ the differential in $\Xi$ then $d(u\wedge v)=\pm [u,v]{\mapsto} \pm [u,v]$, hence we must have $\lambda=0$ in order for the map $\phi$ to be compatible with the differentials. Similarly, $d(\underline{fu})=fu {\mapsto} fu$ so we must have $\mu=\nu=0$ in order for $\phi$ to be compatible with the differentials.
Hence \eqref{equ:cogercomm} commutes.

By applying the cobar construction we hence get a commutative diagram of Gerstenhaber algebras
\[
 \begin{tikzpicture}
  \matrix(m)[diagram]{ \Omega_{\Ger}(B_{\Ger^{\vee}}(\Tpoly^{alg})) & \Omega_{\Ger}(\Xi) & \Omega_{\Ger}(\FF_{\Ger^\vee}(\Dpoly^{alg})) \\};
\draw[->]
 (m-1-1) edge[bend left] node[auto]{$\scriptstyle \Omega_{\Ger^\vee}(\phi)$} (m-1-3)
  (m-1-2) edge (m-1-1)
 (m-1-2) edge (m-1-3);
 \end{tikzpicture}.
\]
Here the lower two morphisms are as in \eqref{equ:DTTdef}, and hence one can conclude that the morphism $\phi$ and \eqref{equ:DTTdef} are indeed homotopic.
\end{proof}

\begin{proof}[Proof of Theorem \ref{thm:kontstamequiv}]
 Fix for now the morphism $\Ger_\infty \to \Br$ constructed in section \ref{sec:extendsproof}. Call it $K: \Ger_\infty\to \Br$, and call the $\Ger_\infty$-formality morphism obtained in that section again $\mU:\Tpoly^{alg}\to \Dpoly^{alg}$.

Next consider the morphism $T: \Ger_\infty \to \Br$ defined as in section \ref{sec:tamreview} by using D. Tamarkin's formality morphism of the little disks operad with the Alekseev-Torossian Drinfeld associator.
In fact, we can assume that the diagram
\[
\begin{tikzpicture}
  \matrix(m)[diagram]{  & \Lie^{(1)}_\infty &  \\ \Ger_\infty & & \Br \\};
\draw[->]
 (m-1-2) edge (m-2-1)
 (m-1-2) edge (m-2-3)
 (m-2-1) edge (m-2-3);
\end{tikzpicture}
\]
commutes. Otherwise, we change the map $\Ger_\infty \to \Br$ to a homotopic morphism so that the diagram commutes.
We then obtain a $\Ger_\infty$-formality morphism $\mU_T: \Tpoly^{alg}\to \Dpoly^{alg}$ implicitly defined by \eqref{equ:DTTdef}. 

Note that the $\Ger_\infty$-structures on $\Dpoly^{alg}$ as occurring in $\mU$ and $\mU_T$ are different, and we will write $\mU:\Tpoly\to \Dpoly^K$, $\mU_T: \Tpoly\to \Dpoly^T$ to emphasize this distinction. We know by the results of \cite{pavol} that the two morphisms $K$ and $T$ are homotopic. Since for both morphisms the above diagram commutes, one can pick the homotopy in such a way that its $\Lie_\infty^{(1)}$-part vanishes. 
It follows that there is a $\Ger_\infty$ quasi-isomorphism 
\[
 \Phi: \Dpoly^{K} \to \Dpoly^{T}
\]
between two copies of $\Dpoly^{alg}$ with the two $\Ger_\infty$-structure coming from $K$ and $T$, such that the $\Lie_\infty^{(1)}$-part of $\Phi$ vanishes.\footnote{The nontrivial part of this statement is that the $\Lie_\infty^{(1)}$-part vanishes.} Furthermore, all components of $\Phi$ may be expressed through $\Br$-operations.
But now the map  $\Phi\circ \mU_K: \Tpoly^{alg}\to \Dpoly^{T}$ satisfies the assumptions of Lemma \ref{ref:lemhomotopy}, and hence the maps $\Phi\circ \mU_K$ and $\mU_T$ are homotopic. In particular this implies that the $\Lie_\infty^{(1)}$-parts are homotopic.  But the $\Lie_\infty^{(1)}$-part of $\Phi\circ \mU_K$ is the same as that of $\mU_K$, and hence the Theorem is shown.

\end{proof}

\makeatletter{}\appendix
\section{Actions of \texorpdfstring{$\Br$}{Br} and \texorpdfstring{$\KS$}{KS}}
\label{app:KSactproof}
This section contains the proofs of Propositions \ref{prop:KSactonDK} and \ref{prop:KSactonDS}.

\subsection{The construction of the action}
\label{sec:bimodgeneral}
Let us start with some Swiss Cheese type operad $\op P$, see definition \ref{def:SCtype}.
Our goal is to construct a $\Br$-$\op P^1$ operadic bimodule structure on the $\mathbb{S}$-module $\op M(\cdot)=\op P^2(\cdot, 0)$. Here, as before, $\op P^\alpha$ is the space of operations with output in color $\alpha\in \{1,2\}$. We proceed in the following steps:
\begin{enumerate}
\item By example \ref{ex:SCtoPTmod} there is a natural $\PT$-$\op P^1$ operadic bimodule structure on 
\[
\op M'(\cdot) = \prod_n \op P^2(\cdot, n)[-n].
\]

\item Suppose we are given a Maurer-Cartan element $\nu \in \op M'(0)$. Then we can twist the left $\PT$-module structure to a $\Tw\PT$- and hence also to a $\Br\subset \Tw\PT$-structure, following Appendix \ref{sec:optwists}. We call the resulting $\Br$-$\op P^1$ operadic bimodule $\op M''$. As an $\mathbb{S}$-module in graded vector spaces it is the same as $\op M'$, but the differential contains an additional term introduced by the twisting.

\item Suppose further that there is a map 
\[
\op P^2(\cdot, 0) \to \op M''(\cdot)
\]
such that (i) the image is an operadic $\Br$-$\op P^1$ sub-module, and (ii) the map is a right inverse to the natural projection $\op M''(\cdot)\to \op P^2(\cdot, 0)$. In this case we endow $\op M=\op P^2(\cdot, 0)$ with the induced operadic $\Br$-$\op P^1$ bimodule structure. Of course, there is a natural inclusion of $\Br$-$\op P^1$ bimodules 
\[
\op M \to \op M''\, .
\]
\end{enumerate}

\subsection{The braces action}
We want to construct an operadic $\Br$-$C(\FM_2)$ bimodule structure on $C(D_K)$, the chains on the space of configurations of points in the upper halfplane. We do this by following the program of the previous subsection. Here the role of the Swiss Cheese type operad $\op P$ is played by the operad of (semi algebraic) chains $C(\SC)$ on the Swiss Cheese operad $\SC$ from section \ref{sec:DKdef}. In particular note that $D_K=\SC^2(\cdot,0)$, $\FM_2=\SC^1$.
Two pieces of data have to be provided, according to the previous section. First, we need a Maurer-Cartan element 
\[
 \nu \in \op M'(0) := \prod_n C(\SC^{2}(0,n))[-n].
\]
\begin{lemma}
\label{lem:nuismumc}
 The element 
\[
 \nu := \sum_{n\geq 2}  \mathit{Fund}(\SC^{2}(0,n)),
\]
where $\mathit{Fund}(\cdot)$ denotes the fundamental chain, is a $\mu$ Maurer-Cartan element, for $\mu$ the element 
\[
 \begin{tikzpicture}[scale=.75,
vert/.style={draw,outer sep=0,inner sep=0,minimum size=5,shape=circle,fill},
helper/.style={outer sep=0,inner sep=0,minimum size=5,shape=coordinate},
default_edge/.style={draw,->},
ext/.style={draw,outer sep=0,inner sep=2,minimum size=5,shape=circle},
every loop/.style={}]

\node (v0) at (5,6.5) [ext] {2};
\node (v1) at (5,7.5) [ext] {1};
\node (v3) at (5,8.4) [helper] {};
\node (v10) at (6.6,7.5) [ext] {2};
\node (v9) at (6.6,6.5) [ext] {1};
\node (v12) at (6.6,8.4) [helper] {};
\node (v14) at (5.6,7.1) [helper,label=0:{$+$}] {};

\draw (v1) -- +(0,.6);
\draw (v10) -- +(0,.6);

\draw[default_edge] (v0) to (v1);
\draw[default_edge] (v9) to (v10);
\end{tikzpicture}
\]
in the total space of $\PT$.
\end{lemma}
To be more concrete, the element $\mu$ defines a map from the (degree shifted) Lie operad into $\PT$. But $\op M'(0)$ is a $\PT$ algebra and hence in particular a Lie algebra. Hence it makes sense to speak about Maurer-Cartan elements in this Lie algebra (the $\mu$-Maurer-Cartan elements).
\begin{proof}
 The spaces $\SC^{2}(0,n)$ form a non symmetric sub-operad which is isomorphic to the nonsymmetric version of $\FM_1$ (Stasheff's associahedra). It has a natural stratification and the chains associated to strata span a suboperad of $C(\SC^{2}(0,n))$ which is isomorphic to the non-symmetric $A_\infty$-operad. The fundamental chains $\mathit{Fund}(\SC^{2}(0,n))$ correspond to the generators $a_n$ ($n=2,3,\dots$). The $\mu$-Maurer Cartan equation for $\nu$ translates into the usual relations expressing the differential of $a_n$ in terms of the $a_j$, $j<n$.
\end{proof}

With this Maurer-Cartan element, we obtain a  $\Br$-$C(\FM_2)$ bimodule structure on
\[
 \op M''(\cdot) := \prod_n C(\SC^{2}(\cdot,n))\, .
\]
The differential has the form
\[
 d = \partial + d_\nu
\]
where $\partial$ is the boundary operator, i.e., the differential on $\op M'$, and $d_\nu$ is the part contributed by the twisting. It in turn has two terms. For a chain $c$
\begin{equation}
\label{equ:dnudef}
 d_\nu c = L_\mu(\nu, c) + L_\mu(c, \nu)
\end{equation}
where $L_\mu$ is the operadic left action of the $\PT$ tree with two vertices.

The second datum we need to provide according to the previous subsection is a map 
\[
F: C(D_K)\to \op M''
\]
satisfying the aforementioned properties.
To construct $F$ consider the forgetful maps 
\[
 \pi_{m,n} \colon \SC^{2}(m,n) \to \SC^{2}(m,0).
\]
On semi-algebraic chains, there is an operator 
\[
 \pi_{m,n}^{-1} \colon C(\SC^{2}(m,0)) \to C(\SC^{2}(m,n))
\]
taking the semi-algebraic fibers, see Appendix \ref{sec:invimgchains} or \cite{HLTV} for details. 
It has the property that for a chain $c\in C(\SC^{2}(m,0))$
\[
 \partial \pi_n^{-1} (c) = \pi_n^{-1} (\partial c) + (-1)^{|c|}  ((\pi_{m,n})^\p)^{-1} (c)
\]
where $(\pi_{m,n})^\p$ is defined as in Appendix \ref{sec:invimgchains}. It takes the fiberwise boundary.
Let us define
\[
 F(c) := \sum_n \pi_{m,n}^{-1} c.
\]
for $c\in C(D_K(m))$. 
It is clear that this map is a right inverse to the projection 
\[
 \pi : \op M'' \to C(\SC^{2}(\cdot,0))=C(D_K)\, .
\]
Note that $\pi_{m,0}^{-1}=\mathit{id}$.
In particular, $F$ is an embedding (of $\mathbb{S}$-modules in graded vector spaces). What remains to be checked is that the image of $F$ 
is an operadic $\Br$-$\FM_2$ sub-bimodule.

\begin{lemma}
\label{lem:Fcompatible}
 The map $F: C(D_K)\to \op M''$ is an embedding of right $C(\FM_2)$-modules.
\end{lemma}
\begin{proof}
 Is is clear that the right action is preserved. The nontrivial part of this statement is that $F$ is compatible with the differentials.
We have, for a chain $c\in C_k(D_K(m))$
\[
 F(\partial c) = \sum_n \pi_{m,n}^{-1}( \partial c) = \sum_n \pi_{m,n}^{-1}( \partial c) =  \sum_n \partial \pi_{m,n}^{-1}(  c) - (-1)^k ((\pi_{m,n})^\p)^{-1}(  c).
\]
where $(\pi_{m,n})^\p$ is the projection of the bundle of fiberwise boundaries, as recalled in section \ref{sec:invimgchains} (or \cite{HLTV}, Propositions 5.17 and 8.2).
Let us look at the fiberwise boundary. Each fiber is, as semi algebraic manifold, a space of configurations of $m$ points on a the real line, possibly with some other (fixed) points infinitely close to the real line. For $m=0$ the statement of the Lemma follows from Lemma \ref{lem:nuismumc}. Assume $m\geq 1$. Then there are two types of boundary strata: (i) Some set of points on the real axis come infinitely close to each other and (ii) some points move to infinity, while zero or more points stay ``at finite distance'' to the $m$ points in the upper halfplane.
These two types of strata contribute the two terms in \eqref{equ:dnudef}.
\end{proof}

The harder part is to check that the image of $F$ is closed under the left $\Br$-action as well.

\begin{prop}
\label{prop:Fimclosed}
 The image of $F$, i.~e., $F(C(D_K))$ is an operadic $\Br-C(\FM_2)$ sub-bimodule.
\end{prop}
\begin{proof}
It is clear by the previous lemma that $F(C(D_K))$ is closed under the differential and the right $C(\FM_2)$ action. What has to be checked is that it is closed under the left braces action. Of course it is sufficient to check this on the generators $T_n$, $T_n'$ of $\Br$, see Lemma \ref{lem:KSgenrel}. We consider here only $T_n$ and leave the simpler proof for $T_n'$ to the reader. So let $c_0\in C(D_K(m_0)),c_1\in C(D_K(m_1)), \dots, c_n\in C(D_K(m_n))$ be chains. It suffices to show that 
\[
 F(\pi ( L_{T_n}(F(c_0),F(c_1),\dots, F(c_n)) ) ) = L_{T_n}(F(c_0),F(c_1),\dots, F(c_n)).
\]
Let us abbreviate $\gamma := \pi (L_{T_n}(F(c_0),F(c_1),\dots, F(c_n)))$. Using the notation of \cite{HLTV}, the left and right hand sides have the form $\gamma \ltimes \Phi$ and $\gamma \ltimes \Phi'$ respectively, for two strongly continuous (families of) chains $\Phi,\Phi'\in \prod_k C^{str}(D_{Ke}(\sum_j m_j, k)\to D_{K}(\sum_j m_j))$. (See \cite{HLTV}, Definition 5.13 for the notation.) Here $\Phi$ is the same as in Appendix \ref{sec:invimgchains}. We have to show that $\Phi=\Phi'$. Clearly it suffices to check this separately for every fiber of the semi-algebraic bundles $D_{Ke}(\sum_j m_j, k)\to D_{K}(\sum_j m_j)$ and for each $k$. The fiber is a space of configurations of points on the real axis, possibly with some fixed points in the upper halfplane, that may be infinitely close to the real axis. 
Let us denote the set of type I vertices involved in configurations in $c_j$ by $S_j$. So $|S_j|=m_j$.
Since by definition $F(c_j)=\sum_{k_j}  \pi_{m_j,k_j}^{-1}( c_j)$, $\Phi'$ can be decomposed as a sum of terms
\[
 \Phi' = \sum_{k_0\geq m} \Phi_{k_0,\cdots, k_m}'
\]
in such a manner that $\Phi_{k_0,\cdots, k_m}$ is the sum of strata of configurations where $k_1$ type II vertices are close to vertices in $S_1$, $k_2$ to vertices in $S_2$ etc.
In a similar manner $\Phi$ can also be decomposed into parts where $k_1$ type II vertices are close to vertices in $S_1$ etc. Each term is again the sum over all such strata and hence $\Phi=\Phi'$ and the Proposition follows.
\end{proof}

To conclude let us re-state the findings of this section.
\begin{prop}
\label{prop:cdk_action}
 There is an operadic $\Br-C(\FM_2)$ bimodule structure on $C(D_K)$ extending the usual right $C(\FM_2)$ module structure.
\end{prop}

A schematic picture of the $\Br$-action can be found in Figure \ref{fig:BractiononDK}. 

\begin{figure}
 \centering
\makeatletter{} \usetikzlibrary{arrows}
\begin{tikzpicture}[
scale=.5,
int/.style={circle, draw, fill, minimum size=5pt, inner sep=0},
ext/.style={circle, draw, fill=white, minimum size=5pt, inner sep=1pt},
helper/.style={coordinate},point/.style={circle, draw, fill, inner sep =1pt},
de/.style={-triangle 60},
point/.style={circle, draw, fill, minimum size=3pt, inner sep=0pt},
xst/.style={cross out, draw, minimum size=5 },
boxed/.style={ draw, inner sep=0.5 },
]

\draw(-6,-1) -- (8,-1);
\node [ext] (v2) at (-13,2.5) {1};
\node [ext] (v1) at (-14,1) {2};
\node [ext] (v3) at (-13,1) {3};
\node [ext] (v4) at (-12,1) {4};
\draw (v1) edge (v2);
\draw (v3) edge (v2);
\draw (v4) edge (v2);
\node at (-9,1.7) {$(c_1,c_2,c_3,c_4) \quad =$};
\draw (v2) -- +(0,.8);

\draw[dashed](-2,-1) arc (0:180:1);
\node[label=90:$c_2$] (c) at (-3,-1) {};
\draw[triangle 60-triangle 60] (c) +(55:1.5)-- +(125:1.5);

\begin{scope}[xshift=4cm]
\draw[dashed](-2,-1) arc (0:180:1);
\node[label=90:$c_3$] (c) at (-3,-1) {};
\draw[triangle 60-triangle 60] (c) +(55:1.5)-- +(125:1.5);
\end{scope}

\begin{scope}[xshift=8cm]
\draw[dashed](-2,-1) arc (0:180:1);
\node[label=90:$c_4$] (cc) at (-3,-1) {};
\draw[triangle 60-triangle 60] (cc) +(55:1.5)-- +(125:1.5);
\end{scope}

\node at (1,4) {\huge $c_1$};
\node [int] at (-2,4) {};
\node [int] at (3,5) {};
\node [int] at (5,3) {};
\node [int] at (1,3) {};

\begin{scope}[yshift=-6cm]
\draw(-6,-1) -- (10,-1);
\node [int] (v2) at (-13,2.5) {};
\node [ext] (v5) at (-14,1) {1};
\node [ext] (v1) at (-13.33,1) {2};
\node [ext] (v3) at (-12.66,1) {3};
\node [ext] (v4) at (-12,1) {4};

\draw (v1) edge (v2);
\draw (v3) edge (v2);
\draw (v4) edge (v2);
\draw (v5) edge (v2);
\node at (-9,1.7) {$(c_1,c_2,c_3,c_4) \quad =$};
\draw (v2) -- +(0,.6);
\begin{scope}[xshift=-1cm]
\draw[dashed](-2,-1) arc (0:180:1);
\node[label=90:$c_1$] (c) at (-3,-1) {};
\end{scope}

\begin{scope}[xshift=3cm]
\draw[dashed](-2,-1) arc (0:180:1);
\node[label=90:$c_2$] (c) at (-3,-1) {};
\draw[triangle 60-triangle 60] (c) +(55:1.5)-- +(125:1.5);
\end{scope}

\begin{scope}[xshift=7cm]
\draw[dashed](-2,-1) arc (0:180:1);
\node[label=90:$c_3$] (cc) at (-3,-1) {};
\draw[triangle 60-triangle 60] (cc) +(55:1.5)-- +(125:1.5);
\end{scope}

\begin{scope}[xshift=11cm]
\draw[dashed](-2,-1) arc (0:180:1);
\node[label=90:$c_4$] (cc) at (-3,-1) {};
\end{scope}
\end{scope}

\end{tikzpicture} 
\caption{\label{fig:BractiononDK} An illustration of the left operadic action of $\Br$ on $C(D_K)$. Here the elements of $\Br$ shown on the left act on chains $c_1,\dots, c_4\in C(D_K)$. The dashed semicircles shall indicate that, e.g., configurations occuring in $c_2$ are placed infinitely close to the real axis. There are some issues due to the compactification that are swept under the rug by this picture. For example in the upper picture $c_1$ itself could contain configurations with points moving infinitely close to the real axis. In this case $c_2,c_3,c_4$ would move even closer to the real axis under all points involved in configuration in $c_1$. For a more precise definition of the action, see the text.}
\end{figure}
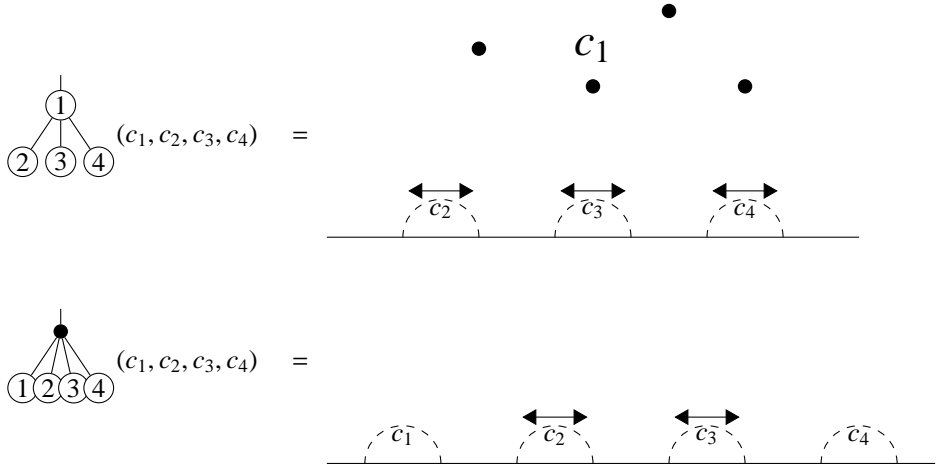

\begin{figure}
 \centering
\makeatletter{}\usetikzlibrary{arrows}
\[
\begin{tikzpicture}[
scale=.6,
int/.style={circle, draw, fill, minimum size=5pt, inner sep=0},
ext/.style={circle, draw, fill=white, minimum size=5pt, inner sep=1pt},
helper/.style={coordinate},point/.style={circle, draw, fill, inner sep =1pt},
de/.style={-triangle 60},
point/.style={circle, draw, fill, minimum size=3pt, inner sep=0pt},
]

\begin{scope}[xshift=22cm]
\draw (-5,0)--(5,0);
\node[point] (e3) at (-1.8,0) {};
\node[point,label=-90:{$k_1$}] (e2) at (-2.2,0) {};
\node[point,label=-90:{$k_{n}$}] (e4) at (2,0) {};
\node[point] (c) at (2.4,0) {};
\node[point] (e5) at (2.8,0) {};
\draw[dashed] (3.3,0) arc (0:180:.9);
\node at (-2,.4) {$S_1$};
\node at (2.4,.5) {$S_{n}$};
\node at (0, .5) {\large$\cdots$};
\node[point] at (-4.5,0) {};
\node[point] at (-4,0) {};
\node[point] at (-3.7,0) {};
\node[point] at (-.9,0) {};
\node[point] at (0,0) {};
\node[point] at (4.5,0) {};
\node[point] at (4.3,0) {};

\draw[dashed] (-1.3,0) arc (0:180:.7);

\node[label=0:{$S_0$}] at (0,3) {};
\node[point] (t1) at (1.8,3) {};
\node[point] (t2) at (-1,2.5) {};
\node[point] (t3) at (0.1,3.8) {};
\end{scope}
\end{tikzpicture}
\] 
\caption{\label{fig:xf} Illustration of a a statement in the proof of Proposition \ref{prop:Fimclosed}. We have to compute the chain $\pi_{\sum m_j, k}^{-1}(\gamma)$. There will be various contributing strata, which can be organized regarding how many type II vertices are infinitely close to points in the upper halfplane in $S_1$, $S_2$ etc. }
\end{figure}
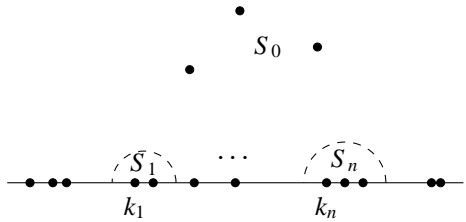

\subsection{An extension}
\label{sec:bimodgeneralex}
Next we want to extend the construction of Appendix \ref{sec:bimodgeneral} to moperadic bimodules. 
We start with an Extended Swiss Cheese type operad $\op Q$, see definition \ref{def:ESCtype}.
By Appendix \ref{sec:bimodgeneral} we can build (-given some extra data-) an operadic $\Br$-$\op Q^1$ bimodule structure on the spaces $\op M(\cdot) = \op Q^2(\cdot, 0,0)$.
Our goal is to construct a moperadic $\Br$-$\op Q^1$-$\KS_1$-$\op Q^3(\cdot, 0,1)$-$\op M$ bimodule structure on the spaces $\op N=\op Q^3(\cdot, 1,0)$. We proceed as follows. In the following let $\op M'$, $\op M''$ and the Maurer-Cartan element $\nu$ be as in Appendix \ref{sec:bimodgeneral}.

\begin{enumerate}
\item By example \ref{ex:ESCtoPT1mod} and the cyclic module structure on $\op Q^3(\cdot, \cdot,0)$, there is a moperadic $\PT$-$\op Q^1$-$\PT_1$-$\op Q^3(\cdot, 0,1)$-$\op M'$ bimodule structure on the spaces
\[
\op N'(\cdot) = \prod_n \op Q^3(\cdot, n,0)[1-n].
\]

\item We suppose that the $\PT_1$ action can be extended to a $\PT_1^{\mathbb{1}}$-action. We obtain a moperadic $\PT$-$\op Q^1$-$\PT_1^{\mathbb{1}}$-$\op Q^3(\cdot, 0,1)$-$\op M'$ bimodule $\op N''$.

\item Using the Maurer-Cartan element $\nu$ from Appendix \ref{sec:bimodgeneral} we can twist the actions to create a $\Tw\PT$-$\op Q^1$-$\Tw\PT_1^{\mathbb{1}}$-$\op Q^3(\cdot, 0,1)$-$\op M''$ bimodule $\op N'''$.

\item We assume that there is a map 
\[
\op N= \op Q^3(\cdot, 1,0) \to \op N'''
\]
(of $\mathbb{S}$-modules in differential graded vector spaces for now) that is right inverse to the natural projection $\op N'''\to \op Q^3(\cdot, 1,0)$. 

\item We next suppose that the moperadic $\Tw\PT$-$\op Q^1$-$\Tw\PT_1^{\mathbb{1}}$-$\op Q^3(\cdot, 0,1)$-$\op M''$ bimodule structure on $\op N'''$ descends to a $\Br$-$\op Q^1$-$\KS_1$-$\op Q^3(\cdot, 0,1)$-$\op M$ bimodule structure on $\op N$.
This means that the following things have to be checked:
\begin{enumerate}
\item The relations of section \ref{sec:KSdef} are respected, so as to obtain an action of the subquotient $\KS_1\subset \Tw\PT_1^{\mathbb{1}}$.
\item The subspaces $\op N$ have to be closed under the actions of $\op Q^1$, $\op Q^3(\cdot, 0,1)$ and under the combined action of $\KS_1$ and $\op M$.
\end{enumerate}

\item We endow $\op N=\op Q^3(\cdot, 1,0)$ with the moperadic sub-bimodule structure. Of course, there is an inclusion of 
moperadic $\Br$-$\op Q^1$-$\KS_1$-$\op Q^3(\cdot, 0,1)$-$\op M$ bimodules
\[
\op N \to \op N'''.
\]
\end{enumerate}

\subsection{The \texorpdfstring{$\KS_1$}{KS1} action}
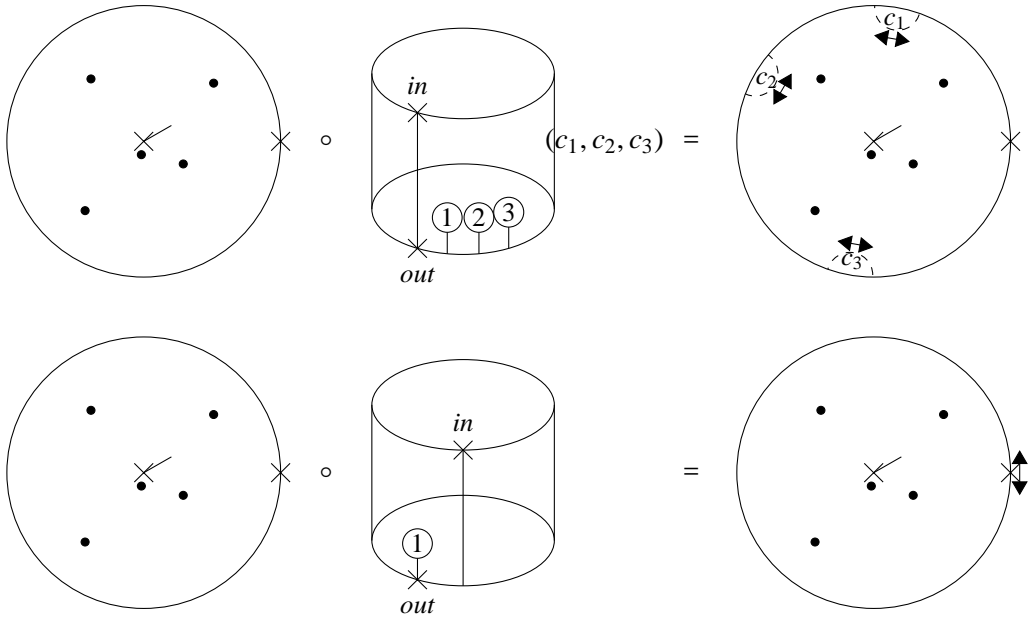
\begin{figure}
 \centering
\makeatletter{}\usetikzlibrary{matrix}
\usetikzlibrary{arrows}
\usetikzlibrary{shapes}
\usetikzlibrary{through}
\usetikzlibrary{calc,3d}
\usetikzlibrary{decorations,decorations.pathmorphing}
\[
\begin{tikzpicture}[
scale=.6,
int/.style={circle, draw, fill, minimum size=5pt, inner sep=0},
ext/.style={circle, draw, fill=white, minimum size=5pt, inner sep=1pt},
helper/.style={coordinate},point/.style={circle, draw, fill, inner sep =1pt},
de/.style={-triangle 60},
point/.style={circle, draw, fill, minimum size=3pt, inner sep=0pt},
xst/.style={cross out, draw, minimum size=5 },
]
\begin{scope}[]
\draw (0,0) node[xst] (c) {} circle (3) -- (30:.7);
\node[xst] at (3,0) {};
\node[point] at (40:2) {};
\node[point] at (-30:1) {};
\node[point] at (-100:.3) {};
\node[point] at (130:1.8) {};
\node[point] at (-130:2) {};

\end{scope}

\begin{scope}[xshift=7cm, yshift=-1.5cm]
\node at (-3,1.5) {$\circ$};
 \draw (0,0) ellipse (2cm and 1cm);
 \draw (0,3) ellipse (2cm and 1cm);
 \draw (-2,0)--(-2,3) (2,0)--(2,3);
\node [xst, label=-90:{$out$}] (out) at ($(0,0)+(-120:2 and 1)$) {};
\node [xst, label=90:{$in$}] (in) at ($(0,3)+(-120:2 and 1)$) {};
\draw ($(0,0)+(-100:2 and 1)$) -- +(0,.8) node[ext] {1};
\draw ($(0,0)+(-80:2 and 1)$) -- +(0,.8) node[ext] {2};
\draw ($(0,0)+(-60:2 and 1)$) -- +(0,.8) node[ext] {3};

\node at (3.1,1.5) {\large ${ (c_1, c_2, c_3)}$};
  \draw (in.base)--(out.base);
\end{scope}
\begin{scope}[xshift=16cm]
\node at (-4,0) {$=$};

\draw (0,0) node[xst] (c) {} circle (3) -- (30:.7);
\node[xst] at (3,0) {};
\node[point] at (40:2) {};
\node[point] at (-30:1) {};
\node[point] at (-100:.3) {};
\node[point] at (130:1.8) {};
\node[point] at (-130:2) {};

\clip (c) circle (3) ;
\node (c1) at (80:3) {};
\draw (c1)+(260:.3) node {$c_1$};
\node (c2) at (150:3) {};
\draw (c2)+(-30:.3) node {$c_2$};
\node (c3) at (-100:3) {};
\draw (c3)+(80:.3) node {$c_3$};
\foreach \x in {c1,c2,c3}
  \draw[dashed] (\x) circle (.5);

\draw[triangle 60-triangle 60] (c1)+(230:.8) --+(290:.8);
\draw[triangle 60-triangle 60] (c2)+(0:.8) --+(-60:.8);
\draw[triangle 60-triangle 60] (c3)+(50:.8) --+(110:.8);
\end{scope}

\end{tikzpicture}
\] 
\makeatletter{}\usetikzlibrary{matrix}
\usetikzlibrary{arrows}
\usetikzlibrary{shapes}
\usetikzlibrary{through}
\usetikzlibrary{calc,3d}
\usetikzlibrary{decorations,decorations.pathmorphing}
\[
\begin{tikzpicture}[
scale=.6,
int/.style={circle, draw, fill, minimum size=5pt, inner sep=0},
ext/.style={circle, draw, fill=white, minimum size=5pt, inner sep=1pt},
helper/.style={coordinate},point/.style={circle, draw, fill, inner sep =1pt},
de/.style={-triangle 60},
point/.style={circle, draw, fill, minimum size=3pt, inner sep=0pt},
xst/.style={cross out, draw, minimum size=5 },
]
\begin{scope}[]
\draw (0,0) node[xst] (c) {} circle (3) -- (30:.7);
\node[xst] at (3,0) {};
\node[point] at (40:2) {};
\node[point] at (-30:1) {};
\node[point] at (-100:.3) {};
\node[point] at (130:1.8) {};
\node[point] at (-130:2) {};

\end{scope}

\begin{scope}[xshift=7cm, yshift=-1.5cm]
\node at (-3,1.5) {$\circ$};
 \draw (0,0) ellipse (2cm and 1cm);
 \draw (0,3) ellipse (2cm and 1cm);
 \draw (-2,0)--(-2,3) (2,0)--(2,3);
\node [xst, label=-90:{$out$}] (out) at ($(0,0)+(-120:2 and 1)$) {};
\node [xst, label=90:{$in$}] (in) at ($(0,3)+(-90:2 and 1)$) {};
\draw (out.base) -- +(0,.8) node[ext] {$\mathrm{1}$};

  \draw (in.base)--($(0,0)+(-90:2 and 1)$);
\end{scope}
\begin{scope}[xshift=16cm]
\node at (-4,0) {$=$};

\draw (0,0) node[xst] (c) {} circle (3) -- (30:.7);
\node[xst] (in) at (3,0) {};
\node[point] at (40:2) {};
\node[point] at (-30:1) {};
\node[point] at (-100:.3) {};
\node[point] at (130:1.8) {};
\node[point] at (-130:2) {};

\draw[triangle 60-triangle 60] (in)+(.2,-.5) --+(.2,.5);
\end{scope}

\end{tikzpicture}
\] 
\caption{\label{fig:ks1ondsaction} Graphical ``definition'' of the moperadic action of $\KS_1$ on $C(D_S)$.}
\end{figure}

In this section we want to define the moperadic action of $\KS_1$ on the space of (semi algebraic) chains $C(D_S)$ on the Shoikhet configuration spaces. Before reading on, the reader should look at Figure \ref{fig:ks1ondsaction} from which the action, up to signs, should be clear.
 However, let us proceeed in a more careful way. We will follow the construction of the previous subsection. 
Consider the three colored extended Swiss cheese operad $\ESC$ from section \ref{sec:DSforgetful} and the three colored operad of semialgebraic chains $C(\ESC)$. Recall also its extension $\EESC$ from section \ref{sec:DSforgetful}, which also incorporates the forgetful maps. The role of the operad $\op Q$ from the previous subsection will be played by $C(\EESC)$. Three things need to be checked. First we need to extend the action of $\PT_1$ on $\op N'$ can be extended to a $\PT_1^{\mathbb{1}}$-action.\footnote{We use the notation of the previous subsection throughout.} To do this we follow example \ref{ex:PT11mod}. The element in $\op M'(0)$ corresponding to the unary operation $\mathbb{1}$ is given by the ``forgetful'' operation we also denoted by $\mathbb{1}$. No relations have to be checked at this point.

Next we twist by the Maurer-Cartan element $\nu$ as in the previous subsection, to obtain a moperadic $\Tw\PT$-$C(\FM_2)$-$\Tw\PT_1^{\mathbb{1}}$-$C(\FM_{2,1})$-$\op M''$ bimodule $\op N'''$. 
The differential on $\op N'''$ has the form $\p + d_\nu$, where $\p$ is the usual boundary operator and $d_\nu$ is the part contributed by the twisting.
Next let us define the embedding $G: C(D_S)\to \op N'''$ that sends $c\in C(D_S(m))$ to
\[
 G(c) := \sum_n \pi_{m,n}^{-1} c.
\]
Here $\pi_{m,n}$ is the forgetful map $\pi_{m,n}: D_{Se}(m,n)\to D_{Se}(m,1)$. Clearly this map is a right inverse to the projection $\op N'''\to C(D_S(\cdot))$. 
Similarly to Lemma \ref{lem:Fcompatible} one proves the following:
\begin{lemma}
\label{lem:Gcompatible}
 The map $G: C(D_S)\to \op N'''$ is an embedding of right $C(\FM_2)$-modules.
\end{lemma}

Proceeding along the lines of the previous subsection, the next thing we have to show is that the moperadic $\Tw\PT$-$C(\FM_2)$-$\Tw\PT_1^{\mathbb{1}}$-$C(\FM_{2,1})$-$\op M''$ bimodule structure on $\op N'''$ descends to a $\Br$-$C(\FM_2)$-$\KS_1$-$C(\FM_{2,1})$-$\op M$ bimodule structure on the image of $G$, i.e., on $\op N$.

\begin{figure}
 \centering
\makeatletter{}\usetikzlibrary{arrows}
\tikzset{ext/.style={circle, draw,inner sep=1pt},int/.style={circle,draw,fill,inner sep=1pt},nil/.style={inner sep=1pt}}
\tikzset{arr/.style={-triangle 60}}
\begin{tikzpicture}[scale=.7]

\draw (-4,1.5) node [int] (v1) {};
\node [ext] (v2) at (-4.5,1) {$\mathbb{1}$};
\node [ext] (v3) at (-3.5,1) {T};
\draw (v1) edge (v2);
\draw (v1) edge (v3);
\node at (-3,1.6) {\large $\Rightarrow$};
\node (v4) at (-4,2.5) {$(\cdots)$};
\draw (v1) edge (v4);

\begin{scope}[shift={(-1.4,0)}]
\draw(-1,0) -- (4,0);
\draw [dashed](3.5,0) arc (0:180:2);
\node at (1.5,3) {$(\cdots)$};
\draw [dashed](1.2,0) arc (0:180:0.6);
\draw [dashed](3.2,0) arc (0:180:0.8);
\node [int] at (0.6,0.2) {};
\node at (2.4,0.2) {$<T>$};
\draw[arr](0,-1.2) node {forget this} -- (0.4,-0.2);
\end{scope}

\node at (3.4,1.4) {\large $\Rightarrow$};

\begin{scope}[shift={(5.2,0)}]
\draw(-1,0) -- (4,0);
\node at (1.5,3) {$(\cdots)$};
\draw [dashed](2.3,0) arc (0:180:0.8);
\node at (1.5,0.2) {$<T>$};
\end{scope}

\begin{scope}[shift={(14.8,-0.4)}]
\node [ext] (v3) at (-4,1) {T};
\node (v4) at (-4,2.5) {$(\cdots)$};
\draw (v3) edge (v4);
\end{scope}

\node at (10,1.4) {\large $\Leftarrow$};
\end{tikzpicture} 
\caption{\label{fig:unitrelgeometry} Illustration of the unit relation in the action of $\widetilde{\PT}_1$ on chains. First from left: Part of some $\widetilde{\PT}_1$ tree $\Gamma$. $T$ stands for some subtree. The $(\cdots)$ stand for the remainder of the tree. Second: Part of configuration space. The action of $\Gamma$ is depicted. The rules for handling a unit vertex say that first it is treated as some auxiliary vertex and then this vertex is forgotten (third picture). Here $<T>$ stands for whatever configurations are produced by the subtree $T$. All these are ``at smaller scale'' and don't interfere with our picture. On the other hand $(\cdots)$ stands for whatever is done by the rest of the tree at a bigger scale, this also does not interfere. Note that the action in effect is the same as that of the tree obtained by removing the internal and the unit vertex (right). This is the unit relation. }
\end{figure}
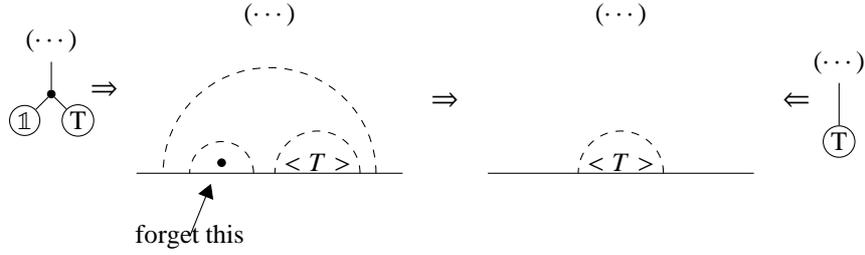

\begin{lemma}
The relations of section \ref{sec:KSdef} are respected by the joint action of $\Tw\PT_1^{\mathbb{1}}$ and $\op M$ on $\op N$.
\end{lemma}
\begin{proof}
We have to check the following 4 relations coming from those of section \ref{sec:KSdef}. 
\begin{enumerate}
\item Graphs containing a unit vertex whose parent is an internal vertex with three or more children act as zero. This is true because the internal vertices represent Maurer-Cartan elements, and the part of the Maurer-Cartan element $\nu$ with $\geq 3$ children is given by chains of degree $\geq 1$. Forgetting one tupe II point of the configuration produces a degenerate, and hence the zero chain.\footnote{ Recall from \cite{HLTV} that semi algebraic chains are by definition currents, representable in a certain way. In particular, there are no degenerate semi algebraic chains (except zero) since the associated currents are automatically zero. Hence, if the forgetful map produces a degenerate chain, it is zero.}

\item A graph containing an internal vertex with two children, one of which is the unit symbol acts in the same way as the graph without the internal vertex and unit vertex. Here the piece of the Maurer-Cartan element $\nu$ corresponding to an internal vertex with two children is (the chain of) a single configuration of two points. Forgetting one, one obtains the operadic unit in $\op Q^2(0,1,0)$. This is illustrated in Figure \ref{fig:unitrelgeometry}.

\item Graphs containing a unit vertex whose parent is an external vertex act as zero. Here it is important that we restrict the action of $\op M''$ to one of $\op M$, otherwise this statement is false.
But the part of a (sum of) chains in $\op M$ that is represented by an 
external vertex with $\geq 1$ children is a chain of degree $\geq 1$. In particular, forgetting one type II vertex yields a degenerate chain, hence zero.

\item Graphs containing a unit vertex whose parent is $\vout$ and the incident edge is not marked act as zero. Here it is important that we restricted the action to the subspace $\op N$, otherwise the statement would be false. But the part of the chain in $\op N$ represented by $\vout$ (which has at least one non-marked edge) has degree $\geq 1$. The forgetful map again produces a degenerate chain, and hence zero.

\end{enumerate}
\end{proof}

We still have to show that $\op N$ is closed under the actions of $C(\FM_{2,1})$ and the combined action of $\KS_1$ and $\op M$. The former statement is trivial. For the latter, there is an analog of Proposition \ref{prop:Fimclosed}.
\begin{prop}
\label{prop:Gimclosed}
 The image of $G$, i.e., $\op N$ is closed under the combined action of $\KS_1$ and $\op M$..
\end{prop}

We hence obtain:
\begin{prop}
\label{prop:ks1_action}
 There is a moperadic $\Br$-$C(\FM_2)$-$\KS_1$-$C(\FM_{2,1})$-$C(D_S)$ bimodule structure on $C(D_S)$ extending the usual right $C(\FM_2)$ (moperadic) module structure.
\end{prop}

\section{The inverse image on semi algebraic chains}
\label{sec:invimgchains}
Let $\pi: Y\to X$ be a semi algebraic (SA) bundle in the sense of \cite[Definition 8.2]{HLTV}, with $l$ dimensional fiber. In this section we want to define (or rather recall the definition of) the map
\[
 \pi^{-1}: C_\bullet(X) \to C_{\bullet+l}(Y).
\]
The relevant examples for us are $Y=D_{Ke}(m,n)$ or $Y=D_{Se}(m,n)$ and $X=D_{K}(m)$ or $X=D_{S}(m)$ with $\pi$ being the forgetful map in each case.
Everything said in this section is already contained in \cite{HLTV}, but a little scattered, so we recall here the relevant statements for the reader's convenience.

In Proposition 8.2 of \cite{HLTV} the authors define the strongly continuous chain $\Phi\in C_l^{str}(Y\to X)$ (see \cite{HLTV}, Definition 5.13) associated to the bundle $\pi:Y\to X$, such that 
\begin{enumerate}
 \item The image of a point $x\in X$ under $\Phi$ is the fundamental chain of the fiber $\pi^{-1}(x)$. 
 \item The boundary of $\Phi$ is the strongly continuous chain associated to the bundle $\pi^\p:Y^\p \to X$, given by the fiberwise boundaries of $Y$.
\end{enumerate}

For a semi algebraic chain $c\in C_k(X)$ we define 
\[
\pi^{-1}(c):= c \ltimes \Phi \in C_{k+l}(Y)
\]
where the operation $\ltimes$ is defined in \cite{HLTV}, Proposition 5.17. From the same Proposition it follows that 
\[
 \p \pi^{-1}(c)= \p c \ltimes \Phi +(-1)^k  c \ltimes \p \Phi = \p c \ltimes \Phi +(-1)^k (\pi^\p)^{-1}c
\]
with $\pi^\p$ as above.

\makeatletter{}\section{Operadic twisting}
\label{sec:optwists}
The notion of ``twisting'' for operads has been described in \cite{megrt}, Appendix G, while a more complete account has been given in \cite{vasilydeligne}. The reader is strongly advised to look at loc. cit. before reading this section. However, let us recall the main properties of operadic twists.
Let $\op P$ be any dg operad. Assume that there is an operad map
\[
 F: \Lie^{(k)}_\infty \to \op P
\]
 where, as before, $\Lie^{(k)}_\infty=(\Lie\{k\})_\infty$ is the minimal resolution of the degree shifted Lie operad. Then one can define the twisted operad $\Tw \op P$, which depends on the chosen map $F$ above, as follows. The underlying $\mathbb{S}$-module is
\begin{equation}
 \label{equ:twdef}
 \Tw \op P(n) = \prod_{j\geq 0} (\op P(n+j) \otimes (\R[k+1])^{\otimes j} )^{S_j}.
\end{equation}
Here $S_j$ acts on $\op P(n+j)$ by permutation of the last $j$ slots and on the factors of $\R[k+1]$ by permutation, i.e., by appropriate Koszul signs. For formulas for the operadic composition and the differential we refer to \cite{megrt}. 
Here we recall several properties:
\begin{itemize}
 \item Suppose $A$ is a $\op P$-algebra. By the map $F$ above, $A$ is also a $\Lie^{(k)}_\infty$-algebra. For a (pro-)nilpotent commutative algebra $\alg{n}$ one can hence define the notion of Maurer-Cartan element in $A\otimes \alg{n}$. This notion depends on $F$. One can twist the differential of $A\otimes \alg{n}$ using such a Maurer-Cartan element $m$. 
Let $p\in \Tw\op P(n+j)$ be some element symmetric (with the right signs) under permutation of the last $j$ slots, and let $\tilde p$ be the corresponding element in $\Tw\op P$.
The operad $\Tw \op P$ is defined such that the formulas
\begin{equation}
 \label{equ:twistedaction}
\tilde p(x_1,\dots, x_n) = \frac{1}{j!} p(x_1,\dots,x_n,m,\dots,m)
\end{equation}
define an action of $\Tw\op P$ on $A\otimes \alg{n}$ (with the $m$-twisted differential). Here $x_1,\dots, x_n\in A$.

 \item There is a natural action of the deformation complex $\Def(\Lie^{(k)}_\infty \to \op P)$ (it is a dg Lie algebra) on the operad $\Tw\op P$.
 \item There is a natural projection $\Tw\op P\to \op P$.
\end{itemize}

\begin{rem}
 The role of the nilpotent algebra $\alg{n}$ above is merely to ensure convergence in formulas like \eqref{equ:twistedaction}. One can define a Maurer-Cartan element directly as an element of $A$, if one imposes as extra condition that the infinite sums in \eqref{equ:twistedaction} are defined. In practice, this is often the case. In particular, it is the case for all left $\op P$-modules we twist in this paper.
\end{rem}

Note that there is always the natural map
\[
 \Tw\op P \to \op P
\]
projecting to the $j=0$ part in \eqref{equ:twdef}. For many operads, this map has a right inverse.
\begin{defi}
 Let $\op P$ be an operad, with an operad map $F: \Lie^{(k)}_\infty \to \op P$ for some $k$. We say that $\op P$ is \emph{natively twistable} if there is a map $\op P\to\Tw\op P$ such that the composition
\[
 \op P\to \Tw\op P \to \op P
\]
is the identity on $\op P$.
\end{defi}
If an operad $\op P$ is natively twistable, it means that $\op P$-algebras can be twisted by Maurer-Cartan elements, without leaving the category of $\op P$-algebras.
\begin{ex}
The operads $\op P =\Lie, \op P=\Lie_\infty, \op P=\Ger, \op P= \Ger_\infty$ are natively twistable. 
\end{ex}

Categorically, the operation of twisting is a functor (or rather a family of functors, one for each $k$)
\[
 \Tw \colon \Lie^{(k)}_\infty \downarrow \mathit{Operads} \to \mathit{Operads}.
\]

Here $\mathit{Operads}$ is the category of dg operads and $\Lie^{(k)}_\infty \downarrow \mathit{Operads}$ denotes the undercategory of $\Lie^{(k)}_\infty$, i.~e., the category of arrows $\Lie^{(k)}_\infty\to (\cdot)$. 
From this it follows that from the commutative triangle 

\begin{center}
  \begin{tikzpicture}[description/.style={fill=white,inner sep=2pt}, scale=.7]
\matrix (m) [matrix of math nodes, row sep=1em,
column sep=1.8em, text height=1.5ex, text depth=0.25ex]
{  &\Lie^{(k)}_\infty &  \\
\Lie^{(k)}_\infty&  & \op P\\ };
\path[->,font=\scriptsize]
(m-1-2) edge (m-2-1) edge (m-2-3)
(m-2-1) edge (m-2-3);
 \end{tikzpicture}
\end{center}
we obtain a map $\Tw \Lie^{(k)}_\infty\to \Tw \op P$. Then, because $\Lie^{(k)}_\infty$ is natively twistable we obtain an arrow $\Lie^{(k)}_\infty\to \Tw \op P$. 
It also follows that $\Tw$ is actually an endofunctor of the undercategory $\Lie^{(k)}_\infty \downarrow \mathit{Operads}$.
But in particular the following Lemma makes sense.

\begin{lemma}
\label{lem:native}
 Let $\op P$ be an operad. Let $F: \Lie^{(k)}_\infty \to \op P$ be an operad map (for some $k$). Then $\Tw \op P$ is natively twistable.
\end{lemma}
On the algebra level, this result says the following. If we have some $\op P$-algebra $A$ and twist it by a Maurer-Cartan element $m\in A$, we obtain some $\Tw \op P$ algebra $A^m$. If then we have another Maurer-Cartan element $m'\in A^m$, we do not need to change the operad again but can twist in the category of $\Tw \op P$-algebras. In fact, the resulting algebra will have the form $A^{m+m'}$.
\begin{proof}[Proof sketch]
 To show the lemma, we have to define a map $f: \Tw \op P\to \Tw \Tw \op P$, such that the composition with $\Tw \Tw \op P\to \Tw \op P$ is the identity. 
Let $\tilde p\in \Tw \op P(n)$ be given, with underlying (partially symmetric) element  $p\in \op P(n+j)$. To construct its image in $\Tw \Tw \op P(n)$ it is sufficient to define the projections to the components $(\op P(n+j_1+j_2)\otimes (\R[k+1])^{\otimes j_1} \otimes (\R[k+1])^{\otimes j_2})^{S_{j_1}\times S_{j_2}}$. We define it to be zero unless $j=j_1+j_2$. If $j=j_1+j_2$ we set the projection of the image equal to the image under the natural inclusion
\[
 (\op P(n+j) \otimes (\R[k+1])^{\otimes j} )^{S_j} \to (\op P(n+j_1+j_2)\otimes (\R[k+1])^{\otimes j_1} \otimes (\R[k+1])^{\otimes j_2})^{S_{j_1}\times S_{j_2}}
\]
\end{proof}

\begin{rem}
 The twisted operad can be introduced in two ways. We define it in \eqref{equ:twdef} as a ``partial deformation complex'' of the map $F: \Lie^{(k)}_\infty \to \op P$. Note that the zero-ary operations are really the deformation complex. Alternatively, one could (essentially) define it as the operad generated by $\op P$ and some zero-ary operation, modulo suitable relations. Essentially this amounts to replacing invariants by coinvariants in \eqref{equ:twdef}. It does not matter too much. We will stick to the ``partial deformation complex'' style in this paper, being consistent with \cite{megrt}.
\end{rem}

\begin{rem}
In fact one can show that the functor $\Tw$ is a co-monad on the undercategory $\Lie^{(k)}_\infty \downarrow \mathit{Operads}$. Its co-algebras are natively twistable operads $\op P$, that satisfy an additional condition saying that twisting a $\op P$ algebra with an MC element $m$ and then again with an MC element $m'$ is the same as twisting only once with an MC element $m+m'$. We will not need this additional condition here. For more details on the categorial properties of $\Tw$ see \cite{vasilydeligne}.
\end{rem}

\subsection{Twisting right modules}
\label{sec:twistrightmodules}
We saw how to twist left $\op P$-modules to left $\Tw\op P$-modules in the last section. Now let us consider twisting of right $\op P$-modules. Here again $\op P$ is some (dg) operad, equipped with a map $F: \Lie^{(k)}_\infty \to \op P$ for some $k$. Everything will depend on $F$ and $k$, though we do not indicate the dependence in the notation.
Let $\op M$ be an operadic right $\op P$-module. As an $\mathbb{S}$-module the twisted module is
\[
 \Tw\op M(n) = \prod_j (\op M(n+j) \otimes (\R[k+1])^{\otimes j} )^{S_j}
\]
with the $S_j$ action similar to the one in the operadic case. The operadic right action is defined by similar formulas as for the operadic composition.
Let $p\in \op P(n_2+j_2)$ be (signed) symmetric in the last $j_2$ slots and let $\tilde p\in \Tw\op P(n_2)$ be the corresponding element in $\Tw\op P$. Similarly, let $m\in \op M(n_1+j_1)$ be symmetric (with the correct signs) in the last $j_1$ slots and let $\tilde m\in \Tw\op M(n_1)$ be the corresponding element. Then, for $1\leq l\leq n_1$ we define the operadic composition so that
the element $\tilde m\circ_l \tilde p$ is described by the following element in $\op M(n_1+n_2-1+j_1+j_2)$, which is symmetric (with signs) under permutations of the last $j_1+j_2$ slots.
\begin{multline}
\label{equ:twrightactiondef}
 (m\circ_l p)(s_1,\dots, s_{n_1+n_2-1}, \bar 1,\dots \overline{j_1+j_2})
= \\ = 
\sum_{I\sqcup J=[j_1+j_2]}
\sgn(I,J)^{k+1}(-1)^{|m|j_1 (k+1)}
m(s_1,\dots,s_{l-1},p(s_l,\dots,s_{l+n_2-1},\bar J), \dots,s_{n_1+n_2-1}, \bar I).
\end{multline}
Here notation from \cite{megrt}, Appendix G is used. The $s_1,s_2,\dots$ are some ``placeholder'' symbols for the input of operations in the operad (or module). The placeholder symbols $\bar 1, \bar 1,\dots$ indicate the slots in which the operations $p, m$ should be symmetric. For the operadic right action a ``functional'' notation is used.

We still need to define the differential on $\Tw\op M$. For this, let temporarily $\widetilde{\Tw\op M}$ be the above operadic right module, with the differential solely that coming from $\op M$.
Let 
\[
 \alg{g} = \Def(\Lie^{(k)}_\infty \stackrel{0}{\to} \op P) = \prod_{j\geq 0} (\op P(j) \otimes (\R[k+1])^{\otimes j} )^{S_j}.
\]
be the deformation complex of the zero map. It is a dgla. Similar to \cite{megrt}, Lemma 9, there is a right action of $\alg{g}$ on $\Tw\op M$. For $m,\tilde m$ as above and $x\in g$ from the $j_3$-th term in the product above, the element $\tilde m \cdot x\in\Tw\op M(n_1)$ is defined by the following (partially symmetric) element $m\cdot x$ of $\op M(n_1+j_1+j_3-1)$
\[
 (m\cdot x)(1,\dots,n_1,\bar 1,\dots \overline{j_1+j_3})
=
\sum_{I\sqcup J=[j_1+j_3]}
\sgn(I,J)^{k+1}
m(1,\dots,n_1, \bar I, x( \bar J)).
\]
Now the Lie algebra $\widetilde{\Tw\op P}(1)$ (notation as in \cite{megrt}) acts on $\widetilde{\Tw\op M}$ from the right, by operations
\[
 \tilde m \cdot q = \sum_{j=1}^{n_1} \tilde m \circ_j q
\]
for $q\in \widetilde{\Tw\op P}(1)$. The Lie algebra $\alg{g}$ also acts on $\widetilde{\Tw\op P}$ by operadic derivations, hence also on $\widetilde{\Tw\op P}(1)$. Both actions on $\widetilde{\Tw\op M}$ can be merged into one right action of the Lie algebra 
\[
 \hat{\alg{g}}=\alg{g}\ltimes \widetilde{\Tw}\op P(1)
\]
by operadic right module derivations. $\hat{\alg{g}}$ also acts on $\widetilde{\Tw\op P}$ from the right by operadic derivations. By multiplying with a sign, we can change the right action to a left action.
Picking any Maurer-Cartan element in $\hat{\alg{g}}$, we can twist simultaneously the operad $\widetilde{\Tw\op P}$ and its module $\widetilde{\Tw\op M}$. The operad map $F$ defines a Maurer-Cartan element $\mu\in \alg{g}$. Then Lemma 10 of \cite{megrt} produces for us a Maurer-Cartan element $\hat{\mu}\in \hat{\alg{g}}$.

\begin{defi}[Definition 2 of \cite{megrt}]
 The twisted operad $\Tw \op P$ is defined to be $\widetilde{\Tw\op P}$ as a graded operad, equipped with differential 
\[
 d_{\op P} + \hat{\mu}\cdot.
\]
where $d_{\op P}$ is the differential coming from $\op P$, $\hat{\mu}$ is as above and $\hat{\mu}\cdot$ denotes its left action.
\end{defi}
\begin{defi}
 The twisted right $\Tw \op P$ module $\Tw \op M$ is defined to be $\widetilde{\Tw\op M}$ as a graded operadic right module, equipped with differential 
\[
 d_{\op M} + \hat{\mu}\cdot.
\]
where $d_{\op M}$ is the differential coming from $\op M$, $\hat{\mu}$ is as above and $\hat{\mu}\cdot$ denotes its left action.
\end{defi}

From this definition it follows that there is naturally an action of the $\mu$-twisted version of the dg Lie algebra $\alg{g}$ on both $\Tw \op P$ and $\Tw \op M$.
This twisted Lie algebra is the deformation complex of the map $F$ from above.

\subsection{Twisting bimodules}
\label{sec:bimodtwist}
Little is to be said about the twisting of bimodules.
Let $\op P$, $\op Q$ be operads and let $\op M$ be a $\op P$-$\op Q$ operadic bimodule. We first twist $\op Q$ to an operad $\Tw \op Q$ (this depends on a chosen map $\Lie^{(k)}_\infty \to \op Q$), and twist $\op M$ to a $\op P$-$\Tw \op Q$ operadic bimodule $\Tw \op M'$. Then we twist $\op P$ to $\Tw \op P$ and twist the operadic left module $\Tw \op M'$ as in indicated in beginning of Appendix \ref{sec:optwists}. I.e., we assume that there is some Maurer-Cartan element $m\in \Tw \op M'(0)$ and give $\Tw \op M'$ a new differential using $m$. The resulting operadic $\op P$-$\op Q$ bimodule we call $\Tw \op M$.
We will furthermore assume that the Maurer-Cartan element $m$ is sent to zero under the projection $\Tw \op M\to\op M$. The reason is that then we again have a map of colored operads
\[
 \bpm\Tw \op P & \Tw \op M & \Tw \op Q \epm \to  \bpm \op P &  \op M & \op Q \epm
\]
similar to the analogous map in the uncolored operad case.

\subsection{Twisting moperads}
\label{sec:moperadtwist}
Let next $\op P$ be an operad and $\op P_1$ be a $\op P$-moperad. I.e., 
\[
 \bpm \op P & \op P_1 \epm
\]
is a two colored operad. Let further $\Lie_1\{k\}$ be the $\Lie\{k\}$ moperad governing Lie algebra modules. It has a minimal resolution, $\hooLie_{k,1}$, so that the two colored operad
\[
 \bpm  \Lie^{(k)}_\infty & \hooLie_{k,1}\epm
\]
governs homotopy $\Lie\{k\}$ algebras (e.g., $\Lie_\infty$ algebras for $k=0$) together with homotopy $\Lie\{k\}$ modules (e.g., $\Lie_\infty$ modules for $k=0$). Suppose that we have a map 
\[
 F\colon \Lie^{(k)}_\infty \to \op P
\]
as before and additionally a map 
\[
 F_1\colon \hooLie_{k,1} \to \op P_1
\]
of moperads, i.e., altogether we have a map of colored operads
\[
 \bpm \Lie^{(k)}_\infty & \hooLie_{k,1}\epm \to \bpm \op P & \op P_1 \epm.
\]
In this situation one can twist the colored operad $\bpm \op P & \op P_1 \epm$ by a similar construction as above. 
Let $\alg{g}$ be the deformation complex as before and
\[
 \alg{h}:= \Def(\hooLie_{k,1} \stackrel{0}{\to} \op P_1) := \prod_{j\geq 0} (\op P_1(j) \otimes (\R[k+1])^{\otimes j} )^{S_j}.
\]
It is a dg Lie algebra and a $\alg{g}$ module, so that together we have a Lie algebra
\[
 \alg{g} \ltimes \alg{h}.
\]
The maps $F$ and $F_1$ provide a Maurer-Cartan element $\mu+\nu$ in this dg Lie algebra, with the part $\mu\in \alg{g}$ contributed by $F$ as above and $\nu\in \alg{h}$ contributed by $F_1$. This Lie algebra acts on $\widetilde{\Tw\op P}$ and $\widetilde{\Tw\op P_1}$ compatibly with (as derivations of) the (m)operadic structures. Here $\alg{h}$ acts trivially on $\widetilde{\Tw\op P}$, the actions of $\alg{g}$ on $\widetilde{\Tw\op P}$ and $\widetilde{\Tw\op P_1}$ we have encountered before and the action of $\alg{h}$ on $\widetilde{\Tw\op P_1}$ is by a similar formula. Again as before, we also have an action of the Lie algebra $\widetilde{\Tw\op P}(0)$ on $\widetilde{\Tw\op P}$ and on $\widetilde{\Tw\op P_1}$. We can form the Lie algebra 
\[
  \alg{g} \ltimes \alg{h}\ltimes \widetilde{\Tw\op P}(0) \cong \alg{g} \ltimes\widetilde{\Tw\op P}(0)\ltimes \alg{h} .
\]
In this Lie algebra we have a Maurer-Cartan element $\hat{\mu}+\nu$ (here $\hat{\mu}$ is as in \ref{sec:twistrightmodules}). Twisting $\widetilde{\Tw\op P}$ and $\widetilde{\Tw\op P_1}$ with this Maurer-Cartan element we obtain the twisted operad $\Tw\op P$ (same as before) and the twisted moperad $\Tw\op P_1$. Concretely, as a graded moperad $\Tw\op P_1$ is the same as $\widetilde{\Tw\op P_1}$, but is has the differential 
\[
 d_{\op P_1} + \hat{\mu}\cdot + \nu\cdot. 
\]

\begin{rem}
 The difference of between twisting $\op P_1$ as a right $\op P$-module as in section \ref{sec:twistrightmodules} and twisting as a moperad is the part $\nu\cdot$ in the above differential. The information contained in the map $F_1$ goes into this part of the differential, while $\hat{\mu}\cdot$ depends only on $F$ and $d_{\op P_1}$ is (of course) independent of both.
\end{rem}

It follows that there is an action of the twisted (by $\mu+\nu$) version of $\alg{g} \ltimes \alg{h}$ on $\Tw\op P$ and $\Tw\op P_1$.

\begin{rem}
 Note however that we do not directly obtain an action of (the twisted version of) $\alg{g}$ alone. E.g., suppose we have a ($\mu$-)closed element $x\in \alg{g}$, i.e., $dx+\co{\mu}{x}=0$. Then to find a ($\mu+\nu$-)closed element $x+y\in \alg{g}\ltimes \alg{h}$, we have to solve the equation 
\[
 D y := dy+\co{\mu}{y}+\co{\nu}{y} = -\co{\nu}{x}
\]
for $y\in \alg{h}$. It is not a priori clear that such a $y$ should always exists. \end{rem}

\subsection{Twisting of moperadic bimodules}
\label{sec:mopbimodtwist}
Finally let us twist moperadic bimodules. We start with operads $\op P$, $\op Q$, a $\op P$-moperad $\op P_1$, a $\op Q$-moperad $\op Q_1$, a $\op P$-$\op Q$ bimodule $\op M$ and the moperadic bimodule $\op M_1$. Concretely, $\op M_1$ is endowed with a right action of $\op Q$, a left action of $\op Q_1$ and a right action of $\op P$ and $\op M$ combined (see Figure \ref{fig:moperadicbimodcomp}). 

\begin{rem}
 To reduce confusion about the many letters, keep in mind that in our situation $\op P$ acts on $\Dpoly$, $\op P_1$ on $C_\bullet$, $\op Q$ acts on $\Tpoly$, $\op Q_1$ on $\Omega_\bullet$. The bimodule $\op M$ controls a map $\Tpoly \to \Dpoly$, the moperadic bimodule $\op M_1$ controls a map $C_\bullet \to \Omega_\bullet$.
\end{rem}
In the preceding sections, we have already seen how to twist $\op P$, $\op Q$, $\op P_1$, $\op Q_1$ and $\op M$ to $\Tw\op P$, $\Tw\op Q$, $\Tw\op P_1$, $\Tw\op Q_1$ and $\Tw\op M$. These twists depend on some choices (concretely maps from $\Lie^{(k)}_\infty$ or $\hooLie_{k,1}$ to the (m)operads and a choice of a Maurer-Cartan element in $\Tw \op M'$), which we assume have been made. 
Let us next consider $\op M_1$. Disregarding the differential, the twisted version of $\op M_1$ is 
\[
 \Tw \op M_1(n) := \prod_{j\geq 0} (\op M_1(n+j) \otimes (\R[k+1])^{\otimes j} )^{S_j}.
\]
Here $S_j$ acts by permuting the last $j$ input slots colored by the color of $\op Q$. So, disregarding the differntial $\Tw \op M_1$ is the same as the twisted version of $\op M_1$, regarded (only) as a right $\op Q$-module. This also defines the right action of $\Tw\op Q$ on $\Tw \op M_1$ (the formula is identical to \eqref{equ:twrightactiondef}).
There is a left moperadic action of $\Tw\op Q_1$. Let $\tilde q\in \Tw\op Q_1(n_1)$, with the underlying partially symmetric element $q\in \op Q_1(n_1+j_1)$. Similarly let $\tilde m_1\in\Tw \op M_1(n_2)$ with the underlying partially symmetric element $m_1\in \op M_1(n_2+j_2)$. Then the composition $\tilde q\circ \tilde m_1\in\Tw \op M_1(n_1+n_2)$ is defined by the partially symmetric element $q\circ m_1\in \op M_1(n_1+n_2+j_1+j_2)$ given by the following formula.
\begin{multline}
 (q\circ m_1)(s_1,\dots, s_{n_1+n_2}, \bar 1,\dots \overline{j_1+j_2}; t)
= \\ = 
\sum_{I\sqcup J=[j_1+j_2]}
\sgn(I,J)^{k+1}(-1)^{|m|j_1 (k+1)}
q(s_1,\dots, \dots,s_{n_1}, \bar I; m_1(s_{n_1+1},\dots,s_{n_1+n_2},\bar J; t)).
\end{multline}
The formula for the missing right action of $\Tw\op P_1$ and $\Tw \op M$ is simialar, but notationally too horrible to display.

Let us finally give the formula for the differential on $\Tw \op M_1$. It has the form
\[
 d_{\op M_1} + \mu_{\op Q}\cdot + \nu_{\op Q_1}\cdot + L_(\nu_{\op P_1}, m).
\]
Let us consider the various parts. First $d_{\op M_1}$ comes from the differential on $\op M_1$. Next, because we have an action of $\Tw \op Q$ on $\op M_1$, we in particular have an action of the dg Lie algebra $\alg{g}_{\op Q}$. Here $\alg{g}_{\op Q}$ is\footnote{Previously we called this dg Lie algebra just $\alg{g}$, but now we need to distinguish two versions of this object, one for $\op P$ and one for $\op Q$. Similarly we will distinguish $\alg{h}_{\op P_1}$ and $\alg{h}_{\op Q_1}$ and $\mu_{\op P}, \nu_{\op P_1}$ and $\mu_{\op Q}, \nu_{\op Q_1}$.}
\[
 \alg{g}_{\op Q} := \Def(\Lie^{(k)}_\infty \stackrel{0}{\to} \op Q).
\]
The Maurer-Cartan element $\mu_{\op Q}\in \alg{g}_{\op Q}$ is the one corresponding to the operad map $\Lie^{(k)}_\infty \to \op Q$ and $\mu_{\op Q}\cdot$ is its action. Similarly, we have an action of 
\[
\alg{h}_{\op Q_1} := \Def(\hooLie_{k,1} \stackrel{0}{\to} \op Q_1).
\]
The term $\nu_{\op Q_1}\cdot$ is the action of the element $\nu_{\op Q_1}\in \alg{h}_{\op Q_1}$. This element has appeared in section \ref{sec:moperadtwist} as $\nu$.
The most difficult term is $L_(\nu_{\op P_1}, m)$. This is the action of the element $\nu_{\op P_1}\in \alg{h}_{\op P_1}$, together with the Maurer-Cartan element $m$ used to twist $\op M$. Here $\nu_{\op P_1}$ and $\alg{h}_{\op P_1}$ are the counterparts for $\op P$, $\op P_1$ of $\nu_{\op Q_1}$ and $\alg{h}_{\op Q_1}$. We do not want to give the lengthy formula. But consider Figure \ref{fig:moperadicbimodcomp} (top part) to see a picture of the right moperadic bimodule action. In our situation the dark grey box represents the element $\nu_{\op P_1}$. It has inputs of two different colors, namely multiple inputs in the color of $\op P$ and exactly one input in the output color of $\op P_1$. In each of the $\op P$ colored inputs, one inserts one copy of $m$ and divides by a factorial. So, in Figure \ref{fig:moperadicbimodcomp}, the white circles have to be filled with copies of $m$. The grey circles are not present (i.~e., are filled by copies of the unit of $\op Q$). 

\begin{rem}
Note that in order to twist the moperadic bimodule, we do not need any additional data, or make additional choices.
\end{rem}

\subsection{Colored case}
\label{sec:coloredtwist}
We have to deal with four-colored operads of the form 
\[
 \op C = \bpm \op P & \op M & \op Q \\ \op P_1 & \op M_1 & \op Q_1 \epm
\]
where $\op P$, $\op Q$, $\op P_1$, $\op Q_1$, $\op M$, $\op M_1$ are as in the previous subsection.
We will write 
\[
 \Tw\op C = \bpm \Tw\op P & \Tw\op M & \Tw\op Q \\ \Tw\op P_1 & \Tw\op M_1 & \Tw\op Q_1 \epm
\]
for its twisted version. The twist depends on various choices as detailed above. We will hide those choices in the notation.
\begin{defi}
 Let $\op C$ be a four colored operad as above and let $\Tw\op C$ be its twist as above. We say that $\op C$ is \emph{natively twistable} if there is a map $\op C\to\Tw\op C$ such that the composition
\[
 \op C\to \Tw\op C \to \op C
\]
is the identity on $\op C$.
\end{defi}

\makeatletter{}\section{Colored operadic twists of several operads }
\label{sec:coloredtwists}
\subsection{Twisting $\bigChains$}
\label{sec:bigchainstwist}
Let us consider the colored operad 
\[
\bigChains
=
\bpm
\Br & C(D_K) & C(\FM_2) \\
\KS_1 & C(D_S) & C(\FM_{2,1})
\epm.
\]
We want to twist it to a colored operad $\Tw\bigChains$. According to our conventions (see section \ref{sec:coloredtwist}) for this we need to pick Maurer-Cartan elements in the operadic bimodule $\Tw C(D_K)(0)$. Note that by definition 
\[
\Tw C(D_K)(0) =  \sum_n C(D_K(n))^{S_n}.
\]
We will take as Maurer-Cartan element
\[
 \nu := \sum_{n\geq 0} \Fund(C(D_K(n))) 
\]
the sum of fundamental chains. 
\subsection{Twisting of $\bigGra$}
\label{sec:biggraphstwist}
Let us consider the colored operad 
\[
\bigGra
=
\bpm
\Br & \SGra & \Gra \\
\KS_1 & \SGra_1 & \Gra_1
\epm.
\]
To twist it, we again need to specify a Maurer-Cartan element $\kappa\in\Tw\SGra(0)$.
However, such a Maurer-Cartan element can be taken to be the image of $\nu$ as defined above under the map $\bigChains\to \bigGra$.
Concretely, the Maurer-Cartan element $\nu$ is the universal Kontsevich star product.

\subsection{\texorpdfstring{$\bigChains$}{bigChains} is natively twistable}
\label{sec:bigchainsnattwist}
In this section we want to construct a map
\[
 F: \bigChains \to \Tw\bigChains
\]
such that the composition with the natural projection 
\[
 \bigChains \to \Tw\bigChains \to \bigChains 
\]
is the identity. The colored operad $\bigChains$ is generated by 6 parts (2 operads. 2 moperads, an operadic and a moperadic bimodule). To construct the map, we will construct the map for each of these 6 parts.
The simplest parts are $\KS$ and $\KS_1$. Recall that $\bpm \KS & \KS_1 \epm$ is defined as a certain suboperad of $\Tw\bpm \PT & \PT_1 \epm$. Since by Lemma \ref{lem:native} (or, more precisely, its colored analogon) every twisted operad is automatically natively twistable, we can map
\[
\Tw \bpm \KS & \KS_1 \epm \to \bpm \PT & \PT_1 \epm \to \left( \Tw\bpm \PT & \PT_1 \epm \right) \subset \Tw\bigChains.
\]
Next consider the part $C(\FM_2)$. Recall that $\FM_2(n)$ consists of configurations of $n$ numbered points in $\R^2$, modulo scaling and translation (up to compactification). The twisted version $\Tw C(\FM_2)(n)$ can be interpreted as chains on the configuration space of $n$ numbered and arbitrarily many unnumbered points, symmetric under permutations of the unnumbered points. Concretely, by definition
\[
\Tw C(\FM_2)(n) := \prod_k C(\FM_2(n+k))^{S_k}.
\]
There are forgetful maps $\pi_{n,k}: \FM_2(n+k)\to \FM_2(n)$. Recall from Appendix \ref{sec:invimgchains} that there is an inverse image map on chains
\[
 \pi_{n,k}^{-1} : C(\FM_2(n))\to C(\FM_2(n+k)).
\]
Chains in the image of $\pi_{n,k}^{-1}$ are invariant under permutation of the $k$ additional points. Hence we can define for a chain $c\in C(\FM_2(n))$
\[
 F(c) := \sum_k \pi_{n,k}^{-1}(c) \in \Tw C(\FM_2)(n).
\]
One can check that this defines a map of operads $C(\FM_2)\to \Tw C(\FM_2)$ and clearly the composition with the natural projection $\Tw C(\FM_2)\to C(\FM_2)$ is the identity.
By analogous formulas one extends this map to $C(\EFM_2)$, the operadic bimodule $C(D_K)$ and the moperadic bimodule $C(D_S)$. In each case the untwisted space is a space of (semi algebraic chains of) configurations of numbered points. The twisted space is the space of configurations of numbered and unnumbered points. The map $F$ is defined by taking the ``inverse image'' on chains of appropriate forgetful maps as above.

\makeatletter{}\section{Proof of Proposition \ref{prop:graphs1cohom}}
\label{sec:graphs1cohomproof}
\subsection{Review of a proof by P. Lambrechts and I. Volic}
\label{sec:LVproof}
Let us recall here the proof of Theorem \ref{thm:graphscohom} (i.~e., the proof that $H(\Graphs)\cong \Ger$) due to Lambrechts and Volic \cite{LV}. We will present it in a way that can be modified with little complication to a proof of Proposition \ref{prop:graphs1cohom} (i.~e., a proof of the fact that $H(\Graphs_1)\cong\calc$).
We will actually (and equivalently) compute the cohomology of the predual $\pdu\Graphs$. This is the space of linear combination of $\Graphs$-graphs, but with the differential $\delta$ being edge contraction instead of vertex splitting. Then $H(\Graphs)=H(\pdu\Graphs)^*$

\begin{lemma}
\label{lem:pduGraphs}
 A basis for $H(\pdu\Graphs)$ is given by the cohomology classes of graphs without internal vertices and such that the external vertex $j$ is directly connected to at most one of the vertices $1,2,\dots,j-1$. In particular $H(\pdu\Graphs)=\pdu\Graphs_0 / \delta \pdu\Graphs_1$, where $\pdu\Graphs_k$ is the subspace spanned by graphs with $k$ internal vertices.
\end{lemma}
\begin{proof}[Proof by Lambrechts and Volic]
 We show the statement for $H(\pdu\Graphs(n))$ by induction on $n$. For $n=0$ there is nothing to be shown. For $n\geq 1$ we can split $\pdu\Graphs(n)$ (as complexes)
\[
 \pdu\Graphs(n) = \pdu\Graphs(n)^0 \oplus \pdu\Graphs(n)^{\geq 1}
\]
where $\pdu\Graphs(n)^0$ is the subscomplex spanned by graphs where the external vertex $n$ has valence 0 and $\pdu\Graphs(n)^{\geq 1}$ is spanned by all other graphs (i.e., those in which vertex $m$ has valence $\geq 1$). Clearly 
\[
H(\pdu\Graphs(n)^0)
\cong
 \begin{cases}
  \R &\text{for $n=1$} \\
  H(\pdu\Graphs(n-1))&\text{otherwise}.
 \end{cases}
\]
Hence, using the induction hypothesis, we see that $H(\pdu\Graphs(n)^0)$ contributes those graphs to the (tentative) basis described in the lemma, for which vertex $n$ has valence 0.
Let us continue to evaluate $H(\pdu\Graphs(n)^{\geq 1})$. Split (as vector spaces):
\[
 \pdu\Graphs(n)^{\geq 1} = \pdu\Graphs(n)^{1} \oplus \pdu\Graphs(n)^{\geq 2}
\]
into parts where vertex $n$ has valence 1 or $\geq 2$ respectively. The direct sum here is not a direct sum of complexes. However, there is an associated spectral sequence whose first differential is the part of $\delta$ mapping $\pdu\Graphs(n)^{1}\to\pdu\Graphs(n)^{\geq 2}$. This map is surjective, its kernel is spanned by graphs in which vertex $n$ connects to another external vertex. Call this space $\pdu\Graphs(n)^{1,e}$. It is the $E^1$ term of our spectral sequence. Furthermore is is a subcomplex of $\pdu\Graphs(n)^{\geq 1}$, and hence
\[
 H(\pdu\Graphs(n)^{\geq 1}) = H(\pdu\Graphs(n)^{1,e}).
\]
However, the latter complex splits into $n-1$ subcomplexes according to which external vertex the vertex $n$ connects to. Each of these is isomorphic to $\pdu\Graphs(n-1)$. Hence, by using the induction hypothesis again, we see that $H(\pdu\Graphs(n)^{\geq 1})$ contributes those elements of the basis announced in the lemma for which $m$ has valence 1.
\end{proof}
Next we want to show that $H(\pdu\Graphs)^*\cong \Ger$. To do that, it is convenient to use a slightly different basis of the cohomology.
\begin{lemma}
\label{lem:anotherbasis}
 Another basis of $H(\pdu\Graphs)^*$ is given by the classes of graphs $\Gamma$ such that 
\begin{itemize}
 \item There are no internal vertices in $\Gamma$.
 \item Each external vertex has at most valence 2.
 \item The external vertex with lowest label in each connected component has valence at most 1.
\end{itemize}
\end{lemma}
This basis consists of graphs formed by several ``strings'' of external vertices, such that the lowest numbered vertex in each string is at the end of the string.
\begin{proof}
 By counting we see that the basis has the same cardinality as the one constructed in the previous lemma. Hence it suffices to show that any graph of the previous basis can be written as a linear combination of graphs in the new (tentative) basis, modulo relations $\delta \pdu\Graphs_1$. This is a simple exercise.
\end{proof}

\begin{cor}
\label{cor:hgraphs}
 $H(\pdu\Graphs)^*\cong \Ger$ and hence $H(\Graphs)\cong \Ger$.
\end{cor}
There is an explicit embedding $\Ger\to \Graphs$ given by the formulas of the remark in section \ref{sec:GrdGr}. We will show the stronger statement that this embedding is a quasi-isomorphism.
\begin{proof}
 A basis of $\Ger(n)$ is given by symmetric products of $\Lie\{1\}$-words, i.e., expressions of the form 
\[
 L_1(X_1,\dots, X_n)\wedge \dots \wedge L_k(X_1,\dots, X_n)
\]
for $k=1,2,\dots$ such that
\begin{itemize}
 \item Each $L_j$ has the form $[X_{\alpha_1}, [X_{\alpha_2}, \cdots, [X_{\alpha_{r-1}},X_{\alpha_{r}}]\cdots]$ where $\alpha_{r}<\alpha_1,\dots,\alpha_{r-1}$.
 \item In the symmetric product of Lie words, each $X_j$, $j=1,\dots, n$ occurs exactly once.
\end{itemize}

We claim that under the map $\Ger\to H(\Graphs)$, this basis is dual to the basis of Lemma \ref{lem:anotherbasis}. This will proof the Corollary.
In fact, mapping a product of Lie words as above by $\Ger\to \Graphs$ one can see that nonzero value is attained on exactly one element of the basis from Lemma \ref{lem:anotherbasis}. Namely, this element has one ``string'' of vertices for each Lie word in the product, and the order of vertices on the string is the same as in the Lie word. 
\end{proof}

\subsection{The proof}
We want to prove Proposition \ref{prop:graphs1cohom} by similar arguments as in the previous subsection. To do this, one first convinces oneself that the graphs of Figure \ref{fig:dliota} are indeed closed and satisfy the $\calc$-relations from the introduction. One concludes that there is a map of moperads
\[
\calc_1 \to \Graphs_1.
\]
We want to show that it is a quasi-isomorphism.  We do this by showing that $H(\pdu\Graphs_1)^*=\calc_1$. In fact, the central part of the proof will be again be the trick due to P. Lambrechts and I. Volic \cite{LV}. The predual $\pdu\Graphs_1$ has the following description:
\begin{itemize}
 \item Elements are linear combinations of $\Graphs_1$-graphs.
 \item The differential $\delta$ is the dual differential. Concretely, it contains two terms:
\begin{enumerate}
 \item Edge contraction: A part contracting each edge, which is incident to at least one internal vertex.
 \item Merging of an internal vertex with $\vin$ or $\vout$: Each internal vertex with incoming edges only, and not connected to $\vin$, is merged with $\vin$ and an edge removed. Each internal vertex with outgoing edges only, and not connected to $\vout$, is merged with $\vout$ and an edge removed.
\end{enumerate}
\end{itemize}
\begin{figure}
 \centering
\makeatletter{}\usetikzlibrary{matrix}
\usetikzlibrary{arrows}
\usetikzlibrary{shapes}
\usetikzlibrary{through}
\usetikzlibrary{calc,3d}
\usetikzlibrary{decorations,decorations.pathmorphing}
\[
\begin{tikzpicture}[
int/.style={circle, draw, fill, minimum size=5pt, inner sep=0},
ext/.style={circle, draw, fill=white, minimum size=5pt, inner sep=1pt},
helper/.style={coordinate},point/.style={circle, draw, fill, inner sep =1pt},
de/.style={-triangle 60},
point/.style={circle, draw, fill, minimum size=3pt, inner sep=0pt},
xst/.style={cross out, draw, minimum size=5 },
boxed/.style={ draw, inner sep=0.5 },
]
\begin{scope}[xshift=2cm, yshift=-1cm]

\begin{scope}[xshift=-4.5cm]
\node at (2.7,0) {$\delta$};
\draw (3.5,.3) node[int] (i2) {} --+(0,-.6) node[int] (i3) {};
\draw[-triangle 60] (i2) edge (i3);
\foreach \x in {-.6,.2,-.2,.6}
{
\draw (i2)--+(\x,.6) ;
\draw (i3)--+(\x,-.6) ;
}
\end{scope}

\node at (0,0) {$=$};
\node[int] (i1) at (1,0) {};
\foreach \y in {.6,-.6}
\foreach \x in {-.6,.2,-.2,.6}
\draw (i1)--+(\x,\y) ;

\end{scope}

\begin{scope}[xshift=7cm,yshift=-1cm]
\node at (0,0) {$\delta$};

\draw (1,.6) node[int] (i1) {} --+(0,-.8) node[ext] (e2) {$j$};
\draw[-triangle 60] (e2) edge (i1);
\foreach \x in {.4,0,-.4}
\draw (i1)--+(\x,.6) ;
\foreach \x in {.7,-.7}
\draw (e2)--+(\x,.6) ;

\node at (5.2,0) {$=$};
\draw (4,.6) node[int] (i1) {} --+(0,-.8) node[ext] (e2) {$j$};
\draw[-triangle 60] (i1) edge (e2);
\foreach \x in {.4,0,-.4}
\draw (i1)--+(\x,.6) ;
\foreach \x in {.7,-.7}
\draw (e2)--+(\x,.6) ;

\node[ext] (e1) at (6.5,0) {$j$};
\foreach \x in {-.8,.4,0,-.4,.8}
  \draw (e1)--+(\x,.7) ;
\node at (2.5,0) {$=\quad \delta$};

\end{scope}

\begin{scope}[yshift=-4cm]
\node at (0,0) {$\delta$};
\node[boxed] (e1) at (3.5,0) {$\vin$};
\foreach \x in {-.8,.4,0,-.4,.8}
\draw[triangle 60-] (e1)--+(\x,.7) ;

\node at (2.25,0) {$=$};
\draw (1,.6) node[int] (i1) {} --+(0,-.8) node[boxed] (e2) {$\vin$};
\draw[triangle 60-] (e2) edge (i1);

\foreach \x in {.4,0,-.4}
\draw[triangle 60-] (i1)--+(\x,.7) ;
\foreach \x in {.7,-.7}
\draw[triangle 60-] (e2)--+(\x,.7) ;

\end{scope}

\begin{scope}[yshift=-7cm]

\node at (0,0) {$\delta$};

\draw (1,.6) node[int] (i1) {} +(0,-.8) node[boxed] (e2) {$\vin$};

\foreach \x in {.4,0,-.4}
\draw[triangle 60-] (i1)--+(\x,.7) ;
\foreach \x in {.7,-.7}
\draw[triangle 60-] (e2)--+(\x,.7) ;
\node at (2.5,0) {$=\quad \sum$};
\node[boxed] (e1) at (4,0) {$\vin$};
\foreach \x in {-.8,.4,0,.8}
\draw[triangle 60-] (e1)--+(\x,.7) ;
\draw[triangle 60-, dotted] (e1)--+(-.4,.7) ;
\end{scope}

\begin{scope}[yshift=-4cm, xshift=7cm]
\node at (0,0) {$\delta$};
\node[boxed] (e1) at (3.5,0) {$\vout$};
\foreach \x in {-.8,.4,0,-.4,.8}
\draw[-triangle 60] (e1)--+(\x,.7) ;

\node at (2.25,0) {$=$};
\draw (1,.6) node[int] (i1) {} --+(0,-.8) node[boxed] (e2) {$\vout$};
\draw[-triangle 60] (e2) edge (i1);

\foreach \x in {.4,0,-.4}
\draw[-triangle 60] (i1)--+(\x,.7) ;
\foreach \x in {.7,-.7}
\draw[-triangle 60] (e2)--+(\x,.7) ;

\end{scope}

\begin{scope}[yshift=-7cm, xshift=7cm]

\node at (0,0) {$\delta$};

\draw (1,.6) node[int] (i1) {} +(0,-.8) node[boxed] (e2) {$\vout$};

\foreach \x in {.4,0,-.4}
\draw[-triangle 60] (i1)--+(\x,.7) ;
\foreach \x in {.7,-.7}
\draw[-triangle 60] (e2)--+(\x,.7) ;
\node at (2.5,0) {$=\quad \sum$};
\node[boxed] (e1) at (4,0) {$\vout$};
\foreach \x in {-.8,.4,0,.8}
  \draw[-triangle 60] (e1)--+(\x,.7) ;
\draw[-triangle 60, dotted] (e1)--+(-.4,.7) ;
\end{scope}

\end{tikzpicture}
\]  
\caption{\label{fig:graphs1dualdiff} The differential on $\pdu\Graphs_1$. A dotted arrow means that the edge is deleted from the graph. If no arrow is drawn on an edge it means that the edge can have either orientation.}
\end{figure}
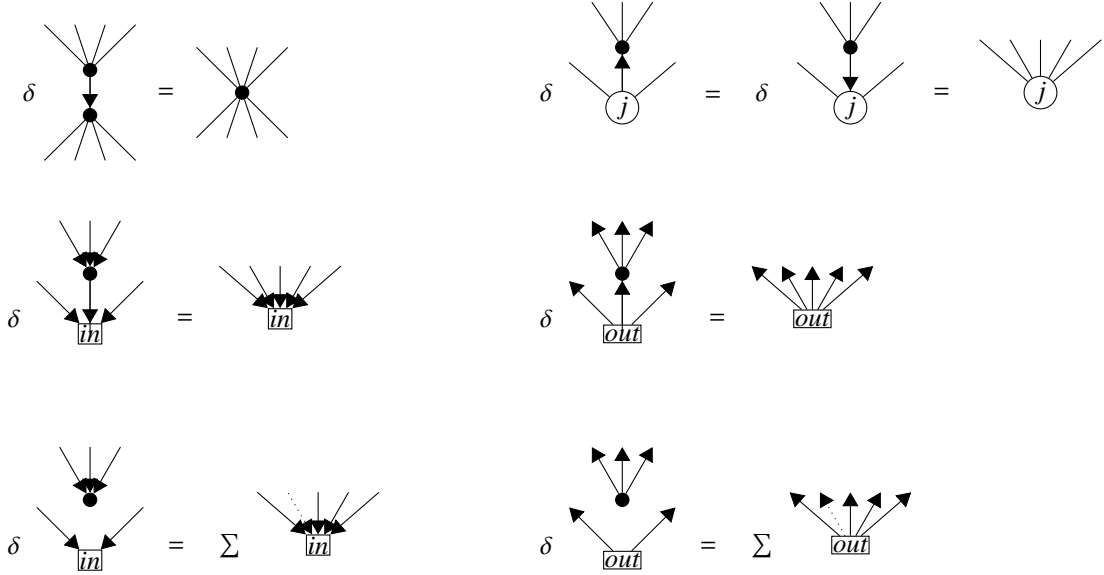

See Figure \ref{fig:graphs1dualdiff} for a graphical description of the differential.
Note that automatically every graph without any internal vertices is a cycle.

\begin{lemma}
\label{lem:g1hom}
 A basis for the homology $H(\pdu\Graphs_1(m))$ is given by the classes represented by graphs without internal vertices of the following form.
\begin{enumerate}
 \item The vertex $\vout$ either has no incident edges or it has exactly one, which connects it to $\vin$.
 \item The external vertex $j$ ($j=1,2,\dots,m)$ either does not have any outgoing edges, or it has exactly one, by which it is connected to one of the vertices $1,2,\dots, j-1,\vin$. 
\end{enumerate}
In particular
\[
H(\pdu\Graphs_1)= (\pdu\Graphs_1)_0 / \delta (\pdu\Graphs_1)_1.
\]
where $(\pdu\Graphs_1)_j$ is the subspace spanned by graphs with $j$ internal vertices. 
\end{lemma}

The proof of this statement is lengthier than that of the analogous statement for $\pdu\Graphs$ (Lemma \ref{lem:pduGraphs}). It will be given below. For now, let us believe the statement. 

\begin{lemma}
\label{lem:anotherbasis1}
 An alternative basis for $H(\pdu\Graphs_1(m))$ is given by classes of graphs of the following form
\begin{itemize}
 \item The vertex $\vout$ either has no incident edges or it has exactly one, which connects it to $\vin$.
 \item Each external vertex has at most one incoming and at most one outgoing edge.
 \item The vertex $\vin$ has at most one incoming edge not connecting to $\vout$.
 \item In each connected component, the lowest labelled vertex has no outgoing edge. If the connected component contains $\vin$, then $\vin$ counts as the lowest labelled vertex.
\end{itemize}

\end{lemma}
This basis consists of graphs formed by several ``strings'' of external vertices, such that the lowest numbered vertex in each string is at the end of the string, and the arrows point to this lowest vertex.

\begin{proof}
 The new (tentative) basis has the same cardinality as that of Lemma \ref{lem:g1hom}, and hence it suffices to check that each element of the basis of Lemma \ref{lem:g1hom} can be expressed as a linear combination of elements in our new basis, modulo $\delta (\pdu\Graphs_1)_1$. This step is tedious, we only sketch it. Let $\Gamma$ be some graph of the basis of Lemma \ref{lem:g1hom}. First, using the relations we can make all connected components not containing $\vin$ into strings with the lowest vertex in each string at one end, similarly to Lemma \ref{lem:anotherbasis}. This uses relations coming from exact elements as in Figure \ref{fig:graphs1dualdiff}, top right. Here care has to be taken that the internal vertex never has all arrows incoming or outgoing for otherwise we will produce more terms through the operations in Figure \ref{fig:graphs1dualdiff}, bottom row. Next we turn around arrows which are pointing in the wrong direction using boundaries of graphs with an internal vertex of valence 2, both edges incoming. As ``side effect'' this may attach some vertices to $\vin$ due to the part of $\delta$ depicted in Figure \ref{fig:graphs1dualdiff}, bottom left. Ignore this for now. At the end we have a sum of graphs, in which each connected component not containing $\vin$ is a string and properly oriented towards the lowest vertex of that string, which is at one end. We still have to take care of the connected component of $\vin$, which is a tree. First, fix one of the graphs produced, and assume all edges in the connected component of $\vin$ were already oriented towards $\vin$. Then, using the same tricks as above, we could rewrite the graph as a linear combination of graphs in which the connected component of $\vin$ is a properly oriented string, and we were done. If not, we can reduce the number of improperly oriented edges by adding the boundary of a graph with a valence 2 internal vertex, as we did before. This might, as side product produce graphs with (i) additional connected components, or (ii) graphs with (unoriented) cycles. Graphs with oriented cycles can be checked to be exact. In graphs with additional connected components, we first make the new connected components into properly oriented strings as we did before and then proceed. In any case, at the end there is either (i) at least one vertex less in the connected component of $\vin$ or (ii) at least one improperly oriented edge less. So the procedure converges.
\end{proof}

\begin{cor}
\label{cor:g1congcalc1}
 $H(\Graphs_1(m)) \cong \calc_1(m)$.
\end{cor}
\begin{proof}
A basis for $\calc_1(m)$ is given by expressions of the form 
\[
 e_{G,k,j_1,\dots, j_k,0}=\iota_{G(X_1, \dots,\hat X_{j_1}, \dots, \hat X_{j_k}, \dots,  X_m)}L_{X_{j_k}} \dots  L_{X_{j_1}}
\]
and 
\[
 e_{G,k,j_1,\dots, j_k,1}=\iota_{G(X_1, \dots,\hat X_{j_1}, \dots, \hat X_{j_k}, \dots,  X_m)}L_{X_{j_k}} \dots  L_{X_{j_1}}d
\]
where $k=0,1,\dots m$, $j_1,\dots,j_k\in [m]$ such that $j_p\neq j_q$ for $p\neq q$, and $G$ ranges over some basis of $\Ger(m-k)$. Let us take a basis of $\Ger(m-k)$ as described in section \ref{sec:LVproof}.
We want to show that the above map $\calc_1(m)\to H(\Graphs_1(m))$ is a bijection. As in corollary \ref{cor:hgraphs}, there is a natural one to one map between the basis of $\calc_1(m)$ above and the basis described in Lemma \ref{lem:anotherbasis1}. Namely, each Lie word in the product of lie words $G(\cdots)$ becomes one string-like connected component. There is a string of vertices $j_k,\dots, j_1$ connecting to $\vin$. If the $d$ is present, there is an additional edge $\vout\to \vin$. 
Unfortunately, the two basis are not dual to each other. in other words, the matrix describing the pairing (call it \emph{pairing matrix}) is not the identity matrix. But one can change the ordering such that the pairing matrix becomes triangular, with nonzero diagonal. This also proves the Corollary. Order the basis vectors of Lemma \ref{lem:anotherbasis1} as follows:
Graphs with $\vout$ of valence 1 are considered higher than those with $\vout$ of valence 0. Amoung both groups, graphs which have more vertices in the connected component of $\vin$ are considered higher. Among the remaining equivalence classes order the graphs arbitrarily.

We claim that with this ordering, the pairing matrix is triangular, with nonzero diagonal.
The verification of this fact is lengthy to write down, but straightforward. So we leave it to the reader.
\end{proof}
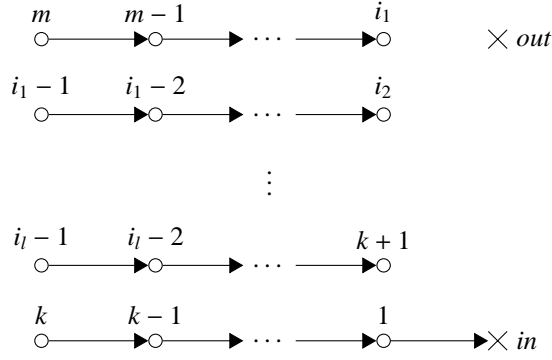
\begin{figure}
 \centering
\makeatletter{}\usetikzlibrary{matrix}
\usetikzlibrary{arrows}
\usetikzlibrary{shapes}
\usetikzlibrary{through}
\usetikzlibrary{calc,3d}
\usetikzlibrary{decorations,decorations.pathmorphing}
\[
\begin{tikzpicture}[
int/.style={circle, draw, fill, minimum size=5pt, inner sep=0},
ext/.style={circle, draw, fill=white, minimum size=5pt, inner sep=1pt},
helper/.style={coordinate},point/.style={circle, draw, fill, inner sep =1pt},
de/.style={-triangle 60},
point/.style={circle, draw, fill, minimum size=3pt, inner sep=0pt},
xst/.style={cross out, draw, minimum size=5 },
]
\node[xst, label=0:{$out$}] (out) at (6,0) {};
\node[ext,label=90:{$m$}] (em) at (0,0) {};
\node[ext, label=90:{$m-1$}] (em1) at (1.5,0) {};
\node[] (emd1) at (3,0) {$\dots$};
\node[ext, label=90:{$i_1$}] (ei1) at (4.5,0) {};
\draw[-triangle 60] (em) edge (em1) (em1) edge (emd1) (emd1) edge (ei1);

\begin{scope}[yshift=-1cm]
\node[ext,label=90:{$i_1-1$}] (em) at (0,0) {};
\node[ext, label=90:{$i_1-2$}] (em1) at (1.5,0) {};
\node[] (emd1) at (3,0) {$\dots$};
\node[ext, label=90:{$i_2$}] (ei1) at (4.5,0) {};
\draw[-triangle 60] (em) edge (em1) (em1) edge (emd1) (emd1) edge (ei1);
\end{scope}

\begin{scope}[yshift=-3cm]
\node at (3,1.2) {$\vdots$};
\node[ext,label=90:{$i_l-1$}] (em) at (0,0) {};
\node[ext, label=90:{$i_l-2$}] (em1) at (1.5,0) {};
\node[] (emd1) at (3,0) {$\dots$};
\node[ext, label=90:{$k+1$}] (ei1) at (4.5,0) {};
\draw[-triangle 60] (em) edge (em1) (em1) edge (emd1) (emd1) edge (ei1);
\end{scope}

\begin{scope}[yshift=-4cm]
\node[ext,label=90:{$k$}] (em) at (0,0) {};
\node[ext, label=90:{$k-1$}] (em1) at (1.5,0) {};
\node[] (emd1) at (3,0) {$\dots$};
\node[ext, label=90:{$1$}] (ei1) at (4.5,0) {};
\node[xst, label=0:{$in$}] (in) at (6,0) {};
\draw[-triangle 60] (em) edge (em1) (em1) edge (emd1) (emd1) edge (ei1) (ei1) edge (in);
\end{scope}

\end{tikzpicture}
\] 
\caption{\label{fig:dualsg1} Drawing of the graph $\Gamma$ constructed in the proof of Corollary \ref{cor:g1congcalc1}.}
\end{figure}

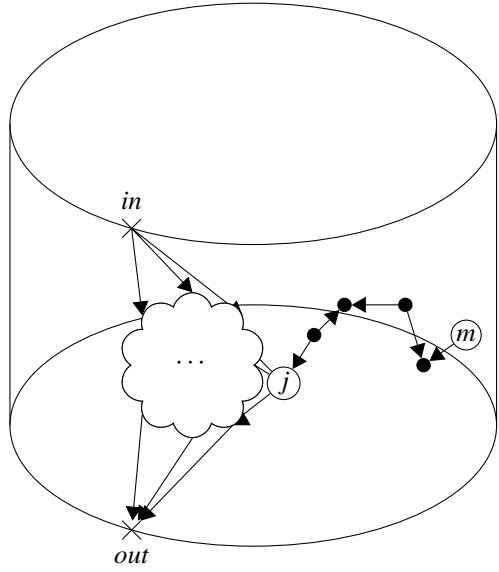
\begin{figure}
 \centering
\makeatletter{}\usetikzlibrary{matrix}
\usetikzlibrary{arrows}
\usetikzlibrary{shapes}
\usetikzlibrary{through}
\usetikzlibrary{calc,3d}
\usetikzlibrary{decorations,decorations.pathmorphing}
\[
\begin{tikzpicture}[
scale=.8,
int/.style={circle, draw, fill, minimum size=5pt, inner sep=0},
ext/.style={circle, draw, fill=white, minimum size=5pt, inner sep=1pt},
helper/.style={coordinate},point/.style={circle, draw, fill, inner sep =1pt},
de/.style={-triangle 60},
point/.style={circle, draw, fill, minimum size=3pt, inner sep=0pt},
xst/.style={cross out, draw, minimum size=5 },
]
\begin{scope}[]
\draw (0,0) ellipse (4cm and 2cm);
\draw (0,5) ellipse (4cm and 2cm);
\draw (-4,0)--(-4,5) (4,0)--(4,5);
\node [xst, label=-90:{$out$}] (out) at ($(0,0)+(-120:4 and 2)$) {};
\node [xst, label=90:{$in$}] (in) at ($(0,5)+(-120:4 and 2)$) {};

\node [cloud, draw, fill=white, minimum size=20, inner sep =10] (cl) at (-1,1) {$\cdots$};
\node[ext] (ej) at (.5,.7) {$j$};
\node[ext] (em) at (3.5,1.5) {$m$};
\node[int] (i1) at (2.8,1) {};
\node[int] (i2) at (2.5,2) {};
\node[int] (i3) at (1.5,2) {};
\node[int] (i4) at (1,1.5) {};
\draw[-triangle 60] (em) edge (i1) (i2) edge (i1) (i2) edge (i3) (i4) edge (i3) (i4) edge (ej) ;

\draw[ -triangle 60]  (in.base) edge (cl.north) 
  (in.base) edge (cl.north west)
  (in.base) edge (cl.north east)
  (cl.south) edge (out)
  (cl.south west) edge (out)
  (cl.south east) edge (out)
;
\draw (ej) edge (cl.east) (ej) edge (cl.north) (ej) -- (cl.south east);

\node [cloud, draw, fill=white, minimum size=20, inner sep =10] (cl) at (-1,1) {$\cdots$};

\end{scope}
\end{tikzpicture}
\] 
\caption{\label{fig:2valstring} Two spectral sequences in the proof of Lemma \ref{lem:g1hom} come from the filtrations by type of and distance to the closest non-bivalent-internal vertex. In the example shown, the type is ``external'' and the distance is 5. The cloud stands for the rest of the graph, which we do not care about.}
\end{figure}

Next let us turn to the proof of Lemma \ref{lem:g1hom}.
The proof will proceed by an induction on the number of external vertices $m$. For $\pdu\Graphs_1(0)=\calc(0)\cong\R \oplus \R[1]$ the statement of the Lemma is obviously true. 
Consider next the case $m>0$. Adapting \cite{LV}, one can decompose 
\[
\pdu\Graphs_1 = C_0 \oplus C_1 \oplus C_{\geq 2}
\]
where the part $C_0$ is spanned by graphs with 0 edges incident at the external vertex $m$, $C_1$ is spanned by graphs with exactly one edge incident at $m$ and $C_{\geq 2}$ is spanned by graphs with two or more edges incident at $m$. There are several components of the differential between these spaces as follows:
\[
\begin{tikzpicture}
  \matrix(m)[diagram]{ C_0 & C_1 &  C_{\geq 2} \\};
\draw[->]
 (m-1-1) edge[loop above]  node[above] {$\delta_{0}$}  ()
 (m-1-2) edge[loop above]  node[above] {$\delta_{1}$}  ()
 (m-1-3) edge[loop above]  node[above] {$\delta_{2}$}  ()
 (m-1-2) edge  node[above] {$\delta_{10}$}  (m-1-1)
 (m-1-2) edge[bend left]  node[above] {$\delta_{12}$} (m-1-3)
 (m-1-3) edge[bend left]  node[below] {$\delta_{21}$} (m-1-2);
\end{tikzpicture}
\]
We take the associated spectral sequence, such that the first differential is $\delta_{12}$. Let us call this spectral sequence ``spectral sequence 1'' to distinguish it from a second one we need below. The differential $\delta_{12}$ contracts the edge incident at $m$, if the vertex the edge connects to is internal and at least trivalent. It is not hard to see that $\delta_{12}$ is surjective. Hence the first convergent of the spectral sequence is
\[
 (C_0, \delta_0) \oplus ( C_{1,cl}, \delta_1)
\]
where $C_{1,cl}\subset C_1$ is the $\delta_{12}$-closed subspace. 
Let us compute the next term in the spectral sequence. The homology of $(C_0, \delta_0)$ is easy to evaluate. Since vertex $m$ is not connected to anything, this complex is isomorphic to the complex $\pdu\Graphs_1(m-1)$. But, by the induction hypothesis, we know its homology.
The homology of $( C_{1,cl}, \delta_1)$ is a bit harder to compute.
The space $C_{1,cl}$ decomposes as follows.
\[
 C_{1,cl} = C_{1,cl}^i \oplus C_{1,cl}^e \oplus C_{1,cl}^{io}
\]
Here the space $C_{1,cl}^i\subset C_{1,cl}$ is spanned by graphs in which the closest vertex to the external vertex $m$, which is not bivalent and internal, is an at least trivalent internal vertex. Similarly, $C_{1,cl}^i\subset C_{1,cl}$ is spanned by graphs in which the closest non-bivalent-internal vertex to the external vertex $m$ is an external vertex, and the space $C_{1,cl}^{io}$ is the spanned by graphs such that the closest non-bivalent-internal vertex is either $in$ or $out$. See Figure \ref{fig:2valstring} for a graphical explanation. The differential $\delta_1$ has the following components.

\[
\begin{tikzpicture}
  \matrix(m)[diagram]{ C_{1,cl}^i & C_{1,cl}^e &  C_{1,cl}^{io} \\};
\draw[->]
 (m-1-1) edge[loop above]  node[above] {$\tilde \delta_{1}$}  ()
 (m-1-2) edge[loop above]  node[above] {$\tilde \delta_{1}$}  ()
 (m-1-3) edge[loop above]  node[above] {$\tilde \delta_{1}$}  ()
 (m-1-1) edge (m-1-2)
 (m-1-2) edge (m-1-3)
 (m-1-1) edge[bend right] (m-1-3);
\end{tikzpicture}
\]
Let us take a spectral sequence, such that the first differential is $\tilde \delta_{1}$. Let us call this spectral sequence ``spectral sequence 2''. Each of the spaces $C_{1,cl}^i$,  $C_{1,cl}^e$, $C_{1,cl}^{io}$ splits further:
\[
 C_{1,cl}^i = C_{1,1}^i \oplus C_{1,2}^i \oplus C_{1,2}^i \oplus \cdots
\]
and similarly for $C_{1,cl}^e$, $C_{1,cl}^{io}$. Here $C_{1,j}^i\subset C_{1,cl}^i$ is spanned by graphs in which vertex $m$ has distance $j$ from the nearest non-bivalent-internal vertex of valence $\geq 3$. Hear the ``distance'' between two vertices is the length, counted in edges, of the shortest path between them. The differential $\tilde \delta_1$ has the following components:
\[
\begin{tikzpicture}
  \matrix(m)[diagram]{ C_{1,1}^? & C_{1,2}^? &  C_{1,3}^? & \cdots \\};
\draw[->]
 (m-1-1) edge[loop above]  node[above] {$\delta_{1}'$}  ()
 (m-1-2) edge[loop above]  node[above] {$\delta_{1}'$}  ()
 (m-1-3) edge[loop above]  node[above] {$\delta_{1}'$}  ()
 (m-1-2) edge node[above] {$\delta_{1}''$}  (m-1-1)
 (m-1-3) edge node[above] {$\delta_{1}''$} (m-1-2)
 (m-1-4) edge node[above] {$\delta_{1}''$} (m-1-3);
\end{tikzpicture}
\]
Here the $?$ can be either $i$, $e$ or $io$.
We take another spectral sequence, call it ``spectral sequence 3'', such that the first differential is $\delta_1''$. The homologies of the three complexes $(C_{1,1}^?, \delta_1'')$, $?=i,e, io$ have to be evaluated separately.
\begin{lemma}
\label{lem:onecomplexcohom}
\[
H(C_{1,cl}^e, \delta_1'') \cong \Graphs_1(m-1)\otimes \R^{m-1}
\]
Representatives of the homology classes are given by graphs which are obtained from a graph $\Gamma\in\Graphs_1(m-1)$ by adding the external vertex $m$ and connecting it with an edge to one of the external vertices $1,\dots, m-1$. \end{lemma}
\begin{proof}
 The complex splits into subcomplexes $C_{\Gamma,j}$ labelled by elements $\Gamma\in\Graphs_1(m-1)$ and a number $j\in [m-1]$. Concretely,  $C_{\Gamma,j}$ is spanned by graphs in which $m$ is connected to $j$ by a chain of bivalent internal vertices. The differential $\delta_1''$ contracts one edge in this chain. It is now helpful to change the basis of $C_{\Gamma,j}$. Instead of taking the basis where the edges in the chain are decorated by a direction, we take a basis where edges are decorated by symbols $a$ or $s$, standing for the symmetric or antisymmetric combination of directions. Then the differential $\delta_1''$ contracts only the edges labelled by $s$. Hence the complex $C_{\Gamma,j}$ splits again into subcomplexes
\[
 C_{\Gamma,j} = \oplus_{k\geq 0} C_{\Gamma,j}^k
\]
where $C_{\Gamma,j}^k$ is the subcomplex spanned by graphs with $k$ edges labelled by $a$ in the chain.
It is not hard to see that $(C_{\Gamma,j}^k, \delta_1'')$ is acyclic for $k>0$, and that $H(C_{\Gamma,j}^k, \delta_1'')$ is one-dimensional, with the cohomology class represented by a graph as in the statement of the Lemma.
\end{proof}

\begin{lemma}
 The complex $(\sum_{j\geq 1} C_{1,j}^i, \delta_1'')$ is acyclic.
\end{lemma}
\begin{proof}
 This proof is a copy of the previous one. The difference here is that $C_{\Gamma,\alpha}^0=0$, where $C_{\Gamma,\alpha}^0$ is defined as in the previous proof. This is because by construction of $(C_1)_{cl}$ the edge at $m$ must not be contractible if it connects to a $\geq 3$-valent internal vertex.
\end{proof}

\begin{lemma}
\[
H(\sum_{j\geq 1} C_{1,j}^{io}, \delta_1'') \cong \Graphs_1(m-1)
\]
Here representatives of the homology classes are given by graphs which are obtained from a graph $\Gamma\in\Graphs_1(m-1)$ by adding the external vertex $m$ and an edge connecting $m$ to $\vin$. \end{lemma}
\begin{proof}
 The proof is a variation the previous proofs. Again by setting up another spectral sequence on the length of the string of vertices connecting $m$ to $\vout$ or $\vin$ (same as above), we can restrict to the part of the differential $\delta''$, that reduces that length by one. Let us call it $\delta_1''$ (same as above). Note that there are no parts of $\delta''$ that reduce the length of the string by more than 1. This is because we had forbidden graphs with components with only internal edges between $\vin$ and $\vout$ in the definition of $\Graphs_1$. The differential $\delta_1''$ concretely does the following (i) contract some edge along the string, which does not connect directly to $\vin$ or $\vout$, or (ii) contract the edge attached to $\vin$ ($\vout$) if the adjacent edge is pointing towards (away from) $\vin$ ($\vout$) or (iii) delete the edge attached to $\vin$ ($\vout$) and reconnect the string to $\vout$ ($\vin$) if the next adjacent edge is pointing away from (towards) $\vin$ ($\vout$). Note that, if we temporarily identify $\vin$ and $\vout$, then operations (ii) and (iii) together become just the contraction of the edge adjacent to $\vin$/$\vout$. Hence our complex is combinatorially the same as the one considered in the proof of Lemma \ref{lem:onecomplexcohom} and the same arguments used there show the present Lemma.
\end{proof}

Let us compute the next term in spectral sequence 3, i.e., take the homology of  
\[
 \Graphs_1(m-1)\otimes \R^m
\]
under $\delta_1'$. The $\delta_1'$ here is just the usual differential on $\Graphs_1(m-1)$, and hence we obtain
\[
 H(\Graphs_1(m-1))\otimes \R^m.
\]
The first term we again know by the induction hypothesis. The spectral sequence 3, and also spectral sequences 2 and 1 terminate at this point, since the differentials always annihilate one internal vertex, and the classes in $H(\Graphs_1(m-1))\otimes \R^m$ can be repesented by graphs without any internal vertices.
Also, from the proof one sees that one can indeed take the representatives in the form stated in Lemma \ref{lem:g1hom}.
 Hence the lemma is proven. 
\hfil\qed
\begin{rem}
 We didn't discuss here the convergence of the spectral sequences. They converge to the homology. This can be seen as follows. The degree is defined as
\[
 \deg = 2\#(\text{internal vertices})-\#(\text{edges}).
\]
Since the differential always annihilates one internal vertex and one edge, we may equivalently take the degree to be
\begin{align*}
 \deg' &= 2\#(\text{internal vertices})-\#(\text{edges}) - \frac{3}{2}(\#(\text{internal vertices})-\#(\text{edges}))
\\&=\frac{1}{2}(\#(\text{internal vertices})+\#(\text{edges}))
\end{align*}
The filtrations leading to the spectral sequences above are compatible with the grading by $\deg'$. Furthermore they are automatically bounded 
since the subspace of graphs of fixed $\deg'$ is finite dimensional. Hence the spctral sequences converge to homology.
\end{rem}

\makeatletter{}\section{The cohomology of \texorpdfstring{$\fSGraphs$ and $\SGraphs$}{fSGraphs and SGraphs} }
\label{sec:sgraphsproofsuppl}
The goal of this section is to compute the cohomology of the operadic bimodules $\fSGraphs$ and $\SGraphs$ introduced in section \ref{sec:SGraphs}. We will proceed in two steps. First we will construct a list of cocycles. Secondly, we show that these cocycles span the cohomology.

\subsection{Construction of cocycles and result}

There is a special element $\mathbf H\in \SGraphs(1)$, which is the image of the fundamental chain of $D_K(1)\cong \{pt\}$ under the map $D_K(1)\to \SGraphs(1)$. It is the
graphical version of the ``twisted'' Hochschild-Kostant-Rosenberg isomorphism. The element $\mathbf H$ is closed since it is the image of a point.
One has a morphism of right $\fGraphs$ modules 
\begin{gather*}
\fGraphs[-1] \to \fSGraphs \\
\Gamma \mapsto \mathbf{H}(\Gamma):= \mathbf{H}\circ \Gamma.
\end{gather*}
Here the right hand ``$\circ$'' comes from the operadic right action of $\fGraphs$ on $\fSGraphs$.
By restricting the same formula, one also has a morphism has of right $\Graphs$-modules
\begin{gather*}
\Graphs \otimes S(\GC_2\oplus \prod_{k=4,8,\dots}\R[-k] )[-1] \to \SGraphs \\
\Gamma \mapsto \mathbf{H}(\Gamma):= \mathbf{H}\circ \Gamma.
\end{gather*}

Here $\GC_2$ is again M. Kontsevich's graph complex. The $\R[-k]$ stand for even wheel graphs of the form 
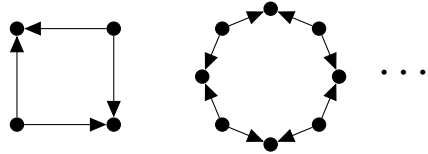
\[
\begin{tikzpicture}[every edge/.style={draw, -triangle 45}
, scale=.9]
\begin{scope}
\node [int] (i1) at (45:1) {};
\node [int] (i2) at (135:1) {};
\node [int] (i3) at (225:1) {};
\node [int] (i4) at (-45:1) {};
\draw (i1) edge (i2) edge (i4)
      (i3) edge (i2) edge (i4);
\end{scope}
\begin{scope}[xshift=3cm]
\node [int] (i1) at (45:1) {};
\node [int] (i2) at (135:1) {};
\node [int] (i3) at (225:1) {};
\node [int] (i4) at (-45:1) {};
\node [int] (i5) at (0:1) {};
\node [int] (i6) at (90:1) {};
\node [int] (i7) at (180:1) {};
\node [int] (i8) at (-90:1) {};
\draw (i1) edge (i5) edge (i6)
      (i2) edge (i6) edge (i7)
      (i3) edge (i7) edge (i8)
      (i4) edge (i5) edge (i8);
\end{scope}

\node at (5,0) {\huge$\cdots$};
\end{tikzpicture}
\]
In $\fGraphs$ these wheels are coboundaries of some odd wheels. However, since any odd wheel has necessarily a valence 2 vertex with one incoming and one outgoing edge, these wheels will actually produce nontrivial cohomology classes in $\SGraphs$. This should be seen as an artifact of our definition of $\SGraphs$. 

Now let us compute the cohomology of $\fSGraphs$ and $\SGraphs$.\footnote{Unfortunately the cohomology if $\SGraphs$ will not be $\Ger$. This indicates that the author's definition of $\SGraphs$ is bad. The author agrees, but currently does not know a better one. The definition given here will suffice to perform the globalization later on.}
Let us first discuss what to expect. First the (total space of) the cohomology of $\fSGraphs$ ( and $\SGraphs$) is a commutative algebra, even a Gerstenhaber algebra because of the left $\Ger=H(\Br)$ action. Furthermore the above embeddings of right modules (in particular of complexes) produce a lot of cohomology classes. There is however one class in $\SGraphs(0)$ we are missing so far. Let $m$ again be the Maurer-Cartan element (the universal star product) from above. The Maurer-Cartan equation can be written in the form
\[
\delta m +\frac{1}{2} \co{m}{m}=0.
\]
Here the differential $\delta$ has two parts, one splitting internal type $I$ vertices into two and one splitting type II vertices.
The bracket formally resembles the Gerstenhaber bracket.
Let $G$ be the gradation operator, multiplying a graph by the total number of internal type I and type II vertices.
We have the following equations:
\begin{align*}
G\delta &= \delta (G+\mathit{id})
& 
G\co{\cdot}{\cdot}
&=
\co{G\cdot}{\cdot}
+
\co{\cdot}{G\cdot}
-\co{\cdot}{\cdot}.
\end{align*}
It follows from these equations and the Maurer-Cartan equation that the element $M:=(G-2)m$ satisfies
\[
\delta M +\co{m}{M} = 0.
\]
Hence $M$ is a cocycle with respect to the twisted differential. The following proposition says that the above cohomology classes and the multiplicative product generate all cohomology of $\fSGraphs$ and $\SGraphs$.

\begin{prop}
\label{prop:SGraphscohom}
 The map $\fGraphs\to \fSGraphs$ induces an embedding on cohomology. Moreover we have 
\[
H(\fSGraphs(n))
\cong
\begin{cases}
S\left(H(\GC_2)[-2] \oplus \R M[-2] \oplus \prod_{k=5,9,\dots} \R[-k]\right)[-1]
& \quad \text{for $n=0$} \\
\Ger(n)
\otimes
S\left(H(\GC_2)[-2] \oplus \R M[-2] \oplus \prod_{k=5,9,\dots} \R[-k]\right)[-1]
& \quad \text{for $n>0$} 
\end{cases}
\] 
In both cases the generator $M$ is as above and the classes $\R[-k]$ are given by odd wheels (cf. proposition \ref{prop:twgracohom}).
 The cohomology of $\Graphs$ is 
 \[
 H(\SGraphs) 
 \cong
  \Ger \otimes S\left( H(\GC_2)[-2] \oplus \prod_{k=4,8,\dots}\R[-k]\right)[-1].
 \]
Here the additional generators $\R[-k]$ stand for even wheels as discussed above.
\end{prop}
We will not need the proposition in this paper, so we only sketch the proof.

\begin{rem}
V. Dolgushev's preprint \cite{vasilystable} contains a computation very similar to the following. In particular, the cohomology of $\fSGraphs(0)$ is computed in loc. cit. 
\end{rem}

\begin{proof}[Sketch of proof]
The proofs for $\SGraphs$ and $\fSGraphs$ are nearly identical. Let us first do the proof for $\fSGraphs$ and then discuss the necessary changes for $\SGraphs$.
Consider the filtration on the number of internal type I vertices, i.e.
\[
 F^p\fSGraphs = \mathit{span}\{\Gamma \mid \Gamma \text{ an $\SGraphs$-graph with $\geq p$ internal type I vertices.}\}.
\]
The filtration is descending and bounded above.
\[
 \fSGraphs = F^0\fSGraphs \supset F^1 \fSGraphs \supset F^2\fSGraphs\supset \cdots
\]
It is furthermore complete since
\[
\fSGraphs = \lim_{\leftarrow}  \fSGraphs / F^p \fSGraphs.
\]
Let us take the associated spectral sequence. 
By what is above the spectral sequence does not necessarily converge to the true cohomology. However, together with the following claim convergence to cohomology follows.

\vspace{3mm}

{\bf Claim 1:} The spectral sequence abuts at the $E^2$ page.

\vspace{3mm}

Note that the differential on $\fSGraphs$ splits as
\[
 \delta = \delta_H + \delta_1 + \delta_{\geq 2}
\]
where $\delta_H$ leaves the number of internal type I vertices constant, $\delta_1$ increases it by one, and $\delta_{\geq 2}$ increases it by two or more. Concretely, $\delta_H$ is the ``Hochschild differential'' given by splitting type II vertices into two type II vertices.
It is the differential on the $E_0$-page of our spectral sequence.
 By arguments similar to those in the proof of the Hochschild-Kostant-Rosenberg Theorem, one can compute $E_1=H(E_0, \delta_H)$ and see that the cohomology classes are represented by graphs with all type II vertices univalent, and antisymmetric under interchange of the order of the type II vertices. The classes of these graphs span the page $E_1$ of the spectral sequence. The differential on $E_1$ is induced by $\delta_1$ from above. This differential splits as 
\[
 \delta_1 = \delta_{1,-1} + \delta_{1,0} + \delta_{1,\geq 1}
\]
where $\delta_{1,-1}$ reduces the number of type II vertices by one, $\delta_{1,0}$ leaves it constant and $\delta_{1,\geq 1}$ increases it by one or more. Concretely, $\delta_{1,-1}$ removes a type II vertex and makes it a univalent type I vertex, with one incoming edge. 

\vspace{3mm}

{\bf Claim 2:} The cohomology $H(E_1, \delta_{1,-1})$ can be identified with the quotient complex 
\[
E_1' = E_1^{(0)} /  \delta_{1,-1} E_1^{(1)}
\]
where $E_1^{(j)}$ is spanned by graphs with $j$ type II vertices.

\vspace{3mm}

Believing the claim, is is not hard to see that the projection 
\[
(E_1, \delta_1) \to (E_1', \delta_{1,0})
\]
is a quasi-isomorphism. Here we denote by $\delta_{1,0}$ the induced differential on $E_1'$, abusing notation.
Let us show Claim 2. A natural basis of $E_1$ is given by graphs without type II vertices, but with an extra number attached to each vertex. The number signifies how many type II vertices are to be attached to that vertex.
For example, the graph 
\[
\begin{tikzpicture}
\node[int, label=180:2] (i1) at (0,0) {}; 
\node[ext, label=0:1] (e1) at (1,0) {$1$};
\draw (i1) edge[bend left, -triangle 45] (e1);
\draw (e1) edge[bend left, -triangle 45] (i1);
\end{tikzpicture}
\]
corresponds to the following antisymmetric linear combination:
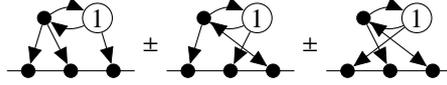
\[
\begin{tikzpicture}[scale=.7,
every edge/.style={draw, -triangle 45}]
\begin{scope}
\draw (-.7,-1) -- (1.7, -1);
\node[int] (i1) at (0,0) {}; 
\node[ext] (e1) at (1,0) {$1$};
\draw (i1) edge[bend left, -triangle 45] (e1);
\draw (e1) edge[bend left, -triangle 45] (i1);
\node[int] (x1) at (-.3, -1) {};
\node[int] (x2) at (.5, -1) {};
\node[int] (x3) at (1.3, -1) {};
\draw (i1) edge (x1) edge (x2);
\draw (e1) edge (x3);
\end{scope}
\begin{scope}[xshift=3cm]
\draw (-.7,-1) -- (1.7, -1);
\node[int] (i1) at (0,0) {}; 
\node[ext] (e1) at (1,0) {$1$};
\draw (i1) edge[bend left, -triangle 45] (e1);
\draw (e1) edge[bend left, -triangle 45] (i1);
\node[int] (x1) at (-.3, -1) {};
\node[int] (x2) at (.5, -1) {};
\node[int] (x3) at (1.3, -1) {};
\draw (i1) edge (x1) edge (x3);
\draw (e1) edge (x2);
\end{scope}
\begin{scope}[xshift=6cm]
\draw (-.7,-1) -- (1.7, -1);
\node[int] (i1) at (0,0) {}; 
\node[ext] (e1) at (1,0) {$1$};
\draw (i1) edge[bend left, -triangle 45] (e1);
\draw (e1) edge[bend left, -triangle 45] (i1);
\node[int] (x1) at (-.3, -1) {};
\node[int] (x2) at (.5, -1) {};
\node[int] (x3) at (1.3, -1) {};
\draw (i1) edge (x3) edge (x2);
\draw (e1) edge (x1);
\end{scope}
\node at (2, -.5) {$\pm$};
\node at (5, -.5) {$\pm$};
\end{tikzpicture}
\]
The differential $\delta_{1,-1}$ acts by decreasing the attached number of some vertex by one and adding a valence 0 vertex:
\[
\begin{tikzpicture}[scale=.7,
every edge/.style={draw, -triangle 45}]
\begin{scope}
\node[int, label=0:{$n$}] (i1) at (0,0) {}; 
\foreach \x in {-.3, 0, .3}
{
\draw (i1) -- +(\x, .5);
\draw (i1) -- +(\x, -.5);
}
\end{scope}
\node at (2, 0) {$\mapsto\quad n \cdot($};
\begin{scope}[xshift=4.2cm]
\node[int] (i2) at (-1,0) {}; 
\node[int, label=0:{$n-1$}] (i1) at (0,0) {}; 
\foreach \x in {-.3, 0, .3}
{
\draw (i1) -- +(\x, .5);
\draw (i1) -- +(\x, -.5);
}
\draw (i2) edge (i1);
\end{scope}
\node at (5.7, 0) {$)$};
\end{tikzpicture}
\]
Here the number attached to the new internal vertex is zero, we do not display it.
Let us define an operator $h$ (a homotopy) on $E_1$ that acts in the reverse way by deleting valence one vertices with one outgoing edge and attached number zero.
\[
\begin{tikzpicture}[scale=.7,
every edge/.style={draw, -triangle 45}]
\begin{scope}[xshift=3cm]
\node[int, label=0:{$n+1$}] (i1) at (0,0) {}; 
\foreach \x in {-.3, 0, .3}
{
\draw (i1) -- +(\x, .5);
\draw (i1) -- +(\x, -.5);
}
\end{scope}
\node at (1.8, 0) {$\mapsto$};
\begin{scope}
\node[int] (i2) at (-1,0) {}; 
\node[int, label=0:{$n$}] (i1) at (0,0) {}; 
\foreach \x in {-.3, 0, .3}
{
\draw (i1) -- +(\x, .5);
\draw (i1) -- +(\x, -.5);
}
\draw (i2) edge (i1);
\end{scope}
\end{tikzpicture}
\]
Being more careful with the signs, one can check that for a graph $\Gamma$
\[
(\delta_{1,-1}h+h \delta_{1,-1} )(\Gamma)
=
(\sum_\alpha n_\alpha + N) \Gamma
\]
where the sum runs over all vertices $\alpha$ of $\Gamma$, $n_\alpha$ is the number attached to $\alpha$ and $N$ is the total number of valence 1 vertices with one incoming edge and attached number 0. From this Claim 2 immediately follows.

Now, given Claim 2, let us compute the cohomology $H(E_1, \delta_1) \cong H(E_1', \delta_{1,0})$. The computation is a variation on the computation of $H(\fGraphs)$ and similar computations in \cite{megrt}. The argument is sketched in Appendix \ref{sec:auxcomputation}.
The result is that the cohomology is given by (the classes of) the elements in the statement of the proposition. I.e., 
\[
H(E_1', \delta_{1,0})
\cong
\begin{cases}
S(H(\GC_2) \oplus \R M[-2] \oplus \prod_{k=5,9,\dots} \R[-k])[2]
& \quad \text{for $n=0$} \\
\Ger(n)
\otimes
S(H(\GC_2) \oplus \R M[-2] \oplus \prod_{k=5,9,\dots} \R[-k])
& \quad \text{for $n>0$} 
\end{cases}
\]

The graphs spanning $E_1'$ may have multiple connected components, some containing external vertices and some not. 
Since the differential acts separately on each connected component, the cohomology will be a symmetric product space with generators the cohomology of connected graphs.  
The connected components with external vertices together produce the $\Ger$-part of the cohomology. The connected components without external vertices produce the remainder of the cohomology, including the wheels and the graph cohomology.

Note that there is one difference to the very similar computations in \cite{megrt}. In our case the graph containing only a single internal vertex is closed. 
This corresponds to the class of $M$.

Now let us return to our spectral sequence.
Since we constructed cocycles in $\fSGraphs$ representing all the classes occurring above,\footnote{To construct representatives of classes of graphs having multiple valence 0 internal vertices, we use the product coming from the left $\Br$-action.} the spectral sequence terminates at this point. All the higher differentials are zero. This shows also Claim 1, and we are done for $\fSGraphs$.

Consider next $\SGraphs$. The proof is formally identical, except that the complex $(E_1', \delta_{1,0})$ is a bit smaller in the case of $\SGraphs$, since certain types of vertices are forbidden to appear. However, almost the same calculation goes through.
\end{proof}

\subsection{Auxiliary computation}\label{sec:auxcomputation}
Let us fill the remaining gap in the proof of Proposition \ref{prop:SGraphscohom} above by computing the cohomology of the complex $(E_1', \delta_{1,0})$ appearing there. We will separately consider the cases $\fSGraphs$ and $\SGraphs$. Very similar computations can be found in \cite{megrt}.

\subsubsection{$\fSGraphs$ case}
Let us compute the cohomology of the complex $(E_1', \delta_{1,0})$ from the proof of proposition \ref{prop:SGraphscohom}. To recall, elements of $E_1'$ are linear combinations of directed graphs with internal and external vertices, modulo graphs with valence 1 internal vertices whose incident edge is incoming. The differential is given by splitting vertices, thus producing one new internal vertex as usual.

Each graph can be decomposed into connected components, and in particular into connected components containing or not containing external vertices.
The complex $E_1'$ can hence be written as
\[
E_1' = E_{1, ext}'\otimes S(E_{1, \conn}').
\]
Here $E_{1, ext}'$ is the subcomplex spanned by graphs all of whose connected components contain an external vertex.
$E_{1, \conn}'$ is the complex formed by connected graphs with only internal vertices.

Let us first consider $E_{1, ext}'$. This complex splits as 
\[
E_{1, ext}' = V_1\oplus V_2
\]
where $V_1$ is spanned by graphs that have at least one univalent internal vertex and $V_2$ is spanned by the graphs that do not contain univalent internal vertices. As in the proof of Proposition 3 of \cite{megrt} one shows that $V_1$ is acyclic. Analogously to \cite[Appendix K]{megrt} one then shows that $V_2$ is quasi-isomorphic to $\Graphs_n$ and hence $H(V_2)\cong e_2(n)$.

Next consider the purely internal components.
In analogy with the proof of Proposition 3 in \cite{megrt}, let us split:
\[
E_{1, \conn}' = C_1 \oplus C_{\geq 2} 
\]
where $C_1$ is the subcomplex of graphs containing at least one valence 1 internal vertex and $C_{\geq 2}$ is the subcomplex of graphs not containing one. 
Let us consider the two parts in turn. $C_1$ splits further
\[
C_{1,2} \oplus C_{1, 3}
\]
where $C_{1,2}$ is the subcomplex of graphs that do not contain an internal vertex of valence $\geq 3$, and $C_{1, 3}$ the subcomplex of graphs that do contain a valence $\geq 3$ internal vertex. Along the lines of \cite{megrt}, Proposition 3, one shows that $C_{1, 3}$ is acyclic.
$C_{1,2}$ is given by ``string-like'' graphs of the form
\[
\begin{tikzpicture}[scale=.7,
every edge/.style={draw, -triangle 45}]
\node[int] (i1) at (0,0) {};
\node[int] (i2) at (1,0) {};
\node[int] (i3) at (2,0) {};
\node (i4) at (3.5,0) {$\cdots$};
\node[int] (i5) at (5,0) {};
\node[int] (i6) at (6,0) {};
\draw (i1) edge (i2) (i2) edge (i3) (i4) edge (i3) (i5) edge (i4) (i6) edge (i5);
\end{tikzpicture}
\]
The first and last edge must be inwards pointing, the orientation of the other edges is arbitrary. Note that for this reason there cannot be a string with two vertices. However, there can be one with only one vertex.
It is not hard to check that
\[
H(C_{1,2}) = \R[-2]
\]
with the single class being represented by the string with one vertex.

Next consider the complex $C_{\geq 2}$. It splits further
\[
C_{\geq 2, 2}\oplus C_{\geq 2, 3}
\]
where $C_{\geq 2, 2}$ is spanned by graphs without at least  trivalent vertices, while $C_{\geq 2, 3}$ is spanned by graphs with at least one $\geq 3$-valent vertex.
The cohomologies have been computed in \cite[Appendix K]{megrt}:
\begin{align*}
H(C_{\geq 2, 2}) &= \prod_{k=5,9,\dots} \R[-k] \\
H(C_{\geq 2, 3}) &= H(\GC_2).
\end{align*}

Hence the result stated in the proof of Proposition \ref{prop:SGraphscohom} follows.

\subsubsection{$\SGraphs$ case}
Next, consider a subcomplex $E_1''\subset E_1'$ spanned by graphs which do not contain vertices of valence $\leq 1$, and that do not contain vertices of valence 2, with one incoming and one outgoing edge. This is the ``version of $E_1'$'' that occurs in the proof of proposition \ref{prop:SGraphscohom} for the case of $\SGraphs$.

The computation is similar to that in the last subsection. Let us merely remark on the differences. First, the complex $V_1$ above obviously does not occur in this case. The complex $V_2$ has to be shrunk so as to be spanned by graphs  
without valence 2 vertices with one incoming and one outgoing edge. Accordingly, we cannot directly use the result of \cite{megrt}, Appendix K, to conclude that the cohomology of the resulting complex, say $V_2'$, is $H(V_2')=e_2(n)$. Copying the trick from Appendix K in \cite{megrt}, one can impose a filtration on the number of valence 2 vertices. The first differential in the associated spectral sequence, say $d$, creates valence 2 vertices. We claim that its cohomology is $\Graphs_2(n)$ (as in loc. cit.). 
To each graph one can associate its \emph{core}, which is the undirected graph obtained by (i) forgetting the orientations of edges and (ii) deleting all valence two internal vertices and joining the two incident edges of each deleted vertex. Then $(V_2', d )$ splits into a direct product of subcomplexes, one for each automorphism class of cores. Each such subcomplex is the space of invariants under the core's automorphism group of a product of complexes, one for each edge. The complex associated to the edge has the form
\[
\R\alpha_1\oplus \R \beta_1
\to 
\R\alpha_2\oplus \R \beta_2
\to 
\R\alpha_3\oplus \R \beta_3
\to \cdots
\]
where $\alpha_j$ stands for a string of $j$ alternatingly oriented edges, starting with an outgoing edge, and $\beta_j$ stands for a string of $j$ alternatingly oriented edges, starting with an incoming edge.
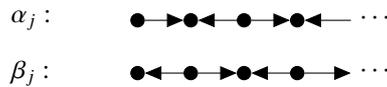
\[
\begin{tikzpicture}[scale=.7,
every edge/.style={draw, -triangle 45}]
\node at (-2,0) {$\alpha_j:$};
\node[int] (i1) at (0,0) {};
\node[int] (i2) at (1,0) {};
\node[int] (i3) at (2,0) {};
\node[int] (i4) at (3,0) {};
\node (i5) at (4.5,0) {$\cdots$};
\draw (i1) edge (i2) (i3) edge (i2) (i3) edge (i4) (i5) edge (i4);
\begin{scope}[yshift=-1cm]
\node at (-2,0) {$\beta_j:$};
\node[int] (i1) at (0,0) {};
\node[int] (i2) at (1,0) {};
\node[int] (i3) at (2,0) {};
\node[int] (i4) at (3,0) {};
\node (i5) at (4.5,0) {$\cdots$};
\draw (i2) edge (i1) (i2) edge (i3) (i4) edge (i3) (i4) edge (i5);
\end{scope}
\end{tikzpicture}
\]

The differential maps 
\begin{align*}
d \alpha_j &= \beta_{j+1} + (-1)^j \alpha_{j+1} \\ 
d \beta_j &= \alpha_{j+1} + (-1)^j \beta_{j+1}.
\end{align*}
It is not hard to see that the cohomology of the resulting complex is one dimensional, and represented by $\alpha_1+\beta_1$. Hence $H(V_2, d)=\Graphs_2(n)$. It follows that $H(V_2) = e_2(n)$.

Next consider the purely internal components. The complex $C_1$ from the previous subsection does not occur in the present case. The complex $C_{\geq 2}$ has to be shrunk so as to be spanned by graphs  
without valence 2 vertices with one incoming and one outgoing edge. Again we can split this subcomplex, say $C_{\geq 2}'$, into 
\[
C_{\geq 2}' = C_{\geq 2, 2}'\oplus C_{\geq 2, 3}'
\] 
as in the previous subsection. $C_{\geq 2, 2}'$ is spanned by wheels of length $4,8,12,\dots$. Note that wheels of odd length necessarily have a vertex with one incoming and one outgoing edge and hence do not occur. Wheels of length $2,6,10,\dots$ are zero by symmetry.
The differential on $C_{\geq 2, 2}'$ acts as zero. Hence
\[
H(C_{\geq 2, 2}') = \prod_{k=4,8,\dots} \R[-k].
\]
Next consider $C_{\geq 2, 3}'$. 
Again, by (almost) the same arguments as before we can show that  $C_{\geq 2, 3}'$ is quasi-isomorphic to its subcomplex of undirected $\geq 3$-valent graphs.\footnote{This subcomplex is embedded into $C_{\geq 2, 3}'$ by mapping an undirected graph to the sum over all graphs obtained by assigning directions on the edges.} This means that $H(C_{\geq 2, 3}', d )\cong \GC_2$, as we wanted to show.

\nocite{*}
\bibliographystyle{plain}
\bibliography{biblio}

\end{document}